\documentclass[11pt,letter]{amsart}
\usepackage{amscd}
\usepackage{amssymb}
\usepackage{graphicx}
\usepackage{color}
\usepackage{framed}
\usepackage[centering,text={15cm,20cm}]{geometry}
\usepackage[all]{xy}
\usepackage{mathrsfs}
\usepackage{bbm}

\definecolor{shadecolor}{rgb}{1,0.9,0.7}

\newtheorem{theorem}{Theorem}[section]
\newtheorem{lemma}[theorem]{Lemma}
\newtheorem{proposition}[theorem]{Proposition}
\newtheorem{corollary}[theorem]{Corollary}
\newtheorem{conjecture}[theorem]{Conjecture}

\theoremstyle{definition}
\newtheorem{definition}[theorem]{Definition}
\newtheorem{construction}[theorem]{Construction}
\newtheorem{convention}[theorem]{Convention}
\newtheorem{example}[theorem]{Example}

\theoremstyle{remark}
\newtheorem{remark}[theorem]{Remark}
\newtheorem{remarks}[theorem]{Remarks}

\numberwithin{equation}{section}
\numberwithin{figure}{section}

\newcommand{\doublebar}[1]{\shortstack{$\overline{\hphantom{#1}}
	$\\[-4pt] $\overline{#1}$}}
\newcommand{\underlinehigh}[1]{\shortstack{#1\\[-5pt]
	\underline{\hphantom{#1}}}}
\newcommand{\lperp}{{}^\perp\!}
\newcommand{\lparallel}{{}^\parallel\!}
\newcommand{\bbfamily}{\fontencoding{U}\fontfamily{bbold}\selectfont}
\newcommand{\textbb}[1]{{\bbfamily#1}}

\newcommand {\lfor} {\mbox{\textbb{[}}}
\newcommand {\rfor} {\mbox{\textbb{]}}}

\newcommand{\NN} {\mathbb{N}}
\newcommand{\ZZ} {\mathbb{Z}}
\newcommand{\QQ} {\mathbb{Q}}
\newcommand{\RR} {\mathbb{R}}
\newcommand{\CC} {\mathbb{C}}

\newcommand{\PP} {\mathbb{P}}
\renewcommand{\AA} {\mathbb{A}}

\newcommand {\shAff} {\mathcal{A}\text{\textit{ff}}}

\newcommand {\shF}  {\mathcal{F}}

\newcommand {\shL}  {\mathcal{L}}

\newcommand {\shN}  {\mathcal{N}}
\newcommand {\shPL} {\mathcal{PL}}
\newcommand {\shQ}  {\mathcal{Q}}

\newcommand {\foB}  {\mathfrak{B}}

\newcommand {\foD}  {\mathfrak{D}}

\newcommand {\foX}  {\mathfrak{X}}

\newcommand {\fob}  {\mathfrak{b}}
\newcommand {\foc}  {\mathfrak{c}}
\newcommand {\fod}  {\mathfrak{d}}
\newcommand {\foe}  {\mathfrak{e}}
\newcommand {\fof}  {\mathfrak{f}}
\newcommand {\fog}  {\mathfrak{g}}
\newcommand {\foh}  {\mathfrak{h}}

\newcommand {\foj}  {\mathfrak{j}}

\newcommand {\fol}  {\mathfrak{l}}

\newcommand {\fop}  {\mathfrak{p}}
\newcommand {\foq}  {\mathfrak{q}}
\newcommand {\forr}  {\mathfrak{r}}

\newcommand {\fou}  {\mathfrak{u}}
\newcommand {\fov}  {\mathfrak{v}}
\newcommand {\fow}  {\mathfrak{w}}

\newcommand {\foz}  {\mathfrak{z}}

\newcommand {\scrA}  {\mathscr{A}}

\newcommand {\Aff}  {\operatorname{Aff}}

\newcommand {\Ad} {\operatorname{Ad}}

\newcommand {\as}  {\mathrm{as}}

\newcommand {\base}  {\operatorname{Base}}

\newcommand {\Chambers} {\operatorname{Chambers}}

\newcommand {\cl}  {\operatorname{cl}}
\newcommand {\codim} {\operatorname{codim}}

\newcommand {\conv} {\operatorname{conv}}

\newcommand {\depth} {\operatorname{depth}}

\newcommand {\dlog} {\operatorname{dlog}}

\newcommand {\eps}  {\varepsilon}

\newcommand {\GL}  {\operatorname{GL}}
\newcommand {\Glue} {\underline{\mathrm{Glue}}}

\newcommand {\gp}  {{\operatorname{gp}}}

\newcommand {\Hom}  {\operatorname{Hom}}
\newcommand {\height}  {\operatorname{ht}}

\newcommand {\id}  {\operatorname{id}}
\newcommand {\im}  {\operatorname{im}}
\newcommand {\inc}  {\mathrm{in}}

\newcommand {\Int}  {\operatorname{Int}}

\newcommand {\Joints} {\operatorname{Joints}}

\newcommand {\kk} {\Bbbk}

\newcommand {\limdir} {\varinjlim}
\newcommand {\liminv} {\varprojlim}

\newcommand {\LogRings} {\underlinehigh{LogRings}}

\newcommand {\LPoly} {\underlinehigh{LPoly}}
\newcommand {\lra}  {\longrightarrow}

\newcommand {\M} {\mathcal{M}}

\newcommand {\Map} {{\operatorname{Map}}}
\renewcommand {\max} {{\operatorname{max}}}

\newcommand {\no}  {\mathrm{no}}

\newcommand {\NS}  {\operatorname{NS}}

\renewcommand{\O}  {\mathcal{O}}

\newcommand {\ord}  {\operatorname{ord}}

\renewcommand{\P}  {\mathscr{P}}

\newcommand {\PM} {\operatorname{PM}}

\newcommand {\pre}  {\mathrm{pre}}

\newcommand {\Proj} {\operatorname{Proj}}

\newcommand {\red}  {{\operatorname{red}}}

\newcommand {\Rings} {\underlinehigh{Rings}}

\newcommand {\scrS}  {\mathscr{S}}

\newcommand {\Sch} {\underline{\mathrm{Sch}}}

\newcommand {\shLS} {\mathcal{LS}}
\newcommand {\sides} {\operatorname{Sides}}
\newcommand {\sing} {\mathrm{sing}}
\newcommand {\Sing} {\operatorname{Sing}}
\newcommand {\SL}  {\operatorname{SL}}
\newcommand {\Spec} {\operatorname{Spec}}

\newcommand {\Strata} {\underline{\mathrm{Strata}}}

\newcommand {\sgn} {\operatorname{sgn}}
\newcommand {\std} {\mathrm{std}}

\newcommand {\tcont} {\operatorname{cont_\mathit{t}}}
\newcommand {\tlog} {\operatorname{tlog}}

\newcommand {\topp}  {\operatorname{Top}}

\newcommand {\trdeg}  {\operatorname{trdeg}}

\newcommand {\D} {\frak D}

\newcommand {\X} {\frak X}

\newcommand {\Z} {\frak Z}

\newcommand {\I} {\mathrm{\,I}}
\newcommand {\II} {\mathrm{\,II}}
\newcommand {\III} {\mathrm{\,III}}

\def\mapright#1{\smash{
  \mathop{\longrightarrow}\limits^{#1}}}

\def\mydate{\ifcase\month \or January\or February\or March\or
April\or May\or June\or July\or August\or September\or October\or 
November\or December\fi \space\number\day,\space\number\year}

\begin{document}

%===========================================================
\title[Affine and complex geometry]{From real affine geometry
to complex geometry}

\author{Mark Gross} \address{UCSD Mathematics,
9500 Gilman Drive, La Jolla, CA 92093-0112, USA}
%\curraddr{}
\email{mgross@math.ucsd.edu}
\thanks{This work was partially supported by NSF grant 0505325 and
DFG priority programs ``Globale Methoden in der komplexen Geometrie''
and ``Globale Differentialgeometrie''.}

% author two information
\author{Bernd Siebert}
\address{Mathematisches Institut,
Albert-Ludwigs-Universit\"at Freiburg, Eckerstra\ss e~1, 79104 Freiburg,
Germany}
\curraddr{Universit\"at Hamburg, FB~Mathematik, Bundesstr.~55,
20146~Hamburg, Germany}
\email{bernd.siebert@math.uni-hamburg.de}

\date{July 20, 2011}

\begin{abstract}
We construct from a real affine manifold with singularities (a
tropical manifold) a degeneration of Calabi-Yau manifolds. This solves
a fundamental problem in mirror symmetry. Furthermore, a striking
feature of our approach is that it yields an explicit and canonical
order-by-order description of the degeneration via families of tropical
trees.

This gives complete control of the $B$-model side of mirror
symmetry in terms of tropical geometry. For example, we expect our
deformation parameter is a canonical coordinate, and expect period
calculations to be expressible in terms of tropical curves. We
anticipate this will lead to a proof of mirror symmetry via tropical
methods. This paper is the key step of the program we initiated in
\cite{logmirror}.
\end{abstract}

\maketitle
\tableofcontents
%\bigskip

%===========================================================
%===========================================================
\section*{Introduction.}

Toric geometry links the integral affine geometry of convex polytopes
to complex geometry. On the complex side this correspondence works via
equivariant partial completions of algebraic tori. It is thus
genuinely linear in nature and deals exclusively with rational
varieties. This paper provides a non-linear extension of this
correspondence producing (degenerations of) varieties with effective
anti-canonical divisor. Among other things we obtain a new and rather
surprising method for the construction of varieties with trivial canonical
bundle by discrete methods. It generalizes the Batyrev-Borisov
construction of Calabi-Yau varieties as complete intersections in
toric varieties and a number of other, non-toric constructions.
One may even hope to obtain all deformation classes of varieties
with trivial canonical bundle which contain maximally unipotent boundary
points.

The data on the affine side consists of a topological manifold $B$
built by gluing integral polytopes in an affine manner along their
boundaries, along with compatible affine charts at the vertices of the
resulting polyhedral decomposition $\P$ of $B$. This endows $B$ with
an integral affine structure on the complement of a codimension two
subset $\Delta\subseteq B$. The subset $\Delta$ is covered by the
cells of the barycentric subdivision of $\P$ neither containing a
vertex of $\P$ nor intersecting the interiors of top-dimensional cells
of $\P$. This notion of affine manifolds with singularities allows
for many more interesting closed examples than without singularities.
For example, the two-torus is the only closed oriented surface with a
non-singular affine structure \cite{benzecri}\cite{milnor}, while
$S^2$ has many affine structures with singularities, for example as
base of an elliptically fibred K3 surface. Furthermore, such
integral affine manifolds arise naturally from boundaries of reflexive
polytopes \cite{GBB}\cite{haasezharkov}.

On the complex side we consider \emph{toric degenerations} $\pi:\X\to
T$ of complex varieties as introduced in \cite{logmirror},
Definition~4.1; the central fibre $X$ is a union of toric varieties
glued torically by identifying pairs of toric prime divisors, and
$\pi$ is \'etale locally toric near the zero-dimensional toric strata
of $X$. Beware that in general $\X$ is not a toric variety.
In the same paper it was shown how such a degeneration gives
rise to an affine manifold with singularities and polyhedral
decomposition $(B,\P)$ as before (\cite{logmirror}, Definition~4.13).
Combinatorially $\P$ is the dual intersection complex of the central
fibre, that is, $k$-cells of $\P$ correspond to codimension-$k$
intersections of irreducible components. The affine structure reflects
the toric nature of the degeneration.

The main result Theorem~\ref{main thm} of this paper deals with the
inverse problem: given $(B,\P)$, find a toric degeneration $\pi:\X\to
T$ with dual intersection complex $(B,\P)$. 

In the case without singularities Mumford already used toric methods
to write down such degenerations, notably for the class of abelian
varieties \cite{mumford}. 

To solve this problem in much greater generality, we need to make
three assumptions. First, we need the existence of the central fibre
of the toric degeneration as a \emph{toric log Calabi-Yau space} as
defined in \cite{logmirror}, Definition~4.3. Constructing such a
space with dual intersection complex $(B,\P)$ was the chief focus of
\cite{logmirror}. This condition is necessary, as $\X\to T$ induces
such a log Calabi-Yau space structure on $X$. Also, Theorem~5.4 in
\cite{logmirror} gives natural sufficient conditions
(\emph{positivity} and \emph{simplicity}) in terms of the local affine
monodromy around $\Delta\subseteq B$. Positivity is a kind of
convexity property that again is necessary, while simplicity should be
viewed as a maximal degeneration property that implies primitivity of
the local monodromy. In instances where $(B,\P)$ is non-simple the
existence of a log Calabi-Yau structure can be explicitly checked by
the results of \S3.3 in \cite{logmirror}, notably by Theorem~3.27.

Second, we assume the existence of a \emph{polarization} for $(B,\P)$.
This is a multi-valued, convex, piecewise affine function on $B$. If
$(B,\P)$ is the dual intersection complex of a toric log Calabi-Yau
space $X$, this condition is closely related to projectivity of
$X$, and is in fact equivalent to it provided $H^1(B,\QQ)=0$, see
\cite{logmirror}, Theorem~2.34. This condition is clearly not
necessary as in some cases, for example in dimension two, log
deformation theory gives the same result without any projectivity
assumptions. While in our algorithm different polarizations lead to
isomorphic families, the polarization is a basic ingredient that
appears to be crucial for globalizing the local deformations
consistently.

Third, we need a condition we term \emph{locally rigid} (see
Definition~\ref{def:locally rigid}) on $X$. This is a somewhat
technical condition, which essentially implies that at each step of
our construction, the choices we need to make are unique. Without this
condition, non-uniqueness can lead to obstructions to solving the
problem. Simplicity implies local rigidity (Remark~\ref{simple implies
locally rigid}). As a result, the Main Theorem implies the following
reconstruction theorem.

\begin{theorem}\label{reconstruction theorem}
Any polarized affine manifold with singularities with positive and
simple monodromy arises as the dual intersection complex of a toric
degeneration.
\end{theorem}

In fact, our Main Theorem (Theorem~\ref{main thm}) applies more
broadly when the dual intersection complex  $(B,\P)$ is non-compact,
corresponding to the case where a general fibre has only effective
anti-canonical class. In principle, one should also be able to deal
with the case when $(B,\P)$ has boundary, where the corresponding
complex manifold is not complete. However, in this situation we
expect a Landau-Ginzburg potential to play an important role, and
there are additional subtleties to the argument. We have chosen not
to deal with these issues here, and will consider this case
elsewhere.

The proof of the Main Theorem gives far more than the existence of a
toric degeneration. It gives a \emph{canonical, explicit} $k$-th order
deformation $X_k\to \Spec \kk[t]/(t^{k+1})$ of $X$ for any $k$.
Furthermore, this degeneration is specified using data of a tropical
nature.

Let us expand on this description. First, we explain the role the
polarization plays. Given the polarization $\varphi$, a piecewise
linear multi-valued function on $B$, one can construct the 
\emph{discrete Legendre transform} of the triple $(B,\P,\varphi)$,
which is another triple $(\check B,\check\P,\check\varphi)$. If
$\varphi$ comes from an ample line bundle on a log Calabi-Yau space
$X$ with dual intersection complex $(B,\P)$, then $(\check B,\check
\P)$ is the \emph{intersection complex}, whose maximal cells are the
Newton polytopes defining the polarized irreducible components of
$X$ (see \cite{logmirror}, \S\S1.5 and 4.2). The data governing the
deformations of $X$ then consist of what we call a \emph{structure},
which is a collection of \emph{slabs and walls}: these are codimension
one polyhedral subsets of $B$, contained locally in affine
hyperplanes, along with some attached data of a ring automorphism
which is used in our gluing construction. This structure has an
important tropical interpretation: morally, a structure can be viewed
as a union of tropical trees in $\check B$ with leaves on $\Delta$. We
will not define the precise notion of tropical curves on affine
manifolds with singularities, as this is not needed for the proof and
the correct general definition is not yet entirely clear, but see
\cite{Gr4} for some further discussion of this. We inductively
construct structures  $\scrS_k$ for $k\ge 0$, with $\scrS_k$ providing
sufficient data to construct $X_k\to\Spec \kk[t]/(t^{k+1})$. Morally,
$\scrS_k$ can be viewed as the union of ``tropical trees of degree
$k$''.

The actual degenerations $X_k\to\Spec\kk[t]/(t^{k+1})$ are
constructed from the structure $\scrS_k$ by gluing together certain
canonical thickenings of affine pieces of irreducible components of
$X$, with the gluings specified by the automorphisms attached to the
slabs and walls of $\scrS_k$. The main difficulty in the inductive
construction of $\scrS_{k+1}$ from $\scrS_k$ is the need to maintain
compatibility of this gluing. For this purpose, we adapt a key lemma
of Kontsevich and Soibelman from \cite{ks}, which essentially
expresses commutators of automorphisms in a standard form as a
product of automorphisms in this standard form.

We end this introduction with a number of remarks, historical and
otherwise.\\[2ex] (1) Our construction can be viewed as a
``non-linear'' generalization of Mumford's and Alexeev and Nakamura's 
construction of degenerations of abelian varieties
\cite{mumford}\cite{AlexNak}\cite{Alex}. If $B=\RR^n/\Gamma$, where
$\Gamma\subseteq\ZZ^n$ is a sublattice, then one obtains a
degeneration of abelian varieties. Here $\Delta$ is empty, and the
structures $\scrS_k$ can be taken to be empty too: there are no
``corrections'' to construct the deformation.\\[1ex] (2) Let
$B_0:=B\setminus\Delta$. Then we can define $X(B_0):=
T_{B_0}/\Lambda$, where $\Lambda$ is the local system of integral flat
vector fields. This is a torus bundle over $B_0$, and it inherits a
natural complex structure from the tangent bundle $T_{B_0}$. One basic
problem that arises in the Strominger-Yau-Zaslow approach to mirror
symmetry \cite{SYZ}  is that one would like to compactify $X(B_0)$ to
a complex manifold $X(B)$. Because of the singularities, the complex
structure is in fact not the correct one, and this compactification
cannot be performed in the complex category, even when it can be
performed in the topological category (see e.g.~\cite{Gr4},
\cite{topology}). One needs to deform the complex structure before
this compactification can be expected to exist. Typically, one
considers an asymptotic version of this problem: consider
$X_{\epsilon}(B_0)= T_{B_0}/\epsilon\Lambda$. Then one expects for
small $\epsilon>0$, there is a small deformation of the complex
structure on $X_{\epsilon}(B_0)$ which can be compactified. This
problem has been discussed already in a number of places, see
\cite{Fuk}\cite{KoSo1}\cite{Gr2}. The results of this paper, combined
with the results of \cite{topology} as described in \cite{Announce},
can be viewed as giving as  complete a solution to this problem as one
could hope for.

This problem was first attacked directly by Fukaya in \cite{Fuk}, in the
two-dimensional case. Fukaya gave heuristic arguments suggesting that
the needed deformation should be concentrated along certain trees
of gradient lines on $B$ with leaves on $\Delta$. This direct
analytic approach seems to be very difficult; nevertheless, it 
gave a hint as to the relevant data for controlling the
deformations.
\\[1ex]
(3) In \cite{ks}, Kontsevich and Soibelman proposed an alternative
approach to the reconstruction problem, suggesting one should
construct a rigid analytic space rather than a complex manifold from
$B$. They showed how to do this in dimension two. Here, the same trees
of gradient flow lines as in \cite{Fuk} emerge, this time with certain
automorphisms attached to the edges of the trees. The proof relies on
a group-theoretic lemma which we also use here. The advantage of using
rigid analytic spaces is that most convergence issues become rather
simple. However, there is one part of their argument which is rather
technical: to prove convergence near the singularities, one has to
control the gradient flow lines to avoid returning to some small
neighbourhood of the singularities. This technical issue, along with
some other points, appears to make this approach rather difficult to
generalize directly to higher dimensions.

In a sense, we surmount these difficulties by discretizing the problem
and passing to the discrete Legendre transform $(\check
B,\check\P)$. The advantage of working with the dual affine
manifold $\check B$ is that the gradient flow lines on $B$ become
straight lines on $\check B$. These are obviously much easier to work
with and control. On the other hand, this leaves us no ability to
avoid a neighbourhood of the singularities. As a result, we have to
deal with some compositions of automorphisms which involve terms of
order $0$; this introduces terms in our expressions with denominators,
which have to be controlled. This is a significant technical problem,
relatively easy in dimension two, but much harder in dimensions three
and higher, and the solution to this problem occupies 
\S4 of this paper. If one were to rewrite this paper in the
dimension two case only, it would be considerably shorter. Given our
current level of understanding, it seems that a price must be paid
somewhere near the singularities, whether it be Kontsevich and
Soibelman's genericity arguments or our algebraic arguments. It would
be nice to find a simpler solution to these problems.

It is also worthwhile making a historical remark here. We had the
original idea of constructing smoothings by gluing thickenings of
affine pieces of irreducible components of $X$ in 2003. It was also
clear to us that the gluing maps should propagate along straight lines
on $\check B$. However, we abandoned this approach for a while,
attempting to find a Bogomolov-Tian-Todorov argument for
smoothability. We returned to the question of explicit smoothings in
2005, and realised the group-theoretic Kontsevich-Soibelman lemma
applied in our situation, thus enabling us to complete the
argument.

We also comment that if one is only interested in the two-dimensional
case and one does not care about the explicit smoothings, but only the
existence of a smoothing, then $X$ can be smoothed using techniques
of \cite{Friedman 1983} or \cite{Kawamata; Namikawa 1994} directly, as
was known to us in 2001.\\[1ex] 
(4) We cannot overemphasize the importance of this result to
understanding mirror symmetry. Our structures, in  a sense, give a
complete description of the $B$-model side of mirror symmetry, at a
much deeper level than the usual description in terms of periods.
Furthermore, our description of the $B$-model side is tropical in
nature. It is well-known \cite{mikhalkin}\cite{nisi} that one should
expect a correspondence between tropical curves on $B$ and families of
holomorphic curves on the corresponding degeneration $\X\to T$. Thus
it is an important point for understanding mirror symmetry that the
construction of the corresponding complex manifold is controlled by
tropical curves on $\check B$, hence by holomorphic curves on the
mirror degeneration $\check\X\to T$ corresponding to $\check B$. This
is what one expects to see in mirror symmetry, and this gives \emph{an
explicit explanation for the connection between deformations and
holomorphic curves} in mirror symmetry.

There remains the question of extracting explicit enumerative
predictions from the structures we use to build our smoothings. We
have performed some calculations in some three-dimensional examples,
and from these, we feel highly confident in the following conjectures,
stated with varying degrees of precision:

\begin{conjecture}
\begin{enumerate}
\item
The coordinate $t$ associated with the canonical $k$-th order
deformations $X_k\to\Spec \kk[t]/(t^{k+1})$ is a canonical
coordinate in the usual sense in mirror symmetry.
\item
The enumerative predictions made by calculating periods of
$X_k\to \Spec\kk[t]/(t^{k+1})$ can be described explicitly in
terms of contributions from each tropical rational curve on $\check B$
of ``degree $\le k$''. The numerical contributions are determined by
the automorphisms appearing in the structure $\scrS_k$.
\item
The automorphisms attached to walls of $\scrS_k$ can be interpreted as
``raw enumerative data'' which morally counts the number of
holomorphic disks with boundary on Lagrangian tori of the mirror dual
Strominger-Yau-Zaslow fibration. The tropical trees arising in the
structures can be viewed as a tropical version of holomorphic disks. 
\end{enumerate}
\end{conjecture}

Ultimately, we believe it will be easier to read off enumerative
information directly from the structures, and that calculation of
periods should be viewed as a crude way of extracting the much more
detailed information present in the structures.\\[1ex] (5) Speculating
further, we expect that our structures will yield a useful description
of (higher) multiplication maps for homological mirror symmetry on the
$B$-model side. Our smoothings $\X\rightarrow T$ come along with
canonical polarizations by an ample line bundle $\shL$. A basis for
the space of sections of $H^0(\X,\shL^{\otimes n})$ as an
$\O_T$-algebra is given by the set of points $B\left({1\over
n}\ZZ\right)$ of points on $B$ whose coordinates are in ${1\over
n}\ZZ$. A structure should then allow us to give explicit descriptions
of the multiplication maps $H^0(\X,\shL^{\otimes n_1})\otimes
H^0(\X,\shL^{\otimes n_2}) \rightarrow H^0(\X, \shL^{\otimes
(n_1+n_2)})$, answering a question of Kontsevich. In discussions with
Mohammed Abouzaid, it has become apparent that it seems likely that
these multiplication maps and higher multiplication maps could be
described in terms of a ``tropical Morse category'' on $B$, once again
making the $B$-model side look very much like the expected structure
of the $A$-model (Fukaya category) side. Hopefully, this approach will
ultimately lead to a proof of Homological Mirror Symmetry.

Further justification of these statements will have to wait for
further work; however, \cite{part II} lays the groundwork for
computation of periods.
\medskip

\noindent
\emph{Conventions.}\ \ We work in the category $\Sch_\kk$ of separated
schemes over an algebraically closed field $\kk$ of characteristic
$0$. A \emph{variety} is a scheme of finite type over $\kk$. All our
toric varieties are normal. A \emph{toric monoid} is a finitely
generated, saturated, integral monoid. These are precisely the monoids
of the form $\sigma^\vee\cap\ZZ^n$ for $\sigma\subseteq \RR^n$ a
strictly convex, rational polyhedral cone. 

%===========================================================
%===========================================================
\section{Fundamentals}

%===========================================================
\subsection{Discrete data}
While our construction is part of the program laid out in
\cite{logmirror} only a fraction of the techniques developed there is
needed for it. To make this paper reasonably self-contained, \S\S1.1
and 1.2 therefore provide the relevant background. At the same
time we discuss a generalization from the projective Calabi-Yau
situation to semi-positive and non-complete cases. For simplicity we
restrict to the case without self-intersecting cells. The treatment of
self-intersections is, however, straightforward; it is merely a matter
of working with morphisms rather than inclusions and with algebraic
spaces rather than schemes, as done consistently in \cite{logmirror}. 

To fix notations recall that a \emph{convex polyhedron} is the
intersection of finitely many closed affine half-spaces in $\RR^n$. As
all our polyhedra are convex we usually drop the attribute ``convex''.
A polyhedron is \emph{rational} if the affine functions defining the
half spaces can be taken with rational coefficients. The dimension of
the smallest affine space containing a polyhedron $\Xi$ is its
\emph{dimension}. Its \emph{relative interior} $\Int\Xi$ is the
interior inside this affine space, and the complement $\Xi\setminus
\Int\Xi$ is called the \emph{relative boundary} $\partial\,\Xi$. If
$\dim\Xi=k$ then $\partial\,\Xi$ is itself a union of polyhedra of
dimension at most $k-1$, called \emph{faces}, obtained by intersection
of $\Xi$ with hyperplanes disjoint from $\Int\Xi$. Faces of dimensions
$k-1$ and $0$ are called \emph{facets} and
\emph{vertices}, respectively. In contrast to \cite{logmirror} our
polyhedra are not necessarily bounded, but we require the existence
of at least one vertex (so half-spaces, for example, are not
allowed). For $y\in \partial\Xi$ the \emph{tangent cone $K_y \Xi$ of
$\Xi$ at $y$} is the cone generated by differences $z-y$ for
$z\in\Xi$. If $\Xi'\subseteq\Xi$ is a face we also define
$K_{\Xi'}\Xi:= K_y\Xi$ for any $y\in \Int\Xi'$. The closure of the
cone $\RR_{\ge0}\cdot \big(\Xi\times\{1\}\big) \subseteq
\RR^n\times\RR=\RR^{n+1}$ is denoted $C(\Xi)$. Any polyhedron $\Xi$
can be written as Minkowski sum $\Xi'+ C$ of a bounded polyhedron
$\Xi'$ and a cone $C$. While $\Xi'$ is not in general unique, $C$ is
determined as the Hausdorff limit $\lim_{\eps\to 0} \eps\Xi$, and is
therefore called the \emph{asymptotic cone of $\Xi$}. Finally, if
$C\subseteq \RR^n$ is a cone then $C^\vee$ denotes
its dual as an additive monoid $\Hom(C,\RR_{\ge 0})$, viewed as a
cone in $\RR^{\dim C}\simeq \Hom(C,\RR)$.

A rational polyhedron is \emph{integral} or a \emph{lattice
polyhedron} if all its vertices are integral. The group of
\emph{integral affine transformations} $\Aff(\ZZ^n)= \ZZ^n\rtimes
\GL(n,\ZZ)$ acts on the set of integral polyhedra. Note that we
require the translational part to also be integral. Finally, if $\Xi$
is an integral polyhedron, $\Lambda_\Xi\simeq \ZZ^{\dim\Xi}$ denotes
the free abelian group of integral tangent vector fields along $\Xi$.
Note that for any $x\in\Int(\Xi)$ there is a canonical injection
$\Lambda_\Xi\to T_{\Xi,x}$ inducing an isomorphism $T_{\Xi,x}\simeq
\Lambda_{\Xi,\RR}:=\Lambda_\Xi\otimes_\ZZ \RR$.
\medskip

We consider topological manifolds with boundary built by gluing
integral convex polyhedra along their faces in an integral affine
manner. To this end consider the category $\LPoly$ with integral,
convex polyhedra as objects and integral affine isomorphisms onto
faces and the identity as morphisms. An \emph{integral polyhedral
complex} is gluing data for a collection of such polyhedra given by a
functor
\[
F: \P\lra \LPoly,
\]
for some category $\P$ such that if $\Xi\in F(\P)$ and $\Xi'\subseteq
\Xi$ is a face then $\Xi'\in F(\P)$. To avoid self-intersections we
also require that for any $\tau,\sigma\in\P$ there is at most one
morphism $e:\tau\to\sigma$. The associated topological space
is the quotient
\[
B=\coprod_{\sigma\in\P} F(\sigma)\Big/\sim,
\]
where two points $p\in F(\sigma)$, $p'\in F(\sigma')$ are equivalent
if there exists $\tau\in\P$, $q\in F(\tau)$ and morphisms $e:\tau\to
\sigma$, $e':\tau \to\sigma'$ with $p= F(e)(q)$, $p'= F(e')(q)$. (Then
$B$ is the colimit of the composition of $F$ with the forgetful
functor from $\LPoly$ to the category of topological spaces.) By abuse
of notation we usually suppress $F$ and consider the elements of $\P$
simply as subsets of $B$, called \emph{cells}, with the structure of
integral convex polyhedra understood. The set of $k$-dimensional cells
is then denoted $\P^{[k]}$, the $k$-skeleton by $\P^{[\le k]}$, and if
$n=\sup \{\dim\sigma\,|\, \sigma\in\P\}$ is finite we write
$\P_\max:=\P^{[n]}$. In this language morphisms are given by
inclusions of subsets of $B$.

While $B$ now has a well-defined affine structure on each cell
our construction also requires affine information in the normal
directions. To add this information recall that the \emph{open star} of
$\tau\in \P$ is the following open neighbourhood of $\Int \tau$:
\[
U_\tau= \bigcup_{\{\sigma\in \P\,|\, \Hom(\tau,\sigma)\neq\emptyset\}}
\Int \sigma.
\]

\begin{definition}\label{fan structure}
Let $\P$ be an integral polyhedral complex. A \emph{fan structure} along
$\tau\in \P$ is a continuous map $S_\tau:U_\tau\to \RR^k$ with
\begin{enumerate}
\item[(i)]
$S_\tau^{-1}(0)=\Int \tau$.
\item[(ii)]
If $e:\tau\to\sigma$ is a morphism then $S_\tau|_{\Int \sigma}$
is an integral affine submersion onto its image, that is, is induced
by an epimorphism $\Lambda_\sigma\to W\cap\ZZ^k$ for some vector
subspace $W\subseteq \RR^k$. 
\item[(iii)]
The collection of cones $K_e:=\RR_{\ge 0}\cdot S_\tau(\sigma\cap
U_\tau)$, $e:\tau\to\sigma$, defines a finite fan $\Sigma_\tau$ in
$\RR^k$.
\end{enumerate}
Two fan structures $S_\tau, S'_\tau: U_\tau\to \RR^k$ are considered
\emph{equivalent} if they differ only by an integral linear
transformation of $\RR^k$.
\qed
\end{definition}

If $S_\tau: U_\tau\to \RR^k$ is a fan structure along $\tau\in\P$ and
$\sigma\supseteq\tau$ then $U_\sigma\subseteq U_\tau$. The \emph{fan
structure along $\sigma$ induced by $S_\tau$} is the composition
\[
U_\sigma\lra U_\tau \stackrel{S_\tau}{\lra} \RR^k
\lra \RR^k/L_\sigma\simeq \RR^l,
\]
where $L_\sigma\subseteq \RR^k$ is the linear span of
$S_\tau(\Int\sigma)$. This is well-defined up to equivalence.

\begin{definition}\label{tropical manifold}
An \emph{integral tropical manifold} of dimension $n$ is an integral
polyhedral complex $\P$, which we assume countable, together with a
fan structure $S_v: U_v\to \RR^n$ at each vertex $v\in \P^{[0]}$ with
the following properties.
\begin{enumerate}
\item[(i)]
For any $v\in\P^{[0]}$ the support $|\Sigma_v|=\bigcup_{C\in \Sigma_v}
C$ is convex with nonempty interior (hence is an $n$-dimensional
topological manifold with boundary).
\item[(ii)]
If $v,w$ are vertices of $\tau\in \P$ then the fan structures along
$\tau$ induced from $S_v$ and $S_w$, respectively, are equivalent.
\qed
\end{enumerate}
\end{definition}

In the case with empty boundary and with all polyhedra bounded this is
what in \cite{logmirror} we called a ``toric polyhedral decomposition
of an integral affine manifold with singularities''.

The underlying topological space $B=\bigcup_{\sigma\in\P}\sigma$ of
an integral tropical manifold carries a well-defined integral affine
structure outside of a closed subset of codimension two. This
\emph{discriminant locus} $\Delta$ can be taken as follows. For each
bounded $\tau\in\P$ let $a_\tau\in\Int\tau$, and for each unbounded
$\tau\in\P$ let $a_\tau\in\Lambda_{\tau,\RR}$ be an element of the
relative interior of the asymptotic cone of $\tau$. The choice of
$a_{\tau}$ in the unbounded case must be made subject to the
constraint that if $\tau'\subseteq\tau$ is a face and $\tau'$ and
$\tau$ have the same asymptotic cone, then $a_{\tau'}=a_{\tau}$. In
this unbounded case, $a_\tau$ should be viewed as a point at
infinity. Then for any chain $\tau_1\subseteq \tau_2\ldots \subseteq
\tau_{n-1}$ with $\dim\tau_i= i$ and $\tau_i$ bounded iff $i\le r$
($r\ge 1$),
\[
\Delta_{\tau_1\ldots\tau_{n-1}}:=
\conv\big\{a_{\tau_i}\,\big|\, 1\le i\le r \big\}+
{\sum}_{i>r} \RR_{\ge0}\cdot a_{\tau_i}\subseteq
\tau_{n-1},
\]
is the Minkowski sum of an $(r-1)$-simplex with a simplicial cone of
dimension at most $n-r-1$.  Define $\Delta$ as the union of all such
polyhedra. In the unbounded case, if some unbounded edges are
parallel, it might happen that these polyhedra are not all
$(n-2)$-dimensional, but  they are always contained in one of this
type which is $(n-2)$-dimensional.

Now if $\rho\in\P^{[n-1]}$ then the connected components of
$\rho\setminus\Delta$ are in one-to-one corespondence with the
vertices of $\rho$. This is clear for bounded cells, while an
unbounded cell together with its discriminant locus retracts onto a
union of bounded cells. Thus the polyhedral structures on interiors
of $n$-cells and the fan structures near the vertices define an
atlas of integral affine charts on $B\setminus\Delta$. For
$x\in\rho\setminus\Delta$ write $v[x]$ for the unique vertex in the
same connected component of $\rho\setminus\Delta$ as $x$.

The fan structures and the polyhedral structure of the cells can be
read off from this affine structure together with the decomposition
into closed subsets provided by $\P$. This motivates the notation
$(B,\P)$ for an integral tropical manifold. In particular, we have a
flat, torsion-free connection on $T_{B\setminus \Delta}$, which we use
to parallel transport tangent vectors along homotopy classes of
paths. By integrality there is also a locally constant sheaf
$\Lambda$ of integral tangent vectors on $B\setminus \Delta$. For
$x\in(\sigma\cap B)\setminus\Delta$ we have a canonical isomorphism
$\Lambda_\sigma\to\Lambda_x\cap T_{\sigma,x}$ that we will use
liberally. If $\omega \subseteq\tau$ then this identification for
$\sigma=\omega,\tau$ and $x\in\omega$ is compatible with the inclusion
$\Lambda_\omega\to \Lambda_\tau$.

For bounded polyhedra a canonical choice of $a_\tau$ is the barycenter
of $\tau$, and this canonical choice, exhibiting $\Delta$ as
a subcomplex of the barycentric subdivision, has been used in
\cite{logmirror}. However, this choice is not in sufficiently general
position for our construction in this paper. We need that the
intersection of any proper \emph{rational} affine subspace of an
$(n-1)$-cell with $\Delta$ is transverse. The following lemma shows
that in the bounded case,
for a sufficiently general choice of the $a_\tau$ the
corresponding  discriminant locus $\Delta= \Delta(\{a_\tau\})$ does
not contain any rational point. This implies that $\Delta$ intersects
any proper rational affine subspace of an $(n-1)$-cell transversally,
for otherwise it would contain a non-empty open subset of it, hence a
rational point. We leave it to the reader to supply the more general
unbounded case.

For the formulation of the lemma note that any $a\in B$ has a
well-defined \emph{field of definition} $\kappa(a) \subseteq\RR$
generated over $\QQ$ by its coordinate entries in any integral affine
chart.

\begin{lemma}\label{perturbation lemma} Given $(B,\P)$ compact,
assume that for all $1\le r\le n-1$, we have
\[
\forall \tau_1\subsetneq\cdots\subsetneq\tau_{r}:\
\trdeg_\QQ\big(\kappa(a_{\tau_1})\cdot\ldots \cdot
\kappa(a_{\tau_r})\big)=\sum_{i=1}^r\dim \tau_i.
\]
Then $\Delta=\Delta(\{a_\tau\})$ contains no rational point. 
\end{lemma}

\proof
It suffices to check the claim on one $(n-2)$-simplex
$\Xi\subseteq\Delta$, say defined by $\tau_1\subsetneq
\ldots\subsetneq\tau_{n-1}$, $\dim\tau_i=i$.  Let
$e_1,\ldots,e_{n-1}$ be a $\QQ$-basis for $\Lambda_{\tau_{n-1},\QQ}$
adapted to the flag of $\QQ$-vector subspaces
$\Lambda_{\tau_1,\QQ}\subsetneq\ldots\subsetneq
\Lambda_{\tau_{n-1},\QQ}$, that is, $\Lambda_{\tau_i,\QQ}=\QQ
e_1+\ldots +\QQ e_i$. Then in the corresponding coordinate system
\[
a_{\tau_i}=\alpha_i^1 e_1+\ldots+\alpha_i^i e_i,\quad i=1,\ldots,n-1,
\]
for $\alpha_i^\mu\in\RR$. Then $\kappa(a_{\tau_i})=
\QQ(\alpha_i^1,\ldots,\alpha_i^i)$, and by assumption the
$\alpha_i^\mu$ are all algebraically independent over $\QQ$.

Now assume the convex hull of the $a_{\tau_i}$ contains a rational
point. Then there exist $\lambda_i\in[0,1]$ with $\sum_i\lambda_i=1$
such that
\[
\sum_i\lambda_i a_{\tau_i}=\left(\begin{matrix}
\lambda_1\alpha_1^1&
+\ldots+&\lambda_{n-1}\alpha_{n-1}^1\\
&\ddots&\vdots\\
&&\lambda_{n-1}\alpha_{n-1}^{n-1}\\
\end{matrix} \right)\ \in\QQ^{n-1}.
\]
Solving inductively expresses each $\lambda_i$ as a rational function
in $\alpha_i^\mu$ with coefficients in $\QQ$, and $\sum_i
\lambda_i=1$ gives an algebraic relation between the $\alpha_i^\mu$.
Now it is not hard to see that in solving inductively for
$\lambda_{n-1}, \ldots,\lambda_1$, any two occurring monomials are
different. Moreover, the coefficients can only all vanish if
$\lambda_i=0$ for all $i$, which is impossible since
$\sum_i\lambda_i=1$. We have thus found a non-trivial algebraic
relation among the $\alpha_i^\mu$, contradicting algebraic
independence.
\qed

\bigskip
To explain the meaning of $\Delta$ we now introduce the concept of
\emph{local monodromy}. Let $\omega\in\P^{[1]}$, $\rho\in\P^{[n-1]}$
with $\omega\subseteq\rho$, $\rho\not\subseteq \partial B$ and
$\omega$ bounded. Then $\rho$ is contained in two $n$-cells
$\sigma^\pm$, and $\omega$ contains two vertices $v^\pm$. Following
the change of affine charts given by (i) the fan structure at $v^+$
(ii) the polyhedral structure of $\sigma^+$ (iii) the fan structure at
$v^-$ (iv) the polyhedral structure of $\sigma^-$ and back to (v) the
fan structure at $v^+$ defines a transformation $T_{\omega\rho}\in
\SL(\Lambda_{v^+})$. It is shown in \cite{logmirror}, \S1.5 that this
transformation has the following form:
\begin{eqnarray}\label{monodromy constant}
T_{\omega\rho}(m)=m+\kappa_{\omega\rho}
\langle m, \check d_\rho\rangle d_\omega.
\end{eqnarray}
Here $d_\omega\in\Lambda_\omega\subseteq \Lambda_{v^+}$ and $\check
d_\rho\in \Lambda_\rho^\perp\subseteq\Lambda_{v^+}^*$ are the primitive
integral vectors pointing from $v^+$ to $v^-$ and, in the chart at
$v^+$, evaluating positively on $\sigma^+$, respectively. The constant
$\kappa_{\omega\rho}\in\ZZ$ is independent of the choices of $v^\pm$
and $\sigma^\pm$. Geometrically meaningful integral tropical manifolds
fulfill $\kappa_{\omega\rho}\ge0$:

\begin{definition}(\cite{logmirror}, Definition~1.54)
\label{def:positivity}
An integral tropical manifold is \emph{positive} if
$\kappa_{\omega\rho}\ge0$ for all $\omega\subseteq\rho$ with
$\omega$ bounded and $\rho\not\subseteq\partial B$.
\end{definition}

Slightly more generally one can consider an analogous sequence of
changes of charts for two arbitrary vertices $v$, $v'$ contained in an
$(n-1)$-cell $\rho\not\subseteq \partial B$. Since $v$ and $v'$ can be
connected by a sequence of $1$-cells contained in $\rho$ the
corresponding monodromy transformation takes the form
\begin{eqnarray}\label{$m^rho$}
m\longmapsto m+\langle m, \check d_\rho\rangle m^\rho_{vv'},
\end{eqnarray}
for a well-defined $m^\rho_{vv'}\in\Lambda_\rho$. In particular,
$m^\rho_{v^+ v^-}=\kappa_{\omega\rho} d_\omega$. In the positive
case this monodromy information can be conveniently gathered in the
\emph{monodromy polytope for $\rho$}
\begin{eqnarray}\label{monodromy polytope}
\Delta(\rho)=\operatorname{conv}\{m^\rho_{vv'}\,|\, v'\in\rho\}.
\end{eqnarray}
Here $v\in\rho$ is a fixed vertex, and a different choice of $v$
leads to a translation of $\Delta(\rho)$. Hence $\Delta(\rho)$ is a
lattice polytope in $\Lambda_\rho\otimes_\ZZ\RR$ that is
well-defined up to translation. Note that $\Delta(\rho)$ can have
any dimension between $0$ and $n-1$, and hence for fixed $v$, the
map from vertices $v'$ of $\rho$ to vertices of $\Delta(\rho)$ needs
not be injective.

\begin{remark}
One can show that the affine structure extends to a neighbourhood of
$\tau\in \P$ if and only if for every $\omega\in\P^{[1]}$,
$\rho\in\P^{[n-1]}$ with $\omega\subseteq \tau\subseteq \rho$ it holds
$\kappa_{\omega\rho}=0$. This has been used in \cite{logmirror},
Proposition~1.27 to find a smaller discriminant locus. In this paper
we choose to work with the larger discriminant locus as it slightly
simplifies the presentation later on.
\qed
\end{remark}

\medskip
Integral tropical manifolds arise algebro-geometrically from certain
degenerations of algebraic varieties whose central fibres are unions
of toric varieties and which are toroidal (``log smooth'') morphisms
near the zero-dimensional toric strata of the central fibre. In the
following we generalize the relevant definitions in \cite{logmirror}
to pairs consisting of a variety and a divisor. Recall that an
algebraic variety is called \emph{algebraically convex} if there
exists a proper map to an affine variety \cite{gola}. A
toric variety is algebraically convex if and only if the defining
fan has convex support.

\begin{definition}\label{totally degenerate varieties} A \emph{totally
degenerate CY-pair} is a reduced variety $X$ together with a reduced
divisor $D\subseteq X$ fulfilling the following conditions: Let
$\nu:\tilde X\to X$ be the normalization and $C\subseteq \tilde X$ its
conductor locus. Then $\tilde X$ is a disjoint union of algebraically
convex toric varieties, and $C$ is a reduced divisor such
that $[C]+\nu^*[D]$ is the sum of all toric prime divisors, $\nu|_C:
C\to \nu(C)$ is unramified and generically two-to-one, and the square
\[\begin{CD}
C@>>> \tilde X\\
@VVV @VV{\nu}V\\
\nu(C)@>>> X
\end{CD}\]
is cartesian and cocartesian.
\qed
\end{definition}

In other words, if $(X,D)$ is a totally degenerate CY-pair then $X$ is
built from a collection of toric varieties by identifying pairs of
toric prime divisors torically. The remaining toric prime divisors define
$D$. Note that by the toric nature of the identification maps it makes
sense to define a \emph{toric stratum} of $X$ as a toric stratum of
any irreducible component.

\begin{definition}\label{log smooth maps}
Let $T$ be the spectrum of a discrete valuation $\kk$-algebra with
closed point $O$ and uniformizing parameter $t\in \O(T)$. Let $\X$ be a
$\kk$-scheme and $\foD,X\subseteq \X$ reduced divisors. A
\emph{log smooth morphism} $\pi:(\X,X;\foD)\to (T,O)$ is a morphism
$\pi:(\X,X)\to (T,O)$ of pairs of $\kk$-schemes with the following
properties: For any $x\in \X$ there exists an \'etale neighbourhood
$U\to \X$ of $x$ such that $\pi|_U$ fits into a commutative diagram of
the following form.
\[\begin{CD}
U@>{\Phi}>> \Spec \kk[P]\\
@V{\pi|_U}VV @VV{G}V\\
T@>{\Psi}>> \Spec \kk[\NN]
\end{CD}\]
Here $P$ is a toric monoid, $\Psi$ and $G$ are defined respectively
by mapping the generator $z^1\in \kk[\NN]$ to $t$ and to a
non-constant monomial $z^{m_0}\in\kk[P]$, and $\Phi$ is
\'etale with preimage of the toric boundary divisor equal to the
pull-back to $U$ of $X\cup \foD$.
\qed
\end{definition}

This definition just rephrases that if we endow $\X$ and $T$ with the
log structures $\M_\X$ and $\M_T$ defined by $X\cup \foD\subseteq \X$
and $O\subseteq T$, respectively, then the map of log spaces
$(\X,\M_X)\to (T,\M_T)$ is (log) smooth and integral. We refer to
\cite{logmirror}, \S3.1 for a quick survey of the relevant log
geometry. However, we will largely avoid the terminology of log
structures here.

\begin{definition}\label{toric degeneration} (cf.\ \cite{logmirror},
Definition~4.1.)
Let $T$ be the spectrum of a discrete valuation $\kk$-algebra and
$O\in T$ its closed point. A \emph{toric degeneration of CY-pairs over
$T$} is a flat morphism $\pi:\X\to T$ together with a
reduced divisor $\D\subseteq \X$, with the following properties:
\begin{enumerate}
\item[(i)]
$\X$ is normal.
\item[(ii)]
The central fibre $X:=\pi^{-1}(O)$ together with $D=\D\cap X$ is a
totally degenerate CY-pair.
\item[(iii)]
Away from a closed subset $\Z\subseteq \X$ of relative codimension two
not containing any toric stratum of $X$, the map $\pi:(\X,X;\D)\to
(T,O)$ is log smooth.
\qed
\end{enumerate}
\end{definition}

In this definition we dropped the requirement that $\pi$ be proper
from \cite{logmirror}. In the non-proper case the deformation theory
of $(X,D)$ does not appear to be very well-behaved, but it still
makes sense to talk about formal toric degenerations of CY-pairs as
in the following definition.

\begin{definition}\label{formal toric degeneration}
Let $T$ be the spectrum of a discrete valuation $\kk$-algebra and
$\hat O$ the completion of $T$ at its closed point.
A \emph{formal toric degeneration of CY-pairs over $\hat O$}
is a flat morphism $\hat \pi: \hat X\to \hat O$ of formal schemes
together with a reduced divisor $\hat D\subseteq \hat X$, with
the following properties:
\begin{enumerate}
\item[(i)]
$\hat X$ is normal.
\item[(ii)]
The central fibre $X=\hat\pi^{-1}(O)\subseteq \hat X$ together with
$D=\hat D\cap X$ is a totally degenerate CY-pair.
\item[(iii)]
Away from a closed subset $Z\subseteq X$ of relative codimension two
and not containing any toric stratum of $X$, the map $\hat\pi: (\hat
X,X;\hat D)\to (\hat O,O)$ is \'etale locally on $\hat X$ isomorphic to the
completion of a log smooth morphism along the central fibre. 
\qed
\end{enumerate}
\end{definition}

Clearly, a toric degeneration of CY-pairs induces a formal toric
degeneration of CY-pairs by completion along the central fibre.

\begin{remark}
While this is not possible with our ad hoc definition of log
(smooth) structures, it does make sense to talk about abstract log
structures also on the codimension two loci $\Z\subseteq \X$ and
$Z\subseteq X$ in Definitions~\ref{toric degeneration}
and~\ref{formal toric degeneration}, respectively. Then $\Z$ or $Z$
contain the locus where this extension fails to be fine, that is,
where the log structure fails to possess a chart locally. By abuse of
notation we therefore refer to $\Z$ or $Z$ as the \emph{singular
locus of the log structure}.
\end{remark}

Before explaining how a toric degeneration defines an integral
tropical manifold we would like to review the basic duality between
convex piecewise linear functions and their Newton polyhedra,
including the unbounded case, see \cite{rockafellar}. Let $\Sigma$ be
a not necessarily complete fan defined on $N_\RR$, where as usual $N$
is a finitely generated, free abelian group. Then as a matter of
convention an (integral) \emph{piecewise linear function} on $\Sigma$
is a map
\[
\varphi:N_\RR\lra \RR\cup\{\infty\}
\]
that is an ordinary (integral) piecewise linear function on $|\Sigma|$
and that takes value $\infty$ everywhere else. The graph
$\Gamma_\varphi\subseteq N_\RR\times\RR$ is the union of the ordinary
graph of $\varphi|_{|\Sigma|}$ with
\[
\big\{(n,h) \in N_\RR \oplus\RR \,\big|\, n\in\partial|\Sigma|, h\ge
\varphi(n)\big\}.
\]
Now given an integral polyhedron $\Xi\subseteq M_\RR$,
$M=\Hom(N,\ZZ)$, define
\[
\varphi:N_\RR\lra \RR\cup\{\infty\},\quad
\varphi(n)=\sup(-n|_\Xi)=-\inf(n|_\Xi).
\]
Then $\varphi$ is a strictly convex piecewise linear function on the
normal fan $\Sigma$ of $\Xi$. For the normal fan we use the convention
that its rays are generated by the inward normals to the facets of
$\Xi$. The signs are chosen in such a way that if $C(\Xi)\subseteq
M_\RR\times\RR$ denotes the closure of the cone generated by
$\Xi\times\{1\}$, then $C(\Xi)^\vee \subseteq N_\RR\times \RR$ is the
convex hull of $\Gamma_\varphi$. This description readily shows that
$|\Sigma|$ is convex.

Conversely, if $\varphi: N_\RR\to \RR\cup\{\infty\}$ is a strictly convex
piecewise linear function on a fan $\Sigma$ on $N_\RR$ with convex
support, then its \emph{Newton polyhedron}
\[
\Xi=\{x\in M_\RR\,|\, \varphi+x\ge 0\}
\]
is an integral polyhedron. It is unbounded if and only if $\Sigma$ is
not complete. An alternative description is
\[
\Xi= C_\varphi^\vee\cap\big(M_\RR\times\{1\}\big),
\]
for $C_\varphi$ the convex hull of $\Gamma_\varphi$.

These two constructions set up a one-to-one correspondence between
integral, strictly convex piecewise linear functions on fans in
$N_\RR$ with convex support on one side and integral polyhedra in
$M_\RR$ on the other side.
\medskip

We are now ready to explain the first method of constructing an
integral tropical manifold out of a toric degeneration of CY-pairs.

\begin{example}\label{fan picture}(\emph{The fan picture.}
\cite{logmirror},~\S4.1)
If $(\pi:\X\to T,\D)$ is a toric degeneration of CY-pairs we can
define an integral tropical manifold as follows. For simplicity we
assume that the irreducible components of $X=\pi^{-1}(O)$ do not
self-intersect, that is, are normal. Let $\Strata(X)$ be the finite
category consisting of toric strata of $X$, with inclusions defining
the morphisms. For $S\in \Strata(X)$ let $\eta_S\in\X$ be the generic
point and let $Y_1,\ldots, Y_r$ be the irreducible components of
$X\cup \D$ containing $S$. Choose the order in such a way that
$Y_i\subseteq X$ iff $i\le s$. Define the monoid
\[
P_S:=\big\{(m_1,\ldots,m_r)\in\NN^r \,\big|\, {\textstyle\sum} m_i [Y_i]
\text{ is a Cartier divisor at }\eta_S\in\X\}.
\]
If $\Phi: U\to\Spec\kk[P]$ is as in Definition~\ref{log smooth
maps} with $\eta_S$ lifting to $\tilde\eta_S\in U$ then $P_S$ is
isomorphic to the monoid localization of $P$ by $A:=\{m\in P\,|\,
\Phi^*(z^m)\in \O^\times_{U,\tilde\eta_S} \}$, that is, the quotient of
the monoid $P-A\subseteq P^\gp$ by its invertible elements.
This shows that $P_S$ is a toric monoid. Hence $P_S$ is the set of
integral points of a rational polyhedral cone in $\RR^r$. Define the
polyhedron $F(S)\subseteq (\RR^r)^*$ by intersecting the dual cone
with the hyperplane $\langle \,.\, ,\rho_S\rangle=1$, where
$\rho_S=(1,\ldots,1,0,\ldots,0)\in P_S$ is the vector with entry $1$
at the first $s$ places:
\[
F(S):=\big\{ \lambda\in\Hom(P_S,\RR_{\ge0})\big|\,
\lambda(\rho_S)=1\big\}.
\]
In this definition $\RR_{\ge0}$ is viewed as additive monoid. The fact
that the central fibre $X$ of the degeneration is reduced  says that
the integral distance of each facet of $P_S$ to $\rho_S$ equals $0$ or
$1$. This implies that the vertices of $F(S)$ are integral. Note also
that $F(S)$ is unbounded iff $\rho_S$ lies in the boundary of the cone
generated by $P_S$, which is the case iff $S\subseteq \D$. Moreover,
if $S_1\subseteq S_2$, generization maps $P_{S_1}$ surjectively to
$P_{S_2}$, and this induces an inclusion of $F(S_2)$ as a face of
$F(S_1)$. Thus
\[
F: \Strata(X)^\text{op}\lra \LPoly,\quad S\longmapsto F(S)
\]
defines an integral polyhedral complex, the \emph{dual intersection
complex} $\check\P$ of $X$ (denoted $\P$ in \cite{logmirror}). Finally, the
toric irreducible components define compatible fan structures at the
vertices, which by assumption on $X$ have convex support. This makes
$\check B=\coprod_{\sigma\in\check\P} F(\sigma)/\sim$ into an integral
tropical manifold.

Because the irreducible components correspond to the fans at the
vertices we refer to this relation of integral tropical manifolds with
toric degenerations as the \emph{fan picture}. In this picture the
maximal cells specify local models for $\X$ at the zero-dimensional toric
strata of $X$.

Note that this construction depends only on $X$ and on the completion
of $\O_\X$ along $X\cup \D$. Hence it works also for formal toric
degenerations.
\qed
\end{example}

\begin{example}
\label{cubicfanexample}
\cite{logmirror}\cite{GBB} give many examples of affine manifolds
obtained from toric degenerations of proper Calabi-Yau varieties. Here
we give an example in the Fano case. Consider the equation
$f_3(u_0,u_1,u_2,u_3) +tu_0u_1u_2=0$ defining
$\foX\subseteq\PP^3\times \Spec \kk\lfor t\rfor$, with $f_3$ a general
choice of cubic form. Then $\foX\rightarrow \Spec\kk\lfor t\rfor$ is a
toric degeneration. The corresponding dual intersection complex
$(B,\P)$ looks like  Figure~\ref{cubicfan}. This picture is slightly
misleading. The three unbounded rays are in fact parallel. The bounded
two-cell is just a standard simplex. The discriminant locus $\Delta$
consists of the three points marked with crosses.
\begin{figure}
\includegraphics{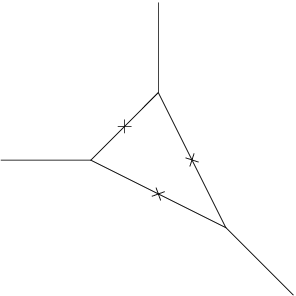}
\caption{}
\label{cubicfan}
\end{figure}
\end{example}
 
The other method for producing an integral tropical manifold, which is
even more relevant to this paper, requires a polarized central fibre.

\begin{example}\label{cone picture}(\emph{The cone picture.}
\cite{logmirror}, \S4.2)
Let $(\pi:\X\to T,\D)$ be a toric degeneration of CY-pairs, and let
$\shL$ be an ample line bundle on $X=\pi^{-1}(O)$. Again we make the
simplifying assumption that no components self-intersect. Then any
$S\in\Strata(X)$ together with the restriction of $\shL$ is a (not
necessarily complete) polarized toric variety. The sections of powers
of the line bundle define the integral points of the cone
$C(\sigma)\subseteq \RR^{n+1}$ over a (not necessarily bounded)
integral polyhedron $\sigma\subseteq\RR^n\times\{1\}$. Now $\sigma$ is
determined uniquely up to integral affine transformations, and
inclusions of toric strata define integral affine inclusions of
polyhedra as faces. Hence we have an integral polyhedral complex
\[
F:\Strata(X)\lra \LPoly,
\]
the \emph{intersection complex} $\P$ of $(X,\shL)$. Note that
the boundary of this polyhedral complex is covered by the cells
corresponding to the toric strata of $\D\cap X$.

The fan structures at the vertices this time come from log smoothness
as follows. At a zero-dimensional toric stratum $\{x\}\subseteq X$ the
degeneration is \'etale locally described by a toric morphism $\Spec
\kk[C\cap \ZZ^{n+1}]\to \Spec \kk[\NN]$ for some rational polyhedral
cone $C\subseteq \RR^{n+1}$. Denote by $\rho_S\in C\cap\ZZ^{n+1}$ the
image of $1\in\NN$. Then the images of the faces of $C$ not containing
$\rho_S$ under the projection $\RR^{n+1}\to \RR^{n+1}/\RR\rho_S\simeq
\RR^n$ define an $n$-dimensional fan. Its support is convex because it
is the image of a convex cone by a linear map. The cones in this fan
are equal to tangent wedges of $F(\{x\})\subseteq F(S)$, for
$S\in\Strata(X)$ containing $x$. This defines the fan structure at
$F(\{x\})$. This construction again works also for formal toric
degenerations of CY-pairs.

In this construction an irreducible component of $X$ is defined by the
cone over a maximal cell $\sigma\subseteq \RR^n$ of $\P$ by
$\Proj\big(\kk[C(\sigma)\cap \ZZ^{n+1}]\big)$. This is why we call this
correspondence the \emph{cone picture}. In contrast to the fan
picture, $(B,\P)$ now carries information about the polarization, but
we have lost some information about the local embedding into $\X$ by
projecting $C$ down to $\RR^{n+1}/\RR\rho_S$. See
Remark~\ref{polarization on (B,P)} for how to keep this information.

Note that in the construction we used a little less than an ample
line bundle on $X$. It suffices to have an ample line bundle on each
irreducible component with isomorphic restrictions on common toric
prime divisors. We call such data a \emph{pre-polarization} of $X$.
If $H^2(B,\ZZ)\neq 0$ a pre-polarization might not arise from a
polarization, cf.\ \cite{logmirror}, Theorem~2.34.
\qed
\end{example}

\begin{example}
\label{cubicconeexample}
Returning to the degeneration of a cubic in Example \ref{cubicfanexample},
polarizing the degeneration with the restriction of $\O_{\PP^3}(1)$,
one obtains the intersection complex $(B,\P)$ which looks like
\begin{center}
\includegraphics{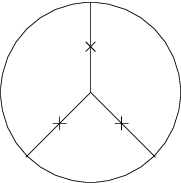}
\end{center}
Again, this figure is misleading: there are three standard simplices,
and the boundary is in fact a straight line with respect to the
affine structure.
\end{example}

Let $(B,\P)$ be an integral tropical manifold. Define an
\emph{(integral) affine function} on an open set $U\subseteq B$ to be
a continuous map $U\to \RR$ that is (integral) affine on $U\setminus
\Delta$. Similarly, an (integral) \emph{PL-function} (``piecewise
linear'') on $U$ is a continuous map $\varphi:U\to \RR$ such that if
$S_\tau: U_\tau\to\RR^k$ is the fan structure along $\tau\in\P$, then
$\varphi|_{U\cap U_{\tau}}= \lambda+S_\tau^*(\bar\varphi)$ for an
(integral) affine function $\lambda:U_\tau\to \RR$ and a function
$\bar\varphi:\RR^k\to \RR$ that is piecewise (integral) linear with
respect to the fan $\Sigma_\tau$ \cite{logmirror}, Definition~1.43.
The integral affine functions and integral PL-functions define sheaves
$\shAff(B,\ZZ)$ and $\shPL_\P(B,\ZZ)$ on $B$.

\begin{remark}\label{polarization on (B,P)}
In the cone picture the boundary of the cone $C$ with $\Spec
\kk[C\cap\ZZ^{n+1}] \to\Spec \kk[\NN]$ describing $\pi:\X\to T$
locally, can be viewed as the graph of an integral PL-function on
$\Sigma_v$, well-defined up to integral affine functions. These glue
to a section $\varphi$ of the sheaf $\shPL_\P(B,\ZZ)/ \shAff(B,\ZZ)$
of multi-valued, integral PL-functions. This is additional information
about \'etale models of the degeneration near the zero-dimensional
strata, or in other words, about the log smooth morphism to $(T,O)$.

Note that a local representative of $\varphi$ as a PL-function is
\emph{strictly convex}, so induces a strictly convex, integral affine
function on each fan $\Sigma_\tau$, $\tau\in \P$. We call such a
$\varphi$ a \emph{polarization} of $(B,\P)$, and $(B,\P,\varphi)$ a
\emph{polarized, integral tropical manifold}.
\qed
\end{remark}

\begin{construction}\label{discrete Legendre transform} 
(\emph{The discrete Legendre transform.} \cite{logmirror}, \S1.4)
There is a duality transformation on the set of polarized, integral
tropical manifolds, the \emph{discrete Legendre transformation}
$(B,\P,\varphi) \mapsto (\check B,\check \P,\check \varphi)$, which is
at the heart of our mirror symmetry construction. It works by defining
$\check \P$ as the opposite category of $\P$. Thus any $\tau\in\P$ is
also an object of $\check \P$, denoted $\check\tau$ for clarity. Then
$\check F: \check \P\to \LPoly$ maps $\check\tau$ to the Newton
polyhedron of $\bar\varphi$, where $\varphi=\lambda+
S_\tau^*(\bar\varphi)$ is as in the definition of PL-functions above.
If $\sigma\in \P$ is a maximal cell then $\check\sigma\in\check\P$ is
a vertex. In this case, the boundary of the dual of the cone over
$\sigma\times\{1\}$ defines the graph of $\check\varphi$ at
$\check\sigma$, hence also the fan structure at $\check\sigma$.
Applying this transformation again retrieves the original polarized
tropical manifold (\cite{logmirror}, Proposition~1.51). Moreover, the
discrete Legendre transformation preserves positivity
(\cite{logmirror}, Proposition~1.55). With the above correspondence
between strictly convex piecewise linear functions on non-complete
fans and unbounded cells the proofs of these facts for closed $B$ in
\cite{logmirror}, \S1.4, extend in a straighforward manner to the
general case.

One important remark is that there is a homeomorphism of
$B\setminus\partial B$ and $\check B\setminus \partial\check B$
mapping the discrimant loci onto each other. If $B$ is closed the
homeomorphism can be constructed by piecewise affine identifications
of the (perturbed) barycentric subdivisions. This is really just a
homeomorphism and certainly does not in general preserve the affine
structures. The general case is a little more subtle, and we leave the
details to the interested reader because we do not need this result
here. What we do use, however, is the consequence that we can view a
local system on $B\setminus \Delta$ as a local system on $\check
B\setminus \check\Delta$ and vice versa. In fact, a local system
$\shF$ on $B\setminus \check\Delta$ is nothing but a collection of
groups $\shF_v$, $\shF_\sigma$, one for each vertex $v$ and maximal cell
$\sigma$ of $\P$, together with a generization isomorphism
\[
\psi_{\sigma v}: \shF_v\lra \shF_\sigma,
\]
whenever $v\in\sigma$. The corresponding local system on $\check B$
then has the generization isomorphisms $\psi_{\check v
\check\sigma}^{-1}$.

By construction it should also be clear that the discrete Legendre
transform of the cone picture $(B,\P,\varphi)$ of a toric degeneration
of CY-pairs with pre-polarized central fibre $X$ leads to the fan
picture $(\check B,\check \P)$ of the same degeneration. The
polarization $\check\varphi$ on $(\check B,\check \P)$ thus obtained
is directly related to the polarization of the irreducible components
of $X$ via the usual description of ample line bundles on a toric
variety by strictly convex, integral PL-functions on the associated
fan, well-defined up to linear functions. This interpretation of a
strictly convex, multi-valued, integral PL-function is what motivates
us to call it a polarization of $(B,\P)$.
\qed
\end{construction}

%===========================================================
\subsection{Algebraic data}\label{subsection algebraic data}
Our aim in this paper is to construct a toric degeneration of CY-pairs
starting from a polarized, integral tropical manifold
$(B,\P,\varphi)$, using the cone picture. This process requires
additional, generally non-discrete input that we now describe.

Let $v\in\P$ be a vertex. Choose a PL-function $\varphi_v$ near
$v$ with $\varphi_v(v)=0$ representing the polarization $\varphi$.
The convex hull of the graph of $\varphi_v$ defines a strictly
convex, rational polyhedral cone $C_v\subseteq T_{B,v}\oplus \RR$, and
the associated toric monoid $P_v=C_v\cap(\Lambda_v\oplus\ZZ)$.
Then according to Example~\ref{cone picture} and
Remark~\ref{polarization on (B,P)}
\[
\kk[t]\lra \kk[P_v],\quad t\longmapsto z^{(0,1)}
\]
is an \'etale local model for any pre-polarized toric degeneration of
CY-pairs with cone picture $(B,\P,\varphi)$ near the zero-dimensional
toric stratum $\{x\}\subseteq \X$ corresponding to $v$. By integrality
of $\varphi$ it follows that the central fibre of this local model is
the union of affine toric varieties
\[
\bigcup_{K} \Spec \kk[K\cap(\Lambda_v\oplus \ZZ)],
\]
where the union runs over all facets $K\subseteq C_v$ not containing
$(0,1)$. These are indexed by maximal cells $\sigma$ containing $v$.
Now the projection $C_v\to T_{B,v}$ defines integral affine
isomorphisms of these facets of $C_v$ with the maximal cones in the
fan $\Sigma_v$. Hence this union depends only on $\Sigma_v$ and we use
the notation $\Spec \kk[\Sigma_v]$ for it. The justification for this
is that if $\Sigma$ is a fan on $M_\RR$ of convex, but not necessarily
strictly convex cones, with $|\Sigma|$ convex, then
\[
m+m':= \left\{\begin{array}{ll}
m+m'&,\, \exists\, K\in \Sigma:\, m,m'\in K\\
\infty&,\,\text{otherwise}.
\end{array} \right.
\]
defines a monoid structure on $(M\cap |\Sigma|)\cup\{\infty\}$. By
formally putting $z^\infty=0$ this yields a monoid $\kk$-algebra
generated by monomials $z^m$, $m\in M\cap |\Sigma|$, that we denote
$\kk[\Sigma]$. In the case of $\Sigma=\Sigma_v$ it is clearly
isomorphic to $\kk[P_v]/(z^{(0,1)})$, and hence $\kk[\Sigma_v]$ is
indeed the coordinate ring of our local model of $X$.

Since the central fibre $X$ of any toric degeneration is a union of
toric varieties, open subspaces isomorphic to $\Spec \kk[\Sigma_v]$
for $v\in \P^{[0]}$ cover $X$, and their mutual intersections have
similar descriptions as unions of affine toric varieties. We thus
arrive at the following gluing construction of $X$.

\begin{construction}\label{open gluing construction} (\cite{logmirror},
\S2.2)
First define certain open subsets of $\Spec \kk[\Sigma_v]$ as
follows. If $\tau\in\P$ and $v\in \tau$ is a vertex,
define the fan of convex, but not necessarily strictly convex, cones
\[
\tau^{-1}\Sigma_v:= \{K_e+\Lambda_{\tau,\RR}\,|\,
K_e\in\Sigma_v,\,e: v\to\sigma \text{ factors through }\tau \}.
\]
Recall $\Lambda_\tau=\Lambda_v\cap T_{\tau,v}$. The quotient of
$\tau^{-1}\Sigma_v$ by the linear space spanned by $\tau$ equals the
fan $\Sigma_\tau$ defining the fan structure along $\tau$. The fan
$\tau^{-1}\Sigma_v$ depends only on $\tau$: for a different choice
of $v'\in\tau$, we can identify $\tau^{-1}\Sigma_v$ and
$\tau^{-1}\Sigma_{v'}$ canonically. This is done via a piecewise
linear identification of $\Lambda_v$ and $\Lambda_{v'}$ which
identifies the cones $K_e+\Lambda_{\tau,\RR}$ and
$K_{e'}+\Lambda_{\tau,\RR}$, for $e:v\rightarrow\sigma$,
$e':v'\rightarrow\sigma$, via parallel  transport between $v$ and
$v'$ through $\sigma\in\P_{\max}$. While the induced bijection
$\Lambda_v\to\Lambda_{v'}$ is not in general linear due to the
effect of monodromy, the scheme
\[
V(\tau):=\Spec \kk[\tau^{-1}\Sigma_v]
\]
is well-defined up to unique isomorphism, independently of the choice
of vertex $v\in\tau$, see \cite{logmirror}, Construction~2.15. Note
that in \cite{logmirror}, which mostly uses the fan picture, this space
was denoted $V(\check\tau)$. The toric strata of $V(\tau)$ are in
bijection with the cones in $\tau^{-1}\Sigma_v$, hence to morphisms
$e:\tau\to\sigma$. The notation is $V_e\subseteq V(\tau)$.

Moreover, if $\omega\subseteq\tau$ there is a well-defined map of fans
\[
\big\{K_e+\Lambda_{\omega,\RR}\in\omega^{-1}\Sigma_v\,\big|\,
e:v\to\sigma \text{ factors through } \omega\to\tau\big\} 
\lra \tau^{-1} \Sigma_v,
\]
which defines an open embedding
\[
V(\tau)=\Spec \kk[\tau^{-1}\Sigma_v]\lra 
\Spec \kk[\omega^{-1}\Sigma_v]=V(\omega).
\]
We can compose this embedding with any toric automorphism of $V(\tau)$.
These are in bijection with maps
\begin{eqnarray}\label{action of PM}
\mu:\Lambda_v \cap|\tau^{-1}\Sigma_v|\lra \kk^\times
\end{eqnarray}
which are \emph{piecewise multiplicative} with respect to
$\tau^{-1}\Sigma_v$, meaning that the restriction to any cone in
$\tau^{-1}\Sigma_v$ is a homomorphism of monoids. In other words, for
each $\sigma\in\P_\max$ with $\tau\subseteq\sigma$ we have a
homomorphism $\mu_\sigma: \Lambda_\sigma \to\kk^\times$ such that for
any $\sigma,\sigma'$ containing $\tau$ the restrictions $\mu_\sigma
|_{\Lambda_{\sigma\cap\sigma'}}$ and $\mu_{\sigma'}
|_{\Lambda_{\sigma\cap\sigma'}}$ coincide. On the irreducible
component $V_e\subseteq V(\tau)$, $e: \tau\to\sigma$, the toric
automorphism of $V(\tau)$ associated to $\mu$ is given by
\[
z^m\longmapsto \mu_\sigma(m)\cdot z^m.
\]
This description also shows that there is a space $\PM(\tau)$ of
piecewise multiplicative functions along $\tau$ that depends only on the
embedding of $\tau$ in $B$ ($\check\PM(\check\tau)$ in the notation of
\cite{logmirror}). A choice of vertex $v\in\tau$ gives a
representation of $\PM(\tau)$ by maps $\Lambda_v
\cap|\tau^{-1}\Sigma_v|\to\kk^\times$ that are piecewise multiplicative with
respect to $\tau^{-1}\Sigma_v$.

To explain how piecewise multiplicative functions can be used to
change the gluing of our affine pieces $V(v)$ we translate
the definition of ``open gluing data  for the fan picture''
(\cite{logmirror}, Definition~2.25) to the cone picture.

\begin{definition}\label{open gluing data}
\emph{Open gluing data} for $(B,\P)$ are data $s=(s_{ e})_{e\in\Hom
\P}$ with the following properties: (1) $s_{ e}\in \PM(\tau)$ for
$e:\omega\to\tau$ (2) $s_{\id_\tau}=1$ for every $\tau\in\P$ (3) if
$e\in \Hom(\tau,\tau')$, $f\in\Hom(\tau',\tau'')$ then $s_{ f\circ e}=
s_{ f} \cdot s_{ e}$ wherever defined:
\[
s_{ f\circ e,\sigma}\ =\
s_{ f,\sigma}\cdot s_{ e,\sigma}\qquad
\text{for all $\sigma\in\P_\max$ with $\sigma\supseteq \tau''$}.
\]

Two open gluing data $s$, $s'$ are \emph{cohomologous} if there exist
$t_\tau\in \PM(\tau)$, $\tau\in\P$ with $s'_{ e}=t_\tau
t_\omega^{-1}\cdot s_{ e}$ for any $e:\omega\to\tau$.
\end{definition}

``Open'' refers to the fact that we glue the open sets $V(
v):=\Spec\kk[\Sigma_v]$ rather than their respective irreducible
components. If $s$ are open gluing data then $s_{ e}$ for
$e:\omega\to\tau$ defines an automorphism of $V(\tau)=\Spec
\kk[\tau^{-1}\Sigma_v]$ that we denote by the same symbol $s_{ e}$.
Thus for any $e:\omega\to\tau$ we obtain an open embedding by
composing $V(\tau)\to V(\omega)$ with $s_{ e}^{-1}$. (For consistency
we need to work with $s_{ e}^{-1}$ instead of $s_{ e}$, see
\cite{logmirror}, proof of Lemma~2.29, for how this arises.) This
yields a category of affine schemes and open embeddings, and saying
that the open sets glue means that there is a scheme $X=X_0(
\check B,\check \P,s)$, together with an open morphism
\[
p: \coprod_{\omega\in\P} V(\omega)\lra X_0(\check B,\check \P,s)
\]
that is a colimit for this category. The existence  of $X_0(\check B,
\check \P, s)$ is shown in \cite{logmirror}, \S2.2. Moreover, two open
gluing data give rise to isomorphic schemes iff they are cohomologous
(\cite{logmirror}, Proposition~2.32).

Conversely, according to \cite{logmirror}, Theorem~4.14, any central
fibre of a toric degeneration of CY-varieties with a pre-polarization
arises in this way. All these results extend in a straightforward
manner to CY-pairs. This ends Construction~\ref{open gluing
construction}.
\end{construction}

\begin{remark}\label{gluing data cone picture remark}
Our notation here differs somewhat from \cite{logmirror} since $s$ are
in fact open gluing data for the fan picture expressed in the cone
picture, rather than gluing data for the cone picture as in
\cite{logmirror}, Definition~2.3. In particular, according to
\cite{logmirror}, Theorem~2.34, there is an obstruction to gluing the
given ample line bundles on the irreducible components of $
X_0(\check B,\check \P,s)$, which in turn needs not be projective.
\end{remark}

The effect of monodromy on open gluing data $s=(s_e)$
(Definition~\ref{open gluing data}) gives rise to a set of
elements in $\kk^\times$, that we now introduce for later use.

\begin{definition}\label{D(lambda,v,rho)}
(Cf.\ \cite{logmirror}, Definition~3.25)
Given $\mu\in\PM(\tau)$ for some $\tau\in\P$, then for any
$\rho\in\P^{[n-1]}$ containing $\tau$ and any vertex $v\in\tau$ we can
measure the \emph{change of $\mu$ along $\rho$ with respect to $v$} as
follows. Let $\sigma,\sigma'$ be the unique maximal cells with
$\rho=\sigma\cap\sigma'$. Let $m\in\Lambda_\sigma$ map to a generator
of $\Lambda_\sigma/\Lambda_\rho\simeq\ZZ$, pointing from $\sigma$
to $\sigma'$. Let $m'\in\Lambda_{\sigma'}$ be obtained by parallel
transport of $m$ through $v$. Then
\begin{eqnarray}\label{D(lambda,rho,v)}
D(\mu,\rho,v):= 
\frac{\mu_{\sigma}(m)}{\mu_{\sigma'}(m')}\in\kk^\times
\end{eqnarray}
does not depend on the choice of $m$, and is also invariant under
changing $\mu$ by a homomorphism $\Lambda_v\to\kk^\times$.
\end{definition}

\begin{remark}\label{D(lambda,rho,v) monodromy}
Formula~\ref{$m^rho$} readily computes the dependence of
$D(\mu,\rho,v)$ on $v$:
\begin{eqnarray}\label{D(lambda,rho,v) v-dependence}
D(\mu,\rho,v')=\mu(m^\rho_{v 'v})^{-1}\cdot D(\mu,\rho,v),
\end{eqnarray}
see \cite{logmirror}, Remark~3.26 for details.
\end{remark}
\medskip

Apart from open gluing data, which specify the central fibre as a
scheme, we need some weak algebraic information about the embedding
into $\X$. This is what the log structure does, and it is more than
just the discrete information retained by the polarization $\varphi$
on $(B,\P)$ in the cone picture. For the purposes of this paper it
seems appropriate to explain this structure in an explicit,
non-abstract form, following \cite{logmirror}, p.~263f. Let $X=
X_0(\check B,\check \P,s)$ be a scheme obtained from open gluing data
in the cone picture as just described. As seen in
Construction~\ref{open gluing construction}, $X$ has a covering by
open sets isomorphic to $V(v)=\Spec \kk[\Sigma_v]$, and $V(v)$ can be
viewed as the toric Cartier divisor $z^{(0,1)}=0$ in the affine toric
variety $\Spec \kk[P_v]$. For $m\in\Lambda_v\cap |\Sigma_v|$, the
closure of the complement of the zero locus of $z^m$,
\[
V_m(v)=\cl\{x\in V( v)\ |\ z^m\in\O^\times_{X,x}\},
\]
is a union of irreducible components. Now a \emph{chart (for a log
smooth structure on $X=X_0(\check B,\check \P,s)$ of type given by
$\varphi$)} is an open set $U\subseteq V(v)$ for some vertex $v$,
together with $h_m\in \Gamma(U\cap V_m(v),\O_{V(v)}^\times)$ for $m\in
\Lambda_v \cap|\Sigma_v|$ that behaves piecewise multiplicatively with
respect to $\Sigma_v$ in the following sense:
\begin{eqnarray}\label{h_m-multiplicativity}
m,m'\in \Lambda_v\cap|\Sigma_v|\quad \Longrightarrow\quad
h_m\cdot h_{m'}= h_{m+m'}\text{ on }V_{m+m'}(v).
\end{eqnarray}
The vertex $v$ is part of the data defining a chart, and
$p(U)\subseteq X$ is the \emph{support} of the chart. Note that if
there is no cone in $\Sigma_v$ containing $m,m'$ then
$V_{m+m'}(v)=\emptyset$ and this condition is empty. Two charts
$(h_m)$, $(h'_m)$ defined on the same open subset $U\subseteq V(v)$
are \emph{equivalent} if there exists a homomorphism
$\lambda:\Lambda_v\to \Gamma(U,\O^\times_X)$ with
\[
h'_m=\lambda(m)\cdot h_m
\]
for any $m\in \Lambda_v\cap|\Sigma_v|$.

In log geometry a chart (for a fine, saturated log structure on $X$)
is just a morphism from an (\'etale) open subset of $X$ to an affine
toric variety. This relates to our definition as follows: For any
$m\in\Lambda_v\cap|\Sigma_v|$ consider $h_m\cdot z^m$ as function on
$U$ by continuation by zero. Then
\[
z^m\longmapsto h_m\cdot z^m
\]
defines an \'etale morphism $U\to \Spec \kk[\Sigma_v]$; the
composition with the closed embedding into $\Spec \kk[P_v]$ is the
associated chart in log geometry. Together with the distingushed
monomial $z^{(0,1)}\in \Spec \kk[P_v]$ this is a chart for a log
smooth morphism to the standard log point. It provides our local model
for a toric degeneration $(\X\to T, \D)$ of CY-pairs with central fibre $X$
as in Definition~\ref{toric degeneration}.

\begin{remark}
There is an important implicit dependence of our log-geometric chart
on $\varphi$, which we suppress in our notion of charts. We can do
this because $\varphi$ fixes the \emph{type of log structure} as
explained in \cite{logmirror}, Definitions~3.15 and 3.16.
\end{remark}

A chart defined on $U\subseteq V(v)$ can be restricted to $U'\subseteq
U$ simply by restricting the $h_m$. To compare arbitrary charts it
remains to explain how to change the reference vertex. So let
$(h_m)_{m\in\Lambda_v \cap|\Sigma_v|}$ be a chart defined on a
non-empty $U\subseteq V(v)$, and assume $v'\in \P$ is another vertex
with $p(U)\subseteq p\big(V(v')\big)$. Let $\sigma$ be a maximal cell
containing $v$ and $v'$. If no such cell exists, $p\big(V(v)\big)\cap
p\big(V(v')\big) =\emptyset$. Otherwise let $\Phi_{v'v}(s): U\to U'$
be the gluing isomorphism, that is, the composition of $p|_U$ with the
inverse of $p|_{V(v')}$. Write $e:v\to\sigma$, $e':v'\to\sigma$, and
denote by $K_e\in\Sigma_v$, $K_{e'}\in\Sigma_{v'}$ the tangent wedges
of $\sigma$ at $v$ and $v'$, respectively. Then for any $m\in
K_e\cap\Lambda_v$ it holds $V_e\subseteq V_m(v)$. Hence
by~(\ref{h_m-multiplicativity}) the map
\[
K_e\cap\Lambda_v\lra \Gamma(U\cap V_e,\O^\times_{V_e}),
\quad m\longmapsto h_m|_{V_e\cap U}
\]
is a homomorphism. Note the collection of these homomorphisms for all
$\sigma$ determine the chart. Now let $\tau\subseteq \sigma$ denote
the minimal cell containing $v$ and $v'$. Then parallel transport
through $\sigma$ gives the identification
\[
K_{e'} +\Lambda_{\tau,\RR}= K_e+\Lambda_{\tau,\RR}.
\]
Moreover, since the tangent wedge to $v$ in $\tau$ is contained in
$K_e$ the above homomorphism extends to $(K_e\cap\Lambda_v)
+\Lambda_\tau$. Let $h^\sigma_m$ denote the image of $m\in
(K_e\cap\Lambda_v) +\Lambda_\tau$ under this extension. We
are then able to define a chart on $U'\subseteq V(v')$ by pulling back
the $h_m^\sigma$ to $U'$:
\[
K_{e'}\cap\Lambda_v\lra \Gamma(U\cap V_{e'},\O^\times_{V_{e'}}),
\quad m\longmapsto
\big(\Phi_{v' v}(s)^{-1}\big)^* (h^\sigma_m).
\]
Finally define two charts, $(h_m)$ on $U\subseteq
V(v)$ and $(h'_m)$ on $U'\subseteq V(v')$, to be \emph{locally
equivalent} if any $x\in U$ has an open neighbourhood $W$ such that
$(h_m|_W)$ is equivalent to the pull-back of $(h'_m|_{W'})$ to
$V(v)$, where $W'=(p|_{V(v')})^{-1}(p(W))$.

\begin{definition}\label{log smooth structure}
An \emph{atlas} (for a log smooth structure on $X=X_0(\check
B,\check \P,s)$, of type defined by $\varphi$) is a system of
locally equivalent charts. A \emph{log smooth structure} on $X$ (of
type defined by $\varphi$) is a maximal atlas on the complement of a
closed subset $Z\subseteq X$ of codimension two which does not
contain any toric stratum.

A \emph{pre-polarized toric log CY-pair} is a scheme of the form
$X=X_0(\check B,\check \P,s)$ together with a polarization $\varphi$
of $(B,\P)$ and a log smooth structure.
\end{definition} 

A (formal) toric degeneration of CY-pairs $(\pi:\X\to T,\D)$ induces a
log smooth structure on the central fibre. Theorem~3.22 and
Definition~4.17 in \cite{logmirror} describe the set of log smooth
structures on $X$ potentially arising in this way, as a quasi-affine
subvariety of the space of sections of a coherent sheaf
$\shLS^+_{\pre,X}$ supported on $X_\sing$. On $V( v)$ this sheaf is
isomorphic to $\bigoplus_e \O_{V_e}$ where the sum runs over all
$e:v\to\rho$ with $\dim\rho=n-1$. A section $(f_e)_e$ of
$\shLS^+_{\pre,V( v)}= \bigoplus_e \O_{V_e}$ over $U\subseteq V(
v)\setminus Z$ defined by a log smooth structure obeys the following
compatibility condition along toric strata of codimension two.
Consider an $(n-2)$-cell $\tau\in\P$ with $v\in\tau$, let
$\rho_1,\ldots,\rho_l$ be a cyclic ordering of the $(n-1)$-cells
containing $\tau$ and write $h:v\to \tau$. Let $\check d_{
\rho_i}\in\Lambda_{\rho_i}^\perp\subseteq \Lambda_v^*$ be generators
compatible with this cyclic ordering. It turns out the following
condition is satisfied by tuples $(f_e)$ defining a log structure:
\begin{eqnarray}\label{multiplicative condition}
\prod_{i=1}^l \check d_{\rho_i}\otimes_\ZZ f_{e_i}|_{V_h} = 0\otimes1
\quad\text{in }\Lambda_v^*\otimes_\ZZ
\Gamma(U,\O_{V_h}^\times).
\end{eqnarray}
Note that the product treats the first factor additively and the
second factor multiplicatively. Conversely, any rational section of
$\shLS^+_{\pre, V( v)}$ with zeros and poles not containing any toric
stratum and fulfilling~(\ref{multiplicative condition}) defines a log
smooth structure on $V(v)$. This is proved by showing that giving
an equivalence
class of charts on an open subset $U\subseteq V(v)$ is equivalent to
giving sections $f_e\in\Gamma(U,\O^{\times}_{V_e})$ fulfilling
(\ref{multiplicative condition}).

Given a chart $(h_m)_{m\in\Lambda_v \cap|\Sigma_v|}$ on $U\subseteq
V(v)$ the associated section $f_e$ of $\O_{V_e}|_U$ for
$e:v\to\rho$ is defined as follows. Since $\rho$ is of
codimension one there exist two unique maximal cells $\sigma^+$,
$\sigma^-$ with $\rho=\sigma^+\cap\sigma^-$. Write $g^\pm:
v\to\sigma^\pm$. Working in an affine chart at $v$ let
$m^+\in\Lambda_v\cap K_{g^+}$ be a generator of $\Lambda_v/
\Lambda_\rho$. For appropriate $m_0\in K_e$ it holds $m_0-m^+\in
K_{g^-}$, and in any case $m_0+m_+\in K_{g^+}$. Then
\[
f_e=\frac{h_{m_0}^2|_{V_e}}{h_{m_0-m^+}|_{V_e}\cdot
h_{m_0+m^+}|_{V_e}}
\in \O^\times_{V_e}(U).
\]
This is independent of the choices of $m^+$ and $m_0$. The meaning of
$f_e$ is that at the generic point of $V_e$ a local model for
a toric degeneration with central fibre $X$ is given by 
\begin{eqnarray}\label{codim one log model}
 V(zw-f_e\cdot t^l)\subseteq \Spec \kk[z,w,t,x_1,\ldots,x_{n-1}].
\end{eqnarray}
Explicitly, $z$ and $w$ may be taken as the continuations by zero of
$z^{m^+}$ and $z^{-m^+}$, and $l\in\NN$ is the integral length of the
one-cell $\check\rho\in\check \P$. From this description it should be
plausible that the geometrically meaningful log smooth structures on
$V(v)$ are defined by sections $(f_e)_e$ of $\bigoplus_e
\O^\times_{V_e}$ over $V( v)\setminus Z$ that extend as sections of
$\shLS^+_{\pre,V( v)}= \bigoplus_e \O_{V_e}$. Such log smooth
structures are called \emph{positive}. (This corresponds to positivity
for integral tropical manifolds.) In fact, the log smooth structure
associated to a toric degeneration of log CY-pairs is positive
(\cite{logmirror}, Proposition~4.20).

The global structure of $\shLS^+_{\pre,X}$ follows from the formula
describing the change of charts. Not surprisingly this depends on the
choice of open gluing data $s$ describing the patching of the
open sets $V( v)$ to yield $X=X_0(\check B,\check \P,s)$. Let
$v,v'\in\P$ be vertices, $U\subseteq V( v)$ with $p(U)\subseteq p\big(
V(v')\big)$ and $(f_e)_{e:v\to \rho}\in\Gamma(U,\shLS^+_{\pre, V(
v)})$. Denote by $\Phi_{v'v}(s)$ the gluing isomorphism from
$U\subseteq V( v)$ to an open subscheme $U'\subseteq V( v')$. Then the
corresponding section of $\shLS^+_{\pre, V( v')}$ on $U'$ is
$(f_{e'})_{e':v'\to\rho}$ with
\begin{eqnarray}\label{shLS change of chart}
\Phi_{v'v}(s)^*(f_{e'})= \frac{D(s_{e'},\rho,v')}{D(s_{ e},\rho,v)}
s_{ e}(m^\rho_{v' v}) z^{m^\rho_{v'v}}
f_e,
\end{eqnarray}
see \cite{logmirror}, Theorem~3.27. If $\tau\in\P$ is a cell
containing $v,v'$ then this equation can be written more symmetrically
by viewing $f_{e'}$ and $f_e$ as functions on an open subset of 
$V(\tau)$ via the canonical open embeddings $V(\tau)\to V( v')$,
$V(\tau)\to V( v)$:
\begin{eqnarray}\label{shLS change of chart, symmetric form}
D(s_{ g'},\rho,v')^{-1} s_{ g'}^{-1}(f_{e'})
= z^{m^\rho_{v' v}} D(s_{ g},\rho,v)^{-1} s_{ g}^{-1}(f_e),
\end{eqnarray}
where $g:v\to\tau$, $g':v'\to\tau$. (Common factors arising from
$\tau\to\rho$ cancel.)

These formulae provide an explicit description of $\shLS^+_{\pre,X}$
as an abstract sheaf: For each $\rho\in\P^{[n-1]}$ fix a vertex
$v\in\rho$. Then the $-m^\rho_{v v'}\in\Lambda_\rho$ for $v'\in\rho$
define a PL-function on the normal fan of $\rho$ (\cite{logmirror},
Remark~1.56), hence an invertible sheaf $\shN_{\rho}$ on the
codimension one stratum $X_{\rho}$ corresponding to $\rho$. The
isomorphism class of $\shN_{ \rho}$ is well-defined since for a
different $v$ the PL-function only changes by a linear function. Now
(\ref{shLS change of chart}) shows
\[
\shLS^+_{\pre,X}\simeq \bigoplus_{\rho\in\P^{n-1}} \shN_{ \rho}.
\]
Summarizing, we now have a complete description of the space of
positive log smooth structures on $X=X_0(\check B,\check \P,s)$
as sections of $\shLS^+_{\pre,X}$ with zeros not containing toric strata
and fulfilling the compatibility condition \eqref{multiplicative condition}
in codimension two. Sections
of this sheaf are given explicitly via tuples $(f_e)_{e:v\to\rho}$ of
regular functions on the codimension one strata of each $V(
v)\subseteq X$, obeying the compatibility condition~(\ref{multiplicative
condition}) and the gluing conditions~(\ref{shLS change of chart}) or
(\ref{shLS change of chart, symmetric form}).

This description also defines uniquely, for each $\rho\in\P^{[n-1]}$,
a codimension two subscheme $Z_\rho \subseteq X$ with preimage in $V(v)$
the zero locus of $f_{v\to\rho}$ in $V_{v\to\rho}$. It is a Cartier
divisor in the $(n-1)$-dimensional toric stratum $X_\rho\subseteq X$
associated to $\rho$. The canonical minimal choice of $Z$ in
Definition~\ref{log smooth structure} is then
\begin{eqnarray}\label{Z}
Z:=\bigcup_{\rho\in\P^{[n-1]}} Z_\rho.
\end{eqnarray}

For given open gluing data the space of sections of $\shLS^+_{\pre,X}$
giving rise to a positive log smooth structure on $X$, can be empty
or complicated. One main result in \cite{logmirror} is, however, that
if $(B,\P)$ is locally sufficiently rigid as an affine manifold
(``\emph{simple}'') and positive then the space of isomorphism classes
of log CY-spaces with dual intersection complex $(B,\P)$ equals
$H^1(B,i_* \Lambda\otimes_\ZZ\kk^\times)$, where $i:B\setminus
\Delta\to B$ is the inclusion (\cite{logmirror}, Theorem~5.4). This
cohomology group is explicitly computable; it is a product of a finite
group and $(\kk^\times)^s$ with $s=\dim_\QQ H^1(B,i_*\Lambda
\otimes_\ZZ \QQ)$. Simplicity of $(B,\P)$ requires certain polytopes
associated to local monodromy to be elementary simplices
(\cite{logmirror}, Definition~1.60). It implies primitivity of local
monodromy in codimension two and can be explained by a complete list
of local models in dimensions up to $3$ (\cite{logmirror},
Example~1.62). In higher dimensions simplicity is a harder to
understand maximal degeneracy condition. In the case of complete
intersections in toric varieties it is related to Batyrev's MPCP
resolutions \cite{GBB}.
\medskip

One final point in this section concerns the interaction of open
gluing data with our chart description of log smooth structures. We
already noted that cohomologous open gluing data lead to isomorphic
$X$. But a change by a ``coboundary''  $(t_\sigma)_{\sigma\in\P}$ has
the effect of composing our charts for the log smooth structure with
the corresponding automorphisms of $V( \sigma)$. Because this changes
the identification of $V(\sigma)$ with an open subset of $X$ this
leads to a different section $(f_e)_{e:v\to\rho}$ of
$\shLS^+_{\pre,X}|_{V( v)}= \shLS^+_{\pre,V(v)}$. We can, however,
partly get rid of this non-uniqueness by requiring that for any
$e:v\to\rho$ and $x\in X$ the zero-dimensional stratum corresponding
to $v$, we have
\begin{eqnarray}\label{normalization condition}
f_e(x)=1.
\end{eqnarray}
In other words, the constant term of $f_e\in \kk[K_e\cap\Lambda_v]$
equals $1$. If this is the case, the log smooth structure and the
corresponding section of $\shLS^+_{\pre,X}$ are called
\emph{normalized} for the given open gluing data (\cite{logmirror},
Definition~4.23).

If a section is not normalized the constant terms of the $f_e$ define
the coboundary of a zero-cycle, whose application to the open gluing
data leads to a normalized section (\cite{logmirror}, p.290). Hence
there is always some open gluing data for which the log smooth
structure of a toric log CY-pair is normalized. The normalization
condition may thus be interpreted as ``gauge fixing'', a process
eliminating infinitesimal automorphisms. A generalization of this
condition will turn up in our deformation process to get rid of pure
$t$-terms (Definition~\ref{def:slab normalization}).

%===========================================================
\subsection{Statement of the Main Theorem}
The input to the Main Theorem is a pre-polarized, positive toric log
CY-pair with intersection complex $(B,\P,\varphi)$. The log smooth
structure needs to satisfy a certain local rigidity condition that we
now explain. If $(B,\P)$ is positive and simple then for any open
gluing data $s$ the set of log smooth structures on $X=X_0(\check B,
\check\P,s)$ is non-empty (\cite{logmirror}, Theorem~5.4) and any choice
fulfills the requested properties, see Remark~\ref{simple implies
locally rigid} below. This implies the Reconstruction Theorem for
integral tropical manifolds stated in the introduction.

The local rigidity condition involves the following notion of
primitivity of Minkowski sums of polyhedra.

\begin{definition}\label{def:minkowski transverse}
Let $\Xi_1,\ldots,\Xi_s\subseteq \RR^n$ be polyhedra, 
$\Xi=\sum_i\Xi_i$ their Minkowski sum and $\Xi_i^{[0]}$, $\Xi^{[0]}$
the respective sets of vertices. For $v\in\Xi^{[0]}$ denote by
$v(i)\in \Xi_i^{[0]}$ the unique solution to the equation
$v=\sum_i v(i)$. Consider the following linear map
\[
F:\prod_i \Map \big(\Xi_i^{[0]},\kk\big)\lra
\Map \big(\Xi^{[0]},\kk\big),\quad
(a_i)_i\longmapsto \Big( v\mapsto {\sum}_i a_i\big(v(i)\big) \Big).
\]
Then $\Xi_1,\ldots,\Xi_s$ are called 
\emph{Minkowski transverse} if $F(a_1,\ldots,a_s)=0$ only has the trivial
solutions
\[
a_i(v)=\alpha_i\in\kk,\quad{\sum}_i\alpha_i=0.
\]
\end{definition}

\begin{remark}\label{rem:Minkowski transversality}
From $F(a_1,\ldots,a_s)=0$ any edge of $\Xi$ leads to a linear equation
of the form
\[
a_1(v_1)+\ldots+a_s(v_s)= a_1(v'_1)+\ldots+a_s(v'_s)
\]
with $v_i,v'_i\in \Xi_i^{[0]}$, and the question is if these impose
enough conditions to imply that $a_i$ is constant. From this
description it is not hard to see that $\Xi_1,\ldots,\Xi_s$ are
Minkowski-transverse if they do not have parallel edges, for then only
one $v_i$ changes along any edge of $\Xi$, and there are enough edges
to compare the values of any two vertices of any $\Xi_i$.

On the other hand, this is not a necessary condition. For example, the
two polygons
\[
\Xi_1=\conv\big\{(0,0),(1,0),(0,1) \big\},\quad
\Xi_2=\conv\big\{(0,0),(1,0),(1,1) \big\},\quad
\]
are Minkowski transverse.
\end{remark}

Suppose now we are given $X=X_0(\check B,\check \P,s)$ along with a
positive log smooth structure. For the following definition recall
from~(\ref{Z}) the components $Z_\rho\subseteq X$,
$\rho\in\P^{[n-1]}$, of the singular locus of the log
structure. For $\tau\subseteq\rho$ and $X_\tau\subseteq X$ the
associated toric stratum, the Newton polytope of $Z_\rho\cap X_\tau$
is
\[
\Delta_\tau(\rho) := \conv\{ m^\rho_{vv'}\,|\, v'\in \tau\},
\]
where $v\in\rho$ is a fixed choice of vertex.
This is naturally a face of the monodromy polytope $\Delta(\rho)$
(\ref{monodromy polytope}) parallel to $\Lambda_\tau$, and is well-defined
up to translation.

\begin{definition}\label{def:locally rigid}
We call a positive, toric log Calabi-Yau space \emph{locally rigid} if
\begin{enumerate}
\item[(i)]
For each $\rho\in\P^{[n-1]}$ and $\tau\in\P^{[n-2]}$, $\tau\subseteq\rho$,
any integral point of $\Delta_\tau(\rho)$ is a vertex of $\Delta_{\tau}(\rho)$.
\item[(ii)]
If $X_\tau\subseteq X$ denotes the toric stratum defined by
$\tau\in\P^{[n-2]}$ then for any $\rho\in\P^{[n-1]}$ containing $\tau$
the intersection $Z_\rho\cap X_\tau$ is reduced and irreducible.
Moreover, no more than three yield the same subset of~$X_\tau$.
\item[(iii)]
For any $\tau\in\P^{[n-2]}$ let $\Xi_i$, $i=1,\ldots,s$, be an
enumeration of $\Delta_\tau(\rho)$ for $\rho\in\P^{[n-1]}$,
$\tau\subseteq\rho$, modulo translation. Then $\Xi_1,\ldots, \Xi_s$
are Minkowski transverse.
\end{enumerate}
\end{definition}

\begin{remark}\label{Xi_i, varphi_i}
Let $\tau\in\P^{[n-2]}$. By (ii) the polynomials
$f_{v\to\rho}|_{V_{v\to\tau}}$ defining $Z_\rho \cap X_\tau$ locally
are irreducible. The compatibility condition~(\ref{multiplicative
condition}) then implies that for any $\rho$ with
$Z_\rho\cap X_\tau\neq\emptyset$
\[
\big\{\check d_{\rho'}\,\big|\, Z_{\rho'}\cap X_\tau= Z_\rho\cap
X_\tau \big\} 
\]
are the edge vectors of a polygon with edges of unit integral length
in the two-dimensional affine space $\Lambda_\tau^\perp\subseteq
\Lambda_{v,\RR}^*$, well-defined up to translation. By the second
requirement in (ii) it has $2$ or $3$ edges. We thus obtain a set
$\{\check\Xi_i\}$ of line segments and triangles and a corresponding set of
convex PL-functions $\varphi_i$ on $\Sigma_\tau$.
\end{remark}

\begin{example}
To illustrate the concept of local rigidity we give three local
examples in dimension $4$. Figure~\ref{fig 1.1} shows the respective
fans $\Sigma_\tau$ together with the values of a PL-function
$\varphi_\tau$ pulling back to $\varphi$ at the generators of the rays
(in square brackets), and the functions $f_{v\to\rho}$.
\begin{figure}[h]
\input{locally_rigid1.pstex_t}
\quad
\input{locally_rigid2.pstex_t}
\quad
\input{locally_rigid3.pstex_t}
\caption{}\label{fig 1.1}
\end{figure}
The five rays in the right-most figure are generated by $(1,0)$,
$(2,5)$, $(-1,4)$, $(-3,-5)$ and $(1,-4)$. In either case $z= z^{(0,0,1,0)},w=
z^{(0,0,0,1)}$ denote monomials generating the coordinate ring of the
maximal torus of $X_\tau$, while $g_1,\ldots,g_4$ are arbitrary
functions vanishing on $X_\tau$. Note that all examples are
normalized and fulfill~(\ref{multiplicative condition}). In the left
figure all four functions $f_{v\to\rho}$ have the same restriction to
$X_\tau$, thus violating condition~(ii). On the other hand, the middle
and right examples are locally rigid. The polytopes according to
Remark~\ref{Xi_i, varphi_i} are $\check\Xi_1=[0,1]\times\{0\}$,
$\check\Xi_2=\{0\}\times [0,1]$ and $\check\Xi_1=
\operatorname{conv}\{0,(4,1)\}$, $\check\Xi_2=
\operatorname{conv}\{(0,0),(0,1),(-5,3)\}$, respectively.
\end{example}

\begin{remark}\label{simple implies locally rigid}
If $(B,\P)$ is simple (\cite{logmirror}, Definition~1.60), then
$X_0(\check B,\check \P,s)$ is locally rigid for any choice of open
gluing data $s$. In fact, (i) follows readily from the fact that in
this case $\Delta(\rho)$ is an elementary simplex for any
$\rho\in\P^{[n-1]}$. This also implies that $Z_\rho\cap X_\tau$ is
reduced and irreducible, as $\Delta_\tau(\rho)$ is the Newton polytope
of $Z_{\rho} \cap X_\tau$. As for the second condition in (ii) let
$v\in\tau$ be a vertex. Simplicity implies the existence of $p\le
\codim\tau=2$ polytopes $\check\Delta_i \subseteq
\Lambda_{\tau,\RR}^\perp \subseteq \Lambda_v^*$ with the following
properties: (1) Each ray of $\Sigma_\tau$ labelled by a codimension
one $\rho$ with $Z_\rho\cap X_\tau\neq\emptyset$ is generated by the
inward normal of some $\check\Delta_i$. (2) The convex hull of
$\bigcup_i \check\Delta_i\times\{e_i\}\subseteq
\Lambda_{\tau,\RR}^\perp \times\RR^p$ is an elementary simplex. By (2)
the tangent spaces of $\check\Delta_i$ are transverse. Thus either
$p=1$ and $\check\Delta_1$ is a triangle or a line segment, or $p=2$
and $\check\Delta_1$, $\check\Delta_2$ are two non-parallel line
segments (cf.\ \cite{logmirror}, p.217 for the latter case). This
implies~(ii). (iii) follows from Remark~\ref{rem:Minkowski
transversality} since $\Xi_i$ are elementary simplices with
$T_{\Xi_1}\oplus\ldots\oplus T_{\Xi_s}$ an internal direct sum.
\end{remark}

We are now in position to state the main result of this paper. The
notions of pre-polarized toric log CY-pair, formal toric degeneration
of CY-pairs, local rigidity and positivity have been introduced in
Definitions~\ref{log smooth structure}, \ref{formal toric
degeneration}, \ref{def:locally rigid} and \ref{def:positivity},
respectively.

\begin{theorem}\label{main thm}
Any locally rigid, positive, pre-polarized  toric log CY-pair with
proper irreducible components arises from a formal toric degeneration
of CY-pairs.
\end{theorem}

Note that the hypothesis of properness is equivalent to the
boundedness of every cell of $\P$ in the intersection complex
$(B,\P)$. We  are confident that this hypothesis is not necessary, but
the unbounded case seems to raise a number of interesting points
involving Landau-Ginzburg potentials, which are better dealt with
elsewhere. However, most of the arguments we give will work  in
general, and we will, in the course of the proof, remark when we are
using this boundedness hypothesis.

If we want actual families we need to restrict to the
projective or compact analytic setting.

\begin{corollary}\label{main thm, polarized version}
Any projective, locally rigid, positive toric log CY-pair $(X,D)$
with $H^1(X,\O_X)= H^2(X,\O_X)=0$  arises from a
projective toric degeneration of CY-pairs over $\kk\lfor t\rfor$.
\end{corollary}

\proof
The cohomological assumptions imply that an ample line bundle on $X$
extends to the formal degeneration. The result then follows from
Grothendieck's Existence Theorem in formal geometry, see
EGA~III, 5.4.5 \cite{EGAIII}.
\qed

\begin{remark}
The assumption on cohomology is indeed superfluous. In fact, one can
show that any ample line bundle $L$ on $X$ extends to the formal
degeneration by applying our construction to the total space of $L$.
Details of this observation will appear elsewhere.
\end{remark}

Note that if $H^2(B,\kk^\times)=0$ for $B$ the dual intersection
complex associated to the toric log-CY pair, then projectivity of
$X_0(\check B, \check \P,s)$ follows from the existence of the
pre-polarization (\cite{logmirror}, Theorem~2.34).

For the analytic formulation we just remark that all the notions we
have introduced so far have straightforward analogues in the
complex-analytic world. In view of the existence of versal
deformations of (pairs of) compact complex spaces
\cite{douady}\cite{grauert}, we obtain the following result.

\begin{corollary}\label{main thm, analytic version}
Any compact, locally rigid, positive toric log CY-pair arises from
a toric degeneration of analytic CY-pairs over $\kk\lfor t\rfor$.
\end{corollary}

%===========================================================
%===========================================================
\section{Main objects of the construction}
\label{section main objects}

The rest of the paper is devoted to the proof of Theorem~\ref{main
thm}. We thus fix, once and for all, a polarized, integral tropical manifold
$(B,\P,\varphi)$, open gluing data $s$ for $(B,\P)$, and a
positive log smooth structure on $X=X_0(\check B,\check\P,s)$
given by a compatible set of sections $(f_e)_e$ of $\shLS^+_{\pre,V(
v)}$ fulfilling the multiplicative condition~(\ref{shLS change of
chart, symmetric form}). We also assume that the discriminant locus
$\Delta=\Delta(\{a_\tau\})$ does not contain any rational points,
as discussed before Lemma~\ref{perturbation lemma}.

In this entire section, we do not assume the cells of $\P$ 
to be bounded.

%===========================================================
\subsection{Exponents, orders, rings}
\label{section rings}
We construct the deformation of $(X,D)$ order by order. In each step
the deformation is a colimit, in the category of separated schemes, of
a system of affine schemes. This system is obtained by chopping $B$
into polyhedral pieces, called \emph{chambers}, of growing number for
higher order; then $\P$ induces a stratification of each chamber, and
there will be one ring for each inclusion of such strata.
Homomorphisms are obtained either by changing strata within one chamber
or by passing from one chamber to a neighbouring one. We start by
explaining how to define the rings.

Recall that for a vertex $v\in B$ we have the local model $V(v)
\subseteq \Spec \kk[P_v]$ for $X\subseteq \X$. Here $P_v$ are the
integral points over the graph in $T_{B,v}\oplus\RR$ of a local
representative $\varphi_v$ of $\varphi$ with $\varphi_v(v)=0$. The
disadvantage of this description is that it depends on the choice of
representative of $\varphi$. To derive a more invariant point of view
recall that the Legendre dual to $v$ is a maximal cell $\check
v\subseteq \Lambda_{v,\RR}^*$ with vertices $\check\sigma=
-\lambda_\sigma \in\Lambda_v^*$, where $\lambda_\sigma$ are the linear
functions defined by $\varphi_v$ for the maximal cells $\sigma$
containing $v$. We may then view $m=(\overline m,h)\in \Lambda_v\oplus\ZZ$
as an affine function on $\check v$ via the sequence of
identifications
\begin{eqnarray*}
\Lambda_v\oplus\ZZ= (\Lambda_v^*)^*\oplus\ZZ
= \Gamma(\Int\check v, \shAff(\check B,\ZZ)).
\end{eqnarray*}
The value of this affine function, denoted by $m$ also, at $\check\sigma$ is
\begin{eqnarray}\label{m_sigma(check sigma)}
m(\check\sigma)= \langle \overline m, -\lambda_\sigma\rangle +h.
\end{eqnarray}
This gives the following description of $P_v$ in terms of affine
functions on $\check v$:

\begin{lemma}\label{P_v}
$P_v=\big\{m=(\overline m,h)\in\Lambda_v\oplus\ZZ \,\big|\,
m|_{\check v}\ge 0 \big\}$.
\end{lemma}

\proof
The condition $m|_{\check v}\ge 0$ is equivalent to requiring
that $m=(\overline m,h)$ lies in the dual of the cone generated by $\check
v\times\{1\} \subseteq \Lambda_{v,\RR}^*\oplus\RR$. By our definition of
Newton polyhedra this agrees with the convex hull of the graph of
$\varphi_v$, whose integral points are $P_v$.
\qed
\smallskip

The preceding discussion motivates the following definition.

\begin{definition}
An \emph{exponent} at a point $x\in B\setminus\Delta$ is an element of
the stalk of $\shAff(\check B,\ZZ)$ at $x$. An exponent on
$\sigma\in\P_\max$ is an exponent at any $x\in\Int\sigma$, that is, an
element of $\shAff(\check B,\ZZ)_{\check\sigma}$. An exponent $m$ on
$\sigma$ defines exponents at any $x\in\sigma\setminus\Delta$ that we
denote by the same symbol.

The image of an exponent $m$ at $x$ (on $\sigma$) under the projection
$\shAff(\check B,\ZZ)_x\to \Lambda_x$ ($\shAff
(\check B,\ZZ)_{\check \sigma} \to \Lambda_\sigma$) is denoted
$\overline m$.
\end{definition}

Here we view $\shAff(\check B,\ZZ)$ as a locally constant sheaf on $B$
as explained in Construction~\ref{discrete Legendre transform}. Note
there is an exact sequence
\[
0\lra \ZZ\lra 
\shAff(\check B,\ZZ)\stackrel{m\mapsto \overline m}{\lra}
i_*\Lambda\lra 0,
\]
see \cite{logmirror}, Definition~1.11. An exponent $m$ at $x\in
B\setminus \Delta$ extends uniquely to a section of $\shAff(\check
B,\ZZ)$ on the interior of each $\sigma\in\P_\max$ containing $x$.
This defines an element of $\shAff(\check B,\ZZ)_{\check\sigma}$ that
we denote $m_\sigma$. Note that if $x\in\sigma\cap\sigma'$ for another
maximal cell $\sigma'$ containing $x$, then parallel transport in
$\check v$ for a vertex $v$ in the same connected component of
$(\sigma\cap\sigma')\setminus\Delta$ as $x$, maps $m_\sigma\in
\shAff(\check B,\ZZ)_{\check\sigma}$ to $m_{\sigma'}\in \shAff(\check
B,\ZZ)_{\check{\sigma}'}$ 
\medskip

To define our rings we need various order functions on exponents. For
$m\in P_v$ and $\sigma\in\P_\max$, $v\in\sigma$, the monomial $z^m$
does not vanish on the irreducible component $V_{v\to\sigma}$ of
$\Spec\kk[\Sigma_v]\subseteq \Spec\kk[P_v]$ if and only if
$h=\lambda_\sigma (\overline m)$. Thus by~(\ref{m_sigma(check
sigma)}), $m_\sigma(\check \sigma)$ equals the order of vanishing of
$z^m$ along $V_{v\to\sigma}$. 

\begin{definition} 
\label{orderdef}
1)\ Let $m$ be an exponent at $x\in B\setminus \Delta$. Then the
\emph{order} of $m$ on $\sigma\in\P_\max$, $x\in\sigma$
is
\[
\ord_{\sigma}(m):=m_\sigma(\check\sigma).
\]
Denote by $\P_{\max}^{\partial}$ the set of codimension one cells of
$\P$ contained in $\partial B$. Then for $x\in\rho$, 
the \emph{order} of $m$ on $\rho\in\P_{\max}^{\partial}$ is
\[
\ord^{\partial}_{\rho}(m):=\langle \overline m,n_{\rho}\rangle,
\]
where $n_{\rho}\in\Lambda^*_x$ is an inward pointing primitive normal to $\rho$.

For $A\subseteq B$ a subset contained in a cell of $\P$ and with
$x\in A$ define
\[
\ord_A(m):=\max\big( 
\big\{ \ord_{\sigma}(m) \,\big|\,
\sigma\in\P_\max,\, A\subseteq\sigma
\big\}\cup
\big\{ \ord^{\partial}_{\rho}(m) \,\big|\,
\rho\in\P_\max^{\partial},\, A\subseteq\rho
\big\}\big).
\]
Note that for $A=\sigma$ a maximal cell this agrees with the
previous definition.\\[1ex]
2)\ Let $\omega\in\P$ be the minimal cell containing $x$. Define
\[
P_x:=\bigg\{m\in\shAff(\check B,\ZZ)_x \,\bigg|\,
\begin{array}{l}
\forall \sigma\in\P_\max, x\in\sigma:\ord_\sigma(m)\ge 0\\
\exists\sigma'\in\P_\max,\omega\subseteq\sigma': \overline m\in
K_\omega\sigma' \end{array}\bigg\}.
\]
\end{definition}

The notion of order is compatible with local monodromy:

\begin{lemma}\label{order well-defined}
Let $\sigma,\sigma'\in\P_\max$ and let $m$ be an exponent on $\sigma$.
If $m'$ is the result of parallel transport of $m$ along a closed loop
inside $(\sigma\cup\sigma')\setminus\Delta$ then
\[
\ord_\sigma(m)=\ord_\sigma(m').
\] 
If $\rho\in\P_{\max}^{\partial}$, $\rho$ is a face of $\sigma$
and $m'$ is the result of parallel transport of $m$ along a closed loop
inside $\Int(\sigma)\cup\Int(\sigma')\cup ((\rho\cap\sigma')\setminus\Delta)$,
then
\[
\ord_{\rho}^\partial(m)=\ord_{\rho}^{\partial}(m').
\]
\end{lemma}

\proof
On $(\sigma\cup\sigma')\setminus\Delta$ the locally constant sheaf
$\shAff(\check B,\ZZ)$ splits non-canonically as $\Lambda\oplus\ZZ$.
This follows from \cite{logmirror}, Proposition~1.12, in connection
with \cite{logmirror}, Proposition~1.29 applied to $\tau=
\sigma\cap\sigma'$. Moreover, by \cite{logmirror}, Proposition~1.29
again, the monodromy for paths in $(\sigma\cup\sigma')\setminus\Delta$
acts trivially on $\Lambda_\tau^\perp\subseteq \Lambda_x^*$, which in
an affine chart at $x\in\tau$ contains $\check\sigma$. Hence
$m_\sigma(\check\sigma)$ remains unchanged under monodromy.

Similarly, for the second statement, let $\tau=\rho\cap\sigma'$.
Then monodromy of loops in $\Int(\sigma)\cup\Int(\sigma')\cup
(\tau\setminus\Delta)$ preserves $\Lambda_{\tau}^{\perp}
\subseteq\Lambda_x^*$. But $\Lambda_{\tau}^{\perp}$ contains the
normal to $\rho$, hence the result.
\qed
\medskip

In view of the lemma it makes sense to define, for $m\in\shAff(\check
B,\ZZ)_{\check \sigma}$, the order on neighbouring maximal cells:

\begin{definition}\label{order adjacent cells}
Let $\sigma,\sigma'\in\P_\max$ with $\sigma\cap\sigma' \neq\emptyset$
and $m\in\shAff(\check B,\ZZ)_{\check\sigma}$ an exponent on $\sigma$.
Define the \emph{order of $m$ on $\sigma'$} as follows. Let
$m'\in\shAff(\check B,\ZZ)_{\check \sigma'}$ be the result of parallel
transport of $m$ inside $\check v$, for any vertex
$v\in\sigma\cap\sigma'$. Then
\[
\ord_{\sigma'}(m):= \ord_{\sigma'}(m').
\]
If in addition $\rho'\in\P_{\max}^{\partial}$ and $\rho'$ is a face of
$\sigma'$, let $m'\in\shAff(\check B,\ZZ)_{\check\sigma'}$ be
the result of parallel transport of $m$ inside $\check v$, for any vertex
$v\in\sigma\cap\rho'$. Then
\[
\ord_{\rho'}^{\partial}(m):=\ord_{\rho'}^{\partial}(m').
\]
\end{definition}
With these definitions it now also makes sense, for $m$ an exponent
on $\sigma\in\P_\max$ and $A\subseteq \sigma$ to define
$\ord_A(m)$ just as in Definition~\ref{orderdef} above for exponents
at a point.

Much of our strategy depends on the idea that if an exponent $m$
is propagated in the direction $-\overline m$, the order of $m$ increases.
This is analogous to the behaviour of the order function $\ord_l$ in
\cite{ks}, \S10.3.

\begin{proposition}\label{exponentprop}
Let $m$ be an exponent at $x\in B\setminus (\Delta\cup\partial B)$,
and let $\tau\in\P$ be the minimal cell containing $x$. If
$\sigma^+,\sigma^-\in\P$ are maximal cells containing $\tau$ such that
the corresponding maximal cones in $\Sigma_\tau$ contain $\overline m$
and $-\overline m$, respectively, then
\begin{eqnarray*}
\ord_{\sigma^-}(m)&=&\max \big\{\ord_\sigma(m)\,\big|\,
\sigma\in\P_\max, \tau\subseteq\sigma\big\},\\
\ord_{\sigma^+}(m)&=&\min \big\{\ord_\sigma(m)\,\big|\,
\sigma\in\P_\max, \tau\subseteq\sigma\big\}.
\end{eqnarray*}
\end{proposition}

\proof
$\Sigma_{\tau}$ is the normal fan of $\check\tau$; given
$e^{\pm}:\tau\to \sigma^{\pm}$, the cones $K_{e^{\pm}}$
of $\Sigma_{\tau}$ are the normal cones to the vertices $\check
\sigma^{\pm}$ of $\check\tau$ (see \cite{logmirror},
Definition~1.38). In particular, on $\check\tau$ an element of
$K_{e^{\pm}}$ achieves its minimal value at $\check\sigma^{\pm}$, from
which the result follows.
\qed
\medskip

Next we construct standard thickenings of the rings describing the
toric strata locally. These will be our basic building blocks. 

\begin{construction}\label{constr rings} (\emph{The rings.})
For $\omega\in\P$ and $\sigma\in\P_\max$ with $\omega\subseteq\sigma$
define the monoid
\[
P_{\omega,\sigma}:= \left\{m\in \shAff(\check B,\ZZ)_{\check \sigma}
\,\left|\,\begin{array}{l}
\forall\, \sigma'\in\P_\max,\, \omega\subseteq\sigma':
\ord_{\sigma'}(m)\ge 0\\
\exists\, \sigma'\in\P_\max,\, \omega\subseteq\sigma': \overline m\in  K_\omega\sigma'
\end{array}\right.\right\}.
\]
The condition $\overline m\in K_\omega\sigma'$ is only relevant if
$\omega\subseteq \partial B$. Note that if $v$ is a vertex then
by Lemma~\ref{P_v}, for any choice of representative $\varphi_v$ of
$\varphi$ at $v$ and any $\sigma\in \P_\max$, it holds
$P_{v,\sigma}=P_v$ canonically.

For any $\sigma'\in\P_\max$ containing $\omega$, parallel transport
through a vertex $v\in \sigma\cap\sigma'$ induces an isomorphism
$P_{\omega,\sigma}\simeq P_{\omega,\sigma'}$. This isomorphism,
however, generally depends on the choice of $v$. Thus $\sigma$ serves
as a reference cell.

Another manifestation of this phenomenon is as follows. If
$x\in\Int(\omega)\setminus\Delta$ there is a canonical isomorphism
$P_{\omega,\sigma}\simeq P_x$. Thus for any $x,x'\in
\Int(\omega)\setminus\Delta$ any choice of maximal cell $\sigma$
containing $\omega$ induces an isomorphism $P_x\simeq P_{x'}$, but this
isomorphism generally depends on the choice of $\sigma$.

If $g:\omega\to\tau\in\Hom(\P)$ and $\sigma\in\P_\max$ with
$\tau\subseteq\sigma$ then for each $k\in \NN$ we have a monoid ideal
\[
P^{>k}_{g,\sigma}:= \big\{m\in P_{\omega,\sigma}
\,\big|\, \ord_{\tau}(m)> k\big\}\subseteq P_{\omega,\sigma}.
\]
Let $I^{>k}_{g,\sigma}$ denote the ideal in $\kk[P_{\omega,\sigma}]$
generated by $P^{>k}_{g,\sigma}$ and define
\[
R_{g,\sigma}^k:=\big(\kk[P_{\omega,\sigma}]/I^{>k}_{g,\sigma}
\big)_{f_{g,\sigma}}.
\]
The function $f_{g,\sigma}$ at which we localize is constructed from
the given section $(f_e)_e$ of $\shLS^+_{\pre,X}$ as follows. Choose a
vertex $v\in\omega$ and write $e:v\to\omega$. Recall that $s_{
e}\in\PM(\omega)$ is a map $\Lambda_v \cap|\omega^{-1}\Sigma_v|\to
\kk^\times$ that is piecewise multiplicative with respect to
$\omega^{-1}\Sigma_v$. Thus the restriction of $s_e$ to the maximal
cone in $\omega^{-1}\Sigma_v$ given by $\sigma$ defines a homomorphism
$\zeta: \Lambda_\sigma\to \kk^\times$. Composing this homomorphism
with $\shAff(\check B,\ZZ)_{\check\sigma}\to \Lambda_\sigma$ leads to
the following ring automorphism of $\kk[P_{\omega,\sigma}]$:
\[
s_{ e,\sigma}: \kk[P_{\omega,\sigma}]\lra
\kk[P_{\omega,\sigma}],\quad
s_{e,\sigma}(z^m)= \zeta(\overline m)\cdot z^m.
\]
Now for any $\rho\in\P^{[n-1]}$ containing $\tau$ denote by $e_\rho$
the composition $v\to \omega\to \tau\to \rho$ and let
$K_{e_\rho}\in\Sigma_v$ be the corresponding cone of codimension
one. Since $V_{v\to\rho}\subseteq V(v)$ equals $\Spec
\kk[K_{e_\rho}\cap \Lambda_v]$ we have the expansion
$f_{e_\rho}=\sum_{\overline m\in  K_{e_\rho}\cap \Lambda_v}
f_{e_\rho,\overline m} z^{\overline m}$. The restriction of this
function to $V_{v\to \tau}$ lifts canonically to $\kk[P_{v,\sigma}]$
as $\sum_{m\in P_{v,\sigma}} f_{e_\rho,m} z^m$ with
\[
f_{e_\rho,m}:= \left\{\begin{array}{ll}
f_{e_\rho,\overline m}&,
\ \overline{m}\in K_{e_\rho},\,\ord_\tau(m)=0,\\
0&,\ \text{otherwise}.
\end{array}\right.
\]
With the inclusion $P_{v,\sigma} \subseteq  P_{\omega,\sigma}$ via
parallel transport inside $\sigma$ understood, we now define
\begin{eqnarray}\label{f_{rho,e,sigma}}
f_{\rho,e,\sigma}:= s_{ e,\sigma}^{-1}\bigg(\sum_{m\in P_{v,\sigma}}
f_{e_\rho,m} z^m\bigg)\in \kk[P_{\omega,\sigma}],
\end{eqnarray}
and the localizing element as
\[
f_{g,\sigma}^v:= \prod_{\rho\supseteq\tau}
f_{\rho,e,\sigma}.
\]
If $\tau\in\P_{\max}$ this product is empty and we take
$f_{g,\sigma}^v=1$. Note that by the normalization condition
$f^v_{g,\sigma}$ has constant term $1$.

For a different choice of vertex $e':v'\to\omega$ equation~(\ref{shLS
change of chart, symmetric form}) implies
\[
s_{ e',\sigma}^{-1}(f_{e'_\rho})=
C\cdot z^{m^\rho_{v'v}}
s_{e,\sigma}^{-1}(f_{e_\rho}),
\]
for some $C\in\kk^\times$. (Here we view $m^\rho_{v'v}\in
\Lambda_\rho\subseteq \Lambda_\sigma$ as an element of
$P_{\omega,\sigma}$ by taking the unique lift $m$ under $P_{\omega,\sigma}
\to \Lambda_\sigma$ with $\ord_\rho(m)=0$. Similar identifications
will occur throughout the text without further notice.) Thus
\begin{eqnarray*}
f_{g,\sigma}^{v'}&=& \prod_{\rho\supseteq\tau}
s_{ e',\sigma}^{-1}\bigg(\sum_{m\in  P_{v',\sigma}}
f_{e'_\rho,m} z^m\bigg)\\
&=& C'\prod_{\rho\supseteq\tau} z^{m^\rho_{v'v}}
s_{ e,\sigma}^{-1}\bigg(\sum_{m\in  P_{v,\sigma}}
f_{e_\rho,m} z^m\bigg)\ =\ C' z^{{\ell m^\rho_{v'v}}} f_{g,\sigma}^v,
\end{eqnarray*}
where $e'_\rho:v'\to\rho$ and $C'\in\kk^\times$ is another constant.
Now the monomials $z^{m^\rho_{v'v}}$ are invertible in
$\kk[P_{\omega,\sigma}]$ since $m^\rho_{v'v}\in\Lambda_\omega$. Hence
the localization of $\kk[P_{\omega,\sigma}]/I^{>k}_{g,\sigma}$ at
$f_{g,\sigma}^v$ does not depend on the choice of $v\in\omega$. We set
$f_{g,\sigma}:= f_{g,\sigma}^v$ for any $v\in\omega$, viewed as
well-defined only up to invertible functions in $\kk[P_{\omega,\sigma}]$.

More generally, if $I\subseteq \kk[P_{\omega,\sigma}]$ is any monomial
ideal with radical $I_{g,\sigma}^{>0}$, set
\[
R_{g,\sigma}^I:= \big(\kk[P_{\omega,\sigma}]/I\big)_{f_{g,\sigma}}.
\]
Any of these rings contains the distinguished monomial
$z^\mathbbm{1}$, where $\mathbbm{1}\in\shAff(\check
B,\ZZ)_{\check\sigma}$ is the constant $1$ function. These monomials
correspond to the deformation parameter $t$ and are preserved by all
our constructions. We therefore write $t=z^\mathbbm{1}$ and keep in
mind that we really work with $\kk[t]$-algebras.
\qed
\end{construction}

\begin{remark}\label{interpretation of rings}
The meaning of the rings $R^k_{\omega\to\tau,\sigma}$ is as follows.
The choice of $e:v\to\omega$ determines the local model $V(v) =
\Spec\kk[\Sigma_v] \subseteq \Spec \kk[P_{v,\sigma}]$ for
$X\subseteq\X$, and the open subsets $\Spec\kk[\omega^{-1}\Sigma_v]
\subseteq \Spec\kk[\Sigma_v]$ and $\Spec\kk[P_{\omega,\sigma}]
\subseteq \Spec \kk[P_{v,\sigma}]$. The open gluing construction of
$X$ yields the open embedding
\[
\Phi_{v,\omega}(s): V(\omega)
\lra V( v)
\]
by twisting $\Spec\kk[\omega^{-1}\Sigma_v] \to \Spec\kk[\Sigma_v]$ by
$s_{ e}^{-1}$. Moreover, $g:\omega\to\tau$ determines the
toric stratum
\[
V_g=\Spec \big(\kk[P_{\omega,\sigma}]/I_{g,\sigma}^{>0}\big)
\subseteq V(\omega).
\]
Thus $\Spec \big(\kk[P_{\omega,\sigma}]/I_{g,\sigma}^{>k}\big)$ is the
$k$-th order thickening of $V_g$ inside
$\Spec \kk[P_{\omega,\sigma}]$.

As for the localization recall that the singular locus $Z$ of the log
structure on $V( v)$ is the union of the zero loci $Z_e$ of $f_e$ for
$e:v\to\rho$, $\rho\in\P^{[n-1]}$. Thus the zero locus of
$f_{g,\sigma}^v$ equals the $\Phi_{v,\omega}(s)$-preimage of
$p^{-1}(Z_\tau)$ for $Z_\tau:=\bigcup_{e:\tau\to\rho} Z_e$. In
summary, $\Spec R^k_{g,\sigma}$ is isomorphic to the $k$-th order
thickening of $V_g\setminus \Phi_{v,\omega}(s)^{-1} (p^{-1}(Z_\tau))$
inside $V(\omega)$.
\qed
\end{remark}

\begin{remark}\label{natural ring homomorphisms}
Let $g:\omega\to\tau$, $g':\omega'\to\tau'$ and assume
$\omega\subseteq \omega'$ and $\tau\supseteq \tau'$. We also fix a
reference cell $\sigma\in\P_\max$ containing $\tau$. Then
$P_{\omega',\sigma}$ differs from $P_{\omega,\sigma}$ by making
invertible those $m\in P_{\omega,\sigma}$ with
$\ord_{\omega'}(m)=0$. Moreover, for any $m\in P_{\omega,\sigma}$ it
holds $\ord_{\tau'}(m)\ge \ord_\tau(m)$ since $\tau'\subseteq\tau$.
Hence $I_{g,\sigma}^{>k}\subseteq I_{g',\sigma}^{>k}$, and we
obtain the canonical homomorphism
\[
\psi_0:\kk[P_{\omega,\sigma}]/ I_{g,\sigma}^{>k} \lra
\kk[P_{\omega',\sigma}]/ I_{g',\sigma}^{>k}.
\]
If $\omega=\omega'$ then $\psi_0(f_{g,\sigma})$ divides $f_{g',\sigma}$,
and hence $\psi_0$ induces a map $R^k_{g,\sigma}\to R^k_{g',\sigma}$.

In the general case we need to take into account the twisting by the
open gluing data as follows. The piecewise multiplicative function
$s_a$, $a:\omega\to\omega'$, coming from the open gluing data,
is given on $\sigma$ by a homomorphism $s_{a,\sigma}:
\Lambda_{\sigma}\to \kk^\times$. This defines a ring automorphism
of $\kk[P_{\omega',\sigma}]$ respecting orders. Hence it induces
an automorphism of $\kk[P_{\omega',\sigma}]/
I_{g',\sigma}^{>k}$ that we also denote  $s_{a,\sigma}$.
The special case where $\omega$ is a vertex is the case used in
Construction~\ref{constr rings}. Now if $e:v\to\omega$ then $s_{a\circ
e,\sigma} = s_{a,\sigma}\cdot s_{e,\sigma}$, which
implies
\[
f_{\rho,a\circ e,\sigma} = (s_{a,\sigma}^{-1}\circ\psi_0)
\big(f_{\rho,e,\sigma} \big).
\]
Thus $s_{a,\sigma}^{-1}\circ\psi_0$ defines a well-defined map
\[
\psi_0(s): R^k_{g,\sigma}\lra R^k_{g',\sigma}.
\]
\vspace{-6ex}

\qed
\end{remark}

\begin{example}\label{2d example, part I}
To illustrate the use of the rings $R^k_{g,\sigma}$ in our
construction let us look at a simple example that captures the
situation in codimension one. Assume that $\rho$ is a one-dimensional
cell in a two-dimensional $B$, with vertices $v_1$, $v_2$, and
monodromy constant $\kappa:=\kappa_{\rho\rho}\ge 0$, see
(\ref{monodromy constant}). For simplicity we assume the open gluing
data to be trivial ($s_e=1$ for all $e$). Let the dual cell $\check
\rho$ have integral length $l$. Let $\sigma_1$, $\sigma_2$ be the two
maximal cells with $\rho=\sigma_1 \cap\sigma_2$. Then $\Spec
\kk[P_{\rho,\sigma_i}]$ is isomorphic to $\AA^1\setminus\{0\}$ times
the two-dimensional $A_{l-1}$-singularity. With $g_i:\rho\to
\sigma_i$, and $\id_\rho:\rho\to\rho$ we have the two maps
$R^k_{g_i,\sigma_i}\to R^k_{\id_\rho,\sigma_i}$, which in appropriate
coordinates are just canonical quotient homomorphisms composed with a
localization:
\begin{eqnarray*}
\kk[w,w^{-1},x_1,y_1,t]/(x_1 y_1-t^l,y_1^{\bar k+1}) 
&\lra&  \kk[w,w^{-1},x_1,y_1,t]_{f_{\id_\rho,\sigma_1}}/
(x_1 y_1-t^l,x_1^{\bar k+1},y_1^{\bar k+1}),\\
\kk[w,w^{-1},x_2,y_2,t]/(x_2 y_2-t^l,x_2^{\bar k+1})
&\lra&  \kk[w,w^{-1},x_2,y_2,t]_{f_{\id_\rho,\sigma_2}}/
(x_2 y_2-t^l,x_2^{\bar k+1},y_2^{\bar k+1}).
\end{eqnarray*}
Here we assumed for simplicity that $k+1=(\bar k+1)\cdot l$ for $\bar
k\in\NN$. Now by parallel transport through $v_i$ we obtain two
isomorphisms
\[
\psi_i:R^k_{\id_\rho,\sigma_1}\to R^k_{\id_\rho,\sigma_2}.
\]
This gives two fibre products $R^k_{g_1,\sigma_1}
\times_{R^k_{\id_\rho,\sigma_2}} R^k_{g_2,\sigma_2}$, and we will
prove in a more general context in Lemma~\ref{model at general points
of Z} that each is isomorphic to $\Spec \kk[w,w^{-1},x,y,t]/(xy-t^l,
t^{k+1})$.

However, for $k+1\ge l$ there exists no isomorphism between these two
fibre products inducing the identity on $R^k_{g_i,\sigma_i}$ unless
$\kappa=0$. In fact, if $w=z^m$ with $\overline m\in\Lambda_\rho$ the
generator pointing from $v_1$ to $v_2$ then $\psi_1$ and $\psi_2$ are
related by the automorphism
\[
\psi_2^{-1}\circ\psi_1: w\mapsto w,\ x_1\mapsto  w^{-\kappa} x_1,\ 
y_1\mapsto w^\kappa y_1,\ t\mapsto t
\]
of $R^k_{\id_\rho, \sigma_1}$,
which is not the identity unless $\kappa=0$ or $k<l$. For a
continuation of this discussion see Example~\ref{2d example, part II}.
\qed
\end{example}

The example illustrates that monodromy yields an obstruction to gluing
the standard $k$-th order deformations of the local models of $X$
consistently. To remedy this we need to compose the maps between rings
by automorphisms. These automorphisms are 
the subject of the next subsection, in a
log setting for our rings that we now discuss.

\sloppy
By construction $R^I_{\omega\to\tau,\sigma}$ comes with a homomorphism
of monoids $(P_{\omega,\sigma},+)\to (R^I_{\omega\to\tau,
\sigma},\cdot)$. This yields a chart for a log structure on $\Spec
R^I_{\omega \to\tau,\sigma}$, and it will be very important in the
algorithm to trace this information. We therefore now introduce a
category of rings with charts.

\fussy
\begin{definition}\label{category of log rings}
A \emph{log ring} is a ring $R$ together with a monoid homomorphism
$\alpha:P\to (R,\cdot)$. A \emph{morphism} of log rings (or \emph{log
morphism}) $(\alpha:P\to R) \to (\alpha':P'\to R')$ consists of a ring
homomorphism $\psi: R\to R'$ and monoid homomorphisms
\[
\beta: P\to P',\quad \theta:P\to (R')^\times,
\]
such that
\begin{equation}\label{log morphism condition}
\psi\circ\alpha= \theta\cdot (\alpha'\circ\beta).
\end{equation}
If $(\beta,\theta,\psi): (P\to R)\to (P'\to R')$ and
$(\beta', \theta',\psi'): (P'\to R')\to (P''\to R'')$ are
log morphisms, their composition is defined as
\[
\big(\beta'\circ\beta,(\psi'\circ\theta)\cdot
(\theta'\circ\beta), \psi'\circ\psi \big).
\]
This is indeed a log morphism from $(P\to R)$ to $(P''\to R'')$ as one
easily checks. Two log morphisms $(\beta_1,\theta_1,\psi_1)$,
$(\beta_2,\theta_2, \psi_2)$ from $\alpha: P\to R$ to $\alpha' :
P'\to R'$ are \emph{equivalent} if there exists a
homomorphism $\eta: P\to (P')^\times$ such that
\[
\beta_2=\beta_1+\eta,\quad
\theta_2=\theta_1\cdot(\alpha'\circ\eta)^{-1},\quad
\psi_2=\psi_1.
\]
Log rings and log morphisms modulo equivalence define the category
$\LogRings$.
\qed
\end{definition}

\begin{remark}\label{convention log morphisms}
1)\ This definition just rephrases the basic notions of log geometry
\cite{kato} on the level of rings. In particular, a log ring
$\alpha:P\to R$ is the same as an affine scheme $X=\Spec R$ with a
chart for a log structure $\gamma:P\to\Gamma(X,\M_X)$ in the Zariski
topology. Note such a chart induces a canonical isomorphism
\[
\big(\O_X^\times\oplus P_X\big)\big/\big\{(h,m)\,\big|\, h\cdot
\alpha(m)=1 \big\} \lra \M_X,
\]
so we can represent elements of $\M_X$ as pairs $(h,m)$,
$h\in\O_X^\times$, $m\in P$.

Similarly, a log morphism $(\beta,\theta,\psi)$ between the log
rings $\alpha:P\to R$ and $\alpha':P'\to R'$ gives rise to a
morphism of the associated affine log schemes as follows. The map
$\underline f:X':=\Spec R'\to X:= \Spec R$ of the underlying schemes
is defined by $\psi$. Then by \eqref{log morphism condition},
\[
\underline f^{-1}(\O_X^\times\oplus P_X)\lra \O_{X'}^\times\oplus
P'_{X'},\quad (h,m)\longmapsto \big(\psi(h)\cdot\theta(m),
\beta(m)\big)
\]
descends to a morphism $\underline f^{-1}\M_X\to \M_{X'}$. Indeed,
$h\cdot \alpha(m)=1$ implies
\[
1=\psi\big (h\cdot\alpha(m)\big ) =\psi(h)\cdot(\psi\circ\alpha)(m)
\stackrel{\eqref{log morphism condition}}{=}
\psi(h)\cdot \theta(m)\cdot\alpha'\big(\beta(m)\big). 
\]
Conversely, under the assumption that for the chart $\gamma':P'\to
\Gamma(X',\M_{X'})$ no non-zero element of $P'$ maps to an
invertible element, any morphism of log schemes
$(\underline f,f^\flat):(X',\M_{X'})\to (X,\M_X)$ arises in this fashion. In
fact, under the stated condition the composition
\[
P\stackrel{\gamma}{\lra} \Gamma(X,\M_X)\stackrel{f^\flat}{\lra}
\Gamma(X',\M_{X'})\stackrel{\kappa}{\lra}
\Gamma(X', \M_{X'}/\O_{X'}^\times),
\]
with $\kappa$ the quotient homomorphism, factors canonically over
$\kappa\circ\gamma'$, thus defining $\beta:P\to P'$. Comparison of
$\gamma'\circ\beta$ with $f^\flat\circ \gamma$ then defines
$\theta$. Note that on the side of log rings the stated condition
translates into the requirement ${\alpha'}^{-1}({R'}^\times)=\{0\}$.

Thus at least for log rings $\alpha:P\to R$ fulfilling
$\alpha^{-1}(R^\times)=\{0\}$ our discussion also shows that the
notion of equivalence is compatible with compositions of
morphisms.\\
2)\ In our case $R=R^k_{g,\sigma}$ is a localization of a quotient of
$\kk[P_{\omega,\sigma}]$, and hence carries canonically the structure
of a log ring via $\alpha: P_{\omega,\sigma}\to R^k_{g,\sigma}$,
$\alpha(m)=z^m$. Because $P_{\omega,\sigma}$ generates
$R^k_{g,\sigma}$ up to localization,~(\ref{log morphism condition})
determines the underlying ring homomorphism $\psi$ of a log morphism
from $\beta$ and $\theta$. Moreover, in the cases we are interested
in, $\beta$ is either canonically given or is fixed in the discussion
and $\theta$ factors through the projection $P_{\omega,\sigma}\to
\Lambda_\sigma$. By abuse of notation we then talk of a group
homomorphism $\theta:\Lambda_\sigma\to (R_{g',\sigma})^\times$ as
being a log morphism. We write $\overline\theta$ for the associated
ring homomorphism, and use $\theta(m)$ and $\theta(\overline m)$
interchangeably. Explicitly, we have
\[
\overline\theta(z^m):=\theta(\overline m)\cdot z^{\beta(m)}
\]
for the underlying ring homomorphism, and the composition of two log
morphisms $\theta_1$, $\theta_2$ reads
\begin{align}\label{composition of log morphisms}
(\theta_1\circ \theta_2)(m) &= \theta_1(m)\cdot
\overline{\theta_1}\big(\theta_2(m)\big).
\end{align}
\nopagebreak
\\[-6ex]
\qed
\end{remark}

%===========================================================
\subsection{Automorphism groups}
\label{subsection automorphism groups}
We will now discuss various groups of log automorphisms of the rings
which appear in our construction. For this subsection fix
$g:\omega\to\tau$, $\sigma\in\P_\max$ with $\tau\subseteq \sigma$,
and a monomial ideal $I\subseteq \kk[P_{\omega,\sigma}]$ with
radical $I_0:=I_{g,\sigma}^{>0}$. Let $f:=f_{g,\sigma}^v\in
\kk[P_{\omega,\sigma}]$, $v\in \omega$ a vertex, be a localizing
element as in Construction~\ref{constr rings}. Write
$P:=P_{\omega,\sigma}$ and $R^I:=R_{g,\sigma}^I=(\kk[P]/I)_f$.
Recall also the projection $P\to\Lambda_\sigma$, $m\mapsto \overline
m$ and the conventions of Remark~\ref{convention log morphisms}. We
are interested in log automorphisms of $P\to R^I$. 

\begin{remark}\label{log automorphism remark}
1) The inverse of a log automorphism $\theta$ is
\[
\theta^{-1}(m)=
\overline\theta^{-1}\left(\frac{1}{\theta(m)}\right).
\]
2) The formula for multiple compositions is
\begin{eqnarray}\label{multiple composition} 
\theta_1\circ\theta_2\circ \ldots\circ\theta_r =
\theta_1\cdot (\overline\theta_1\circ \theta_2)\cdot\ldots\cdot
(\overline\theta_1\circ\overline\theta_2\circ\ldots\circ
\overline\theta_{r-1}\circ\theta_r).
\end{eqnarray}
On the right-hand side the composition symbol denotes ordinary
composition of maps, and the centered dots denote multiplication of
maps with target $R^I$.
\qed
\end{remark}

We will now describe the group of all log automorphisms
$\theta:\Lambda_\sigma\to (R^I)^\times$ with the property that
$\theta(m)=1\mod I_0$, by describing the Lie algebra of this group.
We first consider the \emph{module of log derivations} of $R^I$,
defined by
\[
\Theta(R^I):=R^I\otimes_{\ZZ}\Lambda^*_\sigma=
\Hom(\Lambda_\sigma,R^I).
\]
We view an element $\xi\in\Theta(R^I)$ as an additive map
$\xi:P\to R^I$ factoring through $P\to\Lambda_\sigma$.
In particular, $a\otimes n$ defines the map
\[
P\ni m\longmapsto a\langle\overline m ,n\rangle.
\]
Note $\xi\in\Theta(R^I)$ also induces an ordinary $\kk$-derivation of
$R^I$ via
\[
\overline \xi(z^m):= \xi(m) z^m.
\]
It is then suggestive to write $a\partial_n$ for $a\otimes n\in
\Theta(R^I)$ or its associated ordinary derivation:
\[
(a\partial_n)(z^m)=a\langle\overline m ,n\rangle z^m.
\]

The adjoint action of the group of automorphisms on derivations lifts
to the log setting by defining, for $\theta$ a log automorphism and
$\xi$ a log derivation,
\begin{eqnarray}\label{adjoint action}
\Ad_\theta \xi:=(\overline\theta\circ\overline
\xi\circ\theta^{-1})\cdot \theta + \overline\theta\circ \xi.
\end{eqnarray}

Given $\xi_1,\ldots,\xi_n\in \Theta(R^I)$, we can define a higher
order log differential operator, a map $\xi_1\circ\cdots\circ
\xi_n: \Lambda_\sigma\to R^I$, inductively by the formula
\begin{eqnarray}\label{higher log differentials}
(\xi_1\circ\cdots\circ \xi_n)(m)=\xi_1(m)\cdot(\xi_2\circ\cdots\circ \xi_n)(m)
+\overline \xi_1(\xi_2\circ\cdots\circ \xi_n(m)),
\end{eqnarray}
so that 
\[
z^m(\xi_1\circ\cdots\circ \xi_n)(m)=
(\overline \xi_1\circ\cdots\circ\overline \xi_n)(z^m),
\]
where the composition on the right-hand side is just the composition
of ordinary $\kk$-endomorphisms of $R^I$. The powers of
$\xi\in \Theta(R^I)$ fulfill a higher order Leibniz rule:
\begin{eqnarray}\label{higher order Leibniz rule}
\xi^n(m_1+m_2)=\sum_{i=0}^n \begin{pmatrix}n\\i\end{pmatrix}
\xi^i(m_1)\xi^{n-i}(m_2),\quad m_1,m_2\in \Lambda_\sigma.
\end{eqnarray}

\begin{proposition}
The group
\[
G^I:= \big\{\theta: \Lambda_\sigma\to (R^I)^{\times}\,\big|\,
\theta\ \text{\rm is a log automorphism},
\forall m\in \Lambda_\sigma:\theta(m)=1\mod I_0\big\}
\]
is an algebraic group with Lie algebra $\fog^I:=I_0\cdot\Theta(R^I)$
endowed with the Lie bracket
\[
[\xi_1,\xi_2]:=\xi_1\circ \xi_2-\xi_2\circ \xi_1.
\]
\end{proposition}

\proof
From~(\ref{higher order Leibniz rule}) it follows that if $\xi\in
I_0\cdot\Theta(R^I)$, the formula 
\begin{eqnarray}\label{exp of log derivations}
\exp(\xi)(m):=1+\sum_{i=1}^{\infty} {\xi^i(m)\over i!}
\end{eqnarray}
defines an element $\exp(\xi)\in G^I$ since $\exp(\xi)(m_1+m_2)=
\exp(\xi)(m_1)\cdot \exp(\xi)(m_2)$. Note the sum is finite because
$\sqrt{I}=I_0$. 

Conversely, let $\theta\in G^I$. Define inductively $N_i:
\Lambda_\sigma\to R$ by $N_0:=1$ and
\[
N_i:=\theta\cdot \big(\overline\theta \circ N_{i-1}\big)-N_{i-1}.
\]
The induced map $\overline N_i: z^m\mapsto N_i(m) z^m$ equals
$(\overline\theta-\id)^i$. Note $N_i$ takes values in $I_0^i$. Thus we
can define $\log(\theta):\Lambda_\sigma\to R^I$ by
\[
\log(\theta):=\sum_{i=1}^{\infty}\frac{(-1)^{i+1}}{i} N_i.
\]
Again, this is a finite sum. Noting inductively that
\[
N_n(m_1+m_2)=\sum_{i+j+k=n, i,j,k\ge 0}
{n!\over i!j!k!} N_{i+j}(m_1) N_{i+k}(m_2),
\]
it follows by direct computation that $\log(\theta)$ is additive.
Hence $\log(\theta)\in I_0\cdot\Theta(R^I)$, and then the usual power
series identity implies $\theta=\exp(\log(\theta))$.
\qed
\medskip

On the $\kk$-basis $z^m\partial_n$ of $\fog^I$ the formula for the Lie
bracket is
\begin{align}\label{Lie bracket}
\begin{split}
[z^{m}\partial_n,z^{m'}\partial_{n'}]
&=
(z^{m}\partial_n(z^{m'}))\partial_{n'}
-(z^{m'}\partial_{n'}(z^{m}))\partial_{n}\\[1ex]
&=z^{m+m'}(\langle \overline {m'},n\rangle\partial_{n'}
-\langle\overline m ,n'\rangle\partial_n)
=z^{m+m'}\partial_{\langle \overline {m'},n\rangle n'
-\langle\overline m ,n'\rangle n}.
\end{split}
\end{align}
In particular, $\fog^I$ is a nilpotent Lie algebra.

Later on we will often need to control how the basic elements
$\exp(z^m\partial_n)$ commute with certain more general log
automorphisms. For this we record the following lemma.

\begin{lemma}\label{conjugation lemma}
For ${h}\in(R^I)^\times$ consider the log
automorphism
\[
\theta: m\longmapsto {h}^{-\langle \overline m,n_0\rangle},
\]
of $R^I$, where $n_0\in \Lambda_\sigma^*$ annihilates any
exponent occurring in ${h}$. Then for $m\in \Lambda_\sigma$,
$n\in\Lambda_\sigma^*$
\[
\Ad_{\theta}(z^m\partial_n)= z^m\big(
{h}^{-\langle \overline m,n_0\rangle}\partial_n
+ {h}^{-\langle \overline m,n_0\rangle-1}(\partial_n {h})
\partial_{n_0}\big).
\]
\end{lemma}

\proof
Using the fact that every monomial in ${h}$ is left invariant by
$\overline\theta$ we get $\theta^{-1}(m)= {h}^{\langle
\overline m, n_0\rangle}$ and, with $\xi=z^m\partial_n$,
\begin{eqnarray*}
\Ad_{\theta}(\xi)(m')
&\stackrel{(\mbox{\footnotesize\ref{adjoint action}})}{=}&
(\overline\theta\circ\overline\xi\circ\theta^{-1})(m')
\cdot \theta(m') +(\overline \theta \circ\xi)(m')\\
&=&(\overline\theta\circ \overline\xi) \big({h}^{\langle \overline{m'},n_0
\rangle}\big) \cdot {h}^{-\langle \overline{m'},n_0 \rangle}
+\overline\theta\big(\langle \overline{m'},n\rangle z^m\big)\\
&=& \overline\theta \big(\langle \overline{m'},n_0 \rangle {h}^{\langle
\overline{m'},n_0 \rangle-1} (\partial_n {h}) z^m\big) {h}^{-\langle
\overline{m'},n_0 \rangle}
+ \langle \overline{m'},n \rangle
{h}^{-\langle \overline m,n_0 \rangle}z^m\\
&=&\langle \overline{m'},n_0 \rangle {h}^{-\langle \overline m,n_0
\rangle-1} (\partial_n {h}) z^m
+ \langle \overline{m'},n \rangle
{h}^{-\langle \overline m,n_0 \rangle}z^m.
\end{eqnarray*}
\\[-6ex]
\qed
\medskip

For any sub-Lie algebra $\foh\subseteq \fog^I$, we obtain a subgroup
$H=\exp(\foh)$ of $G^I$ consisting of exponentials of elements of
$\foh$. We shall consider a number of such subgroups.

In what follows, fix a codimension two subspace $T_\foj \subseteq
\Lambda_{\sigma,\RR}$ defined over $\QQ$, and write $\Lambda_{\foj}
=T_\foj \cap \Lambda_\sigma$. Later on $T_\foj$ will be the tangent
space to a polyhedral subset of $\sigma$ of codimension two. Write
$P^{>0}= P_{g,\sigma}^{>0}\subseteq P_{\omega,\sigma}$ and $P^I$ for
the monoid ideals generating $I_0$ and $I$. Then each of the following
subspaces of $\fog^I$ are Lie subalgebras, as is easily checked using
(\ref{Lie bracket}):
\begin{eqnarray*}
\fog^I_{\foj}&:=&
\bigoplus_{m\in P^{>0}\setminus P^I} z^m\big(\kk\otimes
\Lambda_{\foj}^{\perp}\big)\\
\tilde\foh^I_{\foj}&:=&
\bigoplus_{m\in P^{>0}\setminus P^I} z^m\big(\kk\otimes 
(\overline{m}^\perp\cap\Lambda_{\foj}^{\perp})\big)\\
\foh^I_{\foj}&:=&
\bigoplus_{m\in P^{>0}\setminus P^I\atop\overline m \not=0}
z^m\big(\kk\otimes  (\overline{m}^\perp\cap\Lambda_{\foj}^{\perp})\big)\\
{}^\perp\foh^I_{\foj}&:=&
\bigoplus_{m\in P^{>0}\setminus
P^I\atop\overline m \not\in \Lambda_{\foj}}
z^m\big(\kk\otimes (\overline{m}^\perp\cap\Lambda_{\foj}^{\perp})\big)\\
{}^\parallel\foh^I_{\foj}&:=&
\bigoplus_{m\in P^{>0}\setminus P^I\atop
\overline m \in \Lambda_{\foj}\setminus\{0\}} z^m\big(\kk\otimes 
\Lambda_{\foj}^{\perp}\big).
\end{eqnarray*}
The corresponding subgroups of $G^I$ are denoted $G^I_{\foj}$, $\tilde
H^I_{\foj}$, $H^I_{\foj}$, $\lperp H^I_{\foj}$, and $\lparallel
H^I_{\foj}$, respectively. Of these the most essential one for our
construction is $H^I_\foj$ with Lie algebra generated by derivations
$z^m\partial_n$, where $\partial_n$ acts trivially on $z^m$ and $m$
points in a specific direction ($\overline m\neq0$).

\begin{remarks}\label{rem:log automorphisms}
(1) All $\theta\in G^I_{\foj}$ satisfy $\theta(m)=1$ whenever
$\overline m \in\Lambda_{\foj}$.\\[1ex]
(2) The log automorphism associated to an element of $\tilde\foh_\foj^I$
of the form $a\partial_n$ is easy to write down
explicitly:
\[
\exp(a\partial_n)(m)=\exp(\langle \overline m, n\rangle a)
= \exp(a)^{\langle \overline m,n\rangle}.
\]
Here $\exp(a)$ is the usual exponential of a function, which is a
polynomial in $a$ because $a\in I_0$. Indeed, $a$ involves only
monomials $z^m$ with $\langle \overline m, n\rangle=0$, and hence the
composition formula~(\ref{higher log differentials}) inductively shows
\[
(a\partial_n)^i(m)=\big(\langle \overline m,n\rangle a\big)^i.
\]
\\[-2ex]
(3) Denote by
\[
\Omega^p(R^I)=R^I\otimes {\bigwedge}^p\Lambda_\sigma
\]
the module of logarithmic $p$-forms on $R^I$. We write $f\otimes
m=f\dlog m$ for $f\in R^I$, $m\in \bigwedge^p\Lambda_\sigma$. We
have the (ordinary) exterior derivative
\[
d:R^I\lra\Omega^1(R^I)=
\Hom(\Lambda^*_\sigma,R^I),\quad
f\longmapsto \left(n\mapsto \partial_n f\right).
\] 
This gives
\[
d:\Omega^p(R^I)\lra\Omega^{p+1}(R^I),
\]
in the usual way, and $\fog^I$ acts on $\Omega^p(R^I)$ by Lie
derivative. In particular, $\fog^I$ acts on $\Omega^{\dim B}(R^I)$ by
$\xi(\Omega)=\shL_\xi(\Omega)=d(\iota(\xi)\Omega)$ for $\xi\in
\fog^I$, $\Omega\in \Omega^{\dim B}(R^I)$. It is then not difficult to
see that $\tilde\foh^I_{\foj}$ consists of those elements of
$\fog^I_{\foj}$ which preserve
\[
\Omega_{\std}=\dlog(m_1\wedge\ldots\wedge m_n)
=\dlog(m_1)\wedge\ldots\wedge \dlog(m_n),
\]
where $m_1\wedge\ldots\wedge m_n$ is a primitive generator of
$\bigwedge^{\dim B} \Lambda_\sigma$. In fact,
\[
\shL_{z^m\partial_n}\Omega_\std
=d\big( z^m\iota_{\partial_n}\Omega_\std \big)
=\langle m,n\rangle z^m\Omega_\std.
\]

Note also that any log automorphism of $R^I$ acts on $\Omega^p(R^I)$
by
\begin{eqnarray*}
\lefteqn{ \theta
\big(a\dlog(m_1)\wedge\cdots\wedge\dlog(m_p)\big)}\hspace{2cm}&\\
&:=&\overline\theta(a)\Bigl(\dlog
m_1+{d\theta(m_1)\over\theta(m_1)}\Bigr)\wedge
\cdots\wedge\Bigl(\dlog m_p+{d\theta(m_p)\over\theta(m_p)}\Bigr).
\end{eqnarray*}
One can check that whenever $\theta\in G^I$ this agrees with the
exponential of the action of $\fog^I$ on $\Omega^p(R^I)$, that is, for
$\xi\in\fog^I$ and $\alpha\in\Omega^p(R^I)$ it holds
\[
\big(\exp(\xi)\big)(\alpha)=\sum_{i=0}^\infty \frac{1}{i!} \shL_\xi^i (\alpha).
\]
Thus $\tilde H^I_{\foj}$ consists of those log automorphisms in
$G^I_{\foj}$ preserving $\Omega_{\std}$.\\[1ex]
(4) Note that ${}^\parallel\foh^I_{\foj}$ is abelian and
$[{}^\parallel\foh^I_{\foj}, {}^\perp\foh^I_{\foj}]\subseteq
{}^\perp\foh^I_{\foj}$, so we get an exact sequence of Lie algebras
\[
0\lra {}^\perp\foh^I_{\foj}
\lra \foh^I_{\foj}
\lra {}^\parallel\foh^I_{\foj}
\lra 0,
\]
and hence an exact sequence of groups
\[
1\lra \lperp H^I_{\foj}
\lra H^I_{\foj}
\lra \lparallel H^I_{\foj}
\lra 1.
\]
\vspace{-6ex}

\qed
\end{remarks}

%===========================================================
\subsection{Slabs, walls and structures}
Our construction involves splitting $B$ into smaller and
smaller pieces which are separated by \emph{slabs} and \emph{walls}.
We begin with the subdivision of $B$ given by $\P$; the codimension
one elements of $\P$ define slabs. We then proceed to subdivide $B$ through
a scattering process by adding walls,  which are codimension one
polyhedra contained  in maximal elements of $\P$. These walls split
these maximal cells into chambers. The choice of words ``wall'' and
``slab'' is inspired by the first author's house, which is built on a
slab. Just as with this house, over time, the slabs develop cracks,
and here are subdivided, while once a wall is introduced, it remains
unmodified during the process of further subdivisions of $\P$. A slab
also carries additional data, namely the starting data determined by
the log structure and some higher order corrections, while a wall only
carries higher order data. A further difference is that walls, unlike
slabs, have a built-in directionality. Both slabs and walls lead to
log automorphisms of rings $R^k_{g,\sigma}$, which will be used to
glue together these rings to create $k$-th order deformations.

For the following definition recall the open star $U_\tau=
\bigcup_{\sigma\in \P, \sigma\supseteq\tau} \Int \sigma$ of a cell
$\tau$, and the notation $v[x]\in\rho$ for the unique vertex in the
same connected component of $\rho\setminus\Delta$ of some
$x\in\rho\setminus\Delta$, $\rho\in\P^{[n-1]}$. 

\begin{definition}\label{slab}
A \emph{slab} is a convex, rational, $(n-1)$-dimensional polyhedral
subset $\fob$ of a cell $\rho_\fob\in\P^{[n-1]}$ together
with elements
\[
f_{\fob,x} = \sum_{m\in P_x,\, 
\overline m\in\Lambda_{\rho_\fob}} c_m z^m\in \kk[P_x],
\]
one for each $x\in\fob\setminus\Delta$, satisfying the following
properties:
\begin{enumerate}
\item[(i)]
If $x,x'\in\fob\setminus\Delta$, $\Pi:\kk[P_x^\gp] \to
\kk[P_{x'}^\gp]$ is defined by parallel transport along a path inside
$\cl(U_{\rho_\fob})\setminus\Delta$ and $v=v[x], v'=v[x']$ then
\begin{eqnarray}\label{change of vertex - slabs}
D(s_{ e'},\rho_\fob,v')^{-1} s_{ e'}^{-1}(f_{\fob,x'})=
z^{m^{\rho_\fob}_{v'v}}
\Pi\big(D(s_{ e},\rho_\fob,v)^{-1} s_{ e}^{-1}(f_{\fob,x})\big),
\end{eqnarray}
where $e:v\to\rho_\fob$, $e':v'\to\rho_\fob$.
\item[(ii)]
If $e:v\to\rho_\fob$ with $v=v[x]$, and $\Pi:\kk[P_x^\gp]\to
\kk[P_v^\gp]$ is defined by parallel transport from $x$ to $v$ along a
path inside $\rho_\fob\setminus\Delta$ then
\[
f_e=\Pi\Big(\sum_{m\in P_x ,\, \ord_{\fob}(m)=0}c_m z^m \Big),
\]
where $(f_e)$ is the section of $\shLS^+_{\pre,X}$ defining the log
smooth structure on $X$.
\end{enumerate}
\end{definition}

\begin{remarks}
1)\ By condition~(i) the functions $f_{\fob,x}$ determine each other
by parallel transport inside $\cl(U_{\rho_\fob})\setminus\Delta$. In
particular, a slab carries only finitely many non-zero coefficients
as information.\\
2)\ Condition~(ii) says that  $(f_e)$  determines the part of
$f_{\fob,x}$ of order $0$, for every $x\in \fob\setminus\Delta$. Note
that by~(\ref{shLS change of chart, symmetric form}) this is
compatible with condition~(i).
\end{remarks}

\begin{example}\label{2d example, part II}
Continuing on Example~\ref{2d example, part I} let us show how slabs
resolve the problem of incompatible gluings due to monodromy. In this
example add a slab $\fob$ with support the one-dimensional cell
$\rho=\sigma_1\cap\sigma_2$. We view the functions $f_{\fob,v_i}$ as
elements of $R^k_{\id_{\rho},\sigma_2}$ via parallel transport from $v_i$
into $\sigma_2$. Let $\pi_i:\Lambda_{\sigma_i} \to \ZZ$ be the
projection with kernel $\Lambda_\rho$ and which is positive on vectors
pointing from $\sigma_1$ to $\sigma_2$. In going from $\sigma_1$ to
$\sigma_2$  compose the isomorphism $\psi_i:
R^k_{\id_\rho,\sigma_1}\to R^k_{\id_\rho,\sigma_2}$, obtained via
parallel transport through the vertex $v_i\in\rho$, with
\[
z^m\longmapsto (f_{\fob,v_i})^{-\pi_1(\overline m)} \cdot z^m.
\]
In appropriate coordinates these are the homomorphisms of
$\kk[w,w^{-1}]$-algebras sending $x_1,y_1$ to $f_{\fob,v_1}\cdot x_2$,
$f_{\fob,v_1}^{-1}\cdot y_2$ ($i=1$) and to $f_{\fob,v_2}\cdot
w^\kappa x_2$, $f_{\fob,v_2}^{-1}\cdot w^{-\kappa} y_2$ ($i=2$),
respectively. Now~(\ref{change of vertex - slabs}) requires
$f_{\fob,v_2}= w^{-\kappa} f_{\fob,v_1}$, and this is exactly what is
needed to make the two homomorphisms agree. Thus $R^k_{g_1,\sigma_1}
\times_{R^k_{\id_\rho,\sigma_2}} R^k_{g,\sigma_2}$ is well-defined.

Explicitly, computing in the chart at $v_1$, the fibre product is
generated as $\kk[w,w^{-1},t]$-algebra by $X:=(x_1,f_{\fob,v_1} x_2)$,
$Y:=(f_{\fob,v_1} y_1,y_2)$ with
single relation $XY-F_{\fob,v_1}t^l$, where $F_{\fob,v_1}
=(f_{\fob,v_1},f_{\fob,v_1})$. Note how this fits with the
interpretation of the section $f_{v_1\to\rho}$ of $\shLS^+_{\pre,
V(v_1)}|_{V(\rho)}$ in~(\ref{codim one log model}).
\qed 
\end{example}

\begin{definition}\label{wall}
A \emph{wall} is a convex, rational, $(n-1)$-dimensional polyhedral
subset $\fop$ of a maximal cell $\sigma_\fop\in \P^{[n]}$ with
$\fop\cap\Int \sigma_\fop\neq\emptyset$ together with (i)~an
$(n-2)$-face $\foq\subseteq\fop$, $\foq\not\subseteq \partial B$, the
\emph{base} of $\fop$ (ii)~an exponent $m_\fop$ on $\sigma_\fop$ with
$\ord_{\sigma_\fop}(m_\fop)>0$ and $m_{\fop,x}\in P_x$ for every
$x\in\fop\setminus\Delta$, and (iii)~$c_\fop\in\kk$, such that
\[
\fop= (\foq-\RR_{\ge 0}\overline m_\fop)\cap\sigma_{\fop}.
\]
Here we view $\sigma_\fop$ as a polyhedron in
$\Lambda_{\sigma_{\fop},\RR}$ and $m_\fop$ as an element of
$\shAff(\check B,\ZZ)_{\check\sigma_\fop}$. The notation is $(\fop,
m_\fop, c_\fop)$, or simply $\fop$ if $m_\fop$ and $c_\fop$ are
understood.
\begin{figure}[h]
\includegraphics*[scale=0.25]{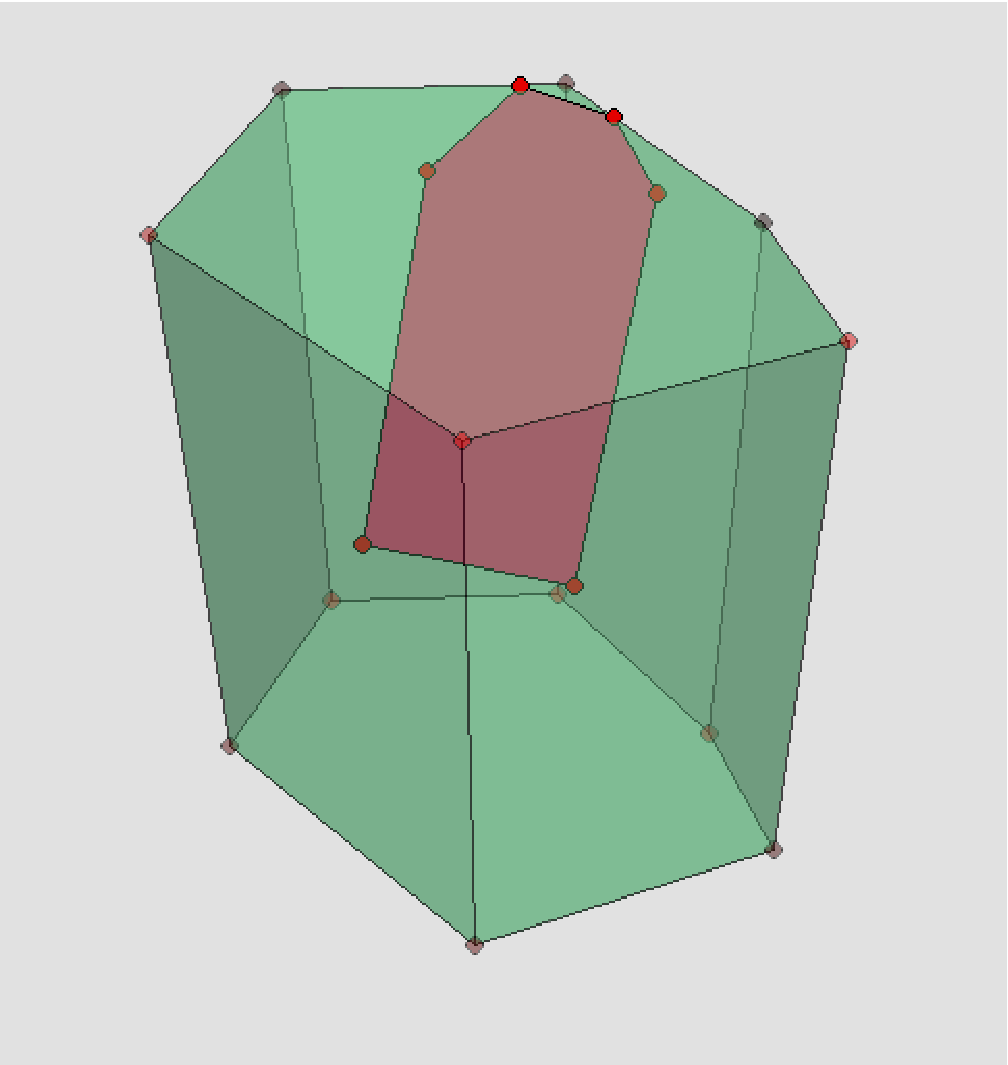}\qquad
\includegraphics*[scale=0.25]{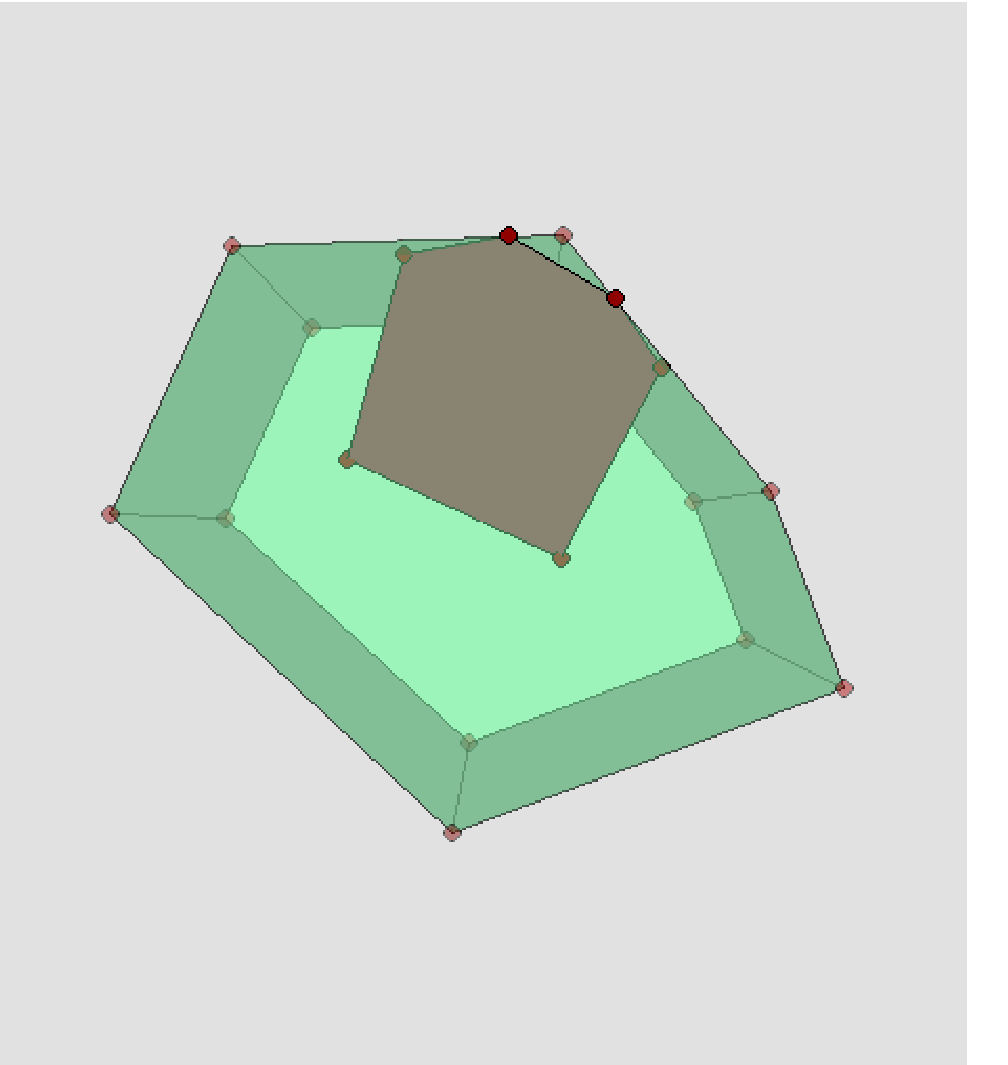}
\caption{Perspective and top views of a wall in a maximal cell ($n=3$)}
\end{figure}
\end{definition}

\begin{remarks}\label{wall remarks}
1)\ Analogous to slabs we have a function for each
$x\in\fop\setminus\Delta$, defined by
\[
f_{\fop,x}:=1+c_\fop z^{m_{\fop,x}}.
\]
The various $f_{\fop,x}$ are transformed to one another via parallel
transport inside $\sigma_\fop \setminus\Delta$. Note that the
$f_{\fop,x}$ are honest functions while the $f_{\fob,x}$ of a slab are
really sections of the line bundle with transition functions
$z^{m^{\rho_b}_{v' v}}$.\\
2)\ We have $\partial\fop=\base(\fop)\cup \sides(\fop)\cup
\topp(\fop)$ with
\begin{eqnarray*}
\base(\fop)&:=&\foq,\\
\sides(\fop)&:=&(\partial\foq-\RR_{\ge0}\overline{m_\fop})\cap \sigma_\fop,\\
\topp(\fop)&:=&\cl\big(\partial\fop\setminus
(\foq\cup\sides(\fop)\big).
\end{eqnarray*}
\end{remarks}
\medskip

In the following we will consider systems of slabs and walls
fulfilling certain additional conditions.

\begin{definition}\label{def:structure}
Let $\scrS=\scrS^b\cup \scrS^p$ with $\scrS^b$ and $\scrS^p$ locally
finite sets of slabs and walls, respectively. Define the
\emph{support} of $\scrS$ as
\[
|\scrS|:=\bigcup_{\fob\in \scrS} \fob.
\]
A \emph{chamber of $\scrS$} is the closure of a connected component of
$B\setminus |\scrS|$.  The set of chambers of $\scrS$ is denoted
$\Chambers(\scrS)$. Two chambers $\fou$, $\fou'$ are \emph{adjacent} if
$\dim\fou\cap\fou'=n-1$.

A \emph{structure} is a locally finite set of slabs and walls $\scrS$ along
with a polyhedral decomposition $\P_\scrS$ of $|\scrS|$, fulfilling
the following conditions.
\begin{enumerate}
\item[(i)]
The map associating to a slab $\fob\in\scrS$ its underlying polyhedral
subset of $B$ defines an injection from $\scrS^b$ to
$\P_\scrS^{[n-1]}$, and any $\rho\in\P^{[n-1]}$ is contained in
$|\scrS^b|$.
\item[(ii)]
Each chamber of $\scrS$ is convex and its interior is disjoint from any
wall.
\item[(iii)]
Any wall in $\scrS$ is a union of elements of $\P_\scrS$.
\item[(iv)]
Any $\sigma\in\P_\max$ contains only finitely many slabs or walls
in $\scrS$.
\qed
\end{enumerate}
\end{definition}

%===========================================================
\subsection{The gluing morphisms}
We assume a structure $\scrS$ to be given. Then for each
chamber $\fou\in\Chambers(\scrS)$ there exists a unique
$\sigma_\fou\in\P_\max$ with $\fou\subseteq \sigma_\fou$. Thus for each
pair $(g,\fou)$ with $(g:\omega\to\tau)\in\Hom(\P)$ and
$\tau\subseteq\sigma_\fou$ we have the rings $R_{g,\sigma_{\fou}}^k$,
$k\in\NN$. These are the rings whose spectra we wish to glue.
Technically this is done by a functor from a ``gluing
category'' to the category of log rings.

\begin{definition}\label{gluing category}
For a structure $\scrS$ define $\Glue(\scrS)$ as the category with
objects $(g,\fou)$ with $(g:\omega\to\tau)\in\Hom(\P)$,
$\fou\in\Chambers(\scrS)$ and $\omega\cap\fou\neq\emptyset$,
$\tau\subseteq \sigma_\fou$. (Then also $\tau\cap\fou\neq\emptyset$,
$\omega\subseteq \sigma_\fou$.) We call $\omega$ and $\tau$ the
\emph{domain} and \emph{target} of $(g,\fou)$, respectively. There is
a (unique) morphism
\[
(g:\omega\to\tau,\fou) \lra(g': \omega'\to\tau',\fou'),
\]
iff $\omega\subseteq\omega'$ and $\tau\supseteq\tau'$.
\end{definition}

Note that each morphism $\foe:(g,\fou)\to (g',\fou')$ in this category
decomposes into a sequence of morphisms of the following two basic
types:
\begin{enumerate}
\item[(I)]
$\omega\subseteq\omega'$, $\tau\supseteq \tau' $,
$\fou=\fou'$ (\emph{change of strata}).
\item[(II)]
$\omega=\omega'$, $\tau=\tau'$, $\dim \fou \cap \fou'=n-1$,
$\omega\cap\fou\cap\fou'\neq\emptyset$ (\emph{change of chamber}).
\end{enumerate}
For these two types of morphisms we now define a morphism of log rings
$R^k_{g:\omega\to\tau,\sigma_\fou}\to R^k_{g':
\omega'\to\tau',\sigma_{\fou'}}$ by specifying homomorphisms of
monoids $\beta:P_{\omega,\sigma_{\fou}}\to P_{\omega',\sigma_{\fou'}}$
and 
\[
\theta: \Lambda_{\sigma_\fou}\lra
\big(R^k_{g',\sigma_{\fou'}}\big)^\times,
\]
following the conventions of Remark~\ref{convention log
morphisms},(2). At this point our definition will still depend on
choices, but this dependence will disappear after imposing the
condition of consistency on $\scrS$ below
(Definition~\ref{consistency}). To put these log morphisms into
context recall from Remark~\ref{interpretation of rings} that if
$(g:\omega\to\tau,\fou) \in \Glue(\scrS)$ then $\Spec
R^k_{g,\sigma_\fou}$ is a $k$-th order thickening of an open subset
of $V_g$. Hence the target $\tau$ of $(g:\omega\to\tau,\fou)$
selects the toric stratum while its domain $\omega$ selects the
affine open subset to consider. Thus a morphism of Type~I in
$\Glue(\scrS)$ should map to a composition of the closed embedding
associated to $\tau'\to\tau$ composed with the open embedding
associated to $\omega\to\omega'$. Changing chambers (II) leads to
the application of log isomorphisms.

\begin{construction} (\emph{The basic gluing morphisms.})\label{gluing
morphisms}\\
I)\ (\emph{Change of strata.}) Let
\[
\foe:(g:\omega\to\tau,\fou)\to (g': \omega'\to\tau',\fou)
\]
be a morphism in $\Glue(\scrS)$ of Type~I and let
$a:\omega\to\omega'$. Denote by $s_{a,\sigma_\fou}:
\Lambda_{\sigma_\fou}\to \kk^\times$ the homomorphism defined by
$s_a\in\PM(\omega)$ for $\sigma_\fou$. Take $\beta:
P_{\omega,\sigma_\fou} \to P_{\omega',\sigma_\fou}$ to be the
canonical map and 
\[
\theta: \Lambda_{\sigma_\fou}\lra (R^k_{g',\sigma_{\fou}})^\times,
\quad
m\longmapsto s_{a,\sigma_\fou}^{-1}(m).
\]
Note that $\overline\theta$ is the canonical ring homomorphism from
Remark~\ref{natural ring homomorphisms}.

\mbox{}
\\[-2ex]
II)\ (\emph{Change of chambers.}) Let
\[
\foe:(g:\omega\to\tau,\fou)\to (g: \omega\to\tau,\fou')
\]
be a morphism in $\Glue(\scrS)$ of Type~II. Then either
$\fou\cap\fou'$ intersects the interior of a maximal cell (that is,
$\sigma_\fou=\sigma_{\fou'}$) or not. This leads to the following two
cases.
\begin{enumerate}
\item \underline{$\sigma_{\fou}=\sigma_{\fou'}$.}
Write $\sigma:=\sigma_\fou=\sigma_{\fou'}$ and $h:\omega\to\sigma$.
The intersection $\fou\cap\fou'$ is an $(n-1)$-dimensional convex
polyhedron not contained in the $(n-1)$-skeleton of $\P$. Since
$\omega\cap\fou\cap\fou'\neq \emptyset$ there exists
$\fov\in\P_\scrS^{[n-1]}$ with $\fov\subseteq\fou\cap \fou'$, $\omega
\cap \fov \neq\emptyset$. Then $\Int \fov\subseteq \Int\sigma$ and any
wall $\fop$ with $\fov\subseteq\fop$ has the property $\omega\cap
\fop\neq\emptyset$. Because $\Delta$ does not contain rational points
then even $\omega\cap (\fop\setminus\Delta)\neq\emptyset$. This is the
first place where we need the perturbation of $\Delta$.

Let $\fop_1,\ldots,\fop_r\in\scrS^p$ be the walls containing
$\fov$. Choose $x\in \Int(\fov)$ and let
$f_{\fop_i,x}=1+c_{\fop_i} z^{m_{\fop_i,x}}$ denote the function
associated to $\fop_i$ at $x$ according to Remark~\ref{wall
remarks},(1). Note that any non-zero exponent $m$ of this function
fulfills $\ord_\tau(m)\ge \ord_\sigma(m)>0$. Denote by $f_i$ the image
of $f_{\fop_i,x}$ in $R^k_{g,\sigma}$.

The tangent space of $\fou\cap\fou'$ defines an $(n-1)$-dimensional
rational subspace $T_{\fou\cap\fou'}\subseteq\Lambda_{\sigma,\RR}$. Let
$\pi: \Lambda_\sigma \to \ZZ$ be the epimorphism which contracts
$T_{\fou\cap\fou'}\cap \Lambda_\sigma$ and which is positive on
vectors pointing from $\fou$ to $\fou'$.

Now define $\beta=\id: P_{\omega,\sigma}\to P_{\omega,\sigma}$ and
\[
\theta=\theta(\fov): \Lambda_{\sigma}\lra (R^k_{g,\sigma})^\times,
\quad
m\longmapsto s_{ h,\sigma} \Big(
\prod_{i=1}^r f_i \Big)^{-\pi(m)}.
\]
This yields a log automorphism of $R^k_{g,\sigma}$ because $f_i=1 \mod
I_{g,\sigma}^{>0}$, and hence the associated automorphism of
$\kk[P_{\omega,\sigma}]/I_{g,\sigma}^{>k}$ changes the localizing
element only by an invertible function.

Without further assumptions our definition of $\theta$ depends on
the choice of $\fov$. We keep this dependence in mind for the time
being by adding $\fov$ to the notation at appropriate places.

\item \underline{$\sigma_{\fou}\neq\sigma_{\fou'}$.}
In this case $\fou\cap\fou'$ is contained in an $(n-1)$-cell
$\rho\in\P$. Let $\fov\in\P_\scrS^{[n-1]}$ be such that $\fov\subseteq
\fou\cap\fou'$ and $\omega\cap \fov\neq \emptyset$. Again we will
emphasize the dependence on $\fov$ in the notation. Since all
polyhedra are rational and $\Delta$ does not contain rational points
it holds $(\omega\cap \fov) \setminus\Delta\neq
\emptyset$. Let $x\in (\omega\cap \fov)\setminus\Delta$ and write
$e:v\to\omega$, $v:=v[x]$. Denote by $\fob$ the unique slab with
underlying polyhedral set $\fov$.

Now define $\beta$ by parallel transport through $v$
\[
\shAff(\check B,\ZZ)_{\check\sigma_\fou}\lra
\shAff(\check B,\ZZ)_v \lra
\shAff(\check B,\ZZ)_{\check\sigma_{\fou'}},
\]
and
\[
\theta=\theta(\fov): \Lambda_{\sigma_\fou}\lra
(R^k_{g,\sigma_{\fou'}})^\times,\quad
m\longmapsto \big(D(s_{  e},\rho,v)^{-1}
s_{ e}^{-1}(f_{\fob,x})\big)^{-\pi(m)}.
\]
Here $\pi:\Lambda_{\sigma_u}\to \ZZ$ is the epimorphism with kernel
$\Lambda_\rho$ which is positive on vectors pointing from $\fou$ to
$\fou'$, and $f_{\fob,x}$ is considered as an element of
$\kk[P_{\omega,\sigma_{\fou'}}]$ via a chart at $v$. Since $\beta$
respects orders (Lemma~\ref{order well-defined}) it identifies
$P_{\omega,\sigma_\fou} \subseteq\shAff(\check
B,\ZZ)_{\check\sigma_\fou}$ with $P_{\omega,\sigma_{\fou'}}$ and
$I_{g,\sigma_\fou}^{>k}$ with $I_{g,\sigma_{\fou'}}^{>k}$. Hence
$\beta$ and $\theta$ define a ring isomorphism
$\kk[P_{\omega,\sigma_\fou}]/ I^{>k}_{g,\sigma_\fou} \to
\kk[P_{\omega,\sigma_{\fou'}}]/ I^{>k}_{g,\sigma_{\fou'}}$. This
isomorphism respects the localizing elements as these only involve
monomials of order zero with respect to $\rho$ and hence are tangent
to $\rho$. This shows that $\theta$ indeed defines a log isomorphism.

Our construction also seems to depend on $x\in \fov\cap
(\omega\setminus\Delta)$. We now show that this is not the case. In
fact, a different choice $x'$ gives a vertex $v'\in\omega$ leading to 
$e':v'\to\omega$, $\theta': \Lambda_{\sigma_{\fou}} \lra
(R^k_{g,\sigma_{\fou'}})^\times$ instead of $e:v\to\omega$,
$\theta:\Lambda_{\sigma_\fou}\lra (R^k_{g,\sigma_{\fou'}})^\times$.
Parallel transport through $v'$ instead of $v$ yields
\[
\beta'(m)=\beta(m)+ \pi(\overline m)\cdot m^\rho_{v' v}.
\]
Then $(\beta,\theta,\overline \theta)$ and $(\beta', \theta',
\overline{\theta'})$ are equivalent via $\eta: \Lambda_{\sigma_\fou}
\to (P_{\omega,\sigma_{\fou'}})^\times$, $\eta(m)= \pi(\overline m)\cdot
m^\rho_{v'v}$:
\begin{eqnarray}\label{basic gluing: change of vertex}
\begin{aligned}
\theta'(m)& = \big(D(s_{  e'},\rho,v')^{-1} s_{ e'}^{-1}
(f_{\fob,x'}) \big)^{-\pi( m)}\\
&=\big(D(s_{  e},\rho,v)^{-1}
s_{ e}^{-1}(f_{\fob,x})z^{m^\rho_{v' v}}\big)^{-\pi( m)}\\
&=\theta(m)\cdot z^{-\pi( m)\cdot m^\rho_{v' v}} \ =\
\theta(m)\cdot (\alpha'\circ\eta(m))^{-1}.
\end{aligned}
\end{eqnarray}
Here we used the change of coordinates formula~(\ref{change of vertex
- slabs}) for $f_{\fob,x}$. This proves independence of the
choice of $x$, up to equivalence.
\qed
\end{enumerate}
\end{construction}

\begin{remark}\label{rem:discussion gluing morphisms}
1)~The consistency check~(\ref{basic gluing: change of vertex}) in
Construction~\ref{gluing morphisms},~II.2 forces the change of
coordinates formula~(\ref{change of vertex - slabs}) in the definition
of slabs.\\[1ex]
2)~While the log isomorphism in Construction~\ref{gluing
morphisms},~II.2 is only defined up to equivalence, a choice of vertex
$v\in\omega$ distinguishes a representative of the equivalence class.
In fact, quite generally for log isomorphisms, equivalent log
isomorphisms can be distinguished by the underlying homomorphism of
monoids, which in the case at hand is given by parallel transport
through~$v$.
\end{remark}

Changing chambers commutes with changing strata:

\begin{lemma}\label{changing strata commutes with changing chambers}
Assume that $g:\omega\to\tau$, $g':\omega'\to\tau'$ and
$\fou,\fou'\in \Chambers (\scrS)$ fulfill $\omega\subseteq\omega'$,
$\tau\supseteq\tau'$,  $\omega\cap\fou\cap\fou'\neq \emptyset$,
$\dim\fou\cap\fou'= n-1$ and $\tau\subseteq
\sigma_\fou\cap\sigma_{\fou'}$. For the morphisms
\begin{alignat*}{2}
\foe_1&:=\big((g,\fou)\lra (g,\fou')\big),\qquad
&\foe_2&:= \big((g,\fou')\lra (g',\fou')\big)\\
\fof_1&:=\big((g,\fou)\lra (g',\fou)\big),
&\fof_2&:= \big((g',\fou)\lra (g',\fou')\big)
\end{alignat*}
in $\Glue(\scrS)$ let $\theta(\foe_i)$, $\theta(\fof_i)$ be the basic
gluing morphisms from Construction~\ref{gluing morphisms}, where
$\theta(\foe_1)$ and $\theta(\fof_2)$ are computed using the same
$\fov\in \P_\scrS^{[n-1]}$. Then
\[
\theta(\foe_2)\circ\theta(\foe_1)= \theta(\fof_2)\circ\theta(\fof_1).
\]
\end{lemma}

\proof
Denote $a:\omega\to\omega'$. Let us first assume $\fou,\fou'$ are
contained in the same maximal cell $\sigma$, that is,
$\sigma_\fou=\sigma_{\fou'}$. With $h:\omega\to\sigma$,
$h':\omega'\to\sigma$ it holds $h=h'\circ a$ and hence
$s_{h,\sigma}=s_{h',\sigma}\cdot s_{a,\sigma}$. Then from
Construction~\ref{gluing morphisms},~I and~II.1 we obtain
\begin{align*}
\theta(\foe_2)\circ\theta(\foe_1)(m)
&= \overline{\theta(\foe_2)} \big(\theta(\foe_1)(m)\big)
\cdot\theta(\foe_2)(m)
=\overline{\theta(\foe_2)}\Big( s_{h,\sigma}\big({\textstyle \prod_i}
f_i \big)^{-\pi(m)}\Big)\cdot s_{a,\sigma}^{-1}(m)\\
&=s_{h',\sigma}\big({\textstyle \prod_i} 
f_i \big)^{-\pi(m)}\cdot s_{a,\sigma}^{-1}(m)
= \theta(\fof_2)(m)\cdot
\overline{\theta(\fof_2)} \big(s_{a,\sigma}^{-1}(m)\big)\\
&=\theta(\fof_2)(m)\cdot
\overline{\theta(\fof_2)} \big(\theta(\fof_1)(m)\big)
=\theta(\fof_2)\circ\theta(\fof_1)(m).
\end{align*} 

Otherwise there exists $\rho\in\P^{[n-1]}$ with
$\fou\cap\fou'\subseteq\rho$. Writing
$e:v\to\omega$, $e':v\to\omega'$ for $v=v[x]$ as in
Construction~\ref{gluing morphisms},~II.2, we have by
the definitions of $D(s_e,\rho,v)$ and $D(s_{e'},\rho,v)$
(Definition~\ref{D(lambda,v,rho)})
\[
\left(\frac{D(s_{e'},\rho,v)}{D(s_e,\rho,v)}\right)^{\pi(m)}
=\frac{s_{e',\sigma_\fou}(m)}{s_{e',\sigma_{\fou'}}(m)}
\cdot \frac{s_{e,\sigma_{\fou'}}(m)}{s_{e,\sigma_\fou}(m)}
= \frac{s_{a,\sigma_\fou}(m)}{s_{a,\sigma_{\fou'}}(m)}.
\]
Hence
\begin{align*}
\theta(\foe_2)\circ\theta(\foe_1)(m)
&=\overline{\theta(\foe_2)} \big(D(s_e,\rho,v)^{-1}
s_e^{-1}(f_{\fob,x})\big)^{-\pi(m)}\cdot \theta(\foe_2)(m)\\
&=s_{a,\sigma_{\fou'}}^{-1} \big(D(s_e,\rho,v)^{-1}
s_e^{-1}(f_{\fob,x})\big)^{-\pi(m)}\cdot s_{a,\sigma_{\fou'}}^{-1}(m)\\
&= \big(D(s_{e'},\rho,v)^{-1}
s_{e'}^{-1}(f_{\fob,x})\big)^{-\pi(m)}\cdot
s_{a,\sigma_{\fou}}^{-1}(m)=\theta(\fof_2)\circ\theta(\fof_1)(m),
\end{align*}
as desired.
\qed

%===========================================================
\subsection{Loops around joints and consistency}\label{par:joints}
To use compositions of basic gluing morphisms to define a functor from
$\Glue(\scrS)$ to $\LogRings$ requires a compatibility condition that
we now discuss.

\begin{definition}\label{joints}
A \emph{joint} $\foj$ of $\scrS$ is an $(n-2)$-cell of $\P_\scrS$ with
$\foj\not\subseteq\partial B$. The set of joints of $\scrS$ is denoted
$\Joints(\scrS)$.
\end{definition}

The minimal cell $\sigma_{\foj}\in\P$ containing a joint $\foj$ has
codimension at most two, and we speak of codimension zero, one and
two joints, respectively. For $v\in\sigma_\foj$ a vertex we
define the normal space of $\foj$ as
\[
\shQ_{\foj,\RR}^v:=\Lambda_{v,\RR}/ \Lambda_{\foj,\RR}.
\]
To denote the image of an object in $\shQ^v_{\foj,\RR}$ we use a
double bar. For example, if $\sigma\in\P_\max$ contains $\foj$ we have
the canonical map
\[
\shAff(\check B,\ZZ)_{\check\sigma}\lra \Lambda_v\lra
\shQ_{\foj,\RR}^v,\quad
m\longmapsto \doublebar m.
\]
Moreover, for a cell inside $B$ containing $\foj$, say $\tau\in\P$,
denote by $\doublebar\tau \subseteq \shQ_{\foj,\RR}^v$ the image of the
tangent wedge to $\tau$ along $\foj$. This is a convex cone, which is
strictly convex iff $\foj$ is contained in an $(n-2)$-dimensional
face of $\tau$.

Note that for any $\sigma\in\P_\max$ containing $v$ parallel transport
defines a canonical isomorphism
\[
\shQ_{\foj,\RR}^v= \Lambda_{\sigma,\RR}/\Lambda_{\foj,\RR},
\]
and local monodromy acts trivially on the right-hand side if
$\codim\sigma_\foj\in\{0,2\}$. Thus in this case $\shQ_{\foj,\RR}^v$
can be defined independently of $v$, while if $\sigma_\foj=\rho
\in\P^{[n-1]}$ we can only define the two half-planes separated by the
line $\doublebar \rho$ invariantly.

Now the $(n-1)$-cells of $\P_\scrS$ containing $\foj$ define a set of
distinct half-lines in $\shQ^v_{\foj,\RR}$. Let $\fov_1,\ldots,
\fov_l,\fov_{l+1}=\fov_1$ be a cyclic numbering of these cells induced
by an orientation on $\shQ^v_{\foj,\RR}$. Then by 
Definition~\ref{def:structure},(ii) for any $i$ there exists a unique
$\fou_i\in\Chambers(\scrS)$ with $\partial \doublebar\fou_i
=\doublebar \fov_i\cup \doublebar\fov_{i+1}$. 

Now let $(g:\omega\to\tau)\in \Hom(\P)$ and assume $\foj\cap
\omega\neq\emptyset$ and $\tau\subseteq\sigma_{\foj}$. Then for
each $i$ we have a morphism
\[
\foe:(g:\omega\to\tau,\fou_i)\to (g: \omega\to\tau,\fou_{i+1})
\]
in $\Glue(\scrS)$ changing chambers. By Construction~\ref{gluing
morphisms} we thus obtain the sequence of log isomorphisms
\[
\theta_i=\theta(\fov_i): \Lambda_{\sigma_{\fou_i}}\to
(R^k_{g,\sigma_{\fou_{i+1}}})^\times,
\]
from $R^k_{g,\sigma_{\fou_i}}$ to $R^k_{g,\sigma_{\fou_{i+1}}}$. Note
that the equivalence class of $\theta_i$ only depends on $g$,
$\fov_i$, $\foj$ and an orientation on $\shQ^v_{\foj,\RR}$.

\begin{definition}\label{consistency}
The structure $\scrS$ is \emph{consistent at the joint $\foj$ to
order $k$} if for any $(g:\omega\to\tau)\in \Hom(\P)$ with
$\foj\cap\omega\neq\emptyset$ and $\tau\subseteq\sigma_\foj$,
the composition
\begin{eqnarray}\label{theta_foj}
\theta_\foj^k:=\theta_l\circ\ldots\circ \theta_1:
\Lambda_{\sigma_{\fou_1}} \lra (R^k_{g,\sigma_{\fou_1}})^\times
\end{eqnarray}
equals $1$. A structure is \emph{consistent to order $k$} if it
is consistent to order $k$ at every joint. 
\end{definition}

\begin{remark}\label{rem:joints}
1) Note that relabelling $\fou_1, \ldots,\fou_l$ to $\fou_i,
\ldots,\fou_l,\fou_1,\ldots,\fou_{i-1}$ only leads to conjugation of
$\theta_\foj^k$ with an isomorphism $R^k_{g,\sigma_{\fou_1}}\to
R^k_{g,\sigma_{\fou_i}}$, and reversing the cyclic order produces the
inverse. Thus consistency around a joint really only depends on
$\scrS$ and $\foj$. Moreover, $\theta_\foj^k$ is well-defined once a
reference chamber $\fou_1$ and orientation of $\shQ_\foj^v$ are
chosen, for a vertex $v\in\sigma_\foj$.\\[1ex]
2) The notion of consistency implicitly depends on the choice of open
gluing data $s$.\\[1ex]
3) Consistency of a structure $\scrS$ does not depend on the choice of
polyhedral decomposition $\P_\scrS$. In fact, if $\scrS$ is consistent
at every joint of $\P_\scrS$ it is also consistent at every joint of
any refinement of $\P_\scrS$.\\[1ex]
4) Assume $\omega\subseteq\omega'\subseteq \tau'\subseteq\tau$ and we
compute~(\ref{theta_foj}) both for $g:\omega\to\tau$ and
$g':\omega'\to\tau'$, with the same sequence of chambers, resulting in
log automorphisms $\theta_\foj^k$ and ${\theta^k_\foj}'$. It then
follows from Lemma~\ref{changing strata commutes with changing
chambers} that ${\theta^k_\foj}'$ equals the composition of
$\theta^k_\foj$ with the canonical homomorphism
\[
(R^k_{g,\sigma_{\fou_1}})^\times\lra (R^k_{g',\sigma_{\fou_1}})^\times.
\]
5) By (4) it suffices to consider the case $\omega=\tau=\sigma_\foj$.
In fact, $R^k_{\omega\to \sigma_\foj,\sigma_{\fou_1}}  \to
R^k_{\id_{\sigma_\foj},\sigma_{\fou_1}}$ is the localization at
non-zero divisors, while $R^k_{\omega\to\sigma_\foj,\sigma_{\fou_1}}
\to R^k_{\omega\to\tau,\sigma_{\fou_1}}$ is a surjection followed by a
localization. Hence a log automorphism of $R^k_{\omega\to\tau,
\sigma_{\fou_1}}$ compatible with the identity on
$R^k_{\id_{\sigma_\foj},\sigma_{\fou_1}}$ is the identity.
\end{remark}

For consistent structures we can define a \emph{gluing functor}
\[
F^k_s: \Glue(\scrS)\lra \LogRings,
\]
mapping $(g,\fou)$ to $R^k_{g,\sigma_\fou}$ and morphisms of
types~I or II to the basic gluing morphisms of the respective types
defined in Construction~\ref{gluing morphisms}. In fact, the following
is true.

\begin{lemma}\label{consistency lemma}
Assume that the structure $\scrS$ is consistent to order $k$. Let
$\foe =\foe_r\circ \ldots\circ\foe_1 = \foe'_{r'}\circ\ldots
\circ\foe'_1$ be two decompositions of $\foe\in\Hom(\Glue(\scrS))$
into basic morphisms. Then if $\theta_r,\ldots,\theta_1$ and
$\theta'_{r'},\ldots,\theta'_1$ are the associated basic log morphisms
it holds
\[
\theta_r\circ\ldots\circ\theta_1=\theta'_{r'}\circ\ldots\circ\theta'_1.
\]
\end{lemma}

\proof The proof proceeds in three steps.
\smallskip

\noindent
\emph{\underline{Step 1.} Independence of choices.} The construction
of the morphism associated to a change of chambers
$(g:\omega\to\tau,\fou)\to (g:\omega\to\tau,\fou')$
(Construction~\ref{gluing morphisms},~II) required a choice of
$\fov\in \P_\scrS^{[n-1]}$. We claim independence of this choice for
consistent structures. We use the notations from
Construction~\ref{gluing morphisms}. If $\fov'\in\P_\scrS^{[n-1]}$ is
another $(n-1)$-cell with $\fov' \subseteq\fou\cap\fou'$,
$\omega\cap \fov'\neq\emptyset$ and adjacent to $\fov$ (i.e.\ $\dim
\fov\cap\fov' =n-2$) then $\foj:=\fov\cap\fov'$ is a joint and the
log isomorphisms $\theta$ and $\theta'$ constructed via $\fov$ and
$\fov'$, respectively, differ by a composition
$\theta_l\circ\ldots\circ \theta_1$ associated to a loop around $\foj$
(Definition~\ref{consistency}). Hence they are equal. The general case
follows since any two $(n-1)$-cells of $\P_\scrS$ contained in
$\fou\cap\fou'$ and intersecting $\omega$ can be connected by a
sequence of adjacent cells with the same properties.
\smallskip

\noindent
\emph{\underline{Step 2.} Reduction to the case $g=g'$.}
By construction the basic gluing morphisms changing strata are
compatible with compositions in $\Hom(\P)$: If
$\foe_i$ and $\foe_{i+1}$ both change strata for a fixed chamber
$\fou=\fou_i=\fou_{i+1}=\fou_{i+2}$ then
\[
F_s^k(\foe_i)\circ F_s^k(\foe_{i+1})=F_s^k(\fof),
\]
where $\fof:(g_i,\fou)\to (g_{i+2},\fou)$. Since by
Lemma~\ref{changing strata commutes with changing chambers} changing
chambers commutes with changing strata we can thus assume that only
$\foe_1$ changes strata. Since the analogous factorization for $\foe'$
leads to the same change of strata this reduces the claim to a
sequence of changes of chambers.
\smallskip

\noindent
\emph{\underline{Step 3.}  Spaces of chambers.}
Given $g:\omega\to\tau$ look at all chambers $\fou$ such that
$(g,\fou)\in\Glue(\scrS)$:
\[
A:=\big\{\fou\in\Chambers(\scrS) \,\big|\,\omega\cap\fou\neq\emptyset,\,
\tau\subseteq \sigma_\fou  \big\}.
\]
Define $\Sigma$ to be the abstract two-dimensional cell complex with
$A$ as set of vertices, edges connecting adjacent $\fou,\fou'\in A$
and a disk glued into any $1$-cycle of chambers $\fou_1,\ldots,\fou_l$
forming a loop around a joint. An edge with vertices $\fou,\fou'$
defines a change of chamber isomorphism (of either type)
$R^k_{g,\sigma_\fou} \to R^k_{g,\sigma_{\fou'}}$. Consistency says
that the composition of these isomorphisms following the boundary of a
$2$-cell is the identity. Thus we obtain the desired independence of
the sequence of adjacent chambers connecting $\fou_1$ with $\fou_r$
once we know $H^1(\Sigma,\ZZ)=0$. This follows from the following
Lemma.

\begin{lemma}\label{pi(Sigma)=0}
$\pi_1(\Sigma)=0$.
\end{lemma}

\proof
Denote by
\[
V:=\bigcup_{\sigma\in\P_\max,\tau\subseteq \sigma} \sigma
\]
the closed star of $\tau$ with respect to $\P$. Then for $\fou\in
\Chambers(\scrS)$ the condition $\tau\subseteq \sigma_\fou$ is
equivalent to $\fou\subseteq V$, and such chambers define a
decomposition of $V$ into closed polyhedra. This is generally not a
proper polyhedral decomposition because the intersection of two chambers
needs not be a face of either chamber. We can however refine the
decomposition of $V$ into chambers into an honest polyhedral
decomposition $\P_V$, for example by replacing each slab or wall by
the hyperplane containing it. Since $\omega$ is topologically a ball
the cells of $\P_V$ intersecting $\omega$ form a polyhedral
decomposition $\P'$ of an $n$-cell, and so does the combinatorial dual
decomposition $\check \P'$. Thus the two-skeleton of $\check \P'$ is
simply connected. Now $\Sigma$ is obtained from $\check \P'$ by
contracting all edges corresponding to adjacent $n$-cells of $\P$
lying in the same chamber. Thus also $\Sigma$ is simply connected.
\qed
\smallskip

This finishes the proof of Lemma~\ref{consistency lemma}.
\qed

%===========================================================
\subsection{Construction of finite order deformation}
Given a structure $\scrS$ that is consistent to order $k$ we are
now able to construct the desired deformation $(X_k,D_k)$ of $(X,D)$
over $\Spec\big(\kk[t]/ (t^{k+1})\big)$, by taking, in a certain
sense, the colimit of $\Spec R^k_{g,\sigma_\fou}$ over $(g,\fou)\in
\Glue(\scrS)$. Taken literally this would lead to a non-separated
scheme because we have removed closed subsets (the zero loci of the
localizing elements $f_{g,\sigma_\fou}$) from lower-dimensional
strata. Instead we proceed as follows. Denote by
\[
\overline{F^k_s}: \Glue(\scrS)\lra \Rings
\]
the composition of $F^k_s$ with the forgetful functor $\LogRings\to
\Rings$. For $(g:\omega\to\tau,\fou)\in \Glue(\scrS)$ the underlying
topological space of $\Spec \overline{F^k_s}(g,\fou)= \Spec
R^k_{g,\sigma_\fou}$ is, according to Remark~\ref{interpretation of
rings}, canonically an open subset of $V_g\subseteq V(\omega)$. Denote
by
\[
i(g,\fou): \big|\Spec R^k_{g,\sigma_\fou}\big|
\lra \big|V(\omega) \big|
\]
the inclusion of the underlying topological spaces. Then
\[
\O_k(g,\fou):= i(g,\fou)_*\O_{\Spec R^k_{g,\sigma_\fou}}
\]
defines a sheaf of $\kk[t]$-algebras on $V(\omega)$, and if $\foe:
(g:\omega\to \tau,\fou)\to (g':\omega'\to\tau',\fou')$ is a morphism
in $\Glue(\scrS)$ and
\[
\Phi_{\omega\omega'}(s): V(\omega')\to V(\omega)
\]
is the open embedding defined from the composition of
$p|_{V(\omega')}$ with the inverse of $p|_{V(\omega)}$, then
$F^k_s(\foe)$ defines a homomorphism of sheaves of $\kk[t]$-algebras
\[
\O_k(g,\fou)\lra \big(\Phi_{\omega\omega'}(s)\big)_*\O_k(g',\fou').
\]
In fact, recall from the discussion following Definition~\ref{open
gluing data} that $\Phi_{\omega\omega'}(s)$ is given by the canonical
embedding $\Spec\kk[{\omega'}^{-1}\Sigma_v]\to \Spec \kk[\omega^{-1}
\Sigma_v]$ for some $v\in\omega$, composed with the automorphism
$s_h^{-1}$ of $\Spec \kk[{\omega'}^{-1}\Sigma]$ coming from open
gluing data, where $h:\omega\to\omega'$. By Remark~\ref{natural ring
homomorphisms} this is compatible with the reduction modulo $t$ of
$\Spec \overline{F^k_s}(\foe)$. In particular, we have
\[
|\Phi_{\omega\omega'}(s)|\circ i(g',\fou')= i(g,\fou)\circ |\Spec
\overline{F^k_s}(\foe)|.
\]
Thus
\begin{eqnarray}\label{ringed space functor}
(g:\omega\to\tau,\fou)\longmapsto \big(|V(\omega)|, \O_k(g,\fou)\big)
\end{eqnarray}
defines a contravariant functor from $\Glue(\scrS)$ to the category of
ringed spaces. The aim of this subsection is to show that the colimit
of this functor defines a deformation of $X$ over $\kk[t]/(t^{k+1})$
of the desired form. This will be achieved in Proposition~\ref{X_k}. 

We first construct deformations $V^k(\omega)$ of the standard affine
sets $V(\omega)$ that cover $X$ and then glue by open embeddings.
Thus for the time being keep $\omega$ fixed and consider
\[
V^k(\omega):=\big(\big|V(\omega)\big|,\liminv \O_k(g,\fou)\big),
\]
where the inverse limit runs over all $(g,\fou)\in\Glue(\scrS)$ with
domain $\omega$.

Let $x\in\Int(\omega)\setminus\Delta$. For each $g:\omega\to\tau$ we
can choose $\fou_g\in\Chambers(\scrS)$ with $x\in\fou_g$ and such that
$(g,\fou_g)\in\Glue(\scrS)$. Then by Lemma~\ref{consistency lemma}, if
$(g,\fou)\in \Glue(\scrS)$ with domain $\omega$ there is a unique
isomorphism $R^k_{g,\sigma_{\fou_g}} \to R^k_{g,\sigma_\fou}$, by a
composition of changes of chambers. This shows that we may replace the
system of rings $R^k_{g,\sigma_\fou}$ and sheaves $\O_k(g,\fou)$,
for $(g,\fou)\in \Glue(\scrS)$ with domain $\omega$, with the system
\[
R_\tau:=R^k_{g:\omega\to\tau,\sigma_{\fou_g}}\quad \text{and}\quad
\O_k(\tau):=\O_k(g:\omega\to\tau,\fou_\tau),
\]
respectively, where now $\tau$ runs over all cells containing $\omega$
and $\fou_\tau:=\fou_{\omega\to\tau}$. In doing this keep in mind that
the homomorphisms
\[
\psi_{\tau'\tau}: R_\tau\to R_{\tau'}
\]
thus obtained now involve compositions of changes of chambers.

Eventually we will argue inductively, noting that reduction modulo
$t^{l+1}$ defines rings $R^l_\tau$, sheaves $\O_l(\tau)$ and ringed
spaces $V^l(\omega)$ for $l<k$. These are related by the complex
\begin{eqnarray}\label{O_{V^k}-sequence}
0\lra \O_{V(\omega)}\stackrel{\cdot t^l}{\lra} \O_{V^l(\omega)}
\lra \O_{V^{l-1}(\omega)}\lra 0,
\end{eqnarray}
of sheaves of $\kk[t]$-algebras on $V(\omega)$, and we have to show
this complex is exact. Then $V^l(\omega)$ is a flat lifting of
$V^{l-1}(\omega)$ from $\kk[t]/(t^l)$ to $\kk[t]/(t^{l+1})$. It then
follows inductively that $V^k(\omega)$ is affine, of finite type over
$\kk$ and flat over $\kk[t]/(t^{k+1})$. The overall strategy for
showing exactness is to write down an isomorphism with standard local
models outside a codimension three locus contained in $Z_\omega$, the
preimage of $Z$ under $p:V(\omega)\to X$, and then to extend. This
will also show log smoothness away from $Z$ in the limit $k\to\infty$.

The first step is the construction of a chart for $V^k(\omega)$ away
from $Z_\omega$. We write $P=P_{\omega,\sigma_{\fou_{\id_\omega}}}
=P_x$. Recall that for any $\tau\supseteq\omega$ we have the
homomorphism
\[
\alpha_\tau:P \lra R_\tau
\]
endowing $R_\tau$ with the structure of a log ring. Whenever
$\tau\supseteq\tau'\supseteq \omega$ there is a log morphism
\[
\theta_{\tau'\tau}: \Lambda_x\simeq
\Lambda_{\sigma_{\fou_\tau}}
\lra (R_{\tau'})^\times
\]
such that $\overline\theta_{\tau'\tau} =\psi_{\tau'\tau}$.

\begin{proposition}\label{charts for V^k(omega)}
Let $p\in|V(\omega)|\setminus Z_{\omega}$. Since $|V(\omega)|$
coincides with the topological space underlying $\Spec\kk[P]/(t^{k+1})$,
there is a prime ideal $\fop\subseteq \kk[P]/(t^{k+1})$ corresponding
to $p$. Furthermore, for $g:\omega\rightarrow\tau$, $|\Spec R_{\tau}|$
is identified with $|V_g|\setminus Z_{\omega}$. Thus if $p\in|V_g|$,
there is a prime ideal $\fop_\tau\subseteq R_{\tau}$ corresponding
to $p$. Then there is an isomorphism
\[
\psi: \big(\kk[P]/(t^{k+1})\big)_{\fop} \lra
\liminv (R_\tau)_{\fop_\tau},
\]
where the inverse limit is over all $g:\omega\rightarrow\tau$ with
$p\in|V_g|$, and the maps of the inverse system are
the localizations of the maps $\psi_{\tau'\tau}$.
\end{proposition}

\proof
For any $\tau'\subseteq\tau$ with $\omega\subseteq\tau'$ it holds
$(\omega\to\tau',\fou_\tau)\in \Glue(\scrS)$. Note that
$\psi_{\tau'\tau}$ is the composition of the canonical homomorphism
\[
R_\tau=R^k_{\omega\to\tau,\sigma_{\fou_\tau}}\lra
R^k_{\omega\to\tau',\sigma_{\fou_\tau}},
\]
changing strata, with a change of chambers isomorphism
\[
R^k_{\omega\to\tau',\sigma_{\fou_\tau}}
\lra R^k_{\omega\to\tau',\sigma_{\fou_{\tau'}}}=R_{\tau'}.
\]
In particular, $\psi_{\tau'\tau}$ needs not be surjective because the
change of strata homomorphism may involve a localization. However, after
localizing at the ideals $\fop_{\tau'}$ and $\fop_\tau$ respectively,
this map becomes surjective.

Note also that if $\omega\subseteq \tau''\subseteq
\tau'\subseteq \tau$ then, by consistency,
\[
\theta_{\tau''\tau}= \theta_{\tau''\tau'}\circ\theta_{\tau'\tau}
=\theta_{\tau''\tau'}\cdot(\psi_{\tau''\tau'}\circ\theta_{\tau'\tau}).
\]
Thus $(\theta_{\tau'\tau} )_{\tau'\subsetneq \tau}$ is a barycentric
$1$-cocycle for the system of groups $\Hom(\Lambda_x,
(R_\tau)_{\fop_\tau}^\times)$, as considered in \cite{logmirror}, A.1.
Since the homomorphisms $(R_\tau)_{\fop_\tau}^\times \to
(R_{\tau'})_{\fop_{\tau'}}^\times$ 
are surjective it is straightforward to check
the exactness criterion~$(*)$ in \cite{logmirror}, Proposition~A.1.
Hence there exist $\theta_\tau: \Lambda_x\to (R_\tau
)_{\fop_\tau}^\times$ such that for any $\tau'\subsetneq\tau$
\[
\theta_{\tau'\tau}=\theta_{\tau'}/ (\psi_{\tau'\tau}\circ
\theta_\tau).
\]
Then $\theta_{\tau}\cdot\alpha_\tau$ defines a compatible system of
homomorphisms
\[
\psi_\tau:\big(\kk[P]/(t^{k+1})\big)_{\fop}\to ( R_\tau)_{\fop_\tau},
\]
hence the desired map $\psi$. Note that $\psi_\tau$ is the composition
of the canonical quotient
\[
\big(\kk[P]/(t^{k+1})\big)_{\fop} \lra
\big(\kk[P]/I_{\omega\to\tau}^{>k}\big)_{\fop}
\]
and an isomorphism of $\big(\kk[P]/I_{\omega\to\tau}^{>k}\big)_{\fop}$ with
$( R_\tau)_{\fop_\tau}$. Now it is easy to check that
$\kk[P]/(t^{k+1})$ is the inverse limit of the rings
$\kk[P]/I_{\omega\to\tau}^{>k}$ with the canonical quotient
homomorphisms between them, and an analogous statement holds for the
localizations. The $\psi_\tau$ induce an isomorphism between this
inverse system and $( R_\tau)_{\fop_\tau}$. This shows that
$\psi$ is indeed an isomorphism.
\qed
\smallskip

\begin{corollary}\label{V^K away from Z}
For any $l\le k$, Sequence~(\ref{O_{V^k}-sequence}) is exact on
$V(\omega)\setminus Z_\omega$.
\qed
\end{corollary}

Following Example~\ref{2d example, part II} we can also check directly
exactness of~(\ref{O_{V^k}-sequence}) at general points of $Z_\omega$.

\begin{lemma}\label{model at general points of Z}
For any $\rho\in\P^{[n-1]}$ containing $\omega$,
Sequence~(\ref{O_{V^k}-sequence}) is also exact at all points of the
maximal torus $\Spec\kk[\Lambda_\rho]\subseteq V_{\omega\to\rho}$.
\end{lemma}

\proof
Let $\sigma^+$, $\sigma^-$ be the maximal cells with
$\rho=\sigma^+\cap\sigma^-$. Only three of the sheaves $\O_k(\tau)$
have support on $\Spec\kk[\Lambda_\rho]$, namely $\O_k(\rho)$ and
$\O_k(\sigma^\pm)$. Moreover, by choosing $\fou_{\sigma^\pm}$ as
adjacent chambers and $\fou_{\rho}= \fou_{\sigma^+}$ we can assume that
$R_{\sigma^+}\to R_\rho$ is the canonical homomorphism while
$R_{\sigma^-}\to R_\rho$ is the canonical homomorphism composed with
a change of chamber automorphism
\[
\psi: R_\rho\lra R_\rho,\quad
z^m\longmapsto f_\rho^{-\pi(\overline m)}\cdot z^m
\]
as defined in Construction~\ref{gluing morphisms},~II.2. Recall that
$\pi:\Lambda_{\sigma^-}\to \ZZ$ maps to $1$ a generator of
$\Lambda_{\sigma^-} /\Lambda_\rho$ pointing from $\sigma^-$ to
$\sigma^+$, and $f_\rho$ has zero locus
$Z_\omega\cap |\Spec\kk[\Lambda_\rho]|\subseteq V(\omega)$. Denote by
$F_\rho\subseteq P$ the face corresponding to $\rho$. Since $\rho$
is codimension one, $P+F_\rho^\gp\simeq \Lambda_\rho\oplus S_e$ for
some $e\ge1$ with $S_e\subseteq\ZZ^2$ the monoid generated by $(1,0)$,
$(e,-1)$, $(0,1)$. Denote by $R_+$, $R_-$ and $R_\cap$ the
localizations at the multiplicative system $\{z^m\}_{m\in F_\rho}$ of
$R_{\sigma^+}$, $R_{\sigma^-}$ and $R_\rho$, respectively. Explicitly,
we may write
\begin{align*}
R_-&=\kk[\Lambda_\rho][x,y,t]/\big(xy-t^e,y^\beta t^\gamma\,\big|\, \beta
e+\gamma\ge k+1\big)\\
R_+&=\kk[\Lambda_\rho][x,y,t]/\big(xy-t^e,x^\alpha t^\gamma\,\big|\,
\alpha e+\gamma\ge k+1\big)\\
R_\cap&=\big(\kk[\Lambda_\rho][x,y,t]/\big(xy-t^e,
x^\alpha y^\beta t^\gamma
\,\big|\,
\max\{\alpha,\beta\}e+\gamma\ge k+1\big)
\big)_{f_\rho},
\end{align*}
and $R_+\to R_\cap$ is the canonical quotient followed by localization
at $f_\rho$, while $R_-\to R_\cap$ is the homomorphism of
$\kk[\Lambda_\rho][t]/(t^{k+1})$-algebras with
\[
x\longmapsto f_\rho x,\quad y\longmapsto f_\rho^{-1} y.
\]
We claim that
\[
R_\cup:=\kk[\Lambda_\rho][X,Y,t]/(XY-f_\rho t^e,
t^{k+1})\lra R_-\times_{R_\cap} R_+,\quad
X\longmapsto (x,f_\rho x),\ Y\longmapsto (f_\rho y, y)
\]
is an isomorphism of $\kk[\Lambda_\rho][t]$-algebras. In fact, the
rings $R_-$, $R_+$ are generated by $1$, $x^i$, $y^j$, $i,j>0$ as
$\kk[\Lambda_\rho][t]/ (t^{k+1})$-modules, and the same monomials
generate $R_\cap$ as $\kk[\Lambda_\rho]_{f_\rho}
[t]/(t^{k+1})$-module. Moreover, the
$\kk[\Lambda_\rho][t]$-submodules of $R_-$ ($R_+$) generated by
$x^i$, $i\ge 0$ ($y^j$, $j\ge0$) is a free direct summand. Thus if
$g_\pm\in R_\pm$ we may write uniquely $g_-=\sum_{i\ge0} a_i x^i+
h_-(y,t)$, $g_+= \sum_{j\ge0} b_j y^j +h_+(x,t)$ with
$h_\pm(0,t)=0$. Thus $(g_-,g_+)\in R_-\times_{R_\cap} R_+$ iff
\[
a_0 = b_0,\quad
h_-(y,t)=\sum_{j> 0} b_j f_\rho^j y^j,\quad
h_+(x,t)=\sum_{i> 0} a_i f_\rho^i x^i,
\]
as elements of $R_\cap$. If this is the case then $(g_-,g_+)$ is the
image of $\sum_{i\ge 0} a_i X^i+\sum_{j>0} b_j Y^j \in R_\cup$. This
shows surjectivity. Injectivity follows along the same lines by noting
that $R_\cup$ is a free $\kk[\Lambda_\rho][t]/(t^{k+1})$-module
with basis $X^i$, $Y^j$, $i\ge0$, $j>0$.

To complete the proof it remains to observe that
$(t^l)\subseteq \kk[\Lambda_\rho][X,Y,t]/(XY-f_\rho t^e, t^{l+1})$ is a
free $\kk[\Lambda_\rho][X,Y]/(XY)$-module, for $0<l\le k$. In fact, by the
same argument as before each element of $(t^l)$ can be uniquely
written as $t^l\big(\sum_{i\ge 0} a_i X^i+ \sum_{j>0} b_j Y^j\big)$.
\qed
\smallskip

We now know that Sequence~(\ref{O_{V^k}-sequence}) is exact on
$V(\omega)\setminus Z'_\omega$, where $Z'_\omega$ is the intersection
of $Z_\omega$ with the union of the codimension two strata of
$V(\omega)$. To extend across $Z'_\omega$ the crucial technical result
is the following.

\begin{lemma}\label{depth lemma}
For $\omega\subseteq\tau\subseteq \sigma$, $\sigma\in\P_\max$, let
$Y=\Spec(\kk[P]/ I^{>k}_{\omega\to\tau,\sigma})$, and let $p\in Y$ be a
scheme-theoretic point contained in a proper toric stratum of $Y$, but
$p$ not the generic point of a toric stratum. Then $\depth\O_{Y,p}\ge
2$. 
\end{lemma}

\proof
There is a $\tau'$ with $\omega\subseteq\tau'\subsetneq\tau$ such that
$V_{\omega\to\tau'}\subseteq Y^\red$ is the smallest toric stratum
containing $p$. Then $p\in V_{\tau'\to\tau'}$, the open torus orbit of
$V_{\omega\to\tau'}$. So $p$ is a point in the open subscheme $U=\Spec
(\kk[P_{\tau',\sigma}]/ I^{>k}_{\tau'\to\tau,\sigma})$ of $Y$.
Note that $P_{\tau',\sigma}$ splits non-canonically, as
\[
P_{\tau',\sigma}=P'\times\ZZ^r,
\]
where $r=\dim\tau'$ and $P'$ is a sharp monoid (i.e. containing
no invertible elements other than $0$). In particular, there is
a monomial ideal $I'\subseteq\kk[P']$ such that
\[
U\simeq\Spec(\kk[P']/I'\otimes_{\kk}\kk[\ZZ^r]).
\]
Let $\fop\subseteq \kk[P']/I'\otimes_{\kk}\kk[\ZZ^r]$ be the prime
ideal corresponding to $p$. We need to find
a regular sequence $a_1,a_2\in\fop$ of length two.
Take $a_1=1\otimes f$, where $f\in\fop\cap\kk[\ZZ^r]$
is a non-zero element, which exists since $p$ is not the generic point
of $V_{\tau'\to\tau}$. Take $a_2=z^m\otimes 1$, where $m\in P'$
is an element in the interior of the face of $P'$ corresponding to
$\tau$. It is then easy to see that $a_1,a_2$ form a regular sequence;
indeed, view $\kk[P']/I'$ and $\kk[\ZZ^r]/(f)$ as $\kk$-vector spaces;
then tensoring the injective map
\[
\kk[\ZZ^r]\mapright{\cdot f}\kk[\ZZ^r]
\]
with $\kk[P']/I'$ shows $a_1$ is not a zero-divisor, and tensoring
the injective map
\[
\kk[P']/I'\mapright {\cdot z^m} \kk[P']/I'
\]
with $\kk[\ZZ^r]/(f)$ shows $a_2$ is not a zero-divisor in
$(\kk[P']/I'\otimes_{\kk}\kk[\ZZ^r])/(a_1)$.
\qed

\begin{remark}
The assumption in Lemma~\ref{depth lemma} that $p$ is not the
generic point of a toric stratum is necessary. In particular, unlike
$V_g$ the thickening $\Spec R^k_{g,\sigma}$ needs not fulfill
Serre's condition $S_2$. As an example, take $\omega$ a point in a
two-dimensional $B$, a polarization $\varphi_\omega$ with Newton
polyhedron the unit square, $\tau\supseteq\omega$ a maximal cell and
$p$ the zero-dimensional toric stratum $V_{\omega\to\omega}\subseteq
\Spec R_\tau$. Then in appropriate coordinates
\[
\kk[P]=\kk[x,y,z,w]/(xy-zw),
\quad I^{>k}_{\omega\to\tau,\tau}=(x,z)^{k+1},
\]
and $y^{-1}z=w^{-1}x$ is a regular function on $\Spec
(R_\tau)\setminus \{p\}$ that does not extend. Such an
extension would be possible by the depth argument if $\depth \O_{\Spec
R_\tau,p}\ge 2$.
\end{remark}

\begin{lemma}\label{extension across Z_omega}
$j_*\O_{V^k(\omega)\setminus Z_\omega}=\O_{V^k(\omega)}$.
\end{lemma}

\begin{proof}
If $U\subseteq V(\omega)$ is an open set then by definition
\begin{eqnarray*}
\O_{V^k(\omega)}(U)&=&\liminv \O_{\Spec R_\tau}(U)\\
j_*\O_{V^k(\omega)\setminus Z_\omega}(U)
&=& \liminv \O_{\Spec R_\tau}(U\setminus Z_\omega).
\end{eqnarray*}
Now for any maximal cell $\sigma\supseteq \omega$, Lemma~\ref{depth
lemma} shows with the usual depth argument that the restriction map
$\O_{\Spec R_\sigma}(U)\to \O_{\Spec R_\sigma}(U\setminus Z_\omega)$
is a bijection. Hence the canonical map $\O_{V^k(\omega)}\to
j_*\O_{V^k(\omega)\setminus Z_\omega}$ is an isomorphism since
membership of a tuple $(f_\sigma)$, $f_\sigma\in \O_{\Spec
R_\sigma}(U)$, in $\liminv \O_{\Spec R_\tau}(U)$ can be checked on
$U\setminus Z_\omega$.
\end{proof}

We are now in position to conclude exactness of
(\ref{O_{V^k}-sequence}) at all points.

\begin{proposition}\label{exactness of O_{V^k}-sequence}
For any $l\le k$, Sequence~(\ref{O_{V^k}-sequence}) is exact.
\end{proposition}

\proof
We know exactness of~(\ref{O_{V^k}-sequence}) on $V(\omega) \setminus
Z'_\omega$, where $Z'_\omega$ is the intersection of $Z_\omega$ with
the union of the codimension two strata of $V(\omega)$. Pushing
forward by $j': V(\omega)\setminus Z'_\omega \to V(\omega)$ yields the
exact sequence
\[
0\lra j'_*\O_{V(\omega)\setminus Z'_\omega} \stackrel{\cdot t^l}{\lra}
j'_*\O_{V^l(\omega)\setminus Z'_\omega}
\lra j'_*\O_{V^{l-1}(\omega)\setminus Z'_\omega}
\lra R^1j'_*\O_{V(\omega)\setminus Z'_\omega}
\]
The term on the right vanishes since $V(\omega)$ is Cohen-Macaulay and
$\codim Z'_\omega\ge 3$.
\qed
\medskip

We have now established that $V^k(\omega)$ is a flat deformation of
$V(\omega)$ over $\Spec \kk[t]/(t^{k+1})$. In particular,
$V^k(\omega)$ is an affine scheme of finite type over $\kk$. Moreover,
whenever $\omega\subseteq \omega'$ the functor~(\ref{ringed space
functor}) induces a map of schemes
\[
\Phi^k_{\omega\omega'}(s): V^k(\omega')\lra V^k(\omega).
\]
These morphisms are compatible with sequences
$\omega\subseteq\omega'\subseteq \omega''$:
\[
\Phi^k_{\omega\omega''}(s)=
\Phi^k_{\omega\omega'}(s)\circ\Phi^k_{\omega'\omega''}(s), 
\]
hence they define a functor from $\P$ to the category of schemes.
Define $X_k$ as the colimit of this functor.

\begin{proposition}\label{X_k}
The maps $\Phi^k_{\omega\omega'}(s)$ are open embeddings. In
particular, $X_k$ is a scheme locally of finite type and flat over
$\kk[t]/(t^{k+1})$. 
\end{proposition}

\proof
Recall that in studying $\O_{V^k(\omega)}$ we reduced the inverse
limit over $(\omega\to\tau,\fou)\in\Glue(\scrS)$ to an inverse limit
over cells $\tau$ containing $\omega$ by choosing one chamber
$\fou_\tau$ for each $\tau$. Now using the same choice of $\fou_\tau$ for
$\omega$ and $\omega'$, for $\tau\supseteq\omega'$, we see that
$V^k(\omega')\to V^k(\omega)$ is defined by
\begin{eqnarray}\label{map of inverse limits}
\liminv_{\tau\supseteq\omega'}
\Phi_{\omega\omega'}(s)^{-1}\O^k(\omega\to\tau,\fou_\tau)
\lra
\liminv_{\tau\supseteq\omega'} \O^k(\omega'\to\tau,\fou_\tau).
\end{eqnarray}
Note that on the left-hand side we dropped the sheaves
$\O^k(\omega\to\tau,\fou_\tau)$ for cells $\tau$ containing $\omega$
but not containing $\omega'$ because they are supported away from
the image of $\Phi_{\omega\omega'}(s)$. Now on the ring level, the
$\tau$-component of~(\ref{map of inverse limits}) is the identity of
$R^k_{\omega'\to\tau,\sigma_{\fou_\tau}}$. Thus~(\ref{map of inverse
limits}) is an isomorphism.

It remains to remark that for vertices $v,v'\in\P$ the open sets
$p\big(V(v)\big)$ and $p\big(V(v')\big)$ intersect in $p \big(
V(\omega)\big)$ for $\omega$ the minimal cell containing $v,v'$. Thus
$\O_{X_k}$, as a sheaf on $|X|$, is isomorphic on $p(V(\omega))$ to
$p_*\O_{V^k(\omega)}$. Hence $X_k$ is a scheme with the claimed
properties. At this point we use crucially that the cells of $\P$ do
not self-intersect; otherwise we would end up with an algebraic space
here.
\qed

\begin{remark}\label{log structure on X_k}
Proposition~\ref{charts for V^k(omega)} also endows $X_k$ with an
abstract log structure, together with a log smooth morphism to
$\Spec\kk[t]/(t^{k+1})$ with the log structure generated by $\NN\to
\kk[t]/(t^{k+1})$, $1\mapsto t$. While this is not relevant to this paper
it is important in order by order computations involving the log
structure, such as in the study of variations of Hodge
structures.
\end{remark}

%===========================================================
\subsection{The limit $k\to\infty$}
So far we have dealt with a fixed structure $\scrS$ that was
consistent to order $k$. We now wish to take the limit $k\to\infty$ by
considering a sequence $\scrS_k$ of structures that are compatible in
the following way.

\begin{definition}\label{compatibility of structures}
Two structures $\scrS$, $\scrS'$ are \emph{compatible to order $k$} if
the following conditions hold.
\begin{enumerate}
\item
If $\fop=(\fop,m, c)\in\scrS$ is a wall with $c\neq 0$ and
$\ord_{\sigma_\fop}(m)\le k$ then $\fop\in \scrS'$, and the analogous
statement holds for $\scrS$ and $\scrS'$ interchanged.
\item
If $x\in \big(\Int(\fob)\cap\Int(\fob')\big)\setminus\Delta$ for slabs
$\fob\in \scrS$, $\fob'\in\scrS'$, then $f_{\fob,x}, f_{\fob',x}
\in\kk[P_x]$ agree modulo $t^{k+1}$.
\end{enumerate}
\end{definition}

If $\scrS$, $\scrS'$ are compatible to order $k$ and $\scrS$ is
consistent to order $k$ then $\scrS'$ is also consistent to order $k$
and the two deformations $X_k$ and $X'_k$ constructed from $\scrS$ and
$\scrS'$, respectively, are canonically isomorphic.

We are now in a position to reduce the Main Theorem to the construction
of a sequence of compatible structures.

\begin{proposition}\label{reduction to sequence of structures}
Assume there is a sequence $(\scrS_k)_{k\ge0}$ of structures on
$(B,\P,\varphi)$ such that for any $k$ (1) $\scrS_k$ is consistent to
order $k$ (2) $\scrS_k$ and $\scrS_{k+1}$ are compatible to order $k$.
Then there exists a formal toric degeneration of CY-pairs $(\hat\pi:
\hat X\to \hat O, \hat D)$ with central fibre $(X,D)$ and intersection
complex $(B,\P,\varphi)$ for the given pre-polarization on $X$. 
\end{proposition}

\proof
By compatibility of $\scrS_k$ and $\scrS_{k+1}$ we have a closed
embedding $X_k\to X_{k+1}$ exhibiting $X_{k+1}$ as flat deformation of
$X_k$, for any $k$. Thus $\hat X:=\limdir_k X_k$ is a formal scheme, flat
over $\kk\lfor t\rfor$. Moreover, the charts
\[
\psi^k: \big(\kk[P]/(t^{k+1})\big)_\fop \lra
\liminv (R_\tau)_{\fop_\tau}
\]
constructed in Proposition~\ref{charts for V^k(omega)} are also
compatible for various $k$ and compatible with the open embeddings and
with other choices of $f$. Hence for any $p\in X\setminus Z$ we obtain
an isomorphism of $\O_{\hat X,p}= \liminv \O_{X_k,p}$ with a
localization of ${\liminv}_k \kk[P]/(t^{k+1})$. Define
the deformation $\hat D\subseteq\hat X$ of $D$ by interpreting
$\psi^k$ as the chart for a log structure. Explicitly, since $D$ is a
toric Cartier divisor there exists $m\in P$, unique up to an
invertible element, such that $(z^m,t)\subseteq \liminv
(R_\tau)_{\fop_\tau}$ is the ideal of $D$. Define $\hat D$ as the closure
of the divisor defined by $z^m$. Note $\hat D$ 
fulfills (iii) in Definition~\ref{formal
toric degeneration} by construction. Since $\codim Z\ge2$ this also
shows regularity of $\hat X$ in codimension one. Furthermore, since
$\liminv$ commutes with push-forward by $j:X\setminus Z\to X$,
Lemma~\ref{extension across Z_omega} implies $\O_{\hat X}= j_*\O_{\hat
X\setminus Z}$. This shows that $\hat X$ is $S_2$, and hence $\hat X$
is normal as required in  Definition~\ref{formal toric
degeneration},(ii). Finally, the central fibre is isomorphic to
$(X,D)$ by construction.
\qed
\medskip

The rest of the paper is devoted to the construction of a sequence of
structures $\scrS_k$ as demanded in Proposition~\ref{reduction to
sequence of structures}.

%===========================================================
%===========================================================
\section{The algorithm}

This section is devoted to the core construction of this paper, the
inductive generation of structures $(\scrS_k)_{k\ge0}$ as required in
Theorem~\ref{reduction to sequence of structures}. We continue with
the polarized tropical manifold $(B,\P,\varphi)$, open gluing data for
the cone picture $s$ and data $(f_e)_e$ defining a positive log smooth
structure on $X=X_0(\check B,\check \P,s)$, as fixed at the
beginning of Section~\ref{section main objects}. We now also assume
local rigidity (Definition~\ref{def:locally rigid}).

\begin{theorem}\label{thm:scrS_k} If all cells of $B$ are bounded
there exists a sequence $(\scrS_k)_{k\ge0}$ of structures on
$(B,\P,\varphi)$ such that for any $k$ (1) $\scrS_k$ is consistent to
order $k$ (2) $\scrS_k$ and $\scrS_{k+1}$ are compatible to order $k$.
\end{theorem}

The proof of this theorem occupies the whole section, with the proof
of one technical result deferred to Section~\ref{section higher
codimension}. As remarked earlier, most of the arguments do not require
bounded cells, so we shall work with the general case, making it
clear where we require the boundedness hypothesis.

%===========================================================
\subsection{The initial structure}
Take $\scrS_{0}$ to consist only of slabs $\fob$ where
$\fob=\rho_{\fob}$ is a codimension one cell of $\P$ and
$f_{\fob,x}=\Pi^{-1}(f_e)$, where $e:v=v[x]\to\rho_{\fob}$ is the
vertex in the connected component of $\rho_{\fob}\setminus\Delta$
containing $x$ and $\Pi: \kk[P^{\gp}_x]\to\kk[P^{\gp}_v]$ is
given by parallel transport from $x$ to $v$ along a path inside
$\rho_{\fob}\setminus\Delta$. Then $\Chambers(\scrS_0)= \P^{[n]}$, and
we can take $\P_\scrS=\P^{[\le n-1]}$.

\begin{proposition}
$\scrS_0$ is consistent to order $0$.
\end{proposition}

\proof
The only joints $\foj$ of $\scrS_0$ are the codimension two cells of
$\P$, so take $\foj=\tau\in\P^{[n-2]}$, $\tau\not\subseteq \partial
B$. Let $\sigma_1,\ldots,\sigma_l$ be the maximal cells of $\P$
containing $\tau$, ordered cyclically, that is,
$\sigma_i\cap\sigma_{i+1}= \rho_i\in\P^{[n-1]}$
with $\sigma_{l+1}:=\sigma_1$, and take $\check
d_{\rho_i}$ (as in~(\ref{monodromy constant})) to be negative on
$\sigma_i$. To check consistency for $g:\omega\to\tau'$ with
$\tau'\subseteq \sigma_i$ for all $i$ it suffices to consider
$\omega=\tau'=\tau$, by  Remark~\ref{rem:joints},(5). Then by
Construction~\ref{gluing morphisms},~II.2, letting 
$x\in\Int(\tau)\setminus\Delta$  and $e:v=v[x]\to \tau$,
\[
\theta_i:=F^0_s(\id_{\tau},\sigma_i): 
m\longmapsto \big(D(s_e,\rho_i,v)^{-1}s_{e}^{-1}
(f_{\fob_{\rho_i},x})\big)^{-\langle m,\check d_{\rho_i}\rangle}.
\]
Thus in $R^0_{\id_\tau,\sigma_1}$, letting $e_{\rho_i}: v\to\rho_i$,
\begin{eqnarray*}
(\theta_l\circ\ldots \circ\theta_1)(m)&=& \Big(\prod_{i=1}^l 
D(s_{e},\rho_i,v)^{\langle  m,\check d_{\rho_i}\rangle}
\Big) s_{e}^{-1}\Big(\prod_{i=1}^l \big(f_{e_{\rho_i}}|_{V_{\id_{\tau}}}\big)^{-\langle
 m,\check d_{\rho_i}\rangle}\Big)\\
&=&\Big( \prod_{i=1}^l {s_{e,\sigma_i}( m)\over
s_{e,\sigma_{i+1}}( m)} \Big)
s_{e}^{-1}\Big(\prod_{i=1}^l\big(f_{e_{\rho_i}}|_{V_{\id_{\tau}}}\big)^{-\langle
 m,\check d_{\rho_i}\rangle}\Big)\ =\ 1,
\end{eqnarray*}
the last equality by~(\ref{multiplicative condition}). This is the
desired consistency.
\qed
\medskip

Note that according to Remark~\ref{log structure on X_k} the
structure $\scrS_0$ defines an abstract log structure on $X$.
Checking consistency to order~$0$ means verifying the multiplicative
condition~(\ref{multiplicative condition}) for the associated
section of $\shLS^+_{\pre, X}$. By construction this is indeed just
the log structure we started with. 

%===========================================================
\subsection{Scattering diagrams}\label{par:scattering diagrams}
Given $\scrS_{k-1}$ the construction of $\scrS_k$ proceeds in
three steps. The first of these introduces new walls, of order $k$,
by a procedure that is strictly local around a joint and is the
subject of this subsection. The second step performs various
semi-global adjustments involving several joints. The remaining
trouble terms are removed in the last step by a normalization
procedure applied to each slab.

Recall from \S\ref{par:joints} the space $\shQ^v_{\foj,\RR}=
\Lambda_{v,\RR}/\Lambda_{\foj,\RR}$ for a joint $\foj$ and the
notations $\doublebar m$, $\doublebar \tau$ etc. We think of
$\shQ^v_{\foj,\RR}$ as being divided by those half-lines $\foc$
emanating from the origin that are contained in $\doublebar \rho$ for
some $\rho= \rho_\foc\in\P^{[n-1]}$, $\rho\supseteq\foj$. We refer to
these half-lines as \emph{cuts}. Observe that if $\codim\sigma_\foj=1$
there are two cuts separating $\shQ^v_{\foj,\RR}$ into two
half-planes, while in the codimension two case the cuts subdivide
$\shQ^v_{\foj,\RR}$ into a number of strictly convex cones. In the
codimension zero case there are no cuts at all.

Once an orientation on $\shQ^v_{\foj,\RR}$ is chosen and $\foc=
\RR_{\ge0}\cdot \doublebar m\subseteq \shQ^v_{\foj,\RR}$,
$m\in\Lambda_\sigma\setminus\{0\}$, is a (rational) half-line
emanating from the origin, the unique generator $n_\foc$ of
$m^\perp\cap\Lambda_\foj^\perp \simeq \ZZ$ with the property that
$\langle m',n_\foc\rangle>0$ for $m,m'$ mapping to an oriented basis of
$\shQ^v_{\foj,\RR}$, is called the \emph{normal vector to $\foc$}.

Consistency at $\foj$ depends only on the local properties around
$\foj$ of slabs and walls containing $\foj$ and hence can be studied
on $\shQ^v_{\foj,\RR}$. The following definition is an abstraction of
the situation.

\begin{definition}\label{def:scattering diagram}
A \emph{ray} in $\shQ^v_{\foj,\RR}$ is a triple $(\forr, m_\forr,
c_\forr)$, where $\forr$ is a one-dimensional, rational cone
$\RR_{\ge0}\cdot \doublebar q$, $q\in \Lambda_v\setminus
\Lambda_\foj$; $m_\forr$ is a nonzero exponent on a maximal cell
$\sigma$ with $\pm\doublebar m_\forr\in \forr \cap\doublebar \sigma$
and such that $m\in P_x$ for all $x\in\foj\setminus\Delta$; $c_\forr$
is a constant in $\kk$. By abuse of notation we often just write
$\forr$ to refer to $(\forr,m_\forr,c_\forr)$. A ray is called
\emph{incoming}, \emph{outgoing} and \emph{undirectional} in the
respective cases $\doublebar m_\forr\in \forr\setminus \{0\}$,
$-\doublebar m_\forr\in \forr\setminus \{0\}$ and $\doublebar
m_\forr=0$. The \emph{order} of a ray $\forr$ is defined as
$\ord_\foj(m_\forr)$.

A \emph{scattering diagram for $\foj$} at a vertex $v\in\sigma_\foj$
consists of (1)~a choice of $\omega\in \P$ with $\foj\cap \Int\omega
\neq \emptyset$, $\omega\subseteq \sigma_\foj$ and $v\in\omega$, (2)~a
finite set of rays $\forr =(\forr, m_\forr, c_\forr)$, (3)~for each
cut $\foc\subseteq \shQ^v_{\foj,\RR}$ and any $x\in (\foj\cap
\Int\omega)\setminus\Delta$ a function $f_{\foc,x}\in \kk[P_x]$ with the
same properties as the functions $f_{\fob,x}$ in
Definition~\ref{slab}, (4)~an orientation of $\shQ^v_{\foj,\RR}$. The
notation is $\foD=\{\forr,f_\foc\}$ with $\omega$ and the orientation
of $\shQ^v_{\foj,\RR}$ understood. For rays $\forr$ and cuts $\foc$ of
a scattering diagram we write $n_\forr$ and $n_\foc$, respectively,
for the now well-defined normal vectors.

Two scattering diagrams $\foD=\{\forr,f_\foc\}$, $\foD'=
\{\forr',f'_\foc\}$ for $\foj$ at $v$ defined with the same $\omega$
are \emph{equivalent modulo a monomial ideal $J\subseteq
\kk[P_{\omega,\sigma}]$}, where $\sigma\supseteq\foj$, if (1)~for any
$m\in P_{\omega,\sigma}$ with $z^m\not\in J$,
and $\eps\in\{-1,1\}$, it holds
\[
\prod_{\{\forr\in \foD\,|\,m_\forr=m,\,
\eps\cdot\doublebar{\scriptstyle m}\in\forr\}}
(1+c_\forr z^{m_\forr})
=\prod_{\{\forr'\in \foD'\,|\,m_{\forr'}=m,\,
\eps\cdot\doublebar{\scriptstyle m}\in\forr'\}}
(1+c_{\forr'} z^{m_{\forr'}})\mod J,
\]
where we use parallel transport through $v$ to interpret $m$ as an
exponent on other maximal cells containing $\foj$; (2)~for any cut
$\foc$ and any $x\in (\foj\cap \Int\omega) \setminus\Delta$ the
functions $f_{\foc,x}, f'_{\foc,x} \in\kk [P_x]$ agree modulo terms of
$\ord_\foj$ at least $k+1$.
\qed
\end{definition}

Given a scattering diagram $\foD=\{\forr_i,f_\foc\}$ and $g:\omega
\to\tau$ with $\foj\cap \Int\omega\neq\emptyset$, $\tau\subseteq
\sigma_\foj$, and $\sigma\in\P_\max$ containing $\foj$, we obtain a log
isomorphism of $R^k_{g,\sigma}$ just as from a loop around a joint.
Specifically, let $\sigma_1,\ldots,\sigma_r=\sigma_0$ be a cyclic
ordering of the maximal cells containing $\foj$ compatible with the
orientation of $\shQ^v_{\foj,\RR}$, and $\rho_j=
\sigma_{j-1}\cap\sigma_j$. This induces a cyclic ordering of the
cuts $\foc_j\subseteq \doublebar \sigma_{j-1}\cap
\doublebar\sigma_j$. In the codimension two case this inclusion
defines $\foc_j$ uniquely, while there are two choices in the case of
codimension one. Assume that the rays are labelled cyclically as
well and in such a way that $\forr_i\subseteq \doublebar \sigma_j$ iff
$i_{j-1}< i\le i_j$. Then for $\forr_i\subseteq\doublebar\sigma_j$ and
for any $k$ we have the log automorphism
\[
\theta_i: \Lambda_{\sigma_j}\lra (R^k_{g,\sigma_j})^\times,\quad
m\longmapsto
\big(s_{\omega\to\sigma_j}(1+c_{\forr_i}z^{m_{\forr_i}})
\big)^{-\langle \overline m, n_{\forr_i}\rangle}
\]
of $R^k_{g,\sigma_j}$, as in Construction~\ref{gluing morphisms},~II.1,
where we think of passing through $\forr_i$ in the sense of the cyclic
ordering of the $\sigma_i$. Note that $\theta_i$ can also be written as
\[
\theta_i= \exp\big(-\log\big( s_{\omega\to\sigma_j}(1+c_{\forr_i}
z^{m_{\forr_i}})\big) \partial_{n_{\forr_i}}\big),
\]
and hence is an element of the group $H^{I_{g,\sigma}^{>k}}_\foj$, acting on $R^k_{g,\sigma_j}$. Similarly,
following Construction~\ref{gluing morphisms},~II.2 with
$x\in\foj \setminus\Delta$ such that $v[x]=v$, the functions
$f_{\foc_j,x}$ define the log isomorphism
\[
\theta_{\foc_j}:P_{\omega,\sigma_{j-1}}\lra
(R^k_{g,\sigma_j})^\times,\quad
m\longmapsto \big(D(s_{ v\to\omega},\rho_j,v)^{-1}
s_{v\to\omega}^{-1}(f_{\foc_j,x})\big)
^{-\langle \overline m, n_{\foc_i}\rangle}
\]
from $R^k_{g,\sigma_{j-1}}$ to $R^k_{g,\sigma_j}$, with monoid
homomorphism defined by parallel transport through $x$. Note that if
we chose an $x$ with $v[x]\neq v$ we would still obtain an equivalent
log isomorphism as verified in~(\ref{basic gluing: change of
vertex}), see also Remark~\ref{rem:discussion gluing morphisms},(2).
Define
\begin{multline}\label{theta_foD}
\theta^k_{\foD,g}:=
\big(\theta_{i_r}\circ\ldots\circ\theta_{i_{r-1}+1}\big)\\
\circ\theta_{\foc_r}
\circ\big(\theta_{i_{r-1}}\circ\ldots\circ\theta_{i_{r-2}+1}\big)
\circ \theta_{\foc_{r-1}}
\circ \ldots \circ
\big(\theta_{i_1}\circ\ldots\circ\theta_{1}\big)
\circ \theta_{\foc_1}.
\end{multline}
After distinguishing $\sigma_1$ this is a well-defined representative
of a log automorphism of $R^k_{g,\sigma_1}$. In fact, any two log
automorphisms associated to rays or slabs in the same direction
commute, so this composition is independent of the chosen indexing.

By definition $\theta^k_{\foD,g}$ depends only on the equivalence
class of $\foD$ to order $k$. Note also that reversing orientations
leads to $(\theta^k_{\foD,g})^{-1}$, while a different choice of
$\sigma_1$ leads to conjugation of $\theta^k_{\foD,g}$ by a log
isomorphism $R^k_{g,\sigma_j} \to R^k_{g,\sigma_1}$ for some $j$. We
are often only interested in properties invariant under these changes
and hence suppress them in the notation for $\theta^k_{\foD,g}$.

\begin{construction}\label{scattering diagrams from structure}
A structure $\scrS$ induces a scattering diagram $\foD_\foj=
\foD_\foj(\scrS,\omega,v)$ for each joint $\foj$, $\omega\in\P$ with
$\foj\cap\Int\omega\neq\emptyset$, $\omega\subseteq\sigma_\foj$ and
vertex $v\in\omega$ as follows. The slabs containing $\foj$ readily
define the functions $f_{\foc,x}$. For a wall $\fop$ containing
$\foj$ there are the following possibilities.
\begin{enumerate}
\item
$\foj\subseteq \partial \fop$. Then add the
ray $(\doublebar \fop, m_\fop, c_\fop)$ to $\foD_\foj$. This ray is
incoming, outgoing or unoriented if $\foj\subseteq \topp(\fop)$,
$\foj\subseteq \base(\fop)$ or $\foj\subseteq \sides(\fop)$,
respectively.
\item
$\foj\cap\Int \fop\neq \emptyset$. Then $\doublebar \fop$ is a line
through the origin, defining two one-dimensional half lines
$\forr$, $\forr'=-\forr \subseteq \shQ_{\foj,\RR}$. Then add the
pair of rays $(\forr, m_\fop,c_\fop)$, $(\forr', m_\fop,c_\fop)$. These
are either both undirectional or a pair of an incoming and an
outgoing ray.
\end{enumerate}
Note that consistency of $\scrS$ around $\foj$ to order $k$ can be
expressed by $\theta^k_{\foD_\foj,\id_{\sigma_\foj}}=1$.
\qed
\end{construction}

\begin{remark}\label{rem: vertices and scattering diagrams}
For a different choice of vertex $v'\in \omega$ there is a
piecewise linear identification of $\shQ^v_{\foj,\RR}$ with
$\shQ^{v'}_{\foj,\RR}$ defined on $\doublebar\sigma$ by parallel
transport from $v$ to $v'$ inside $\sigma\in \P_\max$. This identifies
the scattering diagrams $\foD_\foj=\foD_\foj(\scrS,\omega,v)$ and
$\foD'_\foj=\foD_\foj(\scrS,\omega, v')$.
Note that the respective computations in
(\ref{theta_foD}), for the same maximal
cell $\sigma_1\supseteq\foj$, then only
differ by changing the underlying monoid homomorphisms by parallel
transport from $v$ to $v'$ in $\sigma_1$. In particular,
we have the equality of representatives of log automorphisms of
$R^k_{\omega\to\sigma_\foj,\sigma_1}$
\[
\theta^k_{\foD_\foj}= \theta^k_{\foD'_\foj}.
\]

Similarly, any of the considerations with scattering diagrams
below are independent of the choice of $v$.
\end{remark}

Assuming $\theta^{k-1}_{\foD,\id_{\sigma_\foj}}=1$, we now use the
structure of the group $\lperp H^{I_{g,\sigma}^{>k}}_\foj$ to try to
achieve $\theta^k_{\foD, \id_{\sigma_\foj}}=1$ by adding some rays
and, in the codimension two case, changing the functions $f_\foc$ to
order $k$. The key idea is captured in the following Lemma of
Kontsevich and Soibelman (\cite{ks}, Theorem~6), adapted to our
setting. For the rest of this subsection fix the joint~$\foj$,
$\sigma\in\P_\max$ containing $\foj$, and $g:\omega\to \sigma_\foj$
with $\foj\cap \Int\omega\neq \emptyset$ and $v\in\omega$. Write $I_k=
I_{g,\sigma}^{>k} \subseteq \kk[P_{\omega,\sigma}]$.

\begin{definition}\label{cone groups}
For $K\subseteq \shQ^v_{\foj,\RR}$ a strictly convex cone ($K\cap
-K=\{0\}$), not necessarily closed, and $I\subseteq
\kk[P_{\omega,\sigma}]$ a monomial ideal with radical $I_0$ we define
the following Lie subalgebras of $\fog^{I}_\foj$
\begin{align*}
\fog^{I}_{\foj,K}&:=\bigoplus_{\substack{z^m\in I_0\setminus I\\
-\doublebar{\scriptstyle m}\in K\setminus\{0\}}}
z^m(\kk\otimes\Lambda_{\foj}^{\perp}),\\
\foh^{I}_{\foj,K}& :=\bigoplus_{\substack{z^m\in I_0\setminus I\\
-\doublebar{\scriptstyle m}\in K\setminus\{0\}}}
z^m\big(\kk\otimes(\overline m^\perp\cap \Lambda_{\foj}^{\perp})\big)
=\fog^I_{\foj,K}\cap \foh^I_{\foj}.
\end{align*}
The corresponding Lie groups are denoted $G^I_{\foj,K}$ and
$H^I_{\foj,K}$.
\end{definition}

Note that $\foh^I_{\foj,K}\subseteq {}^\perp\foh^I_\foj$.

\begin{lemma}\label{KSlemma}
Let $\foj$ be a joint with $\sigma_\foj\in\P_\max$ and $K\subseteq
\shQ^v_{\foj,\RR}$ a strictly convex cone. Then for any $\theta\in
H^{I_k}_{\foj,K}$ there exists a scattering diagram $\foD$ for $\foj$
consisting entirely of outgoing rays $\forr$ with $\Int\forr \subseteq
K$ such that $\theta=\theta^k_{\foD,g}$. Moreover, $\foD$ is unique up to
equivalence to order~$k$.
\end{lemma}
\proof
For $k=0$ we may take $\foD=\emptyset$. By induction on $k$ we may
thus assume there exists a unique scattering diagram $\foD'$ with
$\theta^{k-1}_{\foD',g}=\theta \mod I_{k-1}$. Then by the
definition of $\foh_{\foj,K}^{I_k}$ we can write uniquely
\begin{eqnarray}\label{KSsplitting}
\theta^k_{\foD', g} \circ \theta^{-1}=
\exp\big(  \sum_i c_i z^{m_i}\partial_{n_i} \big),
\end{eqnarray}
where $z^{m_i}\in I_{k-1}\setminus I_k$, $-\doublebar m_i \in
K\setminus\{0\}$, $n_i\in \overline{m_i}^\perp \cap
\Lambda_{\foj}^\perp$ and $c_i\in \kk\setminus\{0\}$. By changing
$c_i$ we may assume $n_i= n_{-\RR_{\ge0} \doublebar m_i}$. Define
$\foD$ by adding the outgoing rays $(-\RR_{\ge0} \doublebar m_i, m_i,
c_i)$ to $\foD'$. Noting that $[z^{m_i}\partial_{n_i},\foh_\foj
^{I_k}]=0$, we see that
\[
\theta^k_{\foD,g} \circ \theta^{-1}
= \theta^k_{\foD', g} \circ \theta^{-1}\circ
\prod_i\exp\big(  - c_i z^{m_i}\partial_{n_i} \big)
=\id.
\]
Uniqueness follows from the uniqueness of the expansion in
(\ref{KSsplitting}).
\qed
\medskip

As we will see, the same idea as in the proof of Lemma~\ref{KSlemma}
can be used to add rays to a codimension zero scattering diagram
$\foD'$ with $\theta^{k-1}_{\foD',g}= 1$ to construct a scattering
diagram $\foD$ with $\theta^k_{\foD,g}\in \ker\big( {}^\parallel
H^k_\foj \to {}^\parallel H^{k-1}_\foj \big)$, uniquely up to
equivalence to order $k$. Note that the remaining exponents $m$ with
$\overline m\in\Lambda_\foj$ have to be dealt with by other arguments
since outgoing rays always lead to elements in the subgroup ${}^\perp
H^k_\foj \subseteq H^k_\foj$.

In higher codimension, under the presence of slabs, this is much more
subtle because we have to convert between computations in the various
groups $R^k_{g,\sigma_i}$ using the log isomorphisms associated to
slabs. In particular, it is not clear that the commutation does not
introduce poles in directions different from the cuts $\foc$
corresponding to slabs. We will also have weaker uniqueness properties
because one can always replace a ray in the direction of a cut
$\foc$ with a change of $f_\foc$. The detailed study of this situation
is the subject of the technical last section. Here we content
ourselves with a statement of the results needed for the construction
of $\scrS_k$.

For enhanced readability we introduce the following notations. Recall
we have fixed $\foj$, $\sigma\in\P_\max$ with $\foj\subseteq
\sigma$, and $g:\omega\to\sigma_\foj$ with
$\foj\cap\Int\omega\neq\emptyset$, and we work with various log
automorphisms of $R^k_{g,\sigma} =\big(\kk[P_{\omega,\sigma}]/
I_k\big)_{f_{g,\sigma_\foj}}$, $I_k=I^{>k}_{g,\sigma}$.

\begin{convention}
1)~For a set $V_\mu\subseteq \Lambda_\sigma$, a subspace
$W_\mu\subseteq \Lambda_\sigma^*$ and elements $f_\mu\in
(R^k_{g,\sigma})^\times$  we write
\[
O^k \Big(\sum_\mu \frac{V_\mu}{f_\mu}\otimes W_\mu\Big)
\]
for the set of log automorphisms of $R^k_{g,\sigma}$ of the form
$\theta= \exp\big(\sum_{\mu,i} c_{\mu,i} \frac{z^{m_{\mu,i}}}
{f_\mu}\partial_{n_{\mu,i}}\big)$ with $\overline{m_{\mu,i}}
\in V_\mu$, $n_{\mu,i}\in W_\mu$, $\ord_{\sigma_\foj}(m_{\mu,i})\ge
k$.\\[1ex]
2)~Let $v\in\omega$ be a vertex, $e:v\to\omega$. Then for
$\rho\in\P^{[n-1]}$ containing $\foj$ define
\[
f_\rho= f_{\rho,v}:= D(s_e,\rho,v)^{-1} f_{\rho,e,\sigma}
\in \kk[P_{\omega,\sigma}],
\]
with $f_{\rho,e,\sigma}$ defined in~(\ref{f_{rho,e,sigma}}). Note that
according to~(\ref{shLS change of chart, symmetric form}) a
different choice $e':v'\to\omega$ leads to
\begin{eqnarray}\label{f_rho}
f_{\rho,v'}=  z^{m^\rho_{v'v}}
f_{\rho, v}.
\end{eqnarray}
In particular, $f_\rho$ is well-defined up to multiplication by $z^m$
with $\overline m\in\Lambda_\omega$, $\ord_\omega(m)=0$.
\end{convention}

For example, with this notation $\ker\big( \lparallel H^{I_k}_\foj\to
\lparallel H^{I_{k-1}}_\foj \big) =O^k\big((\Lambda_\foj
\setminus\{0\}) \otimes \Lambda_\foj^\perp\big)$. We are now ready to
state the main result of Section~\ref{section higher codimension}.

\begin{proposition}\label{scattering proposition}
Let $\foD'$ be a scattering diagram for $\foj\in \Joints(\scrS_{k-1})$
with $\theta^{k-1}_{\foD', g}=1$, $g:\omega\to\sigma_\foj$. Then there
exists a scattering diagram $\foD$, equivalent to order $k-1$ to
$\foD'$ and with the sets of rays differing only by outgoing rays
$\forr$ with $\forr\not\subseteq \doublebar\rho$ for any
$\rho\in\P^{[n-1]}$ containing $\foj$, such that
\begin{align}\label{scattering diagram simplified}
\theta^k_{\foD,g}&\in
\left\{\begin{array}{ll}
O^k\big((\Lambda_\foj\setminus\{0\})\otimes
\Lambda_\foj^\perp\big),&
\codim \sigma_\foj=0\\[1ex]
\displaystyle O^k\Big(\Lambda_\foj\otimes \Lambda_\foj^\perp+
\frac{\Lambda_\rho}{f_\rho}\otimes \Lambda_\rho^\perp \Big),&
\codim \sigma_\foj=1\ (\rho=\sigma_\foj)\\
\displaystyle O^k\Big(\Lambda_\foj\otimes \Lambda_\foj^\perp+
{\sum}_{\rho\supseteq \foj}\ \frac{\Lambda_\foj}{f_\rho}
\otimes \Lambda_\rho^\perp \Big),&
\codim \sigma_\foj=2.\\[-1ex]
\end{array}\right.
\end{align}\\[-1ex]
If $\codim \sigma_\foj<2$ the functions $f_{\foc,x}$ of $\foD$ and
$\foD'$ coincide, while if $\codim \sigma_\foj=2$ they may be changed
by adding multiples of $z^m$ with $-\doublebar m\in
\foc\setminus\{0\}$. 

Moreover, up to equivalence $\foD$ is the unique scattering diagram
with these properties.
\end{proposition}

\noindent
\emph{Proof for $\codim \sigma_\foj=0$.}
Arguing similarly to Lemma~\ref{KSlemma} we have the unique
decomposition
\begin{eqnarray}\label{codim 0 splitting}
\theta^{k}_{\foD, g} =
\exp\Big(  \sum_i c_i z^{m_i}\partial_{n_i} \Big),
\end{eqnarray}
but this time only $z^{m_i}\in I_{k-1}\setminus I_k$,
$\overline{m_i}\neq 0$, $n_i\in \overline{m_i}^\perp\cap
\Lambda_\foj^\perp$. Define $\foD$ by adding to $\foD'$ for each $i$
with $\doublebar m_i\neq 0$ the ray $(-\RR_{\ge0} \doublebar m_i,
m_i,c_i)$, assuming without loss of generality $n_i= n_{-\RR_{\ge0}
m_i}$. Since the log automorphisms of these rays are in the center of
$H^{I_k}_\foj$ it holds
\[
\theta^k_{\foD,g}= \theta^{k}_{\foD', g} \circ \prod_{\{i\,|\,
\doublebar{\scriptstyle m_i}\neq 0\}}
\exp\big(-c_i z^{m_i}\partial_{n_i}\big)
= \exp\Big(  \sum_{\{i\,|\, \doublebar{\scriptstyle m_i}=0\}} c_i
z^{m_i}\partial_{n_i} \Big).
\]
This is of the desired form. Finally, uniqueness follows from the
uniqueness statement for~(\ref{codim 0 splitting}).

The proof for $\codim\sigma_\foj>0$ occupies Section~\ref{section higher
codimension}.
\qed
\medskip

It is also important for this section to record the effect on
$\theta^k_{\foD,g}$ of certain simple changes to $\foD$.

\begin{proposition}\label{influence of monomial changes}
Let $\foD$ be a scattering diagram for $\foj$ and assume
$\theta^k_{\foD,g}$ fulfills~(\ref{scattering diagram simplified}) of
Proposition~\ref{scattering proposition}.
\begin{enumerate}
\item
If $\hat\foD$ is obtained from $\foD$ by adding the
term $c z^m$ to some $f_\foc$, with $\ord_\foj m=k$ and $\overline
m\in\Lambda_{\sigma_\foj}$, then
\[
\theta^k_{\hat \foD,g}= \theta^k_{\foD,g}\circ \exp \Big(
-c' \frac{z^m}{f_{\rho_\foc}}\partial_{n_\foc}\Big),
\]
with $c'= D(s_e,\rho_\foc,v)^{-1}\cdot s_e^{-1}(\overline m)\cdot c$,
$e:v\to\omega$.
\item
If $\hat\foD$ is obtained from $\foD$ by adding an
undirectional ray $(\forr,m,c)$ with $\ord_\foj m=k$ then
\[
\theta^k_{\hat \foD,g}=
\left\{\begin{array}{ll}
\theta^k_{\foD,g}\circ \exp \Big(-c' z^m\partial_{n_\forr}\Big),&
\codim \sigma_\foj\neq 1\\[1ex]
\displaystyle \theta^k_{\foD,g}\circ
\exp \Big( -c' z^m\partial_{n_\forr}\Big)\circ
O^k\Big(\frac{\Lambda_{\rho}}{f_{\rho}} \otimes\Lambda_\rho^\perp\Big),
&\codim \sigma_\foj=1\ (\rho=\sigma_\foj),
\end{array}\right.
\]
with $c'= s_h(m)\cdot c$, $h:\omega\to\sigma_j$.
\end{enumerate}
\end{proposition}

\proof
First we observe that the change in (1) has the same effect as
composing the log isomorphism associated to $\foc$ by $\exp(-c'
z^m/f_{\rho_\foc}\partial_{n_\foc})$. Note that since $\overline
m\in\Lambda_{\sigma_\foj}$ it holds $\ord_{\sigma_\foj}(m+m') >k $
whenever $\ord_{\sigma_\foj} (m')>0$. Thus by Lemma~\ref{conjugation
lemma} this log isomorphism commutes with any of the other log
isomorphisms, hence the result. Adding an undirectional ray is
similar, but in the codimension one case $f_\foc$ involves monomials
$z^{m'}$ with $\overline{m'}\in \Lambda_{\rho_\foc} \setminus
\Lambda_\foj$ since $\rho_\foc=\sigma_\foj$. But again by
Lemma~\ref{conjugation lemma} we obtain
\[
\Ad_{\theta_\foc}(-c z^m\partial_{n_\forr})= 
- c z^m\big( \partial_{n_\forr}+ f_{\rho_\foc}^{-1}
(\partial_{n_\forr} f_{\rho_\foc})\partial_{n_\foc}\big),
\]
since $\langle \overline m, n_{\foc} \rangle=0$. The exponential of
this expression is of the form $\exp(-c z^m\partial_{n_\forr})\circ 
O^k\big( (\Lambda_{\rho_\foc}/ f_{\rho_\foc})\otimes
\Lambda_{\rho_\foc}^\perp\big)$.
\qed

%===========================================================
\subsection{Step I: Scattering at joints}
We now begin the algorithm providing the induction step, the
construction of $\scrS_k$ from $\scrS_{k-1}$. Since the various
parts of it are scattered throughout four subsections, text providing
instructions for this process is shaded.

\begin{shaded}
\noindent I.1. \emph{Refinement of slabs.}
The notion of compatibility of structures
(Definition~\ref{compatibility of structures}) allows arbitrary
refinements of slabs. To be able to use local methods we now impose
the following conditions on $\scrS_{k-1}$, which can be achieved by
subdivision of slabs:
\begin{eqnarray}\label{slabs1}
\begin{aligned}
\text{If $\fob\in\scrS_{k-1}$ is a slab and $\fob\cap\partial
\rho_\fob\neq\emptyset$ there exists $\tau\subseteq \partial\rho_\fob
$ with}\\[-1ex]
\fob\cap\Int\tau\neq \emptyset\quad\text{and}\quad
\tau'\in\P^{[\le n-2]},\, \Int \tau'\cap\fob\neq\emptyset\
\Longrightarrow\ \tau'\subseteq\tau.
\end{aligned}
\end{eqnarray}
Of course, this also means refining the polyhedral decomposition
$\P_{\scrS_{k-1}}$.
\end{shaded}

Note that~(\ref{slabs1}) implies that if $\dim \fob\cap\partial
\rho_\fob =n-2$ for a slab $\fob$, then $\fob\cap\partial \rho_\fob$
is a joint.

For each joint $\foj$ of $\scrS_{k-1}$ we now obtain a scattering
diagram $\foD'_\foj$ (Construction~\ref{scattering diagrams from
structure}) to which we can apply Proposition~\ref{scattering
proposition} with $g=\id_{\sigma_\foj}$. Each ray of the scattering
diagram thus obtained defines a new wall with  base $\foj$. However,
if $\codim \sigma_\foj=2$ this also involves a change of slabs by
terms of order $k$ and thus influences the computation at other
joints. We therefore deal with joints of codimension two first.

\begin{shaded}
\noindent I.2. \emph{Adjustments  of slabs from joints of
codimension 2.} For each codimension two joint $\foj$ and any slab
$\fob$ containing $\foj$, Proposition~\ref{scattering proposition},
applied to any $g:\omega\to\sigma_\foj$, defines a change $\tilde
f_{\fob,x}$ of $f_{\fob,x}$, for all $x\in\fob\cap\omega$, by terms
of order $k$ along $\sigma_\foj$. Use \eqref{change of vertex
- slabs} to extend this modification of slab function $f_{\fob,x}$
to all $x\in\fob$.
In view of uniqueness and Remark~\ref{rem:joints},(4)
and~(5) the results for different choices of $\omega$ containing $x$
coincide.

Moreover, by~(\ref{slabs1}) any slab contains at most one joint $\foj$
with $\codim\sigma_\foj=2$. Hence the corrections from different
joints are independent of each other. After applying these changes to
$\scrS_{k-1}$ simultaneously we may therefore assume
Proposition~\ref{scattering proposition} applies without any change of
slabs. Note this replacement does not affect the equivalence class of
$\scrS_{k-1}$ to order $k-1$.
\end{shaded}

The next step produces the new walls. We require the following lemma.

\begin{lemma}\label{walls well-defined}
Let $\foj\subseteq B$, $\foj\not\subseteq\partial B$ 
be an $(n-2)$-dimensional polyhedral subset of
some $\sigma\in\P_\max$ and $m$ a monomial on $\sigma$ with
$\ord_\sigma(m)=k\ge 0$. Assume furthermore $m\in P_x$ for all
$x\in\foj$ and $-\overline m\not\in\Lambda_\foj$ to be contained in
the tangent wedge to $\sigma$ along $\foj$. Then $m\in P_x$ for any
$x\in \fop\setminus\Delta$ with
\[
\fop:=(\foj-\RR_{\ge0}\cdot \overline m)\cap\sigma,
\]
and for any $x\in \topp(\fop):= \cl\big( \partial\fop\setminus\big(\foj
\cup (\partial \foj-\RR_{\ge0}\overline m) \big) \big)$, $x\not\in\partial
B$, there exists
$\sigma'\in\P_\max$ with $x\in\sigma'$ and $\ord_{\sigma'}(m)>k$.
\end{lemma}

\proof
For any $\sigma'\in\P_\max$ with $\sigma'\cap(\fop\setminus
\foj)\neq\emptyset$, Proposition~\ref{exponentprop} with
$\sigma^+=\sigma$ shows
\[
\ord_{\sigma'}(m)\ge  \ord_\sigma(m)=k\ge 0.
\]
Thus $m\in P_x$ for any $x\in\fop$. On the other hand, if
$x\in\topp(\fop)$ then $\overline m$ is not tangent to the minimal cell
$\tau$ containing $x$. Thus $m$, as an affine function on $\check
\tau$, is non-constant. In particular, as it takes its minimal value on
$\check\sigma$, there exists a vertex $\check\sigma'$ of $\check\tau$
such that $\ord_{\sigma'}(m)= m(\check\sigma')> m(\check\sigma)
=\ord_\sigma(m)=k$.
\qed

\begin{shaded}
\noindent I.3. \emph{Scattering at joints.}
For a joint $\foj$ of $\scrS_{k-1}$ let $\foD_\foj$ be the
scattering diagram obtained from $\foD'_\foj:=\foD_{\foj}
(\scrS_{k-1},\sigma_\foj,v)$, for some choice of vertex
$v\in\sigma_\foj$, by the application of Proposition~\ref{scattering
proposition} with $g=\id_{\sigma_\foj}$. By I.2 the functions
$f_{\foc,x}$ remain unchanged, so $\foD_\foj$ differs from
$\foD'_{\foj}$ only by outgoing rays in directions different from
directions of slabs. Moreover, Proposition~\ref{scattering
proposition} applied with various $g:\omega\to\sigma_\foj$, for
$\omega$ with $\foj\cap \Int\omega \neq\emptyset$, implies $m_\forr
\in P_x$ for any $\forr\in \foD_\foj\setminus \foD'_\foj$ and 
$x\in\foj\setminus\Delta$. In fact, by Remark~\ref{rem: vertices and
scattering diagrams} we may assume that the vertex $v\in\sigma_\foj$
lies in $\omega$. Hence $v$ can also be used for the scattering
diagram $\foD_\foj(\omega)$ obtained for $\omega\to\sigma_\foj$.
Then as in Remark~\ref{rem:joints} one sees that uniqueness implies
equivalence of $\foD_\foj$ and $\foD_\foj(\omega)$ to order $k$.
This shows $m_\forr\in P_x$ for $x\in (\foj\cap\Int\omega)
\setminus\Delta$.

Define $\scrS_k^\mathrm{\,I}$ by adding to $\scrS_{k-1}$, for any
joint $\foj$ and any ray $\forr\in \foD_\foj\setminus \foD'_{\foj}$,
the wall $(\fop_\forr, m_\forr,c_\forr)$ with
\[
\fop_\forr:= \big(\foj- \RR_{\ge0}\cdot \overline{m_\forr} \big)
\cap\sigma,
\]
where $\sigma$ is the unique maximal cell with $\forr\subseteq
\doublebar\sigma$. This is indeed a wall by Lemma~\ref{walls
well-defined}. Add some more walls $\fop$ with $c_\fop=0$ to achieve
the requirement of Definition~\ref{def:structure},(ii), for example by
covering $H_\fop\cap\sigma$ by such walls, for each added wall
$\fop\subseteq \sigma$, with $H_\fop\subseteq \Lambda_{\sigma,\RR}$
the affine hyperplane containing $\fop$. (This step is indeed not
necessary as follows by the arguments in \S\ref{par:interstices}, but
we will not prove this.) Choose also a polyhedral decomposition
$\P_{\scrS_k^\I}$ which on $|\scrS_{k-1}|$ refines
$\P_{\scrS_{k-1}}$, and such that each slab or wall of $\scrS_k^\I$ is
a union of $(n-1)$-cells of $\P_{\scrS_k^\I}$. 

We now have produced a new structure $\scrS_k^\I$, with the
superscript ``I'' indicating that it is the result of Step~I of the
algorithm. By subdividing slabs we may assume (i) in the definition of
structures (Definition~\ref{def:structure}) to continue to hold, while
(\ref{slabs1}) from I.1 is true in any case.
\end{shaded}

We now check that the corrections at the joints of $\scrS_{k-1}$ have
the desired effect. In particular, we have to verify that new walls do
not influence the computations at joints different from their bases.

\begin{proposition}\label{step I result}
For any $\foj\in\Joints(\scrS_k^\I)$ the scattering diagram $\foD_\foj=
\foD_\foj(\scrS_k^\I, \sigma_\foj, v)$, $v\in\sigma_\foj$ fulfills
\begin{align}\label{step I reduction}
\theta^k_{\foD_\foj,g}&\in
\left\{\begin{array}{ll}
O^k\big((\Lambda_\foj \setminus\{0\}) \otimes \Lambda_\foj^\perp\big),&
\codim \sigma_\foj=0\\[1ex]
\displaystyle O^k\Big(\Lambda_\foj\otimes \Lambda_\foj^\perp+
\frac{\Lambda_\rho}{f_\rho}\otimes \Lambda_\rho^\perp \Big),&
\codim \sigma_\foj=1\ (\rho =\sigma_\foj)\\
\displaystyle O^k\Big(\Lambda_\foj\otimes \Lambda_\foj^\perp+
{\sum}_{\rho\supseteq \foj}\ \frac{\Lambda_\foj}{f_\rho}
\otimes \Lambda_\rho^\perp \Big),&
\codim \sigma_\foj=2.\\[-1ex]
\end{array}\right.
\end{align}
\end{proposition}

\proof
$\scrS_k^\I$ differs from $\scrS_{k-1}$ effectively by the addition of
walls $\fop$ with base a joint $\foj_\fop$ of $\scrS_{k-1}$, the other
walls having $c_\fop =0$ and thus being irrelevant. We now discuss
the contributions of such $\fop$ to the computation of
$\theta^k_{\foD_\foj,\id_{\sigma_\foj}}\!$. There are the following
possibilities for the relative position of $\foj$ inside $\fop$:
\begin{enumerate}
\item
$\foj\subseteq\base(\fop)=\foj_\fop$.
\item
$\foj\subseteq \topp(\fop)$.
\item
$\foj \not\subseteq\partial\fop$.
\item
$\foj\subseteq\sides(\fop)$.
\end{enumerate}

In (1) $\foj_\fop$ is a joint of $\scrS_{k-1}$ and $\fop$ arose from
an outgoing ray produced in Proposition~\ref{scattering proposition}.
These precisely lead to the desired form~(\ref{scattering diagram simplified}) of
$\theta^k_{\foD_\foj,\id_{\sigma_\foj}}$.

In Case~(2) $\ord_\foj(m_\fop)>k$ by Lemma~\ref{walls well-defined}.
These walls do not make any contribution in order $k$ at $\foj$.

Case~(3) can only happen in the codimension zero case, that is, if
$\sigma_\foj\in\P_\max$. Then $\ord_\foj(m)= \ord_{\sigma_\foj}(m)=k$
and hence $\exp \big(-\log(1+cz^m)\partial_n \big)$ commutes with any
log automorphism of $R^k_{\id_{\sigma_\foj},\sigma_\foj}$, see
(\ref{Lie bracket}). Since the automorphism associated to $\fop$
occurs twice with opposite signs in $\theta^k_{\foD_\foj,
\id_{\sigma_\foj}}$ it makes no contribution.

In (4) with $\sigma_\foj$ maximal, the automorphism associated to
$\fop$ lies in $O^k\big((\Lambda_\foj\setminus \{0\})
\otimes\Lambda_\foj^\perp\big)$. Hence this wall preserves
(\ref{scattering diagram simplified}). Finally, in (4) with
$\codim\sigma_\foj>0$, Proposition~\ref{influence of monomial
changes},(2) shows that the insertion of $\fop$ preserves the form
of~(\ref{scattering diagram simplified}).
\qed

%===========================================================
\subsection{Interstices and consistency in codimension 0}
\label{par:interstices}
The remaining terms in~(\ref{step I reduction}) all involve exponents
tangent to joints or slabs. The topology at intersections of joints
now imply strong compatibility conditions that are the subject of this
subsection. Among other things, these restrictions already
imply consistency at codimension 0 joints.

\begin{definition}\label{interstice}
An \emph{interstice} of a structure $\scrS$ is an $(n-3)$-cell $\fod\in
\P_\scrS$ with $\fod\not\subseteq\partial B$.
\end{definition}

Analogous to the situation for a joint, for a vertex $v\in\sigma_\fod$ we
obtain a normal space $\shQ^v_{\fod,\RR} \simeq \RR^3$, defined as
$\Lambda_{v,\RR}/ \Lambda_{\fod,\RR}$. Again we write $\doublebar m
\in\shQ^v_{\fod,\RR}$ for the image of an exponent on any
$\sigma\in\P_\max$ containing $\fod$, and $\doublebar
\tau\subseteq\shQ^v_{\fod,\RR}$ for the image of the tangent wedge
along $\fod$ of a cell $\tau\subseteq B$, $\tau\supseteq\fod$.

To study the topology of the situation along $\fod$ we look at the
associated $2$-sphere
\[
S_\fod:= \big(\shQ^v_{\fod,\RR}\setminus\{0\}\big)/ \RR_{>0},
\]
which we orient arbitrarily. This $2$-sphere comes with the following
cell-decomposition. For the joints $\foj_1,\dots, \foj_s$ containing
$\fod$ we have $0$-cells $\doublebar{\foj_i}/ \!\sim$, the $1$-cells are
given by $\doublebar\fov/ \!\sim$ for $\fov\in\P_\scrS^{[n-1]}$,
$\fov\supseteq\fod$, and the $2$-cells are $\doublebar \fou/ \!\sim$ for
$\fou\in\Chambers(\scrS)$, $\fou\supseteq\fod$. Let $\Sigma_\fod$ be
the dual cell complex, with $0$, $1$ and $2$-cells $\widehat\fou$,
$\widehat \fov$ and $\widehat{\foj}$ defined by chambers, elements of
$\P_\scrS^{[n-1]}$ and joints containing $\fod$, respectively. Note
this is a subcomplex of the cell-complex $\Sigma$ studied in Step~3 of
the proof of Lemma~\ref{consistency lemma} for the case
$g=\id_{\sigma_\fod}$. It is then clear that for each edge path
$\beta$ in $\Sigma_\fod$ from $\widehat\fou$ to $\widehat\fou'$ we
obtain a log isomorphism
\[
\theta_\beta: P_{\sigma_\fod, \sigma_\fou} \lra
\big(R^k_{\id_{\sigma_\fod},\sigma_{\fou'}} \big)^\times
\]
from $R^k_{\id_{\sigma_\fod},\sigma_\fou}$ to
$R^k_{\id_{\sigma_\fod},\sigma_{\fou'}}$. The underlying monoid
homomorphism is obtained by parallel transport through $v$.

Next choose a base vertex $\widehat{\fou_0}\in \Sigma_\fod$ and let
$\gamma_i$ be a closed loop covering the edges of $\widehat{\foj_i}$
in counterclockwise direction. Because the $1$-skeleton
$\Sigma_\fod^1$ of $\Sigma_\fod$ has the homotopy type of $S^2$ minus
$s$ points, there exist paths $\beta_i$ on $\Sigma_\fod$, with
$\beta_i$ connecting $\widehat{\fou_0}$ with the base point of
$\gamma_i$, such that $\beta_i \gamma_i\beta_i^{-1}$ is a standard
generating set of $\pi_1(\Sigma_\fod^1, \widehat{\fou_0})$:
\begin{eqnarray}\label{pi_1 relation}
\beta_1 \gamma_1\beta_1^{-1}
\beta_2 \gamma_2\beta_2^{-1}
\ldots
\beta_s \gamma_s\beta_s^{-1}&=&1.
\end{eqnarray}
Note that such $\beta_i$ exist regardless of the given order
$\foj_1,\ldots,\foj_s$ of the joints. For the corresponding sequence
of log isomorphisms we conclude
\begin{eqnarray}\label{interstice relation}
\theta_{\beta_s}^{-1}\circ \theta_{\gamma_s}\circ \theta_{\beta_s}
\circ\ldots \circ
\theta_{\beta_1}^{-1}\circ \theta_{\gamma_1}\circ \theta_{\beta_1}
&=&1.
\end{eqnarray}
Note that we may impose additional conditions on the choices of
$\beta_i,\gamma_i$ as long as~(\ref{pi_1 relation}) holds.

For $\scrS=\scrS_k^\I$ constructed in Step~I,~(\ref{interstice
relation}) implies the following result for interstices $\fod$ with
$\codim\sigma_\fod=0$.

\begin{proposition}\label{codim 0 step II lemma}
Assume that $\fod$ is an interstice of $\scrS_k^\I$ with $\codim
\sigma_\fod=0$ and write, using~(\ref{step I reduction}),
\[
\theta_{\gamma_i} =\exp\Big(\sum_m a_{m,i} z^m\partial_{n_i(m)} \Big),
\]
as log automorphism of $R^k_{\id_{\sigma_\fod},\sigma_\fod}$, where
the sum runs over those exponents $m$ on $\sigma_\fod$ with $\overline
m\in\Lambda_{\foj_i}$, $\ord_{\sigma_\fod} (m)=k$, and $a_{m,i}\in\kk$,
$n_i(m)\in \Lambda_{\foj_i}^\perp$. Then for any $m\in
P_{\sigma_\fod,\sigma_\fod}$ with $\ord_{\sigma_\fod}(m) =k$ it holds
in $\Lambda^*_{\sigma_\fod}$
\begin{eqnarray}\label{interior interstice cycle condition}
\sum_i a_{m,i} n_i(m)=0.
\end{eqnarray}
\end{proposition}

\proof
Since $\theta_{\gamma_i}$ only involves monomials of order $k$ it
commutes with $\theta_{\beta_i}$ and any $\theta_{\gamma_j}$. Hence
(\ref{interstice relation}) shows
\[
1=\theta_{\gamma_s}\circ\ldots\circ\theta_{\gamma_1}
=\exp\Big(\sum_{m,i} a_{m,i} z^m\partial_{n_i(m)} \Big),
\]
which readily implies the result.
\qed
\medskip

To deduce analogous restrictions from higher codimension interstices
we need to understand how $\theta_{\gamma_i}$ transforms by
commutation with log isomorphisms changing chambers
$\fou\supseteq\fod$. 

\begin{lemma}\label{interstice interaction}
Let $\theta_\gamma$ be the log isomorphism from
$R^k_{\id_{\sigma_\fod},\sigma_{\fou'}}$ to $R^k_{\id_{\sigma_\fod},
\sigma_\fou}$ associated to an edge path $\gamma$ in $\Sigma_\fod$
connecting $\hat\fou'$ with $\hat\fou$. Then for any $m\in
P_{\sigma_\fod,\sigma_{\fou'}}$ with $z^m\in J_{l-1} := I_0^{l-1}\cdot
I_{k-1}+ I_k$, $n\in\Lambda_\sigma^*$ and $a\in 
\kk(\Lambda_{\sigma_\fod})\cap R^k_{\id_{\sigma_\fod},
\sigma_{\fou'}}$ there exist $b_i\in \kk(\Lambda_{\sigma_\fod})\cap
R^k_{\id_{\sigma_\fod},\sigma_\fou}$ and $n_i\in
\Lambda_{\sigma_\fou}^*$ such that
\[
\theta_\gamma\circ\exp(a z^m\partial_n) \circ\theta_\gamma^{-1}=
\exp\Big( \sum_i b_i z^m\partial_{n_i}\Big)
=\prod_i\exp\big(b_i z^m\partial_{n_i}\big)\mod J_l,
\]
as log automorphisms of $R^k_{\id_{\sigma_\fod},\sigma_\fou}$\!. Here we
identify $\Lambda_{\sigma_\fou}$ and $\Lambda_{\sigma_{\fou'}}$ by
parallel transport through some vertex $v\in\sigma_\fou\cap
\sigma_{\fou'}$.
\end{lemma}

\proof
By induction on the number of edges passed by $\gamma$ it suffices to
consider the case that $\fou$ and $\fou'$ are adjacent chambers and
$\theta_\gamma$ is the associated basic gluing morphism. Thus up to
choosing an isomorphism $R^k_{\id_{\sigma_\fod},\sigma_\fou} \to
R^k_{\id_{\sigma_\fod},\sigma_{\fou'}}$ by parallel transport through
a point in $\Int(\fod)\setminus\Delta$ we are in the situation of
Lemma~\ref{conjugation lemma}. Since $I_0\cdot J_{l-1}\subseteq J_l$
this shows first that we can ignore all expressions in $\theta_\gamma$
involving monomials $z^{m'}$ with $\ord_{\sigma_\fod}(m')>0$. Thus we
may assume $\theta_\gamma$ to be of the form $m'\mapsto f^{\langle
\overline{m'},n_0\rangle}$ with $f\in \kk(\Lambda_{\sigma_\fod})\cap
(R^k_{\id_{\sigma_\fod},\sigma_\fou})^\times$, and then
Lemma~\ref{conjugation lemma} gives the claimed result.
\qed

\begin{remark}\label{note on orders}
The reason for introducing $J_{l-1} =I_0^{l-1}\cdot I_{k-1}+I_k$ here
is that the order function from Definition~\ref{orderdef} is not in
general additive. For example, for a vertex $v$ in a one-dimensonal
$B$ with adjacent maximal cells $\sigma_1,\sigma_2$ and
$\varphi_v|_{\sigma_1}=0$, $\varphi_v|_{\sigma_2}$ having slope $1$,
we have $k[P_{\omega,\sigma}]\simeq \kk[\NN^2]$, $t=z^{(1,1)}$ and
$\ord_v\big((1,0)\big) =1$, $\ord_v\big((0,1)\big)=1$, but also
$\ord_v\big((1,1)\big)=1$. Similar ideals as $J_l$ will occur
repeatedly in the following.

Note that one exception where $\ord_\tau(m+m')=\ord_\tau(m)
+\ord_\tau(m')$ is when $m\in\Lambda_\tau$, for then $\ord_\sigma(m)=
\ord_\tau(m)$ for any $\sigma\in\P_\max$ containing $\tau$.
\qed
\end{remark}

We can now deduce an analogue of Proposition~\ref{codim 0 step II
lemma} for interstices of higher codimension.

\begin{proposition}\label{joints at higher codim interstices}
Assume that $\fod$ is an interstice of $\scrS_k^\I$ with
$\codim\sigma_\fod \ge1$ and let $\sigma\in\P_\max$,
$\sigma\supseteq\fod$. For any joint $\foj\subseteq \sigma$ of
$\scrS_k^I$ containing $\fod$ take the base chamber of the
loop $\gamma_\foj$ around $\foj$ to be contained in $\sigma$ and
oriented according to the chosen orientation of $\Sigma_\fod$ and write
\[
\theta_{\gamma_\foj}=\exp\Big(\sum_{\{(m,\nu)\,|\,m\in A,\,
\ord_{\sigma_\fod}\!(m)=k\}}
a_{\foj,m,\nu} z^m\partial_{n_{\foj,m,\nu}} \Big),
\]
as log automorphism of $R^k_{\id_{\sigma_\fod},\sigma}$, where
$a_{\foj,m,\nu} \in \kk(\Lambda_{\sigma_\fod})
\cap R^k_{\id_{\sigma_\fod},\sigma}$,
$n_{\foj,m,\nu}\in\Lambda_\fod^\perp$ and $A$ is a set of
representatives of $P_{\sigma_\fod,\sigma}/\Lambda_{\sigma_\fod}$.
\begin{enumerate}
\item
If $\sigma_\foj=\sigma$ and $-\overline m\in \Int
K_{\sigma_\fod}\sigma$ then $\sum_{\nu}a_{\foj,m,\nu}
\partial_{n_{\foj,m,\nu}}=0$.
\item Let
$\rho\in\P^{[n-1]}$, $\fod\subseteq\rho\subseteq\sigma$. 
\begin{enumerate}
\item If 
$\codim
\sigma_\fod\ge 2$ and $-\overline m\in \Int K_{\sigma_\fod}\rho$, then
$\sum_{\nu,\{\foj\supseteq\fod\,|\,\sigma_\foj=\rho\}}
a_{\foj,m,\nu}\partial_{n_{\foj,m,\nu}}=0$.
\item If $\codim\sigma_\fod=1$,
assume in addition that $\theta_{\gamma_\foj}=1$ for any joint
$\foj\supseteq\fod$ with $\codim\sigma_\foj =0$.
Then for any $m\in A$,
$\sum_{\nu,\{\foj\supseteq\fod\,|\,\sigma_\foj=\rho\}}
a_{\foj,m,\nu}\partial_{n_{\foj,m,\nu}}=0$.
\end{enumerate}
\item
Assume in addition that $\theta_{\gamma_\foj}=1$ for any joint
$\foj\supseteq\fod$ with $\codim\sigma_\foj \le1$.
\begin{enumerate}
\item
If $\codim \sigma_\fod=2$ and $\foj,\foj'$ are
the unique joints in $\sigma_\fod$ containing $\fod$ then
$\theta_{\gamma_{\foj'}}= \theta_{\gamma_\foj}^{-1}$.
\item
If $\codim \sigma_\fod=3$ and $\foj\supseteq\fod$ is a joint with
$\codim\sigma_\foj=2$ then $\theta_{\gamma_\foj}\in
O^k(\Lambda_\fod\otimes \Lambda_\fod^\perp)$. 
\end{enumerate}
\end{enumerate}
\end{proposition}

\proof 
We proceed inductively, proving the statement for any
$\sigma\in\P_\max$, $\sigma\supseteq\fod$ and those exponents $m$ with
$z^m\in J_l\setminus J_{l-1}$, $J_l= I_0^l\cdot I_{k-1}+ I_k$. For
$l=0$ there is nothing to prove.

The key ingredient is~(\ref{interstice relation}) with a particular
choice of $\gamma_i,\beta_i$. The additional requirement is that for
any $\tau\in\P$ containing $\fod$, the loops around joints $\foj$ with
$\sigma_\foj=\tau$ are numbered consecutively
$\gamma_i,\gamma_{i+1},\ldots,\gamma_{i+r}$, and are based on chambers
contained in the same maximal cell $\sigma(\tau)\supseteq\tau$;
furthermore,
\[
\beta_i \gamma_i\beta_i^{-1}
\ldots
\beta_{i+r} \gamma_{i+r}\beta_{i+r}^{-1}
\]
shall be freely homotopic to an edge path $\gamma_\tau$ passing
along the boundary of $\bigcup_{j=i}^{i+r} \hat \foj_j$ once. There
are two exceptional cases. First, if $\sigma_\fod=\tau=\rho\in
\P^{[n-1]}$ then $\bigcup_{j=i}^{i+r} \hat \foj_j$ is an annulus; in
this case we want a homotopy to an edge path first following one
boundary component, then an edge $\hat\fov$ to the other boundary
component, then following the other boundary component, and finally
back along $\hat\fov$. The other exceptional case occurs for
$\sigma_\fod=\tau\in \P^{[n-2]}$, where $\tau$ contains exactly two
joints $\foj,\foj'$, as in (3)(a). We then take $\gamma_\tau$ to
consist of the composition of two loops with some common base point,
denoted $\hat \fou_\tau$, and which go around $\hat \foj$ and
$\hat{\foj'}$, respectively. In any case, following $\gamma_\tau$
defines a log automorphism $\theta_\tau:=\theta_{\gamma_\tau}$ of
$R^k_{\id_{\sigma_\fod}, \sigma(\tau)}$.

We note at this point that with this selection of paths,
parts (2)(b) and (3)(a) of this proposition follow immediately from
\eqref{interstice relation}, observing that each $\theta_{\gamma_{\foj}}$
commutes with automorphisms attached to any wall containing $\fod$.

Continuing with the other cases, note that for any $\tau$ the
$\theta^k_{\gamma_\foj}$ with $\sigma_\foj=\tau$ commute mutually and with
any automorphism associated to a wall $\fop\subseteq \sigma(\tau)$.
This shows that
\begin{eqnarray}\label{expansion theta_tau}
\theta_\tau=\prod_{\{\foj\,|\, \sigma_\foj=\tau\}}
\theta_{\gamma_\foj}
=\exp\Big(\sum_{\big\{(m,\nu)\,\big|\,\genfrac{}{}{0 pt}{}{m\in A,\,
\ord_{\sigma_\fod}(m)=k}{\ord_\tau(m)=k\hfill}\big\}}
a_{\tau,m,\nu}  z^m\partial_{n_{\tau,m,\nu}} \Big),
\end{eqnarray}
as log automorphism of $R^k_{\id_{\sigma_\fod},\sigma(\tau)}$, with
$a_{\tau,m,\nu} \in \kk(\Lambda_{\sigma_\fod}) \cap
R^k_{\id_{\sigma_\fod}, \sigma(\tau)}$ and
$n_{\tau,m,\nu}\in\Lambda_\fod^\perp$. By~(\ref{step I reduction}) and
Proposition~\ref{exponentprop} the sum runs only over those $m$ with
$-\overline m\in K_{\sigma_\fod}\tau$. Note also that since the
elements of $A$ are not congruent modulo $\Lambda_{\sigma_\fod}$,
$\sum_{\nu}a_{\tau,m,\nu}n_{\tau,m,\nu}$ is uniquely determined.

In the proof special care has to be taken for codimension one joints
$\foj$, because these potentially involve monomials $z^m$ with the
only restriction $-\overline m\in K_{\sigma_\fod}\rho$,
$\rho=\sigma_\foj$. If $\codim\sigma_\fod=3$, these may interact with
terms arising from codimension two joints $\foj'\subseteq \partial
\rho$ in a way spoiling the induction process. We deal with this
problem as follows. If $\codim \sigma_\fod=3$ there are exactly two
cells $\tau_1\neq \tau_2$ of $\P$ of codimension $2$ with
$\fod\subseteq \tau_\mu\subseteq \rho$. Let $\fov_\mu\in
\P_{\scrS_k^\I}^{[n-1]}$ be the unique cell (support of a slab) with
$\fod\subseteq \fov_\mu\subseteq \rho$, and
$\dim\tau_\mu \cap\fov_\mu=n-2$, $\mu=1,2$. Note that~(\ref{slabs1})
from Step~I implies $\fov_1\neq\fov_2$. Thus $\hat \fov_\mu$ separates
the $2$-cell $\hat\foj_\mu\in\Sigma_\fod$, $\foj_\mu$ the unique joint
with $\fod\subseteq \foj_\mu\subseteq \tau_\mu$, from another $2$-cell
$\hat\foj$ with $\foj$ a joint with $\sigma_\foj=\rho$. Now change
$\theta(\fov_\mu)$, the log isomorphism associated to following the
edge $\hat\fov_\mu$ in the same direction as $\gamma_{\foj_\mu}$, by
composition with
\[
\theta_\mu:=\exp\Big(\sum_{\{(m,\nu)\,|\, m\in A,\, 
-\overline m\in \Int(K_{\sigma_\fod} 
\tau_\mu) \}} a_{\rho,m,\nu} z^m\partial_{n_{\rho,m,\nu}}\Big).
\]
This has the effect of composing $\theta_{\tau_\mu} =
\theta_{\foj_\mu}$ with $\theta_\mu$ and $\theta_\rho$ with
$\theta_\mu^{-1}$. Thus this change cancels all terms
$a_{\rho,m,\nu}\partial_{n_{\rho,m,\nu}}$ on the right-hand side of
(\ref{expansion theta_tau}) whenever $-\overline m\in (\partial
K_{\sigma_\fod} \rho)\setminus \Lambda_{\sigma_\fod}$. With this
reinterpreted $\theta(\fov_\mu)$ formula~(\ref{interstice relation}) still
holds, and the conclusions of the proposition remain unchanged. We
henceforth assume these terms do not arise in~(\ref{expansion
theta_tau}) for any $\tau\in\P^{[n-1]}$ in the first place.

After having established this property, if $\codim\sigma_\fod=3$, we
add to the induction hypothesis the following analogue of (1) and (2)(a)
for codimension two joints:
\begin{enumerate}
\item[(4)]
If $\tau\in\P^{[n-2]}$, $\fod\subseteq\tau\subseteq\sigma$ and
$-\overline m\in \Int K_{\sigma_\fod}\tau $ then
$\sum_{\nu} a_{\foj,m,\nu} \partial_{n_{\foj,m,\nu}}=0$ for the unique joint
$\foj\supseteq \fod$ with $\sigma_\foj=\tau$.
\end{enumerate}

Now for any $\tau\supseteq\fod$, if $\foj_i,\ldots, \foj_{i+r}$ are
the joints with $\sigma_\foj=\tau$ then by the above choice of
$\gamma_i,\ldots,\gamma_{i+r}$ there exists an edge path $\beta_\tau$
from $\hat\fou_0$ to the base point $\hat\fou_\tau$ of $\gamma_\tau$
such that
\[
\beta_\tau \gamma_\tau\beta_\tau^{-1}=
\beta_i \gamma_i\beta_i^{-1} \ldots
\beta_{i+r} \gamma_{i+r}\beta_{i+r}^{-1}.
\]
Therefore
\begin{eqnarray}\label{automorphism grouping}
\theta_{\beta_\tau}^{-1}\circ\theta_\tau\circ \theta_{\beta_\tau}
=\theta_{\beta_{i+r}}^{-1}\circ\theta_{\gamma_{i+r}} \circ
\theta_{\beta_{i+r}} \circ\ldots \circ \theta_{\beta_i}^{-1}\circ
\theta_{\gamma_i}\circ\theta_{\beta_i}.
\end{eqnarray}
In particular, we can now rewrite~(\ref{interstice relation}) in the
form
\begin{eqnarray}\label{interstice relation'}
\theta_{\beta'_{s'}}^{-1}\circ\theta_{s'}\circ \theta_{\beta'_{s'}}
\circ\ldots\circ
\theta_{\beta'_1}^{-1}\circ\theta_1\circ \theta_{\beta'_1}
=1,
\end{eqnarray}
where for any $i$ we have $\theta_i=\theta_\tau$,
$\beta'_i=\beta_\tau$ for some $\tau=\tau(i)$, and each
$\tau\supseteq\fod$ occurs exactly once.

We now do the inductive step from $l-1$ to $l$. The plan is to deduce
(1), (2)(a) and (4) for $m$ with $z^m\in J_l\setminus J_{l-1}$ by looking
at~(\ref{interstice relation'}) modulo $J_{l+1}$. For $\tau\supsetneq
\sigma_\fod$ consistency of $\scrS_k^\I$ to order $k-1$ and
the present induction hypothesis show
\begin{eqnarray}\label{theta_tau}
\qquad\theta_\tau=\exp\Big({\sum_{\{(m,\nu)\,|\,m\in A\}}}
a_{\tau,m,\nu}z^m\partial_{n_{\tau,m,\nu}}\Big)
\circ\exp\big(\sum_{\nu'}a_{\tau,\nu'}t^k
\partial_{n_{\tau,\nu'}}\big)\mod J_{l+1}
\end{eqnarray}
with $z^m\in J_l\setminus J_{l-1}$, $-\overline m\in \Int
K_{\sigma_\fod}\tau$ and $a_{\tau,\nu'},a_{\tau,m,\nu} \in
\kk(\Lambda_{\sigma_\fod})\cap R^k_{\id_{\sigma_\fod},\sigma(\tau)}$.
Note that if $\tau\in\P^{[n-2]}$ and $\codim\sigma_\fod=3$
then~(\ref{theta_tau}) follows from~(4). In
(\ref{interstice relation'}), $\theta_\tau$ occurs conjugated by the
log isomorphism associated to the edge path $\beta_\tau$. Now since
$I_0\cdot J_l\subseteq J_{l+1}$ the conjugation by a log automorphism
associated to crossing a wall does not have any effect modulo
$J_{l+1}$ and can thus be ignored in the following. For the
conjugation by a log automorphism associated to crossing a slab,
Lemma~\ref{interstice interaction} shows that likewise
\[
\theta_{\beta_\tau}^{-1}\circ\theta_\tau\circ \theta_{\beta_\tau}
=\exp\Big({\sum_{m,\nu}} a'_{\tau,m,\nu}z^m\partial_{n'_{\tau,m,\nu}}\Big)
\circ\exp\big(\sum_{\nu'}a'_{\tau,\nu'}t^k\partial_{n'_{\tau,\nu'}}\big)
\mod J_{l+1},
\]
with $z^m\in J_l\setminus J_{l-1}$, $-\overline m\in \Int
K_{\sigma_\fod}\tau$ and $a'_{\tau,\nu'}, a'_{\tau,m,\nu} \in
\kk(\Lambda_{\sigma_\fod})\cap R^k_{\id_{\sigma_\fod},\sigma(\tau)}$.
On the other hand, if $\tau=\sigma_\fod$ we readily obtain a similar
expansion without the first factor, that is, with all
$a'_{\tau,m,\nu}=0$. Thus in any case, for any $\tau, m,\nu$ with
$z^m\in J_l\setminus J_{l-1}$ the expression $a'_{\tau,m,\nu}z^m
\partial_{n'_{\tau,m,\nu}}$ can not cancel with any term from
$\theta_{\tau'}$ for any $\tau'\neq\tau$. In view of~(\ref{interstice
relation'}) we therefore conclude
\[
\theta_{\beta_\tau}^{-1}\circ\theta_\tau\circ \theta_{\beta_\tau}
=\exp\big(\sum_{\nu'}a'_{\tau,\nu'}t^k
\partial_{n'_{\tau,\nu'}}\big)\mod J_{l+1},
\]
for any $\tau$, with $a'_{\tau,\nu'}\in \kk(\Lambda_{\sigma_\fod})\cap
R^k_{\id_{\sigma_\fod},\sigma(\tau)}$. Hence also $\theta_\tau=
\exp(\sum_{\nu'}a_{\tau,\nu'}t^k\partial_{n_{\tau,\nu'}})$,  which in
turn gives that modulo $J_{l+1}$, $\sum_{m,\nu}
a_{\tau,m,\nu}z^m\partial_{n_{\tau,m,\nu}}=0$, where the sum is over
those $m$ such that $-\overline m\in\Int K_{\sigma_{\fod}}\tau$.
Expanding out the definition of $\theta_\tau$ now proves the claimed
formulae in (1), (2)(a) and (4) for $m$ with $z^m\in J_l\setminus
J_{l-1}$. This shows (1) and (2)(a) in the statement of the
proposition. (3)(b) follows easily from the fact that under the
additional hypothesis we only have $a_{\tau,m,\nu}\neq0$ if
$\tau\in\P^{[n-2]}$.
\qed
\medskip

We are now in position to prove consistency to order $k$ along
codimension zero joints.

\begin{proposition}\label{codimension 0 consistency}
$\scrS_k^\I$ is consistent to order $k$ at any joint $\foj$ with
$\codim \sigma_\foj=0$ provided $\sigma_{\foj}$ is bounded.
\end{proposition}

\proof
Let $\sigma\in\P_\max$. For a joint $\foj\subseteq\sigma$
intersecting $\Int(\sigma)$ a loop around $\foj$ defines a log
automorphism
\[
\theta^k_\foj= \exp\Big( \sum_i c_i z^{m_i} \partial_{n_i}\Big)
\]
of $R^k_{\id_\sigma, \sigma}$ with $\ord_\sigma(m_i)=k$ and
$\overline{m_i}\in\Lambda_\foj$ for all $i$. This depends only on
the sense of orientation of the loop, since changing the base
chamber leads to conjugation by automorphisms associated to walls,
and these only involve monomials of higher order.

Now fix an exponent $m$ on $\sigma:=\sigma_\foj$. We have to show
that $c_\foj^m:=\sum_{m_i=m} c_i n_i\in
\Lambda^*_{\sigma}\otimes\kk$ vanishes. This is clear if $\overline
m=0$ because $\theta^k_\foj\in H_\foj^{I_k}$. Otherwise $L_x=(x+\RR
\overline m)\cap\sigma$ for $x\in\Int\foj$ is a line segment,
varying in an $(n-3)$-dimensional family with parameter $x$. Since
the boundaries of interstices in $\sigma$ have dimension $n-4$ we
may choose $x\in\foj$ in such a way that the intersection of $L_x$
with any $\fov\in \P_{\scrS_k^\I}^{[n-3]}$ lies in $\Int\fov$. Then
there exist real numbers $0<\lambda_1<\lambda_2<\ldots <\lambda_s <
\infty$ and joints $\foj_1=\foj,\foj_2,\ldots, \foj_s$ with
$x+\lambda \overline m\in\Int \foj_i$ for $\lambda\in
(\lambda_{i-1},\lambda_i)$, and $s$ and $\lambda_s$ maximal with
this property. By the choice of $x$ we see that $x+\lambda_i
\overline m\in \foj_i\cap\foj_{i+1}$, for $1\le i<s$, must be an
interior point of an interstice $\fod_i$ intersecting $\Int\sigma$,
and $\overline m\not\in\Lambda_{\fod_i}$. In this situation only the
two joints $\foj_i$, $\foj_{i+1}$ containing $\overline m$ in their
tangent spaces contribute to the sum in~(\ref{interior interstice
cycle condition}). Thus Proposition~\ref{codim 0 step II lemma}
implies $c_{\foj_i}^m= c_{\foj_{i+1}}^m$, provided we orient the
loops around $\foj_i$ and $\foj_{i+1}$ in the same way.

Inductively we thus see $c_\foj^m=c_{\foj_s}^m$. Now the maximality of
$\lambda_s$ implies that $x+\lambda_s \overline m\in \partial \foj_s$
is an interior point of some $\fov\in\P_{\scrS_k^\I}^{[n-3]}$. 
If $\fov\subseteq\partial B$, noting in any event that $-\overline{m}
\in K_{\sigma_\fov}\sigma$, $\overline{m}$ cannot be in the support
of the fan $\Sigma_{\sigma_{\fov}}$, which is convex. So
$m\not\in P_{\sigma_\fov,\sigma}$, a contradiction.
Thus $\fov=\fod$ is an interstice, and $\overline m\not\in
\Lambda_{\fod}$. If $\fod\not\subseteq\partial \sigma$ we can run
Proposition~\ref{codim 0 step II lemma} again to conclude
$c_\foj^m=c_{\foj_s}^m=0$. If $\fod_s\subseteq\partial\sigma$,
$\fod_s\not\subseteq \partial B$, Proposition~\ref{joints at higher
codim interstices},(1) applies since $-\overline m\in
\Int_{\sigma_{\fod_s}}\sigma$. Hence $c_\foj^m=c_{\foj_s}^m=0$ also in
this case.
\qed

%===========================================================
\subsection{Step II: Homological modification of slabs}
The next step of the algorithm achieves consistency for codimension
one joints. At a single joint this can be done by modifying the
functions associated to the two slabs containing this joint. There is then
a problem of whether or not this can be done consistently, as changes
to slabs dictated by one joint may conflict with changes to slabs
dictated by another joint. Furthermore, monomials in the functions
associated to slabs do not propagate in the same way monomials
associated to walls do, because of~(\ref{change of vertex - slabs}).
Thus it is impossible to emulate the arguments used for codimension
zero joints, and instead, we need to use homological arguments which
will fix the corrections to all slabs in a given codimension one
$\rho\in\P$ simultaneously.

Throughout the following discussion we therefore fix
$\rho\in\P^{[n-1]}$ and a reference cell $\sigma\in\P_\max$,
$\rho\subseteq \sigma$. The discussion in this subsection will only
apply to the case that $\rho$ is bounded. To fix signs orient $\sigma$
and each joint $\foj\subseteq \rho$ arbitrarily. This distinguishes a
sense of orientation of each loop around any joint in $\rho$ that we
tacitly assume in the following. Let $\check d_{\rho}\in
\Lambda_\rho^\perp \simeq\ZZ$ be the generator that is positive on
$\sigma$. If $\foj$ is a joint with $\Int\foj \subseteq \rho$, then
by ~(\ref{step I reduction}) the log automorphism $\theta^k_\foj$ of
$R^k_{\id_\rho,\sigma}$ associated to a loop around $\foj$ lies in
$O^k \big(\Lambda_\foj \otimes\Lambda_\foj^\perp +(\Lambda_\rho/
f_\rho)\otimes \Lambda_\rho^\perp\big)$. Thus for any vertex
$v\in\sigma_\foj=\rho$ we may write
\[
\theta^k_\foj=\exp\Big( {\textstyle\sum_i c_i  z^{m_i}\partial_{n_i}}
+ \frac{\sum_i d_i z^{m'_i}}{f_{\rho,v}}\partial_{\check d_{\rho}}\Big)
\]
with $m_i\in\Lambda_\foj$, $m'_i\in\Lambda_\rho$ and $\ord_\rho(m_i)
=\ord_\rho(m'_i)= k$. The expression $\sum_i d_i z^{m'_i}$ is unique
up to adding multiples of $f_{\rho,v}$. Note that in view of
(\ref{f_rho}) a different choice $v'$ of vertex leads to the
expression $\sum_i d_i z^{m'_i+m^\rho_{v'v}}$.

Thus $\sum_i d_i z^{m'_i}$ defines a well-defined element
$d_\foj=(d_{\foj,v})_v$ in the $\kk$-vector space $W_\rho$ that is
defined as follows. For a vertex $v\in\rho$ let
\[
W_{\rho,v}:=\tilde W_{\rho,v}/ \big(\tilde W_{\rho,v}\cdot f_{\rho,v}\big),
\]
with
\[
\tilde W_{\rho,v}:= \big\{ {\textstyle \sum_i a_i z^{m_i}} \in
\kk[P_{\rho,\sigma}]\,\big|\,
\overline{m_i}\in\Lambda_\rho,\, \ord_\rho(m_i)=k \big\}.
\]
Note that $\tilde W_{\rho,v}\cdot f_{\rho,v}\subseteq \tilde
W_{\rho,v}$ because $f_{\rho,v}$ involves only exponents $m$ with
$\ord_\rho(m)=0$ and $\overline m\in\Lambda_\rho$. For another vertex
$v'\in\rho$ we have the isomorphism
\[
W_{\rho,v}\lra W_{\rho,v'},\quad
h\longmapsto h\cdot z^{m^\rho _{v'v}}.
\]
Then $W_\rho$ is defined as the set of tuples $(h_v)_{v\in\rho}$ with
$h_v\in W_{\rho,v}$ and
\[
h_{v'}= z^{m^\rho_{v'v}}\cdot h_v.
\]

The plan is now to achieve $d_\foj=0$ by changing the slabs contained
in $\rho$. Fix a vertex $v\in\rho$ and let $\foB\subseteq
P_{\rho,\sigma}$ be a set of exponents such that $\big(
z^{m^\rho_{v'v} +m}\big)_{v'}$, $m\in\foB$, forms a basis of $W_\rho$.
For $m\in \foB$ and $x\in \rho\setminus\Delta$ write $m[x]$ for the
parallel transport of $m+m^\rho_{v'v}$ to $x$ inside
$\rho\setminus\Delta$, where $v'=v[x]\in\rho$. Thus we can write
\begin{eqnarray}\label{d_foj^m}
d_{\foj,v'}= \sum_{m\in\foB} d_\foj^m z^{m[v']}
\end{eqnarray}
for some $d_\foj^m\in\kk$. We now follow a procedure to adjust the
functions $f_{\fob,x}$, for $\fob\subseteq \rho$ a slab, by multiples
of $z^{m[x]}$ for a single $m\in\foB$. Write $\P_\rho$ for the polyhedral
decomposition of $\rho$ given by $\P_{\scrS_k^\I}$. For the function
$f_{\fob,x}$ of a slab $\fob\subseteq\rho$ to receive a correction by
a multiple of $z^{m[x]}$ we require that
\begin{enumerate}
\item
$m[x]\in P_x$,
\item
for any joint $\foj\subseteq\fob$ with $x\in\foj$ and any
$\sigma'\in\P_\max$ with $\sigma'\supseteq\foj$
it holds $\ord_{\sigma'}(m[x])\ge k$.
\end{enumerate}
The second condition ensures that such a correction does not influence
lower order computations. Thus we consider only the polyhedral complex
$\P_m\subseteq \P_\rho$ defined as the complement of the open star of
\begin{eqnarray}\label{bad locus in rho}
\Omega:=\big\{ x\in\partial\rho\,\big|\, m[x]\not\in P_x\big\}
\cup\textstyle\bigcup_\foj \Int\foj \subseteq \partial\rho,
\end{eqnarray}
where the union runs over all joints $\foj\subseteq\partial\rho$ such
that $\ord_{\sigma'}(m[x])<k$ for some $\sigma'\in\P_\max$ containing
$\foj$ and $x\in \foj\setminus\Delta$. Said differently, $\P_m$
consists of all cells of $\P_\rho$ not intersecting~$\Omega$.

On the other hand, a codimension two joint $\foj\subseteq\partial\rho$
does not impose any conditions on changing a slab function
$f_{\fob,x}$, $\fob\supseteq\foj$, by $c z^{m[x]}$ if
(i)~$\ord_\foj(m[x])>k$ or (ii)~$\overline{m[x]} \in\Lambda_\foj$,
$x\in (\Int\foj)\setminus\Delta$. In fact, Proposition~\ref{influence
of monomial changes},(1) shows that such changes keep the form
(\ref{step I reduction}) for $\theta^k_{\foD_\foj}$. We therefore work
relative to the subcomplex $\tilde{\scrA_m} \subseteq \P_\rho$
consisting of faces of $(n-2)$-cells $\foj\subseteq\partial \rho$ of
$\P_\rho$ obeying (i) or (ii). ($\tilde\scrA_m$ may include $(n-2)$-cells
$\foj\subseteq\partial B$, which are not joints, but these do not
impose conditions anyway.) Note that by
Proposition~\ref{exponentprop}, an $(n-2)$-cell $\foj\subseteq
\partial\rho$ of $\P_\rho$ is contained in $\tilde A_m$ if and only if
$\overline m$ is contained in the half plane tangent wedge to $\foj$
in $\rho$. In particular, the underlying topological space $\tilde
A_m\subseteq \partial\rho$ of $\tilde \scrA_m$ is a union of facets of
$\rho$, and an alternative description of $\tilde A_m$ is
\begin{eqnarray}\label{tilde A_m}
\tilde A_m=
\cl \big\{x\in\partial\rho\setminus\Delta  \,\big|\,
\overline{m[x]}\in K_x\rho \big\}.
\end{eqnarray}
Finally denote
\[
\scrA_m:= \tilde{\scrA_m}\cap \P_m.
\]
Note that $\tilde A_m\cap \Omega$ is contained in the relative
boundary of $\tilde A_m$, and hence $\scrA_m$ is also obtained by
removing the open star of a subset of $\partial\tilde\scrA_m$. In
fact, if $\fov\in\P_\rho$ is contained in the relative interior of
$\tilde A_m$ then $\overline{m[x]}\in K_x\rho$ for all $x\in
\fov\setminus\Delta$. This implies $\ord_{\sigma'}(m[x])\ge
\ord_\rho(m[x])=k$ for all $x\in\fov\setminus\Delta$ and
$\sigma'\in\P_\max$ containing $x$, and hence $\fov\in \P_m$.

Our interest in $(\P_m,\scrA_m)$ comes from the following result.

\begin{lemma}\label{chain lemma}
The cellular $(n-2)$-chain $(d_\foj^m)_{\foj\in\P_m^{[n-2]}}$ with
$d_\foj^m=0$ for $\foj\subseteq \partial\rho$ and as in
(\ref{d_foj^m}) otherwise,  is a relative cycle for $(\P_m,\scrA_m)$.
\end{lemma}

\proof
Orient each interstice $\fod\subseteq \rho$ arbitrarily. Then for any
$\fod\subseteq \foj \subseteq \rho$ the comparison of the chosen
orientation of $\fod$ with the one induced from $\foj$ defines a sign
$\sgn(\fod,\foj)\in\{\pm1\}$ such that the coefficient of $\fod$ in the
boundary of a cellular $(n-2)$-chain $(c_\foj)_\foj$ is
\[
\sum_{\foj\supseteq\fod} \sgn(\fod,\foj) c_\foj.
\]

Now let $\fov\in\P_m^{[n-3]}$. Then either $\fov\subseteq\partial B$
or $\fov=\fod$ is an interstice. In the first case, if
$\fov\not\in\scrA_m$, then by~(\ref{tilde A_m}),
$\overline{m[x]}\not\in K_x\rho$ for $x$ in a neighbourhood of
$\fod$ in $\partial\rho$.
n particular, $-\overline{m[x]}\in\Int
K_{\sigma_\fov}\rho$, but also $\overline{m[x]}$ maps to
$|\Sigma_{\sigma_{\fov}}|$ since $m[x]\in P_x$. This contradicts
convexity of $|\Sigma_{\sigma_{\fov}}|$. Thus $\fov\in \scrA_m$ and
there is nothing to check.

In the case of an interstice  with $\fod\subseteq\partial\rho$ and
$\fod\not\in \scrA_m$ we have $-\overline{m[x]}\in\Int
K_{\sigma_\fod}\rho$ as before. Proposition~\ref{joints at higher
codim interstices},(2a) now shows $\sum_{\foj\supseteq \fod,\,
\foj\not\in \scrA_m} \sgn(\fod,\foj) d_\foj^m=0$. The sign arises
from the difference in orientation conventions for loops around
joints.

If $\fod\not\subseteq\partial\rho$, Proposition~\ref{joints at higher
codim interstices},(2b) implies $\sum_{\foj\supseteq \fod}
\sgn(\fod,\foj) d_\foj^m=0$,  again observing the different
orientation conventions.

\qed
\medskip

We now prove two lemmas concerning the topology of this
situation.

\begin{lemma}\label{retraction lemma}
The pair $(\P_m, \scrA_m)$ is a deformation retract of
$(\P_\rho,\tilde{\scrA_m})$.
\end{lemma}

\proof
We want to retract the cells in $\P_\rho\setminus\P_m$ successively,
in a way compatible with $\tilde\scrA_m$. For this we use the
following elementary result: If $\Xi\subseteq \RR^k$ is a bounded
convex polytope and $x\in\partial\Xi$ then the projection from a
point $x'\in\RR^k\setminus\Xi$ sufficiently close to $x$ and with
$x-x'\in K_x\Xi$ defines a deformation retraction of $\Xi$ onto the
union of facets of $\Xi$ not containing $x$. Explicitly, for
$y\in\Xi$ define
\[
\alpha(y)=\max\{\alpha\in\RR_{\ge0}\,|\, y+\alpha\cdot(y-x')\in\Xi\}.
\]
Then
\[
[0,1]\times\Xi \lra \Xi,\quad
(\lambda,y)\longmapsto y+\lambda\alpha(y) (y-x'),
\]
is the desired deformation retraction.

We apply this result first to successively retract $(\P_\rho,
\tilde\scrA_m)$ to $(\P_m\cup\tilde \scrA_m, \tilde\scrA_m)$: Let
$\P'\subseteq\P_\rho$ be a subcomplex obtained inductively. Let
$\P'_\partial\subseteq\P'$ consist of subcells of cells $\fow\in
\P'\setminus(\P_m\cup \tilde\scrA_m)$ with the property that there
is a unique $\fov\in\P'\setminus(\P_m\cup \tilde\scrA_m)$ with
$\fow\subsetneq\fov$. This $\P'_\partial$ is the subset of cells
that can be taken as center for the next retraction. We will assume
inductively that $\P'_{\partial} \not=\emptyset$ as long as
$\P'\not=\P_m\cup\tilde\scrA_m$. To see this is true initially, note
first that if  $(\P_m\cup\tilde\scrA_m) \cap\partial\P_\rho
=\partial\P_\rho$ then $\P_m =\P_{\rho}$ and
$\tilde\scrA_m=\partial\P_{\rho}$  anyway and there is nothing to
do. Otherwise, there is a slab $\fob\in\P_{\rho}\setminus
(\P_m\cup\tilde\scrA_m)$  with $\dim\fob\cap\partial \rho=n-1$, and
then $\fob\cap\partial\rho\in\P'_{\partial}$.

Given $\P'_{\partial}\not=\emptyset$, choose a point $x$ in the
interior of a maximal cell $\fow\in\P'_{\partial}$, contained
properly in a unique cell $\fov\in \P'\setminus (\P_m\cup
\tilde\scrA_m)$. Now apply the above deformation retraction of
$\fov$ using the chosen $x$. Since $x$ is disjoint from any cell of
$\P_m\cup\tilde\scrA_m$ this deformation retraction is trivial on
this subcomplex of $\P'$.  We now note that after making this
retraction, we continue to have  $\P'_{\partial}\not=\emptyset$. In
fact, if $\tilde\fow\in\P'\setminus (\P_m\cup \tilde\scrA_m)$, then
$\fow:=\tilde\fow\cap\Omega$ is a non-empty cell of $\P_\rho$.
Moreover, by the inductive construction the link of $\fow$ in $\P'$
is a retraction of the link of $\fow$ in $\P_{\rho}$. From this one
can see that the link of $\fow$ contains a cell in $\P'_{\partial}$,
and hence $\P'_{\partial}$ continues to be non-empty. The process
stops when $\P'= \P_m\cup\tilde\scrA_m$. An analogous argument
then retracts $(\tilde\scrA_m,\scrA_m)$ onto $(\scrA_m,\scrA_m)$,
and hence $(\P_m\cup\tilde\scrA_m, \tilde\scrA_m)$ onto
$(\P_m,\scrA_m)$.
\qed

\begin{lemma}\label{H_{n-2}=0}
$H_{n-2}(\P_m,\scrA_m)=0$.
\end{lemma}

\proof
By Lemma~\ref{retraction lemma} this follows once we prove
$H_{n-2}(\P_\rho,\tilde{\scrA_m})=0$. Let $\Delta(\rho)\subseteq
\Lambda_{\rho,\RR}$ be the convex hull of $\{m^\rho_{vv'}\,|\,
v'\in\rho \text{ vertex}\}$ and $\psi_{\check\rho}$ be the
corresponding PL-function on the normal fan $\check\Sigma_\rho$ of
$\rho$. Recall from~(\ref{tilde A_m}) that $\tilde A_m$ is the union
of facets $\tau\subseteq\rho$ with $\overline{m[x]} \in K_x\rho$ for
some $x\in\Int\tau$. If $n\in\Lambda_\rho^*$ is the inward normal
vector to $\tau$ generating the ray in $\check\Sigma_\rho$ dual to
$\tau$, this is equivalent to
\[
0\le \langle \overline{m[x]},n \rangle = \langle \overline{m},n \rangle +
\langle m^\rho_{v[x]v}, n\rangle = \langle \overline{m},n \rangle
+\psi_{\check\rho}(n).
\]
Thus $\tilde A_m\subseteq\partial\rho$ is dual to the subset of rays
of $\Sigma_{\check\rho}$ on which the convex function $\psi_{\check
\rho}+\overline m$ is non-negative. Thus if
$\overline{m}\in\Delta(\rho)$ we obtain $\tilde A_m =\partial\rho$,
and otherwise $\rho$ deformation retracts to $\tilde A_m$. In any case
it follows $H_{n-2}(\P_\rho,\tilde{\scrA_m}) = H_{n-2}(\rho, \tilde
A_m)=0$.
\qed

\begin{shaded}
\noindent II.1. \emph{First homological modification of slabs.} By
Lemma~\ref{H_{n-2}=0} we can find an $(n-1)$-chain
$(b^m_{\fob})_{\fob\in \P_m^{[n-1]}}$ whose boundary is
$(d^m_{\foj})_\foj$ modulo chains in $\scrA_m$. Then for any slab
$\fob\subseteq\rho$ subtract the term $D(s_{e},\rho,v)
s_{e}(b^m_{\fob} z^{m[x]})$ from $f_{\fob,x}$, where $v=v[x]$ and
$e:v\to\rho$. By construction of $\P_m$ we have $m[x]\in P_x$ for all
$x\in\fob\setminus\Delta$ and the change of vertex formula
(\ref{change of vertex - slabs}) continues to hold, so this makes
sense. Proposition~\ref{influence of monomial changes},(1) shows that
doing so removes the terms $d^m_{\foj}z^{m[v]}/f_{\rho,v}$ from
$\theta^k_\foj$,  whenever $\foj\in \P_m$, $\foj\not\subseteq
\partial\rho$. Furthermore, if $m[v]$ does appear in $\theta^k_\foj$
for some joint $\foj\subseteq\rho$, $\foj\not \subseteq\partial\rho$,
then $\foj\in \P_m$, so this process removes the term involving $m[v]$
from $\theta^k_\foj$ whenever such a term appears.

Repeat this for all exponents $m\in\foB=\foB(\rho)$ and for all
$\rho\in\P^{[n-1]}$.
\end{shaded}

For $\foj\subseteq\rho$ with  $\codim\sigma_{\foj}=1$, we can now
write
\[
\theta^k_\foj=\exp\left(\sum c_iz^{m_i}\partial_{n_i}\right)
\]
with $m_i\in\Lambda_\rho$, $n_i\in\Lambda_\foj^\perp$. Next we would
like to achieve $\theta^k_\foj\in O^k(\Lambda_\foj\otimes
\Lambda_\foj^\perp)$. This is possible by a further, straightforward
modification of slabs.

\begin{shaded}
\noindent II.2. \emph{Further subdivision of slabs to achieve
$\theta^k_\foj\in O^k(\Lambda_\foj\otimes\Lambda_\foj^\perp)$.}
For every $\foj$ with $\codim\sigma_\foj=1$ and every $m_i$ appearing
in $\theta^k_\foj$ with $\overline{m_i}\not\in\Lambda_\foj$, we note $n_i$ must
be proportional to $\check d_{\rho}$, and thus we can assume after
changing $c_i$ that $n_i=\check d_{\rho}$. Viewing $\foj\subseteq\rho$
as a subset of $\Lambda_{\rho,\RR}$, define 
\[
\fob(m_i):=(\foj-\RR_{\ge 0} \overline{m_i})\cap\rho\subseteq\rho,
\]
We then modify slabs contained in $\rho$ by adding  $\pm
c_i  s_e(m_i)z^{m_i}f_e$ to $f_{\fob,y}$ for $y\in
\fob(m_i)\cap\fob\setminus\Delta$, $\fob\subseteq\rho$ a slab and
$e:v[y]\to\sigma_\foj=\rho$. This
of course might mean subdividing the slabs further. By doing so, it
follows again from Proposition~\ref{influence of monomial
changes},(1), that, with proper choice of sign, the term
$c_iz^{m_i}\partial_{n_i}$ disappears from $\theta^k_\foj$. Note this
process might introduce new joints $\foj'$ contained in the sides of
$\fob(m_i)$, but for such $\foj'$, $\overline{m_i'}\in
\Lambda_{\foj'}$, so by a simple calculation analogous to
Proposition~\ref{influence of monomial changes},(2), $\theta^k_{\foj'}$
satisfies
\[
\theta^k_{\foj'}=\exp\left(\sum c_iz^{m_i}\partial_{n_i}\right).
\]
with $\overline{m_i}\in\Lambda_{\foj'}$. In fact, after carrying this
out for every joint $\foj$ with $\codim\sigma_{\foj}=1$, we 
see this now holds for all joints with $\codim\sigma_{\foj}=1$. We
write $\scrS_k^{\II,\pre}$ for the structure thus obtained.
\end{shaded}

The arguments of Proposition~\ref{codimension 0 consistency} for
codimension zero joints now imply that the remaining terms of
$\theta_\foj^k$ are undirectional. Recall that a log automorphism lies
in $O^k(0\otimes \Lambda_\foj^\perp)$ if it is of the form
$\exp(\sum_n a_n\partial_n)$ with $a_n\in\kk[t]$ and
$n\in\Lambda_\foj^\perp$.

\begin{proposition}\label{codimension 0 undirectional}
For $\foj\in \Joints(\scrS_k^\II)$ with $\sigma_\foj=\rho\in
\P^{[n-1]}$ it holds
\[
\theta_\foj^k\in O^k(0\otimes\Lambda_\rho^\perp).
\]
\end{proposition}
\proof
The construction of $\scrS_k^{\II,\pre}$ was designed to achieve
$\theta^k_{\foj} \in O^k (\Lambda_\foj \otimes \Lambda_\foj^\perp)$
for codimension one joints. It then follows exactly as in
Proposition~\ref{codimension 0 consistency} that there are no
contributions of exponents $m$ with $\overline m\neq 0$. In this
argument Proposition~\ref{joints at higher codim interstices},(2)
replaces both Proposition~\ref{codim 0 step II lemma} and
Proposition~\ref{joints at higher codim interstices},(1). Thus
$\theta_\foj^k \in O^k(0\otimes \Lambda_\foj^\perp)$.

To see that then even $\theta_\foj^k\in
O^k(0\otimes\Lambda_\rho^\perp)$ we have to show that
$\theta_\foj^k(m) =1$ for all $m\in\Lambda_\rho$. This follows
easily with the notion of $\tlog$ that comes out
naturally of our discussion of the higher order normalization
procedure in Step~III. We therefore postpone the rest of the proof
to \S\ref{par:Step III}, see after Lemma~\ref{tlog of compositions}.
\qed
\medskip

To remove the remaining undirectional terms we now run a homological
argument again. For $\rho\in\P^{[n-1]}$ and a joint
$\foj\in\scrS_k^{\II,\pre}$ with $\sigma_\foj=\rho$,
Proposition~\ref{codimension 0 undirectional} shows we can write
uniquely
\[
\theta_\foj^k = \exp\big(c_\foj t^k \partial_{\check d_\rho}\big),
\]
for some $c_\foj\in\kk$. Let $\fod$ be an interstice of
$\scrS_k^{\II,\pre}$ with  $\sigma_\fod=\rho$ and
$\foj_1,\ldots,\foj_r\subseteq\rho$ be the codimension one joints
containing $\fod$. Then the $\theta^k_{\foj_i}$ commute and
(\ref{interstice relation}), now interpreted for $\scrS_k^{\II,\pre}$,
implies
\[
1=\theta^k_{\foj_r}\circ\ldots\circ\theta^k_{\foj_1}
=\exp\Big(\sum_{i=1}^r \sgn(\fod,\foj_i) c_{\foj_i}t^k
\partial_{\check d_{\rho}} \Big).
\]
Here we use the signs $\sgn(\fod,\foj_i)$ introduced above. This shows
$\sum_{i=1}^r \sgn(\fod,\foj_i) c_{\foj_i}=0$, and hence
$(c_\foj)_{\foj\in\P_\rho^{n-2}}$ can be viewed as a relative cellular
$(n-2)$-cycle for $(\P_\rho,\partial\P_\rho)$. Again we simply set
$c_\foj=0$ whenever $\foj\subseteq\partial\P_\rho$.

\begin{shaded}
\noindent II.3. \emph{Second homological modification of slabs.}
Since $H_{n-2}(\P_\rho,\partial\P_\rho)=0$ there exists an
$(n-1)$-chain $(b_\fob)_{\fob\in \P_\rho^{n-1}}$ with boundary
$(c_\foj)_\foj$. Then for a slab $\fob\subseteq\rho$ and
$x\in\fob\setminus \Delta$ add $ b_\fob t^k
f_e$ to $f_{\fob,x}$, for $v:=v[x]$ and $e:v[x]\to\rho$.
Proposition~\ref{influence of monomial changes},(1) now shows that
after these changes $\theta_\foj^k=1$ holds for any codimension one
joint $\foj\subseteq \rho$. Repeat this process for all
$\rho\in\P^{[n-1]}$. We denote the structure thus obtained
$\scrS_k^\II$.
\end{shaded}

We have now arrived at a structure $\scrS_k^\II$ that is consistent up
to codimension one. Moreover, for codimension two joints essentially
the same arguments together with local rigidity
(Definition~\ref{def:locally rigid}) give further restrictions.

\begin{proposition}\label{step II reduction}
Let $\foj\in\Joints(\scrS_k^\II)$.
\begin{enumerate}
\item
If $\codim\sigma_\foj\le 1$ then $\theta_\foj^k=1$.
\item
If $\codim\sigma_\foj=2$ let $\tau=\sigma_\foj$ and $v\in\tau$ a
vertex. Then we can write
\[
\theta_\foj^k=\exp\Big( ct^k\partial_n
+{\sum}_{\hspace{-1ex}\begin{array}{ll}
\scriptstyle\rho\supseteq\foj\\[-1ex]
\scriptstyle v'\in\tau
\end{array}}\hspace{-1ex}
c_{\rho,v'} t^k\frac{z^{m^\rho_{vv'}}}{f_{\rho,v}}\partial_{\check d_{\rho}} \Big),
\]
with $c,c_{\rho,v'}\in \kk$ and $n\in \Lambda_\tau^\perp$.
\end{enumerate}
\end{proposition}
\proof
1)\ As we have not changed anything at joints $\foj$ of $\scrS_k^\I$
with $\codim\sigma_\foj=0$ this case follows from
Proposition~\ref{codimension 0 consistency}, while the constructions
in Step~II were designed to achieve $\theta_\foj^k=1$ if
$\codim\sigma_\foj=1$.\\[1ex]
2)\ By the definition of the polyhedral complex $\P_m\subseteq
\P_\rho$ and Proposition~\ref{influence of monomial changes},(1) the
changes from $\scrS_k^\I$ to $\scrS_k^\II$ do not affect the
form~(\ref{step I reduction}) of $\theta^k_\foj$ at codimension two
joints.

Now note that the joints $\foj\subseteq\tau$ are the maximal cells of
the polyhedral decomposition of $\tau$ given by $\P_{\scrS_k^\II}$.
Thus if $\fod$ is an interstice with $\sigma_\fod=\tau$ there are
precisely two codimension two joints $\foj,\foj'\subseteq\tau$
containing $\fod$. In this case Proposition~\ref{joints at higher
codim interstices},(3)(a) shows that $\theta^k_\foj\circ
\big(\theta^k_{\foj'}\big)^{-1}=1$, assuming the normal spaces
$\shQ^v_{\foj,\RR}=\shQ^v_{\tau,\RR} =\shQ^v_{\foj',\RR}$ are oriented
in the same way. Thus  all $\theta^k_\foj$ with $\foj\subseteq\tau$
agree. Thus there exist $c_j, d_{\rho,j}\in \kk$ and $m_j,m'_{\rho,j}$ with
$\overline{m_j}, \overline{m'_{\rho,j}}\in \Lambda_\tau$, $\ord_\tau
(m_j)= \ord_\tau(m'_{\rho,j})=k$, $n_j\in\Lambda_\tau^\perp$, such
that for any $\foj\subseteq\tau$
\[
\theta_\foj^k= \exp\Big( \sum_ j c_j z^{m_j}\partial_{n_j} +
\sum_{j,\,\rho\supseteq\tau} d_{\rho,j} \frac{z^{m'_{\rho,j}}}{f_{\rho,v}}
\partial_{\check d_{\rho}}\Big).
\]
Now by looking at an interstice $\fod\subseteq \partial\tau$
containing $v$, Proposition~\ref{joints at higher codim
interstices},(3)(b) implies that we may assume all $\overline{m_j}$
and $\overline{m'_{\rho,j}}$ to be contained in the half-plane tangent
wedge $K_\fod\tau$. Indeed, otherwise their contribution vanishes by
the proposition. Taking a different vertex $v'$ transforms
$z^{m'_{\rho,j}}/ f_{\rho,v}$ into $z^{m'_{\rho,j}- m^\rho_{v v'}}/
f_{\rho,v'}$. Thus we may take $c_j=0$ for all $j$, and $d_{\rho,j}=0$
unless $\overline{m'_{\rho,j}}- m^\rho_{v v'}\in K_{v'}\tau$ for every
vertex $v'\in\tau$. As in the proof of Lemma~\ref{H_{n-2}=0} one sees
that, in the latter case, $\overline{m'_{\rho,j}}$ must be contained
in the convex hull of $\{m^\rho_{v v'}\,|\, v'\in\tau\text{
vertex}\}$. This is a face of $\Delta(\rho)$, and by
Definition~\ref{def:locally rigid},(i) any integral point of
this face is a vertex. Hence $\overline{m'_{\rho,j}}=
m^\rho_{vv'}$ for some $v'\in\tau$. Because $\ord_\tau(m'_{\rho,j})=k$
it then follows
\[
z^{m'_{\rho,j}}= t^k\cdot z^{m^\rho_{v v'}}.
\]
This proves the claimed formula for $\theta_\foj^k$.
\qed

%===========================================================
\subsection{Step III: Normalization}\label{par:Step III}
For a joint $\foj$ of $\scrS_k^\II$ the remaining terms in
$\theta_\foj^k$ do not propagate --- either they are undirectional
($z^m$ with $\overline m=0$), or they are of the form $z^m/f_{\rho,v}$
and $-\overline m$ points into the tangent wedge of $\rho$ at $v$, for
any choice of vertex $v\in\rho$. Step~III removes these terms by a
normalization procedure.

The normalization condition asks that there be no pure $t$-terms in
the logarithm of the functions $s_e^{-1}(f_{\fob,x})$ that occur in
changing chambers separated by a slab $\fob$, up to order $k$. Because
this expression may contain exponents $m$ with $\ord_\rho(m)=0$ we
need to take appropriate completions of our rings $R^k_{g,\sigma}$ to
make sense of the logarithm.

Again, in this section, we assume all cells of $B$ are bounded.

\begin{construction}
Let $(g:\omega\to\tau)\in\Hom(\P)$, $\sigma\in\P_\max$ a reference
cell containing $\tau$ and $v\in\tau$ a vertex. Because $K_v\tau$ is
a strictly convex cone, the subset
\[
E=\big\{m\in P_{\omega,\sigma} \,\big|\, \overline m\in
K_v\tau\setminus\{0\} \big\},
\]
of $P_{\omega,\sigma}$ is additively closed and
\[
\bigcap_{\nu\ge0} \nu E=\emptyset.
\]
We can thus define a Hausdorff topology on $\kk[P_{\omega,\sigma}]$ by
taking
\[
U_\nu :=\big\{\textstyle\sum_{m\in \nu E} a_m z^m \in 
\kk[P_{\omega,\sigma}]\big\},\qquad \nu\ge 1,
\]
as fundamental neighbourhood system of $0$. Note that if
$\Lambda_\omega\cap(K_v\tau\setminus\{0\})\neq \emptyset$ then $E$
generates the unit ideal, and hence this is not the $I$-adic topology
for any ideal $I$. Denote the completion of this topological ring by
$\kk\lfor P_{\omega,\sigma}\rfor_v$.

Similarly, one defines a topology on $R^k_{g,\sigma}$, with associated
completion ${\vphantom{R^k}}^v\!\widehat{R}^k_{g,\sigma}$. Because the
localizing functions $f_{\rho,v}$ are invertible in $\kk\lfor
P_{\omega,\sigma}\rfor_v/ \hat I_{g,\sigma}^{>k}$, for any
$\rho\supseteq \tau$, we have
\[
{\vphantom{R^k}}^v\!\widehat{R}^k_{g,\sigma} =\kk\lfor P_{\omega,\sigma}\rfor_v
/\hat{I}^{>k}_{g,\sigma},
\]
where $\hat{I}^{>k}_{g,\sigma}\subseteq\kk\lfor
P_{\omega,\sigma}\rfor_v$ is the ideal generated by
$I^{>k}_{g,\sigma}\subseteq \kk[ P_{\omega,\sigma}]$.
\qed
\end{construction}

\begin{definition}
Let $f=\sum_{m\in P_{\omega,\sigma}} a_m z^m \in
\kk\lfor P_{\omega,\sigma}\rfor_v/ \hat I_{g,\sigma}^{>k}$.\\[1ex]
1)~The \emph{$t$-content} of $f$ is defined by
\[
\tcont f:= \sum_{m\in P_{\omega,\sigma}: \overline m=0} a_m z^m
=\sum_{m\in P_{\omega,\sigma}: \overline m=0} a_m t^{\ord_\sigma(m)}
\in \kk\lfor t\rfor/(t^{k+1}).
\]
2)~We say that $f$ \emph{fulfills the cone condition} if $a_m\neq 0$,
$\ord_\omega(m)=0$ implies $\overline m\in K_v\tau$.\\[1ex]
3)~If $a_0\neq0$ and $f$ fulfills the cone condition we define
\[
\tlog f:=\tcont\bigg(-\sum_{i=1}^\infty \frac{1}{i}
\Big(1-\frac{f}{a_0}\Big)^i  \bigg)\in \kk\lfor t\rfor/(t^{k+1}).
\]
If $f\in\kk[P_{\omega,\sigma}]$ we indicate which completion to work
in by writing $\tlog_v f$.
\end{definition}

\begin{remarks}
\label{rem:logadd}
The sum in the definition of $\tlog f$ is the power series
of $\log f$ without the constant term. Thus since $\tcont$ is
additive, the usual power series identity implies
\[
\tlog\big(f_1\cdot f_2\big)= \tlog f_1 +\tlog f_2\mod t^{k+1},
\]
for $f_i\in \kk\lfor P_{\omega,\sigma}\rfor_v/ \hat I_{g,\sigma}^{>k}$
with non-vanishing constant terms and fulfilling the cone condition.
\end{remarks}

We are now in position to formulate the normalization condition.

\begin{definition}\label{def:slab normalization}
A slab $\fob$ is called \emph{normalized to order $k$} if for any
$x\in\fob\setminus \Delta$ and $v'\in\rho_\fob$ a vertex, it holds
\[
\tlog_{v'}\big( z^{m^{\rho_\fob}_{v'v[x]}} s_e^{-1}(f_{\fob,x})
\big)\in (t^{k+1}),
\]
where $e:v[x]\to\rho_\fob$ and we consider $z^{m^{\rho_\fob}_{v'v[x]}}
s_e^{-1}(f_{\fob,x})$ as element of $\kk[P_{\rho_\fob,\sigma}]$. A
structure is \emph{normalized to order $k$} if each of its slabs is
normalized to order $k$. 
\end{definition}

The point of normalization is the following.

\begin{proposition}\label{normalization lemma}
Assume that $\foj$ is a joint of a structure $\scrS$ such that each
slab $\fob\supseteq\foj$ is normalized to order $k$, and let
$\theta^k_\foj$ be the log automorphism associated to a loop around
$\foj$ based on $\sigma\in\P_\max$. Then for any $m\in
P_{\sigma_\foj,\sigma}$ and $v\in\sigma_\foj$ it holds
\[
\tlog_v\big( \theta^k_\foj(m)\big)=0 \mod t^{k+1}.
\]
\end{proposition}

\proof
By parallel transport through $v$ view the basic gluing morphisms
associated to slabs and walls containing $\foj$ as log automorphisms
of $R^k_{\id_{\sigma_\foj},\sigma}$. These all have the form
\[
m\longmapsto f^{\langle \overline m,n\rangle}
\]
for some $n\in\Lambda_{\foj}^\perp\setminus\{0\}$ and $f\in\kk
[P_{\sigma_\foj,\sigma}]$ with $\partial_n f=0$. Moreover, $\tlog_v
f=0$ for automorphisms associated to walls in any case and for slabs
by the normalization condition. Note that the factors
$D(s_e,\rho,v)^{-1}\in\kk\setminus\{0\}$ occurring for slabs have no
influence on $\tlog_v f$. Hence taking into account the composition
formula for log morphisms~(\ref{composition of log morphisms}) and
Remark~\ref{rem:logadd}, the result follows readily from
Lemma~\ref{tlog of compositions} below by induction on the number of
such automorphisms.
\qed 

\begin{lemma}\label{tlog of compositions}
Let $\theta$ be a log automorphism of $R^k_{g,\sigma}$,
$g:\omega\to\tau$, of the form 
\[
m\longmapsto f^{\langle \overline m,n\rangle}
\]
with $n\in\Lambda_\sigma^*\setminus 0$ and $f\in
\kk[P_{\omega,\sigma}]$, $\partial_n f=0$. Assume that $v\in\tau$ is a
vertex such that $a\in R^k_{g,\sigma}$ and $f$, viewed as elements of
$\kk\lfor P_{\omega,\sigma}\rfor_v$, have non-vanishing constant terms
and fulfill the cone condition. Then 
\[
\tlog_v\big( \overline\theta(a)\big)=\tlog_v a \in\kk[t]/(t^{k+1}).
\]
\end{lemma}

\proof
Observe that
\[
\overline \theta(z^m)= f^{\langle \overline m,n\rangle}\cdot z^m
\]
has vanishing $t$-content unless $-m$ occurs as exponent in
$f^{\langle \overline m,n\rangle}$. But then $\partial_n f=0$ implies
$\langle \overline m,n\rangle=0$, and hence $\overline \theta(z^m)
=z^m$. This shows that for any $b\in\kk\lfor P_{\omega,\sigma}\rfor_v$
\[
\tcont(b^i)=\tcont \big( \overline\theta(b^i)\big)
=\tcont\big(\overline\theta(b)^i\big),
\]
and hence, if $b$ fulfills the cone condition and has vanishing
constant term,
\[
\tcont\Big({\sum}_{i\ge1} \frac{1}{i} b^i\Big)
=  \tcont \Big({\sum}_{i\ge1} \frac{1}{i} \overline\theta(b)^i\Big).
\]
The statement follows from this by setting $b=1-(a/a_0)$ for
$a_0\in\kk$ the constant term of $a$.
\qed
\medskip

With the notion of $\tlog$ at hand it is now easy to
complete the proof of Proposition~\ref{codimension 0 undirectional}
left unfinished in Step~II.
\medskip

\noindent\emph{Proof of Proposition~\ref{codimension 0 undirectional}
--- finish.} We have already seen that $\theta_\foj^k\in
O^k(0\otimes\Lambda_\foj^\perp)$, that is,
\[
\theta_\foj^k=\exp\Big(\sum_i c_i t^k\partial_{n_i}\Big)
\]
with $c_i\in\kk$ and $n_i\in\Lambda_\foj^\perp$. If $m\in
P_{\sigma_\foj,\sigma}$ for some $\sigma\in\P_\max$ containing $\foj$
and $v\in\sigma_\foj$ is a vertex, then
\[
\tlog_v\big(\theta_\foj^k(m) \big) = \big\langle
\overline m, \textstyle \sum_i c_i n_i\big\rangle t^k
\quad\mod t^{k+1}.
\]
Thus $\theta_\foj^k\in O^k(0\otimes\Lambda_\rho^\perp)$ if this
expression vanishes for all $m$ with $\overline m\in\Lambda_\rho$.

Now $\theta_\foj^k$ is the composition of log isomorphisms associated
to walls and to two slabs contained in $\rho$. Arguing as in
Proposition~\ref{normalization lemma} we see that the former do not
make any contribution to $\tlog_v\big(\theta_\foj^k(m) \big)$ for any
$m$, while the log isomorphisms associated to the two slabs
are trivial on those $m\in P_{\sigma_\foj,\sigma}$ with $\overline
m\in\Lambda_\rho$. Thus indeed $\tlog_v\big(\theta_\foj^k(m) \big)=0$
for $m$ with $\overline m\in \Lambda_\rho$. 
\qed
\medskip

We will now impose the additional inductive assumption that
$\scrS_{k-1}$ is normalized. This is an empty statement for $k=0$. For
the inductive step from $k-1$ to $k$ note that since no
terms of order $k-1$ have been added to slabs to obtain $\scrS_k^\II$,
all slabs in this latter structure are also normalized to order $k-1$.
It is then easy to normalize $\scrS_k^\II$ to order $k$:

\begin{shaded}
\noindent III. \emph{Normalization of slabs.}
For any slab $\fob\in\scrS_k^\II$ and $x\in\fob\setminus\Delta$,
the inductive assumption shows for any vertex $v'\in\rho_\fob$
\[
\tlog_{v'}\big(z^{m^{\rho_\fob}_{v'v}} s_e^{-1}(f_{\fob,x})
\big)= c_{v'} t^k\quad\mod t^{k+1},
\]
for some $c_{v'}\in\kk$. Here $v=v[x]\in\rho_\fob$ 
and
$e:v\to\rho_\fob$. Fix a set of vertices $\mathscr{V}_{\fob}$ of
$\rho_\fob$ such that
\[
\mathscr{V}_\fob\lra
\big\{ m^{\rho_\fob}_{v'v}\, \big|\, v'\in\rho_\fob\text{ vertex}\big\},
\quad v'\longmapsto m^{\rho_\fob}_{v'v}
\]
is a bijection. Now replace $f_{\fob,x}$ by
\[
f_{\fob,x}-\sum_{v'\in \mathscr{V}_\fob} c_{v'} t^k s_e(z^{-m^{\rho_\fob}_{v'v}}).
\]
Noting that $f_{\fob,x}$ already contains the monomial
$z^{-m_{v'v}^{\rho_\fob}}= z^{m_{vv'}^{\rho_\fob}}$ for each
$v'\in\rho_\fob$, as follows from Equations~(\ref{shLS change of
chart}) and~(\ref{normalization condition}), we must have
$-m_{v'v}^{\rho_\fob}\in P_{x}$. Thus the new collection
$\{f_{\fob,x}\,|\,x\in\fob\setminus\Delta\}$ satisfies the
definition of a slab. Note that $c_{v'}$ for $v'\in \rho_\fob$ depends
only $m^{\rho_\fob}_{v' v}$. This shows that the new $f_{\fob,x}$ is
independent of the particular choice of representative vertices
$\mathscr{V}_{\fob}$. By construction the slab $\fob$ is now
normalized to order $k$.

After modifying each slab in $\scrS_k^\II$ in this way, we obtain
$\scrS_k^\III$.
\end{shaded}

We can now complete the proof of Theorem~\ref{thm:scrS_k}.

\begin{proposition}
The structure $\scrS_k^\III$ is consistent to order $k$.
\end{proposition}

\proof
Since the normalization procedure does not change walls,
Proposition~\ref{step II reduction},(1) still shows consistency to order
$k$ for codimension zero joints.

If $\codim\sigma_\foj=2$, local rigidity (Definition~\ref{def:locally
rigid}) provides a partition of the set of codimension one cells
$\rho\supseteq \foj$ with $Z_\rho\cap X_{\sigma_\foj}\neq \emptyset$
into subsets of cardinalities $2$ and $3$. Let $\rho_i$,
$i=1,\ldots,s$, be a choice of one representative for each such
subset. Let $f_i\in \kk[P_{\sigma_\foj,\sigma}]$ be the sum of the
terms $a_m z^m$ of $f_{\rho_i,v}$ with $\ord_\foj(m)=0$. Note this is
independent of the choice of representative by
Definition~\ref{def:locally rigid},(ii) and the normalization
condition~\eqref{normalization condition}. By local rigidity the
Newton polytope of $f_i$ is
\[
\Xi_i:=\conv\big\{ m^{\rho_i}_{v v'}\,|\, v'\in\sigma_\foj
\big\}\subseteq \Lambda_{\sigma_\foj,\RR},
\]
and any integral point of this polytope is a vertex. Now
Proposition~\ref{step II reduction},(2), which by
Proposition~\ref{influence of monomial changes},(1) continues to hold
after normalization, gives
\begin{align*}
\theta_\foj^k&=\exp\Big( ct^k\partial_n
+{\sum}_{i, v'\in\sigma_\foj}
c_{i,v'} t^k\frac{z^{m^{\rho_i}_{vv'}}}{f_i}\partial_{n_{i,v'}}
\Big)\\
&= \exp\Big( t^k\frac{\big(\prod_j f_j\big)  \partial_{cn}
+\sum_i \big(\prod_{j\neq i} f_j\big)\sum_{m\in\Xi_i}z^m
\partial_{n_i(m)} }{\prod_j f_j}   \Big),
\end{align*}
where $n_i(m)=\sum_{\{v'\in \sigma_\foj\,|\, m^{\rho_i}_{v v'}=m\}}
c_{i,v'} n_{i,v'} \in\kk\otimes\Lambda_\sigma^*$. Writing $f_i=
\sum_{m\in\Xi_i} d_{i,m} z^m$, $d_{i,m}\in \kk\setminus\{0\}$, this
leads to
\begin{align*}
\theta_\foj^k=\exp\Big(\frac{t^k}{\prod_i f_i}
\sum_{m\in\Delta(\sigma_\foj)} z^m\partial_{n(m)}  \Big)
\end{align*}
with $\Delta(\sigma_\foj)=\sum_i\Xi_i$ and
\begin{align}\label{n(m)}
n(m)=\sum_{\sum_j m_j=m,\, m_j \in\Xi_j}
d_{1,m_1}\cdot\ldots\cdot d_{s,m_s} \Big( cn +\sum_i
\frac{n_i(m_i)}{d_{i,m_i}} \Big)
\end{align}
Now Proposition~\ref{normalization lemma} applies and we obtain 
\[
0=\tlog_v\big(\theta_\foj^k(m') \big) = t^k \langle m', n(0)\rangle
\]
Thus $n(0)=0$ by the normalization condition. Expanding at a different
vertex $v'\in\sigma_\foj$ changes $f_i$ to $z^{m^{\rho_i}_{v'v}} f_i$
by~(\ref{f_rho}) and hence $n(0)$ becomes $n(m_{v'})$ with
$m_{v'}=\sum_i m^{\rho_i}_{vv'}$. Because the normal fan of
$\sigma_{\foj}$ is a refinement of the normal fan of
$\Delta(\sigma_{\foj})$, the $m_{v'}$ run over all vertices of
$\Delta(\sigma_\foj)$. Now any vertex $m$ of $\Delta(\sigma_\foj)$ may
uniquely be written $m=\sum_i m_i$ with $m_i\in\Xi_i$. Therefore the
sum in~(\ref{n(m)}) has only one term. This shows that for any vertex
$m=\sum_i m_i\in\Delta(\sigma_{\foj})$
\[
cn+\sum_i \frac{n_i(m_i)}{d_{i,m_i}}=0.
\] 
Now apply Definition~\ref{def:locally rigid},(iii) to the tuple of
functions associating $\frac{cn}{s}+ \frac{n_i(m)}{d_{i,m}}$ to a
vertex $m\in\Xi_i$. It then follows that there are $n_i\in
\Lambda_{\sigma_{\foj}}^\perp\otimes \kk$ with
\[
n_i=cn+\frac{n_i(m)}{d_{i,m}}\quad\text{for all }m\in\Xi_i,
\]
and $\sum_i n_i=0$. Thus $n(m)=0$ for all $m\in\Delta(\sigma_\foj)
\cap\Lambda_{\sigma_\foj}$. Hence $\theta_\foj^k=1$, as desired.

If $\codim\sigma_\foj=1$ we have to show that the normalization
procedure does not spoil the consistency from Proposition~\ref{step II
reduction},(1). Indeed, Proposition~\ref{influence of monomial
changes},(1) shows that we can write
\[
\theta_\foj^k=\exp\Big(\sum_{v'\in\rho} c_{v'} t^k
\frac{z^{m^\rho_{vv'}}}{f_{\rho,v}}\partial_{\check d_\rho} \Big).
\]
Now as in the codimension two case,
\[
0=\tlog_{v'}\big(\theta_\foj^k(m)\big)= t^k\langle
m,c_{v'} \check d_\rho\rangle,
\]
and hence $c_{v'}=0$ for all $v'$.
\qed

%===========================================================
%===========================================================
\section{Higher codimension scattering diagrams}\label{section higher
codimension}
This section fills in the remaining parts of the proof of
Proposition~\ref{scattering proposition}. Nothing here requires
boundedness of cells in $\P$. We first establish a stronger uniqueness
theorem, and then show existence for the two cases $\codim
\sigma_\foj=1,2$ separately. Throughout this section we fix the
following notation. Let $\foj\subseteq B$ be an $(n-2)$-dimensional
polyhedral subset contained in a reference cell $\sigma\in\P_\max$ and
let $g:\omega\to\sigma_\foj$ with $\foj\cap\Int\omega \neq\emptyset$.
Furthermore, choose $x\in (\foj \cap\Int\omega) \setminus\Delta$ and
let $v=v[x]\in\sigma_\foj$ be a vertex in the same connected component
of $\sigma_\foj\setminus\Delta$ as $x$. We work in various rings
$R^I_{g,\sigma}$ for $\sigma\in \P_\max$ containing $\foj$, for ideals
$I$ with radical $I_0=I_{g,\sigma}^{>0}$. We keep the standard
notation $I_l= I_{g,\sigma}^{>l}$ from before. Different choices of
$\sigma$ are identified by parallel transport through $v$ without
further notice. We also fix an orientation of $\shQ^v_{\foj,\RR}=
\Lambda_{v,\RR}/\Lambda_{\foj,\RR}$ and write $\shQ:=
\shQ^v_{\foj,\RR}$ for brevity. Recall that a \emph{cut}
$\foc\subseteq\shQ$ is a one-dimensional cone contained in $\doublebar
\rho\subseteq\shQ$ for some $\rho\in\P^{[n-1]}$, $\rho\supseteq\foj$.
Moreover, as before $m\mapsto \doublebar m$ denotes the quotient maps
$P_{\omega,\sigma} \to \shQ$ and $P_x\to \shQ$.

%===========================================================
\subsection{Uniqueness}

\begin{proposition}\label{processunique}
Let $\foj\in \Joints(\scrS_{k-1})$ and $J,J'\subseteq
\kk[P_{\id_{\sigma_\foj},\sigma}]$ ideals with $I_k\subseteq J
\subseteq J'$ and $I_0\cdot J'\subseteq J$, where $\sigma\in\P_\max$,
$\sigma\supseteq\foj$. Let $\foD,\foD'$ be scattering diagrams for
$\foj$ such that $\theta^k_{\foD, \id_{\sigma_\foj}} = \theta^k_{\foD',
\id_{\sigma_\foj}}=1 \mod J'$, and modulo $J$
\[
\theta^k_{\foD,\id_{\sigma_\foj}},\theta^k_{\foD',\id_{\sigma_\foj}}\in
\left\{\begin{array}{ll}
O^k\big((\Lambda_\foj\setminus\{0\})\otimes \Lambda_\foj^\perp\big),&
\codim \sigma_\foj=0\\[1ex]
\displaystyle O^k\Big(\Lambda_\foj\otimes \Lambda_\foj^\perp+
\frac{\Lambda_\rho}{f_\rho}\otimes \Lambda_\rho^\perp \Big),&
\codim \sigma_\foj=1\quad(\rho=\sigma_\foj)\\
\displaystyle O^k\Big(\Lambda_\foj\otimes \Lambda_\foj^\perp+
{\sum}_{\rho\supseteq \foj}\ \frac{\Lambda_\foj}{f_\rho}
\otimes \Lambda_\rho^\perp \Big),&
\codim \sigma_\foj=2,\\[-1ex]
\end{array}\right.
\]
that is, as log automorphisms of $R^k_{\id_{\sigma_\foj},\sigma}/J$.
Assume that $\foD, \foD'$ only differ by outgoing rays
$(\forr,m_\forr,c_\forr)$ not contained in any cut,  with
$z^{m_\forr}\in J'$ and, if $\codim\sigma_\foj=2$, by changing the
functions $f_{\foc,x}$ by multiples of some $z^m\in J'$ with
$-\doublebar m\in \foc\setminus\{0\}$. Then $\foD$ and $\foD'$ are
equivalent modulo $J$.

Moreover, if $\codim\sigma_\foj=1$ the same conclusion holds if we
also allow adding outgoing rays contained in $\doublebar\rho$, provided
$\theta^k_{\foD, \id_{\sigma_\foj}}= \theta^k_{\foD',
\id_{\sigma_\foj}}=1 \mod J$.
\end{proposition}

\proof
We have to investigate how $\theta^k_{\foD,\id_{\sigma_\foj}}$ changes
when $\foD$ is modified. We first derive formulae for the effect of
adding a single ray $(\forr,m_\forr,c_\forr)$ or of adding $c
z^{m_\foc}$ to $f_{\foc,x}$ for some cut $\foc$, where $z^{m_\forr},
z^{m_\foc}\in J'\setminus J$. In the case of adding a ray,
since $I_0\cdot J'\subseteq J$, the associated log automorphism
\[
\theta_\forr=\exp\big(-\log(1+c_\forr z^{m_\forr})\partial_{n_\forr}\big):
m'\mapsto (1+c_\forr z^{m_\forr})^{-\langle \overline
{m'}, n_\forr\rangle}
\]
commutes with any other log automorphism associated to a ray of
$\foD$, modulo $J$. For the commutation with the log
automorphism $\theta_\foc$ associated to a cut $\foc$,
Lemma~\ref{conjugation lemma} shows
\[
\theta_\foc \circ\theta_\forr\circ \theta_\foc^{-1}
= \exp\big( -c_\forr z^{m_\forr} (f_{\foc,x}^{-\langle
\overline{{m_\forr}}, n_\foc\rangle}  \partial_{n_\forr} +
f_{\foc,x}^{-\langle
\overline{{m_\forr}}, n_\foc\rangle-1}(\partial_{n_\forr}
f_{\foc,x}) \partial_{n_\foc})\big).
\]
If $\codim\sigma_\foj=2$ any monomial $z^m$ in $\partial_{n_\forr}
f_{\foc,x}$ fulfills $\ord_\foj(m) >0$, and hence we may write
\[
\theta_\foc \circ\theta_\forr\circ \theta_\foc^{-1}
= \exp \big(a_\forr z^{m_\forr} \partial_{n_\forr}\big)
\]
for some $a_\forr\in \kk(\Lambda_\foj)\cap
(R^k_{\id_{\sigma_\foj},\sigma'})^\times$, $\sigma'\in\P_\max$ the
relevant reference cell. Note that $a_\forr$ depends, up
to a constant factor, only on
$\doublebar{m_\forr}$. If $\codim\sigma_\foj=1$ this
needs not be true because $\partial_{n_\forr} f_{\foc,x}$ may contain
monomials $z^m$ with $\ord_\rho(m)=0$. This term, however, becomes
irrelevant after restriction to $P_\rho:=\big\{m\in
P_{\id_{\sigma_\foj},\sigma'} \,\big|\, \overline
m\in\Lambda_\rho\big\}$ because it occurs in combination with $
\partial_{\check d_{\rho}}$. We can thus nevertheless write
\[
\theta_\foc \circ\theta_\forr\circ \theta_\foc^{-1}|_{P_\rho}
= \exp \big(a_\forr z^{m_\forr} \partial_{n_\forr}\big)|_{P_\rho},
\]
with $a_\forr\in \kk(\Lambda_{\sigma_\foj})\cap
(R^k_{\id_{\sigma_\foj},\sigma'})^\times$.

As for the change of $f_{\foc,x}$ for $\codim\sigma_\foj=2$ by adding
$c_\foc z^{m_\foc} \in J'$ to $f_{\foc,x}$, we note
that modulo $J$ this is equivalent to composing the log
isomorphism $\theta_\foc$ by $\exp \big(-(c_\foc
z^{m_{\foc}}/f_{\foc,x}) \partial_{n_\foc}\big)$. Then analogous
arguments show
\begin{align*}
\theta_{\foc'} \circ\exp \big(- (c_\foc z^{m_{\foc}}/f_{\foc,x} )
\partial_{n_\foc}\big)  \circ \theta_{\foc'}^{-1}
&= \exp \big(a_\foc z^{m_\foc} \partial_{n_\foc}\big),
\end{align*}
for some $a_\foc\in \kk(\Lambda_\foj)\cap
(R^k_{\id_{\sigma_\foj},\sigma'})^\times$ depending only on
$\doublebar{m_\foc}$ and $c_\foc$.

Now assume without loss of generality that $\foD'$ is obtained from
$\foD$ by adding rays $(\forr_i,m_{\forr_i},c_{\forr_i})$ and addition
of $c_j z^{m_j}$ to $f_{\foc(j),x}$ with $c_{\forr_i}
z^{m_{\forr_i}},\, c_j z^{m_j}\in J'$, $-\doublebar{m_j }\in
\foc(j)\setminus\{0\}$ and all $m_{\forr_i}, m_{\foc_j}$ pairwise
distinct. Then if $\codim\sigma_\foj\neq 1$ the above computations
show
\[
\theta^k_{\foD',\id_{\sigma_\foj}}= \theta^k_{\foD,\id_{\sigma_\foj}}\circ
\exp\Big( {\sum}_i a_{\forr_i} z^{m_{\forr_i}} \partial_{n_{\forr_i}}
+{\sum}_j a_{\foc_j} z^{m_j}\partial_{n_{\foc(j)}} \Big)\mod J.
\]
Under the hypotheses on $\theta_{\foD,\id_{\sigma_\foj}},
\theta_{\foD',\id_{\sigma_\foj}}$ this is only possible if both sums
are empty. Indeed, as the $a_{\forr_i}$'s and $a_{\foc_j}$'s
are determined, up to constant factors, by the
$\doublebar{m_{\forr_i}}$'s and $\doublebar{m_{\foc_j}}$'s,
there is no way non-zero terms in these sums can cancel.

If $\codim\sigma_\foj=1$ the line $\doublebar\rho$ separates $\shQ$
into two half-planes. By symmetry it suffices to show that $\foD'$
differs from $\foD$ at most by adding rays in the half-plane not
containing $\doublebar\sigma$. Letting $\forr_1,\ldots,\forr_s$ be
the rays with $\forr_i\subseteq \doublebar\sigma$ we obtain
\[
\theta^k_{\foD',\id_{\sigma_\foj}}|_{P_\rho}
= \theta^k_{\foD,\id_{\sigma_\foj}}|_{P_\rho}\circ
\exp\Big(- {\sum}_{i\le s} c_{\forr_i} z^{m_{\forr_i}} \partial_{n_{\forr_i}}
+{\sum}_{i> s} a_{\forr_i} z^{m_{\forr_i}} \partial_{n_{\forr_i}}
\Big)\Big|_{P_\rho} \mod J.
\]
Now for any monomial $z^m$ occurring in the second sum, $\doublebar m$
is contained in the interior of the half-plane not containing
$\doublebar \sigma$. So these can not cancel with any term from the
first sum. As before we can thus conclude that the first sum must be
empty. This finishes the proof of the first paragraph of the
proposition.

To prove the second paragraph we just need to add that if
$\sigma_\foj=\rho\in \P^{[n-1]}$, then the log automorphism
$\theta_\forr$ for an outgoing ray $\forr$ with
$\forr\subseteq\doublebar\rho$ commutes with $\theta_\foc$ for the two
cuts $\foc$ present in this case. One then sees easily that adding
such rays destroys the condition $\theta^k_{\foD,
\id_{\sigma_\foj}}=1$ unless the change leaves the equivalence class
of $\foD$ unchanged. 
\qed

%===========================================================
\subsection{Infinitesimal scattering diagrams}
One basic idea in the existence proofs of the next two paragraphs is
to ``perturb'' a scattering diagram in order to simplify the type of
scatterings to be considered. For the case of codimension two we also
need to consider more general log automorphisms and more general
functions asssociated to cuts than before. This leads to a decoration
of the elements of the deformed scattering diagram by group elements.
We obtain the following notion of ``infinitesimal scattering
diagram''.

\begin{definition} 
\label{def:infinitesimal scattering diagram}
A \emph{squiggly ray} or \emph{s-ray} $\fol$ in $\shQ$ is the image of
a $C^{\infty}$-embedding $i:[0,+\infty)\to\shQ$, such that for
$t\gg0$, $i(t)=(t-t_0)\doublebar m +i(t_0)$ for some $m\in
\Lambda_\sigma \setminus\Lambda_\foj$. We call $i(0)$ its
\emph{endpoint}, and denote by $\forr(\fol):=\RR_{\ge0}\cdot
\doublebar m$ the associated \emph{asymptotic half-line}. A
\emph{segment} $\fol$ in $\shQ$ is the image of a $C^{\infty}$-embedding
$i:[0,1]\to\shQ$ with distinguished \emph{initial} and \emph{final
endpoint} $i(0)$ and $i(1)$, respectively. An \emph{orientation} of an
s-ray or segment is an orientation of its tangent bundle.

An \emph{infinitesimal scattering diagram} for a group $G$ of log
automorphisms of $R^I_{g,\sigma}$ is a collection
$\foD=\{(\fol,\theta_{\fol}), f_\foc\}$ where (1)~$\fol\subseteq \shQ$
is an s-ray or segment, either oriented or unoriented, (2)~$\theta_{\fol}$
is an element of $G$ of the form $\exp \big(\sum_i c_{\fol,i}
z^{m_{\fol,i}}\partial_{n_{\fol,i}}\big)$, (3)~for each cut
$\foc\subseteq\shQ$ and any $p\in\foc\setminus \{0\}$ not
contained in any $\fol$ with $\dim\fol\cap\foc=0$, we have a
polynomial $f_{\foc,p} =\sum_{m\in P_x,\, \overline m\in
\Lambda_{\rho_\foc}} c_m z^m\in \kk[P_x]$ defining an invertible
element of $R^I_{g,\sigma}$.

We have the following additional conditions imposed on this data:
\begin{enumerate}
\item[(i)]
If $\fol$ is an oriented s-ray or segment then $-\RR_{\ge0}\cdot\doublebar
{m_{\fol,i}} \subseteq \shQ$ is independent of $i$ and, in
the case of an s-ray, is parallel to $\fol$ outside a compact subset,
extending in the direction of the orientation of $\fol$.
\item[(ii)]
If $\fol$ is an unoriented s-ray or segment, then for any $i$,
$\doublebar{m_{\fol,i}}=0$ and, in the case of an s-ray, $n_{\fol,i}\in
\Lambda_\foj^\perp\cap \forr(\fol)^\perp$.
\item[(iii)]
If $G\subseteq\tilde H^I_{\foj}$, then in addition we will assume that
for any s-ray or segment $(\fol,\theta_{\fol})\in\foD$, $\theta_{\fol}=
\exp\big(-\log(1+cz^{m_{\fol}})\partial_{n_\fol}\big)$ for some
$m_{\fol}\in P_x$, $n_\fol\in m_\fol^\perp$, $c\in\kk$.
\item[(iv)]
For s-rays or segments $(\fol,\theta_\fol)$,
$(\fol',\theta_{\fol'})\in\foD$, either $\fol\cap\fol'$ is a finite
set of points, or $\theta_\fol\circ\theta_{\fol'} =
\theta_{\fol'}\circ\theta_\fol$ and $\fol\subseteq\fol'$ or
$\fol'\subseteq\fol$.
\item[(v)]
If $p,p'\in\foc\setminus\{0\}$ lie in the same connected component of
$\foc\setminus\bigcup_{\dim\fol\cap \foc=0}\fol$ then
$f_{\foc,p}=f_{\foc,{p'}}$.
\end{enumerate}

An oriented s-ray is called \emph{incoming} if it is oriented towards
its endpoint; otherwise it is called \emph{outgoing}. An oriented segment
is \emph{outgoing} from its initial endpoint and \emph{incoming} into its
final endpoint.
\qed
\end{definition}

\begin{remark}\label{rem:normal vector}
Recall from \S\ref{par:scattering diagrams} that for any rational
half-line $\fol\subseteq \shQ$ the chosen orientation of $\shQ$
determines uniquely a primitive \emph{normal vector} $n_\fol\in
\Lambda_\sigma^*$. A polynomial $f_{\foc,p}$ of an infinitesimal
scattering diagram thus defines unambigously the log isomorphism
\[
\theta_{\foc,p}: m\longmapsto f_{\foc,p}^{-\langle \overline m,
n_\foc\rangle}
\]
from $R^I_{g,\sigma_-}$ to $R^I_{g,\sigma_+}$, if $\sigma_\pm$ are
the maximal cells with $\doublebar
\sigma_+\cap\doublebar\sigma_-=\foc$, ordered appropriately.
Note this differs slightly from the convention for ordinary
scattering diagrams since here we have already taken into account
the effect of the open gluing data $(s_e)$. Conversely, this log
isomorphism determines $f_{\foc,p}$ uniquely modulo $I$. For
uniformity of notation we will thus describe both s-rays and the
polynomials $f_{\foc,p}$ for cuts by pairs $(\foz,\theta_\foz)$
consisting of a locally closed submanifold ($=\fol$ or a connected
component of $\foc\setminus \bigcup_{\dim\fol\cap \foc=0} \fol$) of
$\shQ$ and a log isomorphism between rings $R^I_{g,\sigma}$. We then
call $(\foz,\theta_\foz)$ \emph{foundational} if it comes from a
cut, and \emph{non-foundational} otherwise.

\sloppy
We will also sometimes confuse an element of a scattering diagram
$(\foz,\theta_{\foz})$ with its\hspace{12ex} \emph{support}, $\foz$.

\fussy
An element $\foz\in\foD$ comes with an orientation of the normal
bundle: For s-rays or cuts it is defined by the normal vector
$n_\foz$ (for s-rays, of the asymptotic half-line). For a segment
$\foz=\im\big( i:[0,1]\to \shQ\big)$ take the orientation in such a
way that if $b\in\foz$ and $\xi\in T_{\shQ,b}$ maps to a positive
normal vector, then $i_*\partial_t,\xi$ forms an oriented basis of
$T_{\shQ,b}$.
\qed
\end{remark}

\begin{construction}
Given a smooth immersion $\gamma:[0,1]\to\shQ$ which
intersects elements of an infinitesimal scattering diagram $\foD$ for
a group $G$ transversally, with endpoints disjoint
from any element of $\foD$, and which does not pass through
any point of
\[
\Sing(\foD):=\bigcup_{\foz\in\foD}\partial\foz\cup\bigcup_{\foz_1,\foz_2
\in\foD\atop\dim\foz_1\cap\foz_2=0}\foz_1\cap\foz_2,
\] 
we now define $\theta_{\gamma,\foD}\in G$, the $\gamma$-ordered
product of those $\theta_{\foz}$ with $\foz$ crossed by $\gamma$.
Explicitly, we can find numbers 
\[
0<t_1\le t_2\le\cdots\le t_s<1
\]
and elements $\foz_i\in\foD$ such that $\gamma(t_i)\in\foz_i$, and
$\foz_i\neq\foz_j$ if $t_i=t_j$, $i\neq j$, with
$s$ taken to be as large as possible. Then we set
\[
\theta_{\gamma,\foD}=\theta_{\foz_s}^{\eps_s}\circ\cdots\circ
\theta_{\foz_1}^{\eps_1}
\]
with the sign $\eps_i=\pm 1$ positive if and only if
$\gamma_*\partial_t|_{t_i}$ maps to a positive normal vector
along~$\foz$.

Note if $t_i=t_{i+1}$ then $\dim \foz_i\cap\foz_{i+1}=1$ and hence
$\theta_{\foz_i}$ and $\theta_{\foz_{i+1}}$ commute according to (iv)
in Definition~\ref{def:infinitesimal scattering diagram}. Thus the
$\gamma$-ordered product is well-defined.
\qed
\end{construction}

\begin{construction}\label{asympscatterdiagram}
An infinitesimal scattering diagram $\foD$ for $H^I_\foj$ has an
associated \emph{asymp\-totic scattering diagram} $\foD_\as$ 
constructed as follows. Take $\omega\in\P$ in
Definition~\ref{def:scattering diagram},(1) as fixed throughout this
section. For each s-ray $(\fol, \theta_\fol)$, according to
Definition~\ref{def:infinitesimal scattering diagram},(iii), we can
write uniquely $\theta_\fol= \exp\big( -\log(1+c_\fol
z^{m_\fol})\partial_{n_\fol}\big)$. Define $\foD_\as$ as the
collection of rays $(\forr(\fol),m_\fol, s_{\omega\to\sigma_\fol}^{-1}
(m_\fol)\cdot c_\fol)$, where $\fol\in\foD$ is an s-ray and
$\sigma_\fol\in \P_\max$ is such that $\forr(\fol)\subseteq
\doublebar{\sigma_\fol}$, together with the functions
\[
f_{\foc,y}:= D(s_{e_y},\rho_\foc, v[y])\cdot s_{e_y}\big(
z^{m^{\rho_\foc}_{v[y]v}} f_{\foc,p} \big)
\]
for cuts $\foc$ and $y\in(\foj\cap \Int \omega) \setminus\Delta$,
$e_y:v[y]\to\omega$. Here $p$ is any point in the unbounded connected
component of $\foc\setminus\Sing(\foD)$, and we use parallel transport
through a maximal cell 
$\sigma$ containing $\rho_{\foc}$ 
to interpret $f_{\foc,p}\in \kk[P_x]$ as an element
of $\kk[P_y]$. 
\end{construction}

\begin{remark}\label{rem:infinitesimal scattering diagrams}
Conversely, an ordinary scattering diagram $\foD=\big\{(\forr,m_\forr,
c_\forr),f_{\foc,x} \big\}$ gives rise to an infinitesimal scattering
diagram for $H_\foj^I$ for any ideal $I$, with one straight s-ray
$\fol$ with endpoint the origin for each ray $\forr$ and $f_{\foc,p}:=
D(s_e,\rho_\foc,v)^{-1}\cdot s_e^{-1}\big(f_{\foc,x}\big)$,
$e:v\to\omega$, for all $p\in\foc\setminus\{0\}$. The associated
asymptotic scattering diagram is equivalent to $\foD$ modulo $I$.
\end{remark}

The point of the definition of infinitesimal scattering diagrams is of
course that if $\theta_{\gamma,\foD}\in H$ for all small loops
$\gamma$ around points of $\Sing(\foD)$, where $H\subseteq G$ is a
normal subgroup, say of log automorphisms of $R^k_{g,\sigma}$,
then also $\theta^k_{\foD_\as,g }\in H$.

%===========================================================
\subsection{Existence in codimension one}
\label{subsection:codimen 1}
This subsection is devoted to the proof of the existence statement in
Proposition~\ref{scattering proposition} in the case
$\codim\sigma_\foj=1$, that is, $\rho:=\sigma_\foj\in \P^{[n-1]}$. In
this case $\shQ$ has two cuts that we denote by $\foc_+,\foc_-$, and
we take $\check d_{\rho}=n_{\foc_+}$. Denote by $K_+$, $K_-$  the two
connected components of $\shQ\setminus\doublebar\rho$ such that
$\check d_{ \rho}$ is positive on $K_+$. We are given a scattering
diagram $\foD'=\{\forr, f_{\foc,y}\}$ for $\foj$
with $\theta_{\foD',g}^{k-1}=1$.

Interpreting $\foD'$ as an infinitesimal scattering diagram for
$H_\foj^k$ according to Remark~\ref{rem:infinitesimal scattering
diagrams}, let $\foD_\inc$ be the set of all elements of $\foD'$ which
are not outgoing s-rays. We will construct an infinitesimal scattering
diagram  $\foD$ for $H^k_\foj$ containing $\foD_\inc$ such that
$\foD_\as$ is obtained from $\foD_\inc$ by adding only outgoing
s-rays. This $\foD$ will be such that $\theta_{\gamma,\foD} =
\theta^k_{\foD_\as,g}$, for $\gamma$ a large counterclockwise loop
around the origin, and it satisfies the conditions of the
Proposition, except that a priori it is not constructed to contain
$\foD'$. But since in particular $\theta_{\gamma,\foD} =1\mod
I_{k-1}$, it then follows by the uniqueness result of
Proposition~\ref{processunique},  used inductively, that in fact
$\foD'$ and $\foD_\as$ are equivalent to order $k-1$. Thus
$\foD_\as$ will be the desired extension of $\foD'$.

Writing $\foD_\inc=\foD_\no\cup \foD_{\inc,+}\cup\foD_{\inc,-}$ where
$\foD_\no$ consists of the two foundational elements and the
non-oriented s-rays of $\foD'$, and $\foD_{\inc,\pm}$ consists of
those incoming s-rays contained in $\cl(K_\pm)$. As starting point for
the construction of $\foD$ deform $\foD_\inc$ to an infinitesimal
\begin{figure}
\input{deformed.pstex_t}
\caption{The deformation $\foD_0$ of $\foD_\inc$.}\label{fig:foD_0}
\end{figure}
scattering diagram $\foD_0$ with $(\foD_0)_\as= \foD_\inc$, by moving
the s-rays in $\foD_{\inc,\pm}$ so that their endpoints are points
$p_\pm\in\Int(\foc_\pm)$, as illustrated in Figure~\ref{fig:foD_0}.

Set
\[
J_l:=I_0^{l+1}+I_k,
\]
and let $s$ be the smallest integer such that $J_s=I_k$. We then
construct infinitesimal scattering diagrams $\foD_1,\ldots,\foD_s$ for
$H^k_\foj$, always with the same functions associated to cuts
as $\foD'$, that is
\[
f_{\foc_\pm}=f_{\foc_\pm,p}:=
D(s_e,\rho_{\foc_\pm},v)^{-1}\cdot s_e^{-1}\big(f_{\foc_\pm,x}\big),
\]
$e:v\to\omega$, for any $p\in\foc\setminus\Sing(\foD_l)$, enjoying the
following properties.
\begin{enumerate}
\item
\[
\theta_{\delta,\foD_l}=1\mod J_l
\]
for each loop $\delta$ around a singular point of $\foD_l$
\emph{except} the origin, and
\[
\theta_{\delta_0,\foD_l}=1\mod J_l+I_{k-1},
\]
for $\delta_0$ a small loop around the origin.
\item
Each $\foz\in\foD_{l+1}\setminus\foD_l$ is either
oriented and fulfills $\Int(\foz)\cap \Sing(\foD_l)=\emptyset$ and
$\Int(\foz)\cap \doublebar \rho=\emptyset$, or has
support equal to $\foc_-$ or $\foc_+$.
\item
If $\foz\in\foD_{l+1}\setminus\foD_l$  is a segment, then the final
endpoint of $\foz$ is in $\doublebar\rho$ and $\doublebar{m_{\foz}}$
is in the same half-plane $K_{\pm}$ as $\Int(\foz)$.
\item
If $\foz\in\foD_\no\setminus \{\foc_+,\foc_-\}$, and $p\in \Int(\foz)$
is a point of $\Sing(\foD_l)$, then there is exactly one oriented
segment or s-ray $\foz'\in\foD_l$ with $p\in\Int(\foz')$. Furthermore
$\foz'$ intersects $\foz$ transversally at $p$ and there is an open
neighbourhood $U$ of $p$ such that, if $U_+$ and $U_-$ are the
connected components of $U\setminus\foz$ such that $U_+\cap \foz'$ is
oriented away from $p$ and $U_-\cap\foz'$ is oriented towards $p$,
then for any other $\foz''\in\foD_l$ containing $p$, $p$ is an
endpoint of $\foz''$, and $\Int(\foz'')\cap U\subseteq U_+$.
\item
If $p\in\doublebar\rho\cap \Sing(\foD_l)\setminus\{0\}$, then for any
non-foundational $\foz\in\foD_l$ with endpoint $p$, $\langle \overline
m_{\foz},\check d_{\rho}\rangle<0$ or $\langle\overline  
m_{\foz},\check d_{\rho} \rangle>0$, independently of $\foz$.
\end{enumerate}

\begin{figure}
\input{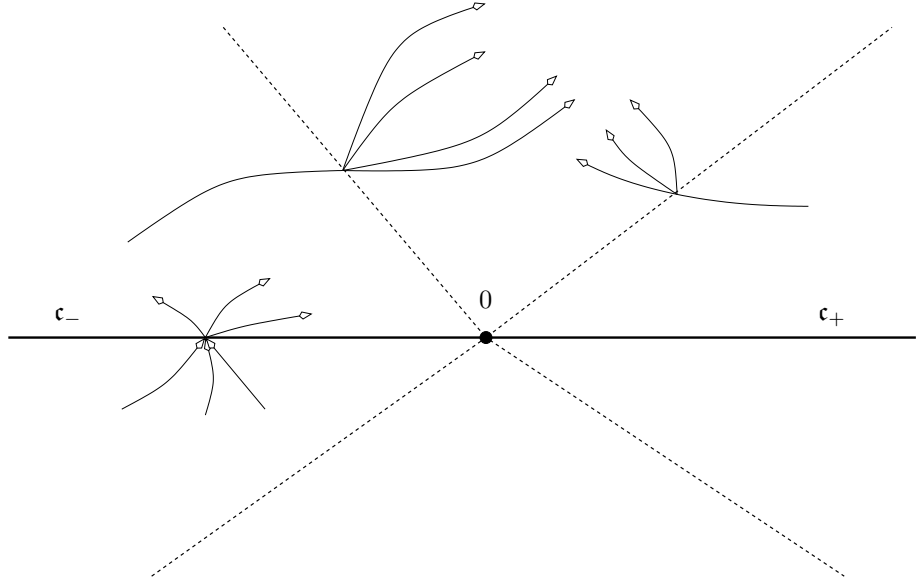}
\caption{Behaviour of $\foD_l$ at non-oriented
s-rays.}\label{fig:nonoriented rays}
\end{figure}

Note here in general the elements of $\foD_l$ are not straight! We
need this so that we do not get ``non-general'' collision points.
Figure~\ref{fig:nonoriented rays} demonstrates allowable behaviour at
non-oriented s-rays and along $\doublebar \rho$, illustrating (4) and
(5).

Note $\foD_0$ satisfies properties (1)--(5).
We shall now construct $\foD_{l+1}$ from $\foD_l$ with these 
properties in several steps, adding new s-rays and segments for each
singular point of $\foD_l$.

\begin{construction}
\emph{Step 1}. If $p\in\Sing(\foD_l)\setminus \doublebar\rho$,  there
are two cases. If $p$ is not contained in a non-oriented s-ray, we
follow essentially the same process as in the proof of the codimension
zero case: If $\delta$ is a small counterclockwise loop around $p$,
then since all $\foz\in\foD_l$ containing $p$ have $\theta_{\foz}\in
\lperp H^{k}_{\foj}$, we can write
\[
\theta_{\delta,\foD_l}=\exp \big(\sum c_iz^{m_i}\partial_{n_i}\big)\mod
J_{l+1}
\]
with $z^{m_i}\in J_l\setminus J_{l+1}$, $\doublebar{m_i}\not=0$ and
$\langle\overline m_i,n_i\rangle=0$. We can then take 
\[
\foD(p)=\big\{(\foz_i,\exp\big(- \log(1+c_iz^{m_i})\partial_{n_i}))\big\},
\]
where: 
\begin{enumerate}
\item
If $-\doublebar{m_i}$ is in the closure of the same connected
component $K_\pm$ of $\shQ\setminus \doublebar\rho$ as $p$, then we
take $\foz_i$ to be an outgoing s-ray (not necessarily straight!)
disjoint from $\doublebar\rho$, with endpoint $p$.
\item
If $-\doublebar{m_i}$ is in the other connected component, then we
take $\foz_i$ to be a segment with initial endpoint $p$ and final
endpoint on $\doublebar\rho$, but not a singular point of $\foD_l$.
\end{enumerate}
In either event, we choose $\foz_i$ so it satisfies the relevant
constraints (2)--(5) listed above. It now holds $\theta_{\delta,
\foD_l\cup\foD(p)} =1\mod J_{l+1}$, just as in the proof of the
codimension zero case of the Proposition.

If $p$ is contained in a non-oriented s-ray $(\foz,\theta_{\foz})$, 
we apply the same process, but now have to argue that
$\theta_{\delta,\foD_l}\in \lperp H^{k}_{\foj}$ rather than just
$H^{k}_{\foj}$. In fact, it follows from constraint (4) that there is
exactly one other s-ray or segment $\foz'\in\foD_l$ with
$p\in\Int(\foz')$, and any other $\foz''\in\foD_l$ containing  $p$ has
endpoint $p$ and is initially contained on the same side of $\foz$ as
the outgoing part of $\foz'$. By abuse of notation denote by
$\theta_\foz$ the composition of all log automorphisms associated to
non-oriented s-rays with support $\foz$. It then suffices to check
$\theta_{\foz}\circ \theta_{\foz'}^{\pm 1} \circ \theta_{\foz}^{-1}\in
\lperp H^{k}_{\foj}$. By Definition~\ref{def:infinitesimal
scattering diagram},(ii) and (iii) we can write
\[
\theta_\foz=\exp \Big( \sum c_j z^{m_j}\partial_{n_\foz}\Big)
\]
with $\doublebar{m_j}=0$, and $\theta_{\foz'}=\exp\big( -\log(1+c
z^{m_{\foz'}}) \partial_{n_{\foz'}}) = \exp\big( \sum_j c'_j z^{j
m_{\foz'}} \partial_{n_{\foz'}}\big)$. It then follows from
Lemma~\ref{conjugation lemma} with $h=\exp(-\sum c_j z^{m_j})$ and
$n_0=n_\foz$ that
\[
\theta_{\foz}\circ \theta_{\foz'}^{\pm 1} \circ
\theta_{\foz}^{-1}
= \exp\Big(\pm{\sum}_j \Ad_{\theta_{\foz}}(c'_j z^{j m_{\foz'}}
\partial_{n_{\foz'}}) \Big)
\]
does not contain any monomials $m'$ with
$\overline{m'}\in\Lambda_\foj$.

We can thus obtain $\foD(p)=\big\{(\foz_i,\exp\big(-\log(1+
c_iz^{m_i})
\partial_{n_i}))\big\}$ as before, with the further property that each
$\foz_i$ lies on the same side of $\foz$ as all the other outgoing
s-rays or segments with endpoint $p$, giving the inductive requirement
(4).

Now set
\[
\foD^{(1)}_l=\foD_l\cup\bigcup_{p\in \Sing(\foD_l)\setminus
\doublebar{\scriptstyle\rho}} \foD(p).
\]
We can make all the choices of s-rays and segments so that $\foD_l^{(1)}$
satisfies the constraints (2)--(5).\\[-3ex]
\smallskip

\noindent
\emph{Step 2}. If $p\in \doublebar\rho\cap\Sing(\foD_l^{(1)})
\setminus\{0\}$ we can construct $\foD(p)$ consisting of outgoing
s-rays with endpoint $p$ such that $\theta_{\gamma,\foD_l^{(1)}
\cup\foD(p)}=1\mod J_{l+1}$. Indeed, by constraint~(5), all incoming
s-rays or line segments at $p$ not contained in $\doublebar\rho$ are
contained in, without loss of generality, $K_+$, and outgoing s-rays
are contained in $K_-$. We claim that  $\theta_{\delta,\foD_l}\in
H^{k}_{\foj,K_-\cup\{0\}}$, the group defined in Definition~\ref{cone groups}.
In fact, it follows again from Lemma~\ref{conjugation lemma} that any
monomial $z^m$ occurring in $\theta_{\delta,\foD_l}$ obeys
$-\doublebar m\in K_-$. To check the claim it remains to
verify that the commutation of an element $\exp\big( c z^m \partial_n
\big)$ of $H_\foj^k$ with the log automorphism associated to crossing
$\foc_-$ preserves $\Omega_\std$, see Remark~\ref{rem:log
automorphisms},(3). If one writes $\Omega_\std= \alpha\cdot
\bigwedge_i \dlog(m_i)$ with $\alpha\in\kk$, $m_i\in\Lambda_\foj$ for
$i\ge 3$, and $\langle m_1, \check d_{\rho} \rangle= \langle
m_2,n\rangle=0$, this follows by a straightforward
computation from Lemma~\ref{conjugation lemma}.

Thus we can now follow the same procedure as in Step~1,
defining $\foD(p)$ to consist of a finite number of outgoing s-rays
with endpoint $p$ and with interior contained in $K_-$. (No segments are
necessary since these s-rays need never cross $\doublebar\rho$.) This
can be done so that $\theta_{\gamma,\foD_l\cup\foD(p)}=1\mod
J_{l+1}$. Doing this for each such $p$, we set
\[
\foD_l^{(2)}=\foD_l^{(1)}\cup\bigcup_{p\in
\Sing(\foD_l^{(1)})\cap \doublebar{\scriptstyle\rho},\, p\neq0} \foD(p).
\]
Again, we can make these choices so that $\foD_l^{(2)}$
satisfies the constraints (2)--(5).
One then sees easily that
\[
\theta_{\delta,\foD_l^{(2)}}=1\mod J_{l+1}
\]
for $\delta$ a loop around any point of $\Sing(\foD_l^{(2)})
\setminus\{0\}$, as follows from the fact that $\ker\big(
H_\foj^{J_{l+1}} \to H_\foj^{J_l}\big)$ is contained in the center of
$H_\foj^{J_{l+1}}$.
\smallskip

\noindent
\emph{Step 3}. 
It remains to address the situation at the origin. Without loss of
generality, we will take $\delta_0$ to have a base point immediately
\begin{figure}[h]
\input{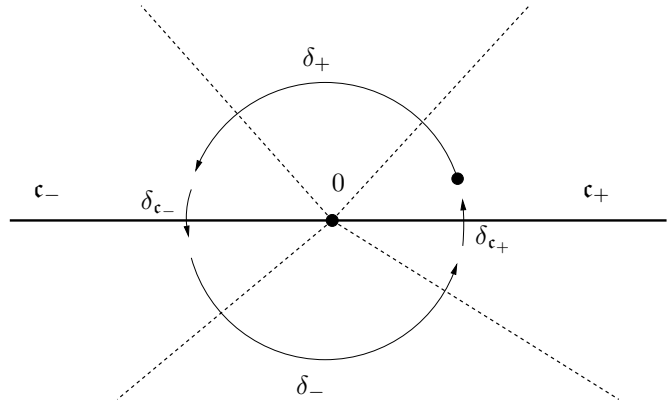}
\caption{Decomposition of $\delta_0$.}\label{fig;semi-circular arcs}
\end{figure}
to one side of $\foc_+$, and split $\delta_0$ up into four
semi-circular arcs in such a way that $\delta_{\foc_-}$ and
$\delta_{\foc_+}$ only cross s-rays with support $\foc_+$ or $\foc_-$,
see Figure~\ref{fig;semi-circular arcs}. Then
\[
\theta_{\delta_0,\foD_l^{(2)}}
=\theta_{\foc_+}\circ\theta_-\circ
\theta_{\foc_-}\circ\theta_+,
\]
where we have written
$\theta_{\foc_{\pm}}:=\theta_{\delta_{\foc_{\pm}},\foD_l^{(2)}}$,
$\theta_{\pm}:=\theta_{\delta_{\pm},\foD_l^{(2)}}$.
Using~(\ref{multiple composition}) we calculate
\begin{eqnarray*}
\theta_{\delta_0,\foD_l^{(2)}}(m)
&=&\theta_{\foc_+}(m)
\cdot\overline\theta_{\foc_+}(\theta_-
(m))
\cdot\overline\theta_{\foc_+}(\overline\theta_-
(\theta_{\foc_-}(m)))
\cdot\overline\theta_{\foc_+}(\overline\theta_-
(\overline\theta_{\foc_-}(\theta_+
(m)))).
\end{eqnarray*}
Now $\theta_{\pm}$ only involves monomials $z^{m'}$ with
$\doublebar{m'}=0$, and these are left invariant by $\overline
\theta_{\foc_{\pm}}$. In addition, $\theta_{\foc_{\pm}}$ take the form
\[
m\mapsto (f_{\pm})^{-\langle\overline m,\check d_{\rho}\rangle}.
\]
(Note that we may not have $f_{\pm}=f_{\foc_{\pm}}$ since it is
possible that some outgoing rays with support $\foc_{\pm}$ may have
been added in this Step for a smaller $l$.) This allows us to simplify
to get 
\begin{eqnarray}\label{codimonestep3comp}
\begin{aligned}
\theta_{\delta_0,\foD_l^{(2)}}(m)&=
(f_{\foc_+})^{-\langle\overline m,\check d_{\rho}\rangle}
\cdot\theta_-(m)
\cdot\overline\theta_-
(f_{\foc_-})^{\langle\overline m,\check d_{\rho}\rangle}
\cdot\overline\theta_-
(\theta_+(m))\\
&=(\theta_-\circ\theta_+)(m)\cdot
\left(\frac{\overline\theta_-(f_{\foc_-})}{f_{\foc_+}}
\right)^{\langle\overline m,\check d_{\rho}\rangle}.
\end{aligned}
\end{eqnarray}
Assuming, say, that $\langle\overline m,\check d_{\rho}\rangle=1$, we
can write this as
\begin{equation}
\label{foundationalmvalue}
\theta_{\delta_0,\foD_l^{(2)}}(m)={f_{\rho}+\cdots\over f_{\rho}+\cdots}
\end{equation}
where $\cdots$ denotes expressions only involving monomials $z^{m'}\in
I_0$ with $\langle\overline{m'},\check d_{\rho}\rangle=0$. Indeed,
$f_{\foc_\pm}= f_\rho\mod I_0$ and $\theta_-\circ\theta_+=1 \mod I_0$.
Since by induction $\theta_{\delta_0,\foD_l^{(2)}} =1 \mod J_l+
I_{k-1}$ the numerator and denominator must agree up to terms in
$J_l+I_{k-1}$. Thus modulo $J_{l+1}+ I_{k-1}$, $\theta_{\delta_0,
\foD_l^{(2)}}(m) = 1+\sum (a_j z^{m_j'}/ f_{\rho})$ with $z^{m'_j}\in
J_l$, $\overline{m'_j}\in\Lambda_\rho$. On the other hand, if $\langle
\overline m,\check d_{\rho}\rangle=0$, then $\theta_{\delta_0,
\foD_l^{(2)}}(m) = (\theta_-\circ \theta_+)(m)$, and we can write
$\theta_{\delta_0, \foD_l^{(2)}}(m)= 1+\sum b_jz^{m_j}$ with
$z^{m_j}\in I_0$ and $\doublebar{m_j}=0$. Thus
\begin{eqnarray}\label{theta_delta0}
\theta_{\delta_0,\foD_l^{(2)}}=\exp \Big(\sum c_jz^{m_j}\partial_{n_j}
+\sum d_j{z^{m_j'}\over f_{\rho}}\partial_{\check d_{\rho}} \Big)
\mod J_{l+1}+ I_{k-1},
\end{eqnarray}
for some coefficients $c_j, d_j\in\kk$ and exponents $m_j,m'_j$ with
$\overline{m_j}\in\Lambda_\foj$, $\overline{m_j'}\in \Lambda_\rho$
(possibly with a different set of $m_j$'s and $m_j'$'s) and
$z^{m_j},z^{m'_j}\in J_l$.

Now taking $\gamma$ to be a large loop around the origin, enclosing
\begin{figure}[t]
\input{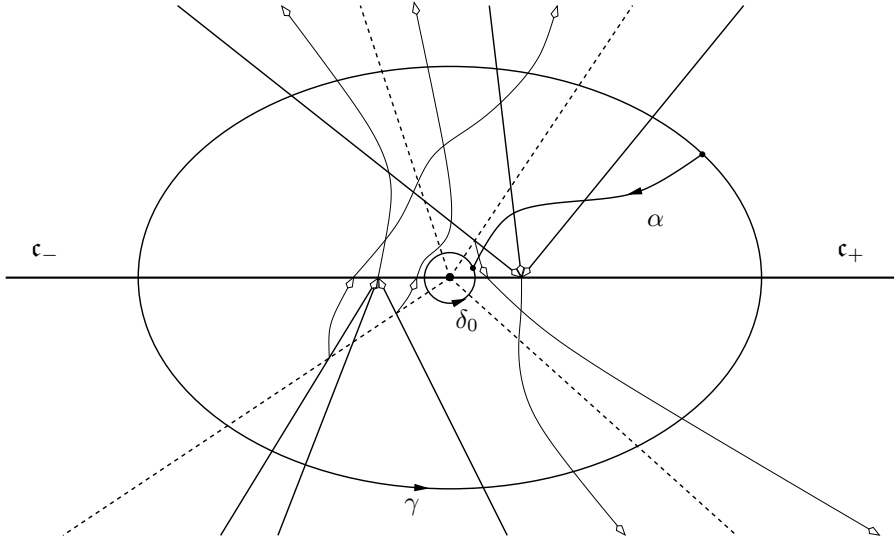}
\caption{Illustration of $\alpha$,$\gamma$ and  $\delta_0$.}
\end{figure}
all singular points and segments of $\foD_l^{(2)}$ and with base
point in the same connected component of
$\shQ\setminus\doublebar\rho$ as $\delta_0$, let $\alpha$ be a path
disjoint from $\doublebar\rho$ joining the base point of $\gamma$ to
the base point of $\delta_0$, so that $\alpha\delta_0\alpha^{-1}$ is
homotopic to $\gamma$ in $\shQ\setminus\{0\}$. Then modulo $J_{l+1}$
\[
\theta_{\gamma,\foD_l^{(2)}}=\theta_{\alpha,\foD_l^{(2)}}^{-1}
\circ\theta_{\delta_0,\foD_l^{(2)}}\circ\theta_{\alpha,\foD_l^{(2)}}.
\]
Since $\theta_{\delta_0,\foD_l^{(2)}}$ only involves monomials $z^m$
in $J_l+I_{k-1}$, one sees that modulo $J_{l+1}+I_{k-1}$,
$\theta_{\alpha,\foD_l^{(2)}}$ commutes with
$\theta_{\delta_0,\foD_l^{(2)}}$. Thus if we let
$\foD_l'=(\foD_l^{(2)})_\as$, then
\[
\theta^k_{\foD_l',g}=\theta_{\gamma,\foD_l^{(2)}}
=\theta_{\delta_0,\foD_l^{(2)}}
\mod J_{l+1}+ I_{k-1}.
\]
This shows that $\theta^k_{\foD'_l,g}$ equals the right-hand side of
(\ref{theta_delta0}) modulo $J_{l+1}+I_{k-1}$. Note that $\foD'_l$
has the same incoming and non-directional rays as $\foD'$ and the
same functions $f_\foc$. Since $\theta^k_{\foD',g}=1 \mod I_{k-1}$
and $\theta^k_{\foD_l',g}=1 \mod J_l+I_{k-1}$, and furthermore
$\foD'$ and $\foD'_l$ agree on cuts, incoming rays, and non-oriented
rays, it follows from the last sentence of
Proposition~\ref{processunique} applied inductively that $\foD'$ and
$\foD'_l$ are equivalent modulo $J_l+I_{k-1}$. To compare $\foD'$
and $\foD'_l$ modulo $J_{l+1}+ I_{k-1}$, let $\foD''_l$ be obtained
from $\foD'_l$ by removing all outgoing rays $\forr$ contained in
$\doublebar\rho$ with $z^{m_{\forr}} \in J_l$. Then
$\theta_{\foD''_l,g}=1\mod J_l+I_{k-1}$ and $\foD'$, $\foD''_l$
satisfy the hypothesis of Proposition~\ref{processunique} modulo
$J_{l+1}+I_{k-1}$. Thus $\foD''_l$ and $\foD'$ are equivalent modulo
$J_{l+1}+I_{k-1}$. Thus by adding a number of outgoing rays
contained in $\doublebar\rho$ to $\foD'_l$ observing
Proposition~\ref{influence of monomial changes},(1) once more, we
can insure that $\theta^k_{\foD'_l,g}=1 \mod J_{l+1}+I_{k-1}$.
Adding these same s-rays to $\foD_l^{(2)}$ to obtain $\foD_{l+1}$
completes the third step of the construction, so that
$\theta_{\delta_0, \foD_{l+1}}=1\mod J_{l+1}+I_{k-1}$. Note this
step does not destroy the result of Step~2.
\qed
\end{construction}

We have now obtained $\foD_1,\ldots,\foD_s$, and from $\foD_s$ we take
$\overline\foD=(\foD_s)_{\as}$, which we modify by throwing out all
outgoing rays in $\overline\foD$ contained in $\doublebar\rho$ to get 
a scattering diagram $\foD$. Now compare $\foD'$, $\overline\foD$ and
$\foD$. All three scattering diagrams share the same incoming and
undirectional rays by construction, and the functions $f_{\foc_{\pm}}$
are the same. In addition, neither $\foD'$ nor $\foD$ have outgoing
rays contained in $\doublebar\rho$. We will use uniqueness
(Proposition~\ref{processunique}) to argue inductively that  $\foD'$,
$\overline\foD$ and $\foD$ are equivalent modulo $J_l+I_{k-1}$ for
every $l$, hence modulo $I_{k-1}$. This is trivially true for $l=0$.
Assume it is true for a given $l$. Now by assumption on $\foD'$ and
by construction of $\overline\foD$
\[
\theta^k_{\foD',g}=\theta^k_{\overline\foD,g}=1
\mod J_{l+1}+I_{k-1}.
\]
Thus by Proposition~\ref{processunique}, $\foD'$ and
$\overline\foD$ are equivalent modulo $J_{l+1}+I_{k-1}$. On the other
hand, inductively $\overline\foD$ and $\foD$ are equivalent modulo
$J_l+I_{k-1}$, so up to equivalence, these two scattering diagrams
only disagree by outgoing rays in $\doublebar\rho$ with attached
monomial in $J_l+I_{k-1}$. It then follows that modulo
$J_{l+1}+I_{k-1}$, $\theta^k_{\foD,g}$ takes the form given in the
codimension one case of Proposition~\ref{processunique}, and hence by the uniqueness
of the first paragraph of that proposition, $\foD'$ and $\foD$ are
equivalent modulo $J_{l+1}+I_{k-1}$.

Now since $\theta_{\delta_0,\foD_s}=1 \mod I_{k-1}$,
(\ref{theta_delta0}) shows
\[
\theta_{\delta_0,\foD_s}
\in O^k \Big(\Lambda_\foj\otimes \Lambda_\foj^\perp +
\frac{\Lambda_\rho}{f_\rho}
\otimes\Lambda_\rho^\perp \Big).
\]
It then follows as before that $\theta^k_{\foD,g}$ has the same
form (but is not equal to $\theta_{\delta_0,\foD_s}$ as we threw out
s-rays contained in $\doublebar\rho$). This finishes the proof of the
existence of scattering diagrams in codimension one.
\qed

%===========================================================
\subsection{Existence in codimension two}
\subsubsection{The denominator problem}
We now want to prove the existence part of Proposition~\ref{scattering
proposition} for the case that $\tau:=\sigma_\foj$ is of codimension
two. Unlike the cases of lower codimension this requires the use of
our hypothesis of local rigidity, specifically, (ii) in
Definition~\ref{def:locally rigid}. We continue to use the notation
set up at the beginning of this section. In addition we write
$P=P_x$, which we also identify with $P_{\omega,\sigma}$ via parallel
transport, for any $\sigma\in \P_\max$ containing $\tau$. Using this
convention we drop the reference cell $\sigma$ in the notation for the
rings, so $R^k_g$ means $R^k_{g,\sigma}$ for any appropriate $\sigma$.

Before we embark on the proof, let us explain the basic difficulty,
why we call this the denominator problem, and why we need some
hypotheses. 

\sloppy
The naive approach to proving this result is to simply emulate the
argument of Lemma~\ref{KSlemma}, proceeding inductively and at each
stage calculating $\theta^k_{\foD,g}$. When doing so, because of the
automorphisms associated to slabs, we obtain terms of the form
$\left(\prod_\mu f_{\foc_\mu,x}^{a_\mu} \right)z^m \partial_n$ for
some $a_\mu\in \ZZ$; in particular, the $f_{\foc_\mu}$'s can appear as
denominators in this expression. We can attempt to get rid of these
terms by adding outgoing rays with support $-\RR_{\ge 0}\overline{m}$
or modifying cuts so as to produce automorphisms of the form 
$\left(\prod_\mu f_{\foc_\mu,x}^{a'_\mu} \right)z^m \partial_n$. Here
$a'_{\mu}$ and $a_{\mu}$ may differ because this automorphism might
need to be commuted past some cuts before it can be
cancelled with the original troubling term. However, we are not
allowed to have terms with denominators appearing in the automorphisms
in $\foD$: if $-\RR_{\ge 0}\overline{m}$ does not coincide with a cut,
we need to have $a_{\mu}'\ge 0$ for all $\mu$, while if $-\RR_{\ge
0}\overline{m}$ does coincide with a cut $\foc_{\mu}$, then we need
$a'_{\mu}\ge -1$ and $a'_{\mu'}\ge 0$ for all $\mu'\not=\mu$. This
means we need to carefully control the powers $a'_\mu$ which appear.
This is what we call the \emph{denominator problem}. If one attempts
this direct approach, then in sufficiently complex examples, the
absence of denominators seems like a miracle. This direct approach can
be carried out using a computer algebra package, and when item (ii) of
the definition of local rigidity is satisfied, this naive algorithm
works. Unfortunately, we have been unable to prove this directly, and
instead have to resort to a more indirect solution to this problem.
However, the fact that Proposition~\ref{scattering proposition} holds,
as proved in this section, implies this naive algorithm does work.

It does not always work if item (ii) of local rigidity fails:

\fussy
\begin{example}
\label{denomproblemcounterexample}
Take $\dim B=3$, $\tau\in\P^{[1]}$, $g=\id_\tau$, and write
$\shQ=\ZZ^2$ with cuts $\foc_i\subseteq \shQ_{\RR}$ 
generated by $(-1,0)$, $(0,-1)$ and $(1,2)$.
Suppose the polarization $\varphi$ is given by
$\varphi(-1,0)=\varphi(0,-1)=0$ and $\varphi(1,2)=2$. Writing
$\Lambda_x=\shQ\oplus\ZZ$ and $\shAff(\check B, 
\ZZ)_x =\shQ\oplus\ZZ\oplus\ZZ$ we obtain
\[
P=P_x=\{(a_1,a_2,a_3,a_4)\in\ZZ^4\,|\,\varphi(a_1,a_2)\le a_4\}.
\]
We then consider the (infinitesimal) scattering diagram $\foD'$
without any rays, thus given by the functions $f_{\foc_\mu} \in
\kk[P_x]$, which we take as follows.
\begin{eqnarray*}
f_{\foc_1}&=& 1+z^{(0,0,1,0)},\\
f_{\foc_2}&=& (1+z^{(0,0,1,0)})^2+ z^{(0,-1,0,0)},\\
f_{\foc_3}&=& 1+z^{(0,0,1,0)}.
\end{eqnarray*}
Noting that $\ord_\tau (0,-1,0,0)=1$, we see $z^{(0,-1,0,0)}\in I_0$
and $\theta^0_{\foD',g}=1$. On the other hand, an elementary
calculation shows that the unique choice of lifting $\foD$ of $\foD'$
so that $\theta^1_{\foD,g}=1$, is obtained by adding the
outgoing s-ray $(\RR_{\ge 0}\cdot (0,1),\theta)$ with
\[
\theta(m)=\Big (1+{z^{(0,-1,0,0)}\over 1+z^{(0,0,1,0)}}
\Big)^{\langle \overline m,(1,0,0)\rangle}
\]
Of course, this is not permitted, as we can not allow denominators in
our automorphisms attached to s-rays. 

Note that the restriction of $f_{\foc_2}$ to $X_\tau$ is (necessarily)
not reduced, and hence this example does not fulfill
Definition~\ref{def:locally rigid},(ii).
\qed
\end{example}

Of course, the foundational elements in this example are completely
determined by the initial data defining the log smooth structure on
$X$, and the point is that there may be local obstructions to
smoothability. In this particular example, one can check there is no
suitable smoothing. Thus additional hypotheses are required. The key
point of Definition~\ref{def:locally rigid},(ii) is that if this
hypothesis holds, there exists a sufficiently rigid local model for
a smoothing, which we now describe.

%=======================

\subsubsection{Model deformations}
Let $\rho_1,\ldots,\rho_r=\rho_0$ be a cyclic numbering of the
codimension one cells containing $\tau$, inducing a counterclockwise
ordering of the corresponding cuts $\foc_\mu= \doublebar{\rho_\mu}$ in
$\shQ$. Let $\check d_{\rho_\mu}=n_{\foc_\mu}$ be as defined in
Remark~\ref{rem:normal vector}. Recall from
Definition~\ref{def:locally rigid} that equality of the subschemes
$Z_{\rho_\mu}\cap X_\tau$ distinguishes subsets of $\{\rho_\mu\}$ of
cardinalities~$2$ and $3$. Moreover, for each such subset we have a
convex PL-function $\varphi$ on $\Sigma_\tau$ with Newton polygon
$\check\Xi$ an integral line segment or triangle with edges of integral
length one, see Remark~\ref{Xi_i, varphi_i}. We number these functions
arbitrarily $\varphi_1,\ldots,\varphi_s$, and write $i(\rho)$ for the
index given by $\rho\supseteq\tau$, provided $Z_\rho\cap
X_\tau\neq\emptyset$. We call such $\rho$ \emph{singular}, while if
$Z_\rho\cap X_\tau=\emptyset$ we put $i(\rho)=0$ and say $\rho$
is \emph{non-singular}. Furthermore, let $\varphi_0$ be the pull-back
to $\Lambda_x$ of a convex PL-function on $\Sigma_\tau$ defining the
polarization. Similarly, we view $\varphi_i$ as PL-functions on
$\Lambda_x$ via composition with $\Lambda_{x}\to \shQ =
|\Sigma_\tau|$.

For $\sigma\in\P_\max$ containing $\tau$ denote by
$n_{i,\sigma}\in\Lambda_x^*$, $1\le i\le s$, the linear function
defining $\varphi_i$ on $K_\tau\sigma$. We also need the partial
linear extension of $\varphi_i$ along $\rho\in\P^{[n-1]}$,
$\rho\supseteq \tau$, defined as follows:
\[
\varphi_{i,\rho}(\overline m):= \max\big\{ \langle
\overline m, n_{i,\sigma}\rangle\,\big|\,
\sigma\in\P_\max,\, \sigma\supseteq\rho \big\}.
\]

Let $e_1,\ldots,e_s$ be the standard generating set
for $\ZZ^s$. Now define $\tilde P\subseteq\shAff(\check B ,\ZZ)_{x}
\oplus\ZZ^s$ to be the monoid
\[
\tilde P=\Big\{m+{\sum}_{i=1}^s a_i e_i\,\Big|\, m\in P,\,
a_i\ge \varphi_i(\overline m),\, i=1,\ldots,s\Big\}.
\]
Then $k[\tilde P]$ is a $\kk[t]$-algebra by setting
$t=z^{\mathbbm{1}}$ for $\mathbbm{1}\in P$ the distinguished element
as before.

To define the ideal of our local model in $\kk[\tilde P]$ we
furthermore use the following notion of $t$-divisibility of a monomial
$z^m\in\kk[P]$ along the $\tau$-stratum: For $m\in P$ define
\[
\height(m)=\min\big\{ \ord_{\sigma}(m)\,\big|\, \sigma\in\P_\max,\,
\sigma\supseteq\tau\big\}.
\]
Note that $\height(m)$ is the integral height of $m$ above the graph of
$\varphi_0$ for an identification $\shAff(\check B,\ZZ)_x=\Lambda_x
\oplus\ZZ$. Note also that by Proposition~\ref{exponentprop},
$\height(m)=\ord_\sigma(m)$ iff $\doublebar m\in\doublebar \sigma$.

For the moment, we will assume we are given functions
$\tilde f_i\in \kk[\tilde P]$ of the form
\begin{eqnarray}\label{tilde f_i}
\tilde f_i=\sum_{m}
c_{i,m} z^{m+\sum_j \varphi_j(\overline m)e_j}
\end{eqnarray}
such that if $\height(m)=0$ and $\overline m\not\in K_\tau\rho$ for
some $\rho$ with $i(\rho)=i$ then $c_{i,m}=0$, and 
\[
f_{\rho}=\sum_{\overline m\in K_\tau \rho,\, \height(m)=0} c_{i,m} z^m.
\]
In particular, for any $i$ we require the functions $f_\rho$ with
$i(\rho)=i$ to have a common extension to $V(\omega)$. This is not
always possible. Our construction of $\tilde f_i$  in \S\ref{par:tilde
f_i} indeed requires the hypothesis of Definition~\ref{def:locally
rigid},(ii). For example, it is impossible to achieve this in
Example~\ref{denomproblemcounterexample}.

Now write $t_i=z^{e_i}\in \kk[\tilde P]$, $1\le i\le s$, and consider
the ideal
\[
\tilde J^{>k}=\big( t_1-\tilde f_1,\ldots,t_s-\tilde f_s \big)
+\big( z^{m+\sum_i a_i e_i} \in\kk[\tilde P]\,\big|\,
\height(m)> k \big)
\]
in $\kk[\tilde P]$. For later use it will be convenient to formally define
$t_0:=1$.

We will also use the following related ideals, with
$A=\sigma,\rho,\tau$ and $l\ge0$:
\begin{eqnarray}\label{strata tilde J}
\begin{aligned}
\tilde J^{>k}_l&= \tilde J^{>k}+
\big( z^{m+\sum_i a_i e_i}\in \kk[\tilde P]\,\big|\, 
z^m\in I_0^{l+1}+I_k \big)\subseteq \kk[\tilde P],\\
\tilde J^{>k}_A&= \tilde J^{>k}+
\big( z^{m+\sum_i a_i e_i}\in \kk[\tilde P]\,\big|\, 
\ord_A(m)>k \big)\subseteq \kk[\tilde P].
\end{aligned}
\end{eqnarray}
The ideal $\tilde J^{>k}$ defines the local model, while $\tilde
J^{>k}_l$, $\tilde J^{>k}_A$ provide the reduction modulo
$I_0^{l+1}+I_k$ and the various primary components, respectively. To
study the situation along a codimension one cell $\rho\supseteq\tau$
it is natural to forget the non-standard behaviour on all other
codimension one cells and work in $\kk[\tilde P_\rho]$ with
\[
\tilde P_\rho=\Big\{m+\sum a_ie_i\,\Big|\, m\in P,\,
a_i\ge \varphi_{i,\rho}(\overline m),\, i=1,\ldots,s \Big\}.
\]
The analogue of $\tilde J^{>k}$ in $\kk[\tilde P_\rho]$ is
\[
{}_\rho\tilde J^{>k}=
\big( t_1-\tilde f_1,\ldots,t_s-\tilde f_s \big)
+\big( z^{m+\sum_i a_i e_i} \in\kk[\tilde P_\rho]\,\big|\,
\height_\rho(m)> k \big),
\]
where
\[
\height_\rho(m)=\min\big\{ \ord_{\sigma}(m)\,\big|\, \sigma\in\P_\max,\,
\sigma\supseteq\rho\big\}.
\]
The formulae for the analogues ${}_\rho\tilde J^{>k}_l$,
${}_\rho\tilde J^{>k}_A$ of $\tilde J^{>k}_l$, $\tilde J^{>k}_A$ are
defined as in (\ref{strata tilde J}) with ${}_\rho\tilde J^{>k}$
replacing $\tilde J^{>k}$.

\begin{remark}
The basic idea of what we are going to do is that 
\[
\Spec \kk[\tilde P]/\tilde J^{>k}\lra \Spec \kk[t]/(t^{k+1})
\]
gives a good local model for the $k$-th order deformation of
$V(\tau)$: Clearly, the reduction of the central fibre is canonically
isomorphic to the closure of $V(\tau)$ in $V(\omega)= \Spec
\kk[P]/(t)$, and the central fibre is reduced on $V(\tau)$. To
describe the map at the generic point of a codimension one stratum
$V_{\tau\to\rho}\subseteq V(\tau)$, find $m,m'$ in a localization  of
$P$ along the face corresponding to $\rho$ such that $\height_\rho(m)=
\height_\rho(m')=0$ and $\overline m=-\overline{m'}$ is a generator of
$\Lambda_x/\Lambda_\rho$. Then $x=z^m,\, y=z^{m'}$ fulfill the
relation $xy=t^e$, while their canonical lifts $\tilde x=z^{m+\sum_i
\varphi_{i,\rho}(\overline m)e_i}$, $\tilde y=z^{m'+\sum_i
\varphi_{i,\rho}(\overline{m'})e_i}$  fulfill
\[
\tilde x\tilde y= z^{e_{i(\rho)}}\cdot t^e
=t_{i(\rho)}\cdot t^e= \tilde f_i\cdot t^e\mod \tilde J^{>k}.
\]
This is just a slight perturbation of the standard local model as
derived in the proof of Lemma~\ref{model at general points of Z}. As a
result, we can essentially describe this local model by gluings of the
standard thickenings via certain automorphisms, which can be written
down explicitly. Some effort is then required to massage these
automorphisms into a standard form.
\qed
\end{remark}

\sloppy
\begin{remark}
\label{iteration}
We will use repeatedly the following observation: \emph{Any monomial in
$\kk[\tilde P]/\tilde J^{>k}$ is equal to an expression only involving
monomials of the form $z^{m+\sum_i \varphi_i(\overline m)e_i}$.}
Indeed, we can show this by downward induction on $\height (m)$,
starting with those monomials with $\height=k+1$, in which case such a
monomial is already in $\tilde J^{>k}$. Now suppose the result is true
for all monomials with $\height>k'$. Then for a monomial $z^{m+\sum_i
a_ie_i}$ with $\height(m)=k'$, we can write in $\kk[\tilde P]$
\begin{equation}
\label{iterateeq}
z^{m+\sum_i a_ie_i}=z^{m+\sum_i\varphi_i(\overline m)e_i}
\prod_{i=1}^s t_i^{a_i-\varphi_i(\overline m)}.
\end{equation}
We can then substitute $t_i=\tilde f_i$ for each $i$. Now by
design $\tilde f_i$ only contains terms of the desired form
$z^{m'+\sum_j \varphi_j(\overline {m'})e_j}$; however, in making the
substitution, cross terms will arise which are not of this form. These
cross terms  in $\prod_{i=1}^s \tilde f_i^{a_i-\varphi_i(\overline
m)}$ are of the form
\[
\prod_j z^{m_j+\sum_i\varphi_i(\overline{m_j})e_i}
=z^{\sum_j m_j+\sum_{i,j}\varphi_i(\overline{m_j})e_i}.
\]
By convexity of $\varphi_i$, $\sum_j\varphi_i(\overline{m_j})\ge
\varphi_i(\sum_j \overline{m_j})$, and if we have inequality, then
there is no $\sigma\in\P_\max$ with $\overline {m_j}\in K_\tau\sigma$
for all $j$; strict convexity of $\varphi_0$ on the fan
$\Sigma_{\tau}$ then implies $\sum_j\varphi_0 (\overline{m_j})>
\varphi_0(\sum_j \overline{m_j})$. Thus any term arising in the
expansion of~(\ref{iterateeq}) of the form $z^{m'+\sum_i a'_ie_i}$
with $a'_i>\varphi_i \big(\overline m\big)$ for some $1\le i\le s$,
must have $\height(m')>k'$. By the induction hypothesis, these terms
can be written in the desired form. We call this process
\emph{reduction} and say a monomial in $\kk[\tilde P]/\tilde J^{>k}$ 
is in \emph{reduced form} if it is of the form $z^{m+\sum_{i=1}^s
\varphi_i(\overline m) e_i}$.

\fussy
The same argument works in the ring $\kk[\tilde P_{\rho}]/
{}_\rho\tilde J^{>k}$ if $\varphi_i$ is replaced by
$\varphi_{i,\rho}$, for $0\le i\le s$. There is one slight difference
here: the terms $z^{m'+\sum_i a_ie_i}$ appearing in $\tilde f_j$ may
not satisfy $a_i=\varphi_{i,\rho}(\overline{m'})$, but if they do not,
then we also have $\height_\rho(m')>0$, allowing the induction process.
Again, we call a monomial in this ring in \emph{reduced form} if it is
of the form $z^{m+\sum\varphi_{i, \rho}(\overline m)e_i}$.
\qed
\end{remark}

\subsubsection{Comparison with the standard model}
We now want to decompose $\Spec \kk[\tilde P]/\tilde J^{>k}$ into
standard pieces. To make the comparison with the standard piece
$\Spec R^k_{\omega\to \sigma}$ for some $\sigma\in\P_\max$ containing
$\tau$, we need to localize at the product of
\[
f_{i,\sigma}:=\sum_{\{m\,|\, \height(m)=0,\,  \overline m\in
K_\omega\sigma \}} c_{i,m}z^{m} \in \kk[P].
\]
Let $\sigma_1,\ldots,\sigma_r$ be the maximal cells containing $\tau$,
labelled modulo~$r$ in such a way that $\rho_\mu=\sigma_\mu\cap
\sigma_{\mu+1}$. Since we often need to consider neighouring cells we
write $\mu^-:=\mu-1$ , $\mu^+:=\mu$ for $1\le \mu\le r$, interpreted
modulo $r$.

\begin{proposition}
\label{saturationprop}
For $1\le \mu\le r$
\begin{equation}
\label{saturation}
z^{m}\longmapsto z^{m+\sum_i \langle \overline m ,
n_{i,\sigma_\mu} \rangle e_i}
\end{equation}
induces ring isomorphisms
\begin{eqnarray*}
\beta_{\mu}^\pm:(R^{k}_{\omega\to\sigma_\mu})_{\prod_i \! f_{i,\sigma_\mu}}
&\lra& \big(\kk[\tilde
P_{\rho_{\mu^\pm}}]/{}_{\rho_{\mu^\pm}}\hspace{-1pt} \tilde
J^{>k}_{\sigma_\mu} \big)_{\prod_i t_i}\\
\kappa_\mu: \big(R^{k}_{\omega\to\tau})_{\prod_i \! f_{i,\sigma_\mu}}
&\lra& \big(\kk[\tilde P]/\tilde J_\tau^{>k}
\big)_{\prod_i t_i}
\simeq \big(\kk[\tilde P_{\rho_{\mu^\pm}}]/
{}_{\rho_{\mu^\pm}}\hspace{-1pt}\tilde J_\tau^{>k}
\big)_{\prod_i t_i}.
\end{eqnarray*}
\end{proposition}

\proof
Let us first consider the case of $\beta_\mu^+$. First we note that
\[
\beta^+_\mu(f_{i,\sigma_\mu})=\sum_{\{m\,|\, \height(m)=0, \,
\overline m\in K_\omega\sigma_\mu \}}
c_{i,m}z^{m+\sum_j\langle \overline m,n_{j,\sigma_\mu} \rangle e_j}
=\tilde f_i\mod {}_{\rho_\mu}\hspace{-1pt}\tilde J_{\sigma_\mu}^{>0}.
\]
Because ${}_{\rho_\mu}\hspace{-1pt}\tilde J_{\sigma_\mu}^{>0}$ is
nilpotent in $\kk[\tilde P_{\rho_\mu}]/{}_{\rho_\mu}\hspace{-1pt}
\tilde  J_{\sigma_\mu}^{>k}$,  we see that if we can invert
$t_i=\tilde f_i$,  we can invert $\beta^+_\mu(f_{i,\sigma_\mu})$, and
vice versa. Thus $\beta^+_\mu$ is defined.
Now set
\[
\tilde R^k_{\mu}:=\big(\kk[P]\otimes_{\kk} 
\kk[t_1^{\pm 1},\ldots,t_s^{\pm 1}]\big) 
\big/ \big( z^m\otimes 1\,\big|\,m\in P,\,\ord_{\sigma_{\mu}}
(m)>k \big).
\]
The formula for $\beta_{\mu}^+$ induces an obvious identification of
$\tilde R^k_{\mu}$ with 
\[
\big(\kk[\tilde P_{\rho_{\mu}}]/\langle 
z^{m+\sum_j a_je_j}\,|\,\ord_{\sigma_{\mu}} (m)>k\rangle\big)_{\prod t_i}.
\]
This identification yields an isomorphism
\[
\tilde R^k_{\mu}/(t_1-\tilde f_1,\ldots,t_s-\tilde f_s)
\simeq \big(\kk[\tilde P_{\rho_{\mu}}]/
{}_{\rho_\mu}\hspace{-1pt}\tilde J^{>k}_{\sigma_{\mu}} \big)_{\prod t_i}.
\]
Now using the same reduction process as in Remark~\ref{iteration},
$\tilde f_i\in \tilde R^k_{\mu}$ is equivalent modulo $(t_1-\tilde f_1,
\ldots,t_s-\tilde f_s)$ to a $\tilde f_i'$ only containing monomials of
the form $z^m\otimes 1$. In particular, in $\tilde R^k_{\mu}$,
we can write
\[
t_i-\tilde f_i'=\sum g_{ij}(t_j-\tilde f_j),
\]
and we can assume the image of $g_{ij}$ in $\tilde R^0_{\mu}$ is 
$\delta_{ij}$. Thus the matrix $(g_{ij})$ is invertible in $\tilde R^k_{\mu}$,
and the ideals $(t_1-\tilde f_1,\ldots,t_s-\tilde f_s)$
and $(t_1-\tilde f_1',\ldots,t_s-\tilde f_s')$ coincide. Thus
\begin{eqnarray*}
(R^k_{\omega\rightarrow\sigma_{\mu}})_{\prod f_{i,\sigma_{\mu}}}
&\simeq&\big(\kk[P]/I^{>k}_{\omega\rightarrow
\sigma_{\mu}}\big)_{\prod \tilde f_i'}\\
&\simeq&\tilde R^k_{\mu}/(t_1-\tilde f_1',\ldots,t_s-\tilde f_s')
\ \simeq\ \big(\kk[\tilde P_{\rho_{\mu}}]/
{}_{\rho_\mu}\hspace{-1pt}\tilde J^{>k}_{\sigma_{\mu}} \big)_{\prod t_i},
\end{eqnarray*}
the first isomorphism since the localizing elements only differ by nilpotent
monomials. Furthermore, this isomorphism is induced by $\beta^+_{\mu}$,
giving the result.

The proofs for $\beta_\mu^-$ and $\kappa_\mu$ run 
identically. For the target of $\kappa_\mu$ note that the inclusion
$\kk[\tilde P] \to \kk[\tilde P_{\rho_{\mu^\pm}}]$ is an isomorphism
after localizing at $\prod t_i$.
\qed
\medskip

For the following proposition recall that we defined $t_{i(\rho)}=t_0=1$
whenever $\rho$ is non-singular.

\begin{proposition} 
\label{logautoconstruction}
For each $1\le \mu\le r$ define the log automorphism
\[
\theta_\mu:P\lra (R^{k}_{g})^{\times},\quad
\theta_\mu(m)=
\kappa_{\mu+1}^{-1}(t_{i(\rho_\mu)})^{-\langle \overline m,
\check d_{\rho_\mu}\rangle}.
\]
Then
\[
\theta_r\circ\cdots\circ\theta_1=1.
\]
\end{proposition}

\proof First note that from the definition of the convex PL-functions
$\varphi_j$ and their defining linear functions $n_{j,\sigma}$
\[
t_{i(\rho_\mu)}^{-\langle \overline m,\check d_{\rho_\mu}\rangle}
=\prod_{j=1}^s t_j^{\langle \overline m,n_{j,\sigma_\mu}-n_{j,\sigma_{\mu+1}}\rangle}.
\]
Thus
\begin{eqnarray*}
\kappa_{\mu+1}(\overline\theta_\mu(z^{m}))&=&
\kappa_{\mu+1}(\theta_\mu(m))\kappa_{\mu+1}(z^{m})
\ =\ \bigg( \prod_{i=1}^s 
t_i^{\langle \overline m,n_{i,\sigma_\mu}-n_{i,\sigma_{\mu+1}}\rangle}
\bigg)
z^{m+\sum_i \langle \overline m,n_{i,\sigma_{\mu+1}}\rangle e_i}\\
&=&z^{m+\sum \langle \overline m,n_{i,\sigma_\mu}\rangle e_i}
\ =\ \kappa_\mu(z^m),
\end{eqnarray*}
and hence
\[
\overline\theta_\mu= \kappa_{\mu+1}^{-1}\circ\kappa_{\mu}.
\]
From this one easily sees that
\begin{eqnarray*}
\theta_r\circ\cdots\circ\theta_1(m)
&=& \prod_{\mu=1}^r \big(\kappa_1^{-1}\circ\kappa_{\mu+1}
\circ\theta_\mu\big)(m)\ =\ 
\kappa_1^{-1}\Big( \prod_{\mu=1}^r
t_{i(\rho_\mu)}^{-\langle\overline m, \check d_{ \rho_\mu}\rangle} \Big)\\
&=&\kappa_1^{-1}\Big(\prod_{i=1}^s \prod_{\{\rho\,|\, i(\rho)=i\}}
t_{i}^{-\langle \overline m, \check d_{\rho}\rangle} \Big)
=1,
\end{eqnarray*}
the last equality by Definition~\ref{def:locally rigid},(ii) and
Remark~\ref{Xi_i, varphi_i}.
\qed

%=======================

\subsubsection{The factorization lemma}
The problem now is that when $\rho_\mu$ is singular, $\theta_\mu$ is
not a very well-behaved log automorphism. Our next task is to factor
it into manageable log automorphisms. The first step achieves a
product decomposition of $\tilde f_{i(\rho_\mu)}$ of the form
$t_i(1+ t_i^{-1}\tilde G_\mu^+)(1+t_i^{-1}\tilde G_\mu^-)$, with
$\tilde G_\mu^\pm$ gathering all monomials propagating into
$\sigma_{\mu^\pm+1}$ and $i=i(\rho_\mu)$. There are two problems
with this. First, this is generally only possible modulo monomials
of higher order or propagating along $\rho$ away from $\tau$. And
second, $\tilde G_\mu^\pm$ will only become divisible by $t_i$ after
pulling back by $\beta_\mu^\pm$. From this factorization
Proposition~\ref{logpunchline} constructs the desired factorization
of $\theta_\mu$.

Recall that $c_{i,m}$ denotes a coefficient of $\tilde f_i$
(\ref{tilde f_i}).

\begin{lemma}
\label{thing2prop}
Let $\rho=\rho_\mu$ be singular and $i=i(\rho_\mu)$.\\[1ex]
(1)~There exist $\tilde F_\mu,\tilde G_\mu^\pm\in \kk[\tilde
P_{\rho}]/{}_{\rho}\tilde J^{>k}$ with the
following properties.
\begin{enumerate}
\item[(i)]
\[
\tilde F_\mu=\sum_{m\in P,\, \overline m\in \Lambda_{\rho}}
d_{\mu,m} z^{m+\sum_j \varphi_{j,\rho}(\overline m)e_j}
\]
with $d_{\mu,m}=c_{i,m}$ if $m\in
K_\omega\rho$ and $\ord_{\rho} (m)=0$.
\item[(ii)]\sloppy
$\tilde G_\mu^\pm$ is a linear combination of monomials
of the form  $z^{m+\sum_j \varphi_{j,\rho}(\overline m)e_j}$ with
$\pm\langle \overline m,\check d_{\rho}\rangle<0$.
\item[(iii)]\fussy
$\tilde G_\mu^-\tilde G_\mu^+$ is divisible by $t_i$ in  $\kk[\tilde
P_{\rho}] /{}_{\rho}\tilde J^{>k}$, and in this
ring,
\begin{equation}
\label{tfimage}
t_i+\tilde G_\mu^-+\tilde G_\mu^++t_i^{-1}\tilde
G_\mu^-\tilde G_\mu^+-\tilde F_\mu=0.
\end{equation}
\end{enumerate}
\smallskip

\noindent
(2)~Suppose for some $i'$, $\tilde f_{i'}$ is replaced by $\tilde
f_{i'}+cz^{m+\sum_j\varphi_j(\overline m)e_j}$, where  $z^m\in
I_0^{l+1}+I_k$.
\begin{enumerate}
\item[(i)]
If $i'\not=i$,  then $\tilde F_\mu$ and $\tilde G_\mu^\pm$
are unchanged modulo ${}_\rho\tilde J_{l+1}^{>k}$.
\item[(ii)]
If $i'=i$ and $\doublebar m\in \doublebar\rho$, then $\tilde
F_\mu$ is replaced by $\tilde F_\mu+cz^{m+\sum_j \varphi_{j,
\rho}(\overline m) e_j}$ modulo ${}_\rho\tilde J^{>k}_{l+1}$,  while $\tilde
G_\mu^\pm$ are unchanged modulo ${}_\rho\tilde J^{>k}_{l+1}$.
\item[(iii)]
If we are not in case (i) or (ii), then modulo ${}_\rho\tilde J^{>k}_{l+1}$,
$\tilde F_\mu$ and $\tilde G_\mu^\pm$ are modified by expresssions of
the form  $az^{m+\sum_j \varphi_{j,\rho}(\overline m)e_j}$, where
$a\in \kk[\tilde P_{\rho}]\setminus{}_\rho\tilde J^{>k}_0$.
\end{enumerate}
\end{lemma}

\proof
We first note that if $\tilde G_\mu^\pm$ satisfy condition (ii) in (1),
then $\tilde G_\mu^-\tilde G_\mu^+$ consists entirely of cross-terms of
the form
\[
z^{m+\sum_j\varphi_{j,\rho}(\overline m)
e_j}\cdot
z^{m'+\sum_j\varphi_{j,\rho}(\overline{m'})
e_j}
\]
with $\langle \overline m, \check d_{\rho}\rangle >0$, $\langle
\overline {m'},\check d_{\rho}\rangle<0$. Then
$\varphi_{i,\rho}(\overline m)+\varphi_{i,\rho}(\overline{m'})
>\varphi_{i,\rho}(\overline m+\overline{m'})$, so $\tilde G_\mu^-\tilde G_\mu^+$ is
divisible by $t_i$.

We will construct $\tilde G_\mu^\pm$ and $\tilde F_\mu$ by induction
on $k$. Begin initially with $\tilde F_\mu= \tilde G_\mu^\pm=0$ for
$k=-1$. 

Now assume $\tilde G_\mu^\pm$ and $\tilde F_\mu$  have been
constructed so that~(\ref{tfimage}) holds in  $\kk[\tilde
P_{\rho}]/ {}_{\rho} \tilde J^{>k-1}$. Now as
in Remark~\ref{iteration}, the left-hand side of
(\ref{tfimage}) can be rewritten, modulo
${}_{\rho} \tilde J^{>k}$, entirely in terms of
monomials of the form $cz^{m+\sum_j \varphi_{j,\rho} (m)e_j}$,
necessarily with $\height_\rho(m)=k$, by the induction hypothesis. 
For each such term, we have three cases:
\begin{enumerate}
\item 
If $\langle \overline m,\check d_{\rho}\rangle>0$,
we subtract this term from $\tilde G_\mu^-$. 
\item 
If $\langle \overline m,\check d_{\rho}\rangle<0$, 
we subtract this term
from $\tilde G_\mu^+$. 
\item
If $\langle \overline m,\check d_{\rho}\rangle=0$, 
we add this term to $\tilde F_\mu$.
\end{enumerate}

After making these changes to $\tilde G_\mu^\pm$ and  $\tilde F_\mu$,
it is then clear that~(\ref{tfimage}) holds modulo ${}_{\rho} \tilde
J^{>k}$. Proceeding inductively, we construct $\tilde F_\mu$ and
$\tilde G_\mu^\pm$, such that conditions~(ii) and (iii) hold. 

Note that when this procedure is carried out for $k=0$, we get for the
left-hand side of ~(\ref{tfimage}) just $t_i$, which is equivalent to
$\tilde f_i$. The only terms of the form $z^{m+\sum_j a_je_j}$
appearing in $\tilde f_i$ with $\langle \overline m, \check d_{\rho}
\rangle=0$ that are not in ${}_{\rho} \tilde J^{>0}$
must have $\ord_\rho(m)=0$, hence $\overline m\in
K_\omega \rho$ and $a_j=\varphi_j(\overline m)
=\varphi_{j,\rho}(\overline m)$. This makes it clear that
condition~(i) holds.

For (2), we just need to look at the terms in~(\ref{tfimage}) and
investigate what effect the reduction process of
Remark~\ref{iteration} has on these. Terms in $\tilde G_\mu^\pm$ and
$\tilde F_\mu$ are already in reduced form by conditions~(i) and (ii)
in (1). On the other hand, any cross-term $z^{m'+\sum_j a_j  e_j}$ in
$t_i^{-1}\tilde G_\mu^-\tilde G_\mu^+$ is necessarily in ${}_\rho \tilde
J_{\tau}^{>0}$, so the effect of the change to $\tilde f_{i'}$ to the
reduction process only affects this term by something in ${}_\rho\tilde
J_{\tau}^{>0}\cdot{}_\rho\tilde J_l^{>k}\subseteq{}_\rho\tilde J_{l+1}^{>k}$. 
Finally, reducing $t_i=\tilde f_i$, we see that if $i\not=i'$, any
term appearing in $\tilde f_i$ which is not in ${}_\rho\tilde J^{>0}_{\tau}$
is already reduced, so similarly the change to this term from the
change in  $\tilde f_{i'}$ is in ${}_\rho\tilde J^{>k}_{l+1}$. This gives
Case (i).

Now if $i=i'$, then $t_i=\tilde f_i$ acquires an additional term
$cz^{m+\sum_j\varphi_j(\overline m)e_j}$. If $\doublebar m\in
\doublebar\rho$, then this is already in reduced form as
$\varphi_j(\overline m)=\varphi_{j,\rho}(\overline m)$, giving
Case~(ii). Otherwise, the reduction process will replace this term
modulo ${}_\rho\tilde J^{>k}_{l+1}$ with an expression 
\[
a z^{m+\sum_j\varphi_{j,\rho}(\overline m)e_j}
\]
with $a\not\in{}_\rho\tilde J^{>k}_0$, giving Case (iii).
\qed

\medskip

The following Lemma is the key point for showing that
no denominators will occur in our factorization.

\begin{lemma}
\label{nodenominators}
For $1\le\mu\le r$ and $\rho=\rho_\mu$ singular there are chains of
surjections and inclusions
\[
\begin{array}{rclcl}
\kk[\tilde P_{\rho}]/{}_{\rho}\tilde J^{>k}
&\twoheadrightarrow
&\kk[\tilde P_{\rho}]/{}_{\rho}
\tilde J^{>k}_{\sigma_{\mu}}
&\stackrel{(\beta_{\mu}^+)^{-1}}{\hookrightarrow}
&R^{k}_{\omega\to\sigma_{\mu}}\\
\kk[\tilde P_{\rho}]/{}_{\rho}\tilde J^{>k}
&\twoheadrightarrow
&\kk[\tilde P_{\rho}]/{}_{\rho}
\tilde J^{>k}_{\sigma_{\mu+1}}
&\stackrel{(\beta_{\mu+1}^-)^{-1}}{\hookrightarrow}
&R^{k}_{\omega\to\sigma_{\mu+1}},
\end{array}
\]
where in each instance the first map is the obvious one.
Furthermore, there exist elements $G_\mu^-\in R^{k}_{\omega\to
\sigma_{\mu}}$, $G_\mu^+\in R^{k}_{\omega\to \sigma_{\mu+1}}$ such
that
\begin{eqnarray*}
\beta_\mu^+(G_\mu^-)&=& t_i\cdot\tilde G_\mu^-\\
\beta_{\mu+1}^-(G_\mu^+)&=& t_i\cdot\tilde G_\mu^+,
\end{eqnarray*}
where $i=i(\rho_\mu)$. In addition, the image of $G_\mu^\pm$ in
$R^{k}_{g}$ is in $I_0$.
\end{lemma}

\proof 
We give the proof for $G_\mu^-$ and $\beta_\mu^+$, the other
case being completely analogous.

The first point is to show that $(\beta_\mu^+)^{-1}$ takes
$\kk[\tilde P_{\rho}]/ {}_{\rho} \tilde
J_{\sigma_{\mu}}^{>k}$ into $R^{k}_{\omega\to\sigma_\mu}$ rather
than, a priori, into $(R^{k}_{\omega\to\sigma_\mu})_{\prod_i
f_{i,\sigma_{\mu}}}$. To see this, let $\tilde
P_{\sigma_{\mu}}$ be defined by
\[
\tilde P_{\sigma_{\mu}}=\big\{m+{\textstyle\sum}_i a_ie_i\,\big|\,
m\in P_{\id_{\sigma_\mu}},\,
a_i\ge \langle \overline m,n_{i,\sigma_{\mu}} \rangle,\, i=1,\ldots,s\big\}.
\]
This can be viewed as the localization of $\tilde P$ along the face
corresponding to $\sigma_{\mu}$. Beware this notation is not in strict
analogy with $\tilde P_\rho$ because now we localize at elements of
$P$. If we denote
\[
{}_{\sigma_\mu} \hspace{-0.5pt} \tilde J^{>k}
:=\big( t_1-\tilde f_1,\ldots,t_s-\tilde f_s \big)
+\big( z^{m+\sum_i a_i e_i} \in\kk[\tilde P_{\sigma_\mu}]\,\big|\,
\ord_{\sigma_\mu}(m)> k \big),
\]
then $z^{m}\mapsto z^{m+\sum_j \langle \overline
m,n_{j,\sigma_{\mu}} \rangle e_j}$ also defines a map
\[
\alpha'_{\mu}:R^{k}_{\id_{\sigma_{\mu}}}\lra
\kk[\tilde P_{\sigma_{\mu}}]/ {}_{\sigma_\mu}
\hspace{-0.5pt} \tilde J^{>k}.
\]
This is an isomorphism. Indeed, $\kk[\tilde P_{\sigma_{\mu}}] \simeq
\kk[P_{\id_{\sigma_{\mu}}}]\otimes_\kk \kk[t_1,\ldots,t_s]$, and as in
Remark~\ref{iteration}, modulo ${}_{\sigma_\mu} \hspace{-0.5pt}\tilde
J^{>k}$, every $t_i$ is equivalent to an element of
$\kk[P_{\id_{\sigma_{\mu}}}]$ under this isomorphism. Thus
\[
\kk[\tilde P_{\sigma_{\mu}}]/{}_{\sigma_\mu}
\hspace{-0.5pt} \tilde J^{>k}
\simeq \kk[P_{\id_{\sigma_{\mu}}}]/(t^{k+1})
=R_{\id_{\sigma_{\mu}}}^{k}.
\]

Now $R^{k}_{\id_{\sigma_{\mu}}}$ is the localization of
$R^{k}_{\omega\to\sigma_\mu}$ at any element $z^{m}$ with
$\ord_{\sigma_\mu}(m)=0$ and $\overline m\in \Int
K_\omega\sigma_\mu$. Thus  $U:=\Spec R^{k}_{\id_{\sigma_{\mu}}}$
and $V:=\Spec  (R^{k}_{\omega\to \sigma_\mu})_{\prod
f_{i,\sigma_{\mu}}}$ are both open subsets of $X=\Spec
R_{\omega\to\sigma_\mu}^{k}$, and $X\setminus (U\cup V)$ is a closed
subset of $X$ of codimension at least $2$ and not contained in a toric
stratum. Hence the restriction map $\Gamma(X,\O_X)\to \Gamma(U\cup V,
\O_X)$ is an isomorphism by Lemma~\ref{depth lemma}.

Similarly, $\tilde U=\Spec \kk[\tilde P_{\sigma_{\mu}}]/
{}_{\sigma_\mu} \hspace{-0.5pt} \tilde J^{>k}$ and
$\tilde V=\Spec \big(\kk[\tilde P_{\rho}]/ {}_{\rho}
\tilde J_{\sigma_{\mu}}^{>k}\big)_{\prod t_i}$ are open subschemes of
$\tilde X=\Spec \kk[\tilde P_{\rho}]/ {}_{\rho} \tilde
J_{\sigma_{\mu}}^{>k}$. The maps $\beta_\mu^+$ and
$\alpha'_{\mu}$ induce isomorphisms $\tilde V\to V$ and $\tilde
U\to U$ respectively, which from their explicit form are clearly
compatible on overlaps, defining an isomorphism $\alpha:\tilde U\cup
\tilde V \to U\cup V$. Thus any $\xi\in \kk[\tilde P_{\rho}]/
{}_{\rho} \hspace{-0.5pt}\tilde J_{\sigma_{\mu}}^{>k}$ defines a
regular function on $\tilde U\cup\tilde V$ by restriction from $\tilde
X$,  and $(\alpha^{-1})^*(\xi) \in\Gamma(U\cup
V,\O_X)=\Gamma(X,\O_X)$. As $(\alpha^{-1})^*$ coincides with
$(\beta_\mu^+)^{-1}$ for functions on $\tilde V$, this shows that
$(\beta_\mu^+)^{-1}$ maps $\kk[\tilde P_{\rho}]/ {}_{\rho}
\tilde J_{\sigma_{\mu}}^{>k}$ into
$R^{k}_{\omega\to\sigma_\mu}$.

Now (a) obviously $(\beta_\mu^+)^{-1} (t_i^{-1} \tilde G_\mu^-) \in
(R^{k}_{\omega\to\sigma_\mu})_{\prod f_{i,\sigma_\mu}}$; (b) any
monomial in $\tilde G_\mu^-$ is of the form  $z^{m+\sum_j
\varphi_{j,\rho}(\overline m) e_j}$ with  $\langle \overline m,\check
d_{\rho}\rangle>0$, so $\varphi_{j,\rho}(\overline m)=\langle
\overline m,n_{j,\sigma_{\mu}}\rangle$ for $j\not=i$, but
$\varphi_{i,\rho}(\overline m)>\langle \overline
m,n_{i,\sigma_{\mu}}\rangle$. This shows 
\[
t_i^{-1}z^{m+\sum_j\varphi_{j,\rho}(\overline m)e_j}
=z^{m-e_i+\sum_j\varphi_{j,\rho}(\overline m)e_j}
\in \kk[\tilde P_{\sigma_{\mu}}],
\]
so $(\beta_{\mu}^+)^{-1} \big(t_i^{-1}\tilde G_\mu^-\big)\in
R^{k}_{\id_{\sigma_{\mu}}}$.

Thus from (a) and (b), we see that
\[
G_\mu^-:=(\beta_{\mu}^+)^{-1} \big(t_i^{-1}\tilde G_\mu^-\big)
\in R^{k}_{\omega\to\sigma_\mu},
\]
as desired. For the statement that $G_\mu^-\in I_0$, we note that if
$z^m \not\in I_0$ then $\overline m\in \Lambda_\tau$, but such
monomials do not occur in $G_\mu^-$ by the properties of $\tilde
G_\mu^-$.
\qed

\begin{proposition}
\label{logpunchline}
For singular $\rho=\rho_\mu$ define log automorphisms
\begin{eqnarray*}
\theta_\mu^\pm: P\lra (R^{k}_{g})^{\times},&&
\theta_\mu^\pm(m)=(1+G_\mu^\pm)^{\pm
\langle \overline m,\check d_{\rho}\rangle}\\
\theta'_{\mu}: P\lra (R^{k}_{g})^{\times},&&
\theta'_\mu(m)=F_\mu^{-\langle \overline m,\check d_{\rho}
\rangle}
\end{eqnarray*}
where 
\[
F_\mu=\kappa_{\mu}^{-1}(\tilde F_\mu)
=\kappa_{\mu+1}^{-1}(\tilde F_\mu).
\]
Then $\theta_\mu^\pm\in G_{\foj}^{k}$ and
\begin{eqnarray}\label{theta_mu factorization}
\theta_{\mu}^+\circ\theta_\mu'\circ (\theta_\mu^-)^{-1}=\theta_\mu.
\end{eqnarray}
\end{proposition}

\proof 
The fact that $\theta_\mu^\pm\in G_{\foj}^{k}$ follows
from the fact that $G_\mu^\pm$ are in $I_0$ and are induced by
elements of $\kk[P]$, by Lemma~\ref{nodenominators}. Similarly,
by condition~(i) of Lemma~\ref{thing2prop},(1), $F_\mu\in
(R^k_g)^\times$.

We are going to verify~(\ref{theta_mu factorization}) in the form
\begin{eqnarray}\label{factorization to be checked}
\theta_{\mu}^+\circ\theta_{\mu}'
=\theta_\mu\circ \theta_\mu^-.
\end{eqnarray}
The proof relies on~(\ref{tfimage}) and the fact that $\overline{
\theta_\mu}$ transforms $\kappa_{\mu+1}$ into $\kappa_\mu$:
\begin{eqnarray}\label{kappa_mu theta_mu}
\begin{aligned}
\kappa_{\mu+1}\big(\overline{\theta_\mu}(z^m)\big)
&= \kappa_{\mu+1}\big( \kappa_{\mu+1}^{-1}\big( t_i^{-\langle
\overline m, \check d_{ \rho}\rangle}\big) z^m\big)
\ =\ t_i^{-\langle \overline m,\check d_{\rho}\rangle}
\cdot z^{m+\sum_j \langle \overline m,n_{j,\sigma_{\mu+1}} \rangle}\\
&= z^{m+\sum_j \langle \overline m,n_{j,\sigma_{\mu}} \rangle}
\ =\ \kappa_\mu(z^m).
\end{aligned}
\end{eqnarray}
Here $i=i(\rho_\mu)$. Now if $\overline m\in\Lambda_{\rho}$ then
$\big(\theta_{\mu}^+\circ\theta_{\mu}'\big)(m) =1=
\big(\theta_\mu\circ \theta_\mu^-\big)(m)$. It thus suffices to
evaluate both sides of~(\ref{factorization to be checked}) at one
$m\in P$ with $\langle \overline m, \check d_{\rho}\rangle =1$. Note
that by~(\ref{tfimage})
\begin{eqnarray}\label{tfimage'}
\tilde F_\mu= t_i(1+ t_i^{-1}\tilde
G_\mu^+)(1+t_i^{-1}\tilde G_\mu^-)
\end{eqnarray}
holds in $\big( \kk[\tilde P_\rho]/ {}_\rho\tilde J^{>k}\big)_{\prod_i t_i}$.
Moreover, since $F_\mu$ only involves monomials $z^{m'}$ with
$\langle \overline{m'}, \check d_{\rho}\rangle=0$, multiplication with
$F_\mu^{-1}$ commutes with $\overline{\theta_\mu^+}$. Thus
\[
\overline{\theta_\mu^+}\big( \theta'_\mu(m)\big)
= F_\mu^{-1}.
\]
Using the composition formula for log automorphisms~(\ref{composition
of log morphisms}) we now compute:
\begin{eqnarray*}
\kappa_{\mu+1} \big(\big(\theta_{\mu}^+\circ\theta_{\mu}'\big)(m) \big)
&=& \kappa_{\mu+1} \big(\overline{\theta_\mu^+}
\big( \theta'_\mu(m)\big)
\cdot \theta_\mu^+(m)\big)
\ =\ \kappa_{\mu+1}\big( F_\mu^{-1} \cdot(1+G_\mu^+)\big)\\
&=& \tilde F_\mu^{-1}\cdot \big(1+t_i^{-1}\tilde G_\mu^+\big)
\ \stackrel{(\ref{tfimage'})}{=}\ 
\big(1+t_i^{-1}\tilde G_\mu^-\big)^{-1}\cdot t_i^{-1}\\
&=& \kappa_\mu\big(1+G_\mu^-\big)^{-1}\cdot t_i^{-1}
\ \stackrel{(\ref{kappa_mu theta_mu})}{=}\ 
\kappa_{\mu+1}\big( \overline{\theta_\mu}\big( \theta_\mu^-(m) \big)\big)
\cdot \kappa_{\mu+1}\big(\theta_\mu(m)\big)\\
&=& \kappa_{\mu+1} \big(\big(\theta_{\mu}
\circ\theta_{\mu}^-\big)(m) \big).
\end{eqnarray*}
By Proposition~\ref{saturationprop} this implies~(\ref{factorization
to be checked}) after localization at $\prod_i f_{i,\sigma_{\mu+1}}$,
which is enough because $\Spec R^k_g$ has no embedded components.
\qed
\medskip

Letting $\theta_\mu^\pm=\theta'_\mu=1$ if $\rho$ is non-singular,
Proposition~\ref{logautoconstruction} and
Proposition~\ref{logpunchline} now show
\begin{eqnarray}\label{theta product}
\theta_{r}^+\circ\theta'_r\circ (\theta_r^-)^{-1}\circ\cdots\circ
\theta_{2}^+\circ\theta_2'\circ (\theta_{2}^-)^{-1}\circ
\theta_{1}^+\circ\theta_1'\circ (\theta_{1}^-)^{-1}=1.
\end{eqnarray}

%=======================

\subsubsection{Construction of the scattering diagram}
\label{par:construction of foD}
Recall from the hypothesis of Proposition~\ref{scattering proposition}
that we are given a scattering diagram $\foD'$ for $\foj$ with
$\theta^{k-1}_{\foD',g}=1$. Following Remark~\ref{rem:infinitesimal
scattering diagrams} view $\foD'$ as an infinitesimal scattering
diagram for $G_\foj^k$ and write $\foD'_\inc$ and $\foD'_\no$ for the
sets of incoming and unoriented rays of $\foD'$, respectively. In this
step we will explain how to produce a scattering diagram $\foD$ with
certain properties given $\foD'_\inc$ and $\foD'_\no$ and an
additional choice of elements $\tilde f_1,\ldots,\tilde f_s\in
\kk[\tilde P]$ as was considered above. This choice of functions will
determine  the functions associated to singular $\rho$ containing
$\tau$. To determine the functions for non-singular $\rho$ let
\[
\NS(\tau)=\{\mu\,|\,\hbox{$\rho_{\mu}$ is non-singular}\},
\]
and assume also given a collection of functions $\{F_{\mu}\in
\kk[P]\,|\,\mu\in\NS(\tau)\}$ with $F_\mu\in 1+I_0$. In the final step
below we will then show how to choose $\tilde f_i$ and
$\{F_{\mu}\,|\,\mu\in\NS(\tau)\}$ in such a way that $\foD$ is
equivalent to $\foD'$ to order $k-1$. In doing this it turns out that 
we need more flexibility for terms $cz^m$ with $\doublebar m=0$.
Hence we also add, as auxiliary input, polynomials $h_1,\ldots, h_r
\in \kk[P]$ with all occurring exponents $m$ fulfilling
$\ord_\tau(m)>0$ and $\doublebar m=0$. These should be thought of as
potential perturbations of $f_{\foc_\mu}$ by undirectional terms. As
compatibility condition for the $h_\mu$ we require
\begin{eqnarray}\label{h_mu condition}
\prod_{\{\rho_\mu\,|\, i(\rho_\mu)=i\}} (1+h_\mu)
\otimes \check d_{\rho_\mu}=1
\end{eqnarray}
in $(R^k_g)^\times\otimes \Lambda_x^*$, for $i=1,\ldots,s$.

Let $\foD_f$ be the infinitesimal scattering diagram
without s-rays or segments, and with
\[
f_{\foc_{\mu},p}=F_{\mu}
\]
for any $p\in\foc_\mu\setminus \{0\}$. Here if $\mu\in\NS(\tau)$, then
$F_{\mu}$ is the chosen element of $\kk[P]$, and if $\mu\not\in\NS(\tau)$,
then $F_{\mu}$ is as constructed in
Proposition~\ref{logpunchline} from the data $\tilde f_1,\ldots,\tilde f_s$.

\begin{figure}[h]
\input{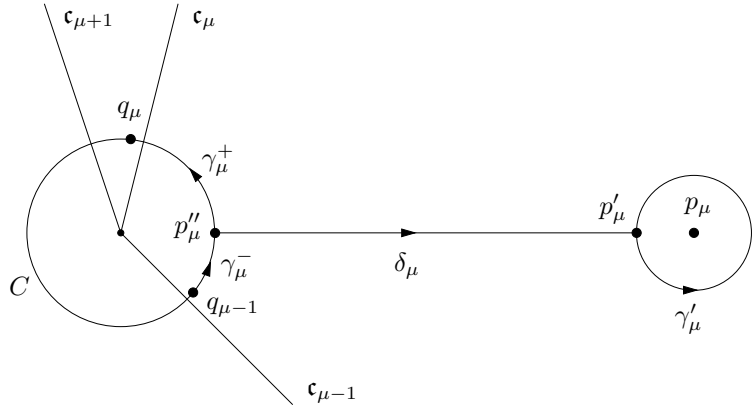}
\caption{Choice of paths. Only a part of the scattering diagram is
shown.}\label{fig:paths}
\end{figure}

For $\mu=1,\ldots,r$ pick a point $p_\mu\in \Int
\doublebar{\sigma_\mu}$. As illustrated in Figure~\ref{fig:paths} let
$p'_\mu$ be a point very close to $p_\mu$ on the line joining the
origin and $p_\mu$, and let $p''_\mu$ be a point on this same line,
but very close to the origin. We can assume $p''_1,\ldots,p''_n$ are
on a small circle $C$ centered at the origin. Let $q_{\mu-1}\in
C\cap\Int \doublebar{\sigma_\mu}$ and close to 
$\foc_{\mu-1}$, so that there are no s-rays of $\foD'_\no$
intersecting $C$ between $C\cap \foc_{\mu-1}$ and $q_{\mu-1}$. Let
$\gamma_\mu^-$,  $\gamma_\mu^+$ be the arcs  of $C$ running from
$q_{\mu-1}$ to $p''_{\mu}$ and from $p''_\mu$ to $q_\mu$,
respectively. Let $\delta_\mu$ be a path connecting $p''_\mu$ to
$p'_\mu$ along the line joining them, and let $\gamma'_\mu$ be a small
loop  around $p_\mu$ based at $p'_\mu$, oriented in the same direction
as $\gamma$, which we take here to be a big loop around the origin. 
The point of these choices is that $\gamma$ is freely homotopic inside
$\shQ\setminus\{0,p_1,\ldots,p_r\}$ to
\begin{eqnarray}\label{tilde gamma}
\tilde\gamma:=\prod_{\mu=1}^r \gamma_\mu^-\delta_\mu\gamma_\mu'\delta_\mu^{-1}
\gamma_\mu^+. 
\end{eqnarray}
Let furthermore $q\in\shQ\setminus \big(\{p_1,\ldots, p_r\}\cup
\bigcup_\mu \foc_\mu\big)$ be a point not on any of the chosen paths
and encircled by $\gamma$ but not by $\tilde\gamma$.

For each $\foz\in\foD'_\inc$, choose a point
$q_{\foz}\in\Int(\foc_{\mu})$, for $\mu$ such that $\foz\subseteq
\doublebar{\sigma_\mu}$. Taking all $q_{\foz}$'s distinct and outside
$C$, we set
\[
\foD_\inc:=\{(\foz+q_{\foz},\theta_{\foz})\,|\,(\foz,
\theta_{\foz})\in\foD'_\inc\}.
\]
This has the effect of translating each incoming s-ray so its endpoint
lies on a cut, but not the origin. Choose the $q_{\foz}$'s so no
element of $\foD_\inc$ passes through $q$ or one of the $p_\mu$'s.
Also write $\foD_\no:=\foD'_\no$.

Set
\[
J_l:=I_0^{l+1}+I_{k},
\]
so for sufficiently large $l$, $J_l=I_{k}$. We now construct
inductively infinitesimal scattering diagrams $\foD_0\subseteq
\foD_1\subseteq\cdots$ for $G_{\foj}^{k}$ with the following
properties.
\begin{enumerate}
\item Modulo $J_l$,
\[
\theta_{\gamma'_\mu,\foD_l}=
\theta_{\delta_\mu,\foD_l}\circ
\theta_{\gamma_\mu^+,\foD_\no}^{-1}\circ
(\theta^{\rm ns}_{\mu})^{-1}\circ
(\theta_\mu^-)^{-1}\circ\theta_{\mu-1}^+\circ
\theta_{\gamma_\mu^-,\foD_\no}^{-1}
\circ\theta_{\delta_\mu,\foD_l}^{-1},
\]
where $\theta_\mu^\pm$ are the log automorphisms of $R^k_g$
constructed in Proposition~\ref{logpunchline}, (the identity if $\mu\in\NS(\tau)$) and
\[
\theta^{\rm ns}_{\mu}=\begin{cases} 
m\mapsto F_{\mu}^{-\langle \overline m,\check d_{\rho_{\mu}}\rangle}&
\mu\in\NS(\tau),\\
1&\mu\not\in\NS(\tau).
\end{cases}
\]
\item
For $\gamma'$ any loop around a singular point
of $\foD_l$ \emph{other} than the origin or any $p_\mu$,
$\theta_{\gamma',\foD_l}=1\mod J_l$.
\item 
No non-foundational element of $\foD_l\setminus \foD_\no$ intersects
$C$ or its interior. 
\item
For each $\foz\in\foD_{l+1}\setminus\foD_l$, $\theta_\foz\equiv 1\mod
J_l$.
\item
If $p\in\Sing(\foD_l)\cap(\foc_\mu\setminus\{0\})$ for some
$\mu$, then either (a) there is an undirectional s-ray $\foz$ with
$p\in\Int(\foz)$ and $\foz$ is the only non-foundational element
containing $p$, or (b) there is exactly one incoming segment or s-ray
with endpoint $p$, and all other non-foundational elements
$\foz\in\foD$ containing $p$ are oriented s-rays or segments with
endpoint $p$ which lie on the other side of $\foc_\mu$, with
$\doublebar{m_{\foz}}$ in the same connected component of
$\shQ\setminus\RR\foc_\mu$ as $\Int \foz$.
\item
Given $\foz\in\foD_\inc$, all s-rays in $\foD_l\setminus \foD_\inc$
asymptotically parallel to $\foz$ are encountered by $\gamma$ (the
large loop around the origin) before encountering $\foz$.
\item
The only elements $\foz$ of $\foD_l$ containing $q$ are unoriented
s-rays $\fol_\mu$, $\mu=1,\ldots,r$, with endpoint $q$,
$\forr(\fol_\mu)=\foc_\mu$ and $\theta_{\fol_\mu}= \exp\big(
-\log(1+h_\mu)\partial_{\check d_{ \rho_\mu}} \big)$.
\end{enumerate}
Note that~(\ref{h_mu condition}) together with~(7) implies (2) for
small loops around $q$.

\begin{construction}\label{construction foD_l}
To start the inductive construction of $\foD_l$ consider $\foD_f\cup
\foD_\no\cup\foD_\inc$. This infinitesimal scattering diagram fulfills
(1)--(6) for $l=0$. To achieve (7) insert unoriented s-rays with
endpoint $q$ as demanded, observing (3), and crossing any cut at most
once. Note that we can not in general avoid crossing cuts, so we need
to adjust the functions $f_{\foc_\mu,p}$ to achieve (2) for small
loops around such intersection points. The adjustments are made
inductively along each cut $\foc_\mu$, treating the points in
$\foc_\mu\cap \bigcup_{\mu'} \fol_{\mu'}$ in the order encountered
along $\foc_\mu$ starting from $0$.

At such a point $p\in\foc_\mu \cap \fol_{\mu'}$ denote by $\foc_\pm$
the connected components of $\foc_\mu\setminus \bigcup_{\mu''}
\fol_{\mu''}$ adjacent to $p$  such that $\foc_+$ is contained in the
unbounded part of $\foc_\mu\setminus \{p\}$, and let $f_{\foc_\pm}$
and $\theta_{\foc_\pm}$ be the associated polynomials and log
automorphisms, respectively. Then for a small counterclockwise loop
$\gamma'$ around $p$ with appropriate base point we find
\[
\theta_{\gamma',\foD_l} =\theta_{\foc_+}\circ\theta_{\fol_{\mu'}}^{\mp 1}\circ
\theta_{\foc_-}^{-1}\circ \theta_{\fol_{\mu'}}^{\pm 1}.
\]
The signs for $\theta_{\fol_{\mu'}}$ are determined by the orientation
of the normal bundle of $\fol_{\mu'}$. Then the same arguments as in
Step~3 of the codimension one case (\S\ref{subsection:codimen 1}) give
\begin{eqnarray}\label{constr:f_{foc_mu}}
\theta_{\gamma',\foD_l}(m) =
\bigg(\frac{(\overline{\theta_{\fol_{\mu'}} })^{\mp 1} 
(f_{\foc_-})}{f_{\foc_+}} \bigg)^{\langle \overline m,
\check d_{\rho_\mu}\rangle},
\end{eqnarray}
see~(\ref{codimonestep3comp}) with $\theta_{\pm}
=\theta_{\fol_{\mu'}}^{\pm 1}$. Now replace all polynomials
$f_{\foc_\mu, p'}$ for $p'$ in the unbounded part of
$\foc_\mu\setminus\{p\}$ by $(\overline{\theta_{l_{\mu'}}})^{\mp 1} 
(f_{\foc_-})$. Continuing in this fashion along $\foc_\mu$ and for all
$\mu$ defines $\foD_0$. Note $\foD_0$ fulfills (2) for loops around
$\foc_\mu\cap \fol_{\mu'}$ for any $l$.
\medskip

We now construct $\foD_{l+1}$ from $\foD_l$ by adding new s-rays and
segments at the singular points of $\foD_l$ and the $p_\mu$.
\smallskip

\noindent
\emph{Step 1.} Let $p$ either be a singular point of $\foD_l$ with $p\in
\Int\doublebar{\sigma_\mu}$ or $p=p_\mu$ for some $\mu$. Let $v^-,
v^+$ be primitive generators of $\foc_{\mu-1}$ and $\foc_\mu$,
respectively. Let $\gamma_p$ be a small loop around $p$, oriented in
the same direction as $\gamma$. Write
\[
\theta_p:=\begin{cases}
\theta_{\delta_\mu,\foD_l}\circ
\theta_{\gamma_\mu^+,\foD_\no}^{-1}\circ
(\theta^{\rm ns}_{\mu})^{-1}\circ
(\theta_\mu^-)^{-1}\circ\theta_{\mu-1}^+\circ
\theta_{\gamma_\mu^-,\foD_\no}^{-1}
\circ\theta_{\delta_\mu,\foD_l}^{-1}&p=p_\mu\\
1&p\not=p_\mu.
\end{cases}
\]
By the inductive assumption~(1),
\[
\theta_p\circ\theta^{-1}_{\gamma_p,\foD_l}=1\mod J_l,
\]
so we can write
\[
\theta_p\circ
\theta^{-1}_{\gamma_p,\foD_l}=\exp \big(\sum c_iz^{m_i}\partial_{n_i}
\big) \mod J_{l+1},
\]
with $\sum c_iz^{m_i}\partial_{n_i}\in \ker\big( \fog_{\foj}^{J_{l+1}}
\to  \fog_{\foj}^{J_l} \big)$. Let $S$ be the set of indices
$S=\{i\,|\,\doublebar{m_i}\not=0\}$. For each $i\in S$, take $\foz_i$
to be a suitably chosen outgoing s-ray with endpoint $p$,
asymptotically parallel to $-\RR_{\ge 0}\doublebar{m_i}$. It will need
to be chosen so $\foz\cap\Sing(\foD_l)=\{p\}$, and so that it passes
through the maximal cones $\doublebar\sigma_\mu$ of $\Sigma_\tau$ in
the same order $p-\RR_{\ge 0}\doublebar{m_i}$ passes through these
cones. In addition, for  $i\not\in S$, we can write
\[
c_iz^{m_i}\partial_{n_i}=c_{i}^- z^{m_i}\partial_{\check v^-}
+c_{i}^+z^{m_i}\partial_{\check v^+},
\]
where $\check v^\pm\in\Lambda_{\foj}^{\perp}\otimes\QQ$ are
the dual basis to $v^\pm$. Define \emph{unoriented} s-rays
$\foz^\pm \subseteq \Int \doublebar{\sigma_\mu}$ with
endpoint $p$, asymptotically parallel to $\foc_{{\mu-1}}$ and
$\foc_{\mu}$, respectively, and with
\[
\theta_{\foz^\pm}:= \exp\Big( \sum_{i\not\in S} c_{i}^\pm z^{m_i}
\partial_{\check v^\pm}\Big).
\]
Setting
\[
\foD(p)= \big\{ \big(\foz_i,\exp(c_iz^{m_i}\partial_{n_i})\big)
\,\big|\,i\in S\big\}\cup\big\{\big(\foz^\pm,
\theta_{\foz^\pm} \big)\big\},
\]
we then have $\theta_{\gamma_p, \foD_l\cup\foD(p)}=\theta_p \mod
J_{l+1}$, as desired.

However, we cannot just add $\foD(p)$ to $\foD_l$, because while
the s-rays of $\foD(p)$, being the identity modulo $J_l$,
commute  modulo $J_{l+1}$ with all non-foundational elements of
$\foD_l$, they do not commute with foundational elements.

To rectify this, we modify $\foD(p)$ as follows. Replace an s-ray
$(\foz,\theta_{\foz})\in\foD(p)$ with $\foz_1,\ldots,\foz_b$, where
$\foz_1,\ldots,\foz_b$ are the closures of the connected components of
$\foz\setminus \bigcup_\mu \foc_{\mu}$. Of course there are a
finite number of such components, and if they are ordered so that
$p\in\foz_1$ and $\foz_i \cap\foz_{i+1}\neq \emptyset$, then $\foz_i$
is a segment for $i<b$ while $\foz_b$ is an s-ray. Note that if $\foz$ was
unoriented, then by the construction of $\foD(p)$, $\foz$ is in fact
contained in the interior of $\doublebar{\sigma_\mu}$ anyway, and no
modification of $\foz$ is necessary. Otherwise, $\foz$ is an outgoing
s-ray, and we define $\theta_{\foz_i}$ inductively as follows,
starting with $\theta_{\foz_1}=\theta_{\foz}$. Let $\mu'$ be chosen so
that $\foz_i\cap\foz_{i+1}\in \foc_{\mu'}$. If $\theta_{\foz_i}=
\exp\big(\sum_{m,n} c_m z^m\partial_n\big)$, then by our choice of $\foz$,
$-\doublebar{m}$ is contained in the connected component of
$\shQ\setminus\RR\foc_{\mu'}$ containing $\foz_{i+1}$. Then we take
\begin{equation*}
\theta_{\foz_{i+1}}=\exp\Big(\sum_{m,n}
c_mz^m f_{\rho_{\mu'}}^{|\langle \overline m, \check d_{
\rho_{\mu'}}\rangle|}\partial_n\Big).
\end{equation*}
Note that $f_{\foc_{{\mu'}},p'}=f_{\rho_{\mu'}} \mod I_0$ for any
$p'\in \foc_{\mu'}\setminus\Sing(\foD_l)$ and $\partial_n
f_{\rho_{\mu'}}=0 \mod I_0$ for $n\in\Lambda_{\foj}^{\perp}$, so that by
Lemma~\ref{conjugation lemma}, for $z^m\in J_l$,
\[
\Ad_{\theta_{\foc_{\mu'}}^{-1}}(z^m\partial_n)
=z^m f_{\rho_{\mu'}}^{\langle \overline m,
\check d_{\rho_{\mu'}}\rangle}\partial_n
\]
for $\langle \overline m,\check d_{\rho_{\mu'}}\rangle >0$ and
\[
\Ad_{\theta_{\foc_{\mu'}}}(z^m\partial_n)
=z^m f_{\rho_{\mu'}}^{-\langle \overline
m,\check d_{\rho_{\mu'}}\rangle}\partial_n
\]
for $\langle \overline m,\check d_{\rho_{\mu'}}\rangle <0$. From this
we conclude $\theta_{\gamma', \{\foc_{{\mu'}},
\foz_i,\foz_{i+1}\}}=1\mod J_{l+1}$ where $\gamma'$ is a loop
around the point $\foz\cap \foc_{\mu'}$. Applying this procedure to
each s-ray in $\foD(p)$, we get a modified $\foD(p)$,  consisting of a
collection of segments and s-rays.
\smallskip

\noindent
\emph{Step 2.}
Assume that $p\in\Sing(\foD_l)\cap (\foc_\mu\setminus\{0\})$. Since 
singular points on cuts contained in an unoriented s-ray or segment
lie on some $\fol_\mu$, which we have already discussed for all $l$,
it suffices to consider the case that $p$ is contained in an
\emph{oriented} s-ray or segment. We then follow a similar but simpler
procedure than in Step~1. By condition~(5)(b), there is precisely one
incoming segment or s-ray $\foz$ with endpoint $p$. Then by
Lemma~\ref{conjugation lemma}, $\theta_{\gamma_p, \{\foz,
f_{\foc_\mu}\}}\in G^{J_{l+1}}_{\foj,K\cup\{0\}}$ (Definition~\ref{cone
groups}), where $K$ is the connected component of $\shQ\setminus\RR
\foc_{\mu}$ disjoint from $\foz$. Thus, using the same technique as
that of the proof of Lemma~\ref{KSlemma} and the previous case, one
can construct a collection of outgoing s-rays $\foD(p)$ with endpoint
$p$ and interiors contained in $K$ such that $\theta_{\gamma_p,
\foD_l\cup\foD(p)}=1\mod J_{l+1}$. We then modify $\foD(p)$ by
subdividing the s-rays as we did above.
\smallskip

We now see if we take 
\[
\foD_{l+1}=\foD_l\cup\bigcup_p\foD(p),
\]
$\foD_{l+1}$ satisfies the inductive properties (1), (2) and (4). For
(1) note that because all monomials occuring in $\foD_\no$ and
$\theta_\mu^\pm$ are in $I_0$, we may replace
$\theta_{\delta_\mu,\foD_{l+1}}$ by $\theta_{\delta_\mu,\foD_l}$ in
the required equation for $\theta_{\gamma'_\mu,\foD_{l+1}}$. Thus (1)
follows from the definition of $\theta_{p_\mu}$ in Step~1 above. With
a further bit of care in making the choices of s-rays above
sufficiently general, all other conditions can be satisfied also. This
completes the construction of the $\foD_l$'s.
\qed
\end{construction}

\begin{lemma}\label{theta(big loop) trivial}
For any $l$, $\theta_{\gamma,\foD_l}=1 \mod J_l$.
\end{lemma}

\proof
By (3) the only elements of $\foD_l$ that $\gamma_\mu^\pm$ crosses are
foundational or in $\foD_\no$. Hence
\[
\theta_{\gamma_\mu^-,\foD_l} = \theta_{\gamma_\mu^-,\foD_\no},\qquad
\theta_{\gamma_\mu^+,\foD_l} = \begin{cases}\theta'_\mu\circ 
\theta_{\gamma_\mu^+,\foD_\no}&\mu\not\in\NS(\tau)\\
\theta^{\rm ns}_{\mu}\circ\theta_{\gamma_\mu^+,\foD_\no}&\mu\in\NS(\tau)
\end{cases}
\]
where $\theta'_\mu= \theta_{\foc_\mu,p}$ is the log automorphism from
Proposition~\ref{logpunchline}, $p$ the intersection point of
$\gamma_\mu^+$ with $\foc_\mu$. In view of property (1) of $\foD_l$,
the definition of $\tilde \gamma$~(\ref{tilde gamma}), and~(\ref{theta
product}) this shows
\begin{eqnarray*}
\theta_{\tilde\gamma,\foD_l}
&=& \prod_{\mu=r}^1 \theta_{\gamma_\mu^+,\foD_l}
\circ\theta_{\delta_\mu,\foD_l}^{-1}
\circ \theta_{\gamma'_\mu,\foD_l}
\circ \theta_{\delta_\mu,\foD_l}
\circ\theta_{\gamma_\mu^-,\foD_l}\\
&=& \prod_{\mu=r}^1 \theta_{\gamma_\mu^+,\foD_l}
\circ \theta_{\gamma_\mu^+,\foD_\no}^{-1}
\circ (\theta^{\rm ns}_\mu)^{-1}
\circ \big(\theta_\mu^-\big)^{-1}
\circ \theta_{\mu-1}^+
\circ \theta_{\gamma_\mu^-,\foD_\no}^{-1}
\circ \theta_{\gamma_\mu^-,\foD_l}\\
&=& \prod_{\mu=r}^1 \theta'_\mu\circ \big(\theta_\mu^-\big)^{-1}
\circ\theta_{\mu-1}^+
\ =\ 1.
\end{eqnarray*}
By property~(2) we conclude
\[
\theta_{\gamma,\foD_l}=\theta_{\tilde\gamma,\foD_l}=1\mod J_l,
\]
because $\gamma$, being a big loop around the origin, is freely
homotopic to $\tilde\gamma$ in $\shQ\setminus\{0,p_1,\ldots,p_r\}$.
\qed

\begin{lemma}\label{theta_forr is in group}
For a rational half-line $\forr\subseteq\shQ$, let $\theta_\forr$ be
the contribution to $\theta_{\gamma,\foD_l}$ from outgoing s-rays
asymptotically parallel to $\forr$. Then $\theta_\forr\in \tilde
H_{\foj}^{J_l}$, and if $\forr$ is not a cut then $\theta_\forr\in
\lperp H_{\foj}^{J_l}$.
\end{lemma}

\proof 
Note automorphisms attached to non-oriented or incoming s-rays are in 
$H_{\foj}^{J_l}$ in any event, and hence preserve
$\Omega_{\std}$. Automorphisms attached to foundational elements
preserve $\Omega_\std$, as follows from their explicit form. Together
with Proposition~\ref{theta(big loop) trivial} this shows
\begin{equation}
\label{omegapreserving}
\Omega_\std=\theta_{\gamma,\foD_l}(\Omega_\std)
=(\theta_{\forr_t}\circ\cdots\circ\theta_{\forr_1})(\Omega_\std)\mod J_l,
\end{equation}
where $\{\forr_1,\ldots,\forr_t\}$ is the finite set of asymptotic directions
of outgoing rays in $\foD_l$ in the order encountered by $\gamma$.
Assume we have shown inductively that each $\theta_{\forr_i}$ preserves
$\Omega_\std$ modulo $J_l$, the base case $l=0$ being trivial. Then
modulo $J_{l+1}$, we can factor  $\theta_{\forr_i}=\theta_{i,1}
\circ\theta_{i,2}$, where $\theta_{i,1}\in
\tilde H_{\foj}^{J_{l+1}}$, and $\theta_{i,2}=\exp \big(\sum_m
z^{m}\partial_{n_m} \big)$ with $n_m\in \Lambda_\foj^\perp
\otimes\kk$, $z^{m}\in J_l$ and $-\doublebar{m} \in \forr_i$.
Then by Remark~\ref{rem:log automorphisms},(3)
\[
\theta_{\forr_i}(\Omega_\std)=\theta_{i,2}(\Omega_\std)
= \Big(1+{\sum}_m \langle \overline{m}, n_m\rangle z^{m} \Big)
\Omega_\std\mod J_{l+1}.
\]
However, monomials $z^{m'},\,z^{m''}$ with
$-\RR_{\ge0}\doublebar{m'}\neq -\RR_{\ge0}\doublebar{m''}$ can never
cancel, so in order for~(\ref{omegapreserving}) to hold modulo
$J_{l+1}$, we in fact must have $\langle \overline m,n_m\rangle=0$ for
each $m$. Thus $\theta_{\forr_\mu} \in \tilde H_{\foj}^{J_{l+1}}$.
Furthermore, if $\forr_\mu$ is not a cut, then $\theta_{\forr_\mu}\in
\lperp H^{J_{l+1}}_\foj$ since only outgoing s-rays contribute to
$\theta_{\forr_\mu}$, and these only involve monomials $z^m$
with $\doublebar m\neq0$.
\qed
\medskip

Now take $l$ sufficiently large so that $J_l=I_{k}$. Similar to
Construction~\ref{asympscatterdiagram} we now construct an asymptotic
scattering diagram  $\foD$ by following a procedure for each rational
half-line $\forr\subseteq \shQ$.

First suppose $\forr$ is not a cut. Consider the contribution $\theta_\forr$ to
$\theta_{\gamma,\foD_l}$ from all s-rays asymptotically parallel to
$\forr$. By property~(6) of $\foD_l$ we can write
\[
\theta_\forr=\theta_{\inc,\forr}\circ\theta_{\no,\forr}\circ\theta'_\forr,
\]
where $\theta_{\inc,\forr}, \theta_{\no,\forr}$ and $\theta'_\forr$
are the contributions from elements of $\foD_\inc$ and from
non-oriented and outgoing s-rays, respectively. Note that by the
explicit commutator formula~(\ref{Lie bracket}), automorphisms
attached to non-oriented s-rays asymptotically parallel to $\forr$
commute with automorphisms attached to oriented s-rays asymptotically
parallel to $\forr$. Now by Lemma~\ref{theta_forr is in group},
$\theta'_\forr\in \lperp H^{k}_{\foj,\forr}$, so by
Lemma~\ref{KSlemma} we can write $\theta'_\forr= \prod_{\forr'\in
\foD_\forr} \theta_{\forr'}$ for $\foD_\forr$ a set of rays with
support~$\forr$ (in particular $\theta_{\forr'} \in
{}^{\perp}H^k_{\foj,\forr}$).

If $\forr=\foc_{\mu}$ for some $\mu$ we are not allowed to add
outgoing rays with support $\forr$ and rather need to modify
$f_{\foc_\mu}$. In this case, the contribution to
$\theta_{\gamma,\foD_l}$ from s-rays asymptotically parallel to
$\forr$ takes the form $\theta_1\circ \theta_{\foc_{\mu}}
\circ\theta_2$, where $\theta_1,\theta_2\in G_{\foj}^{k+1}$ and
$\theta_1 \circ\theta_2 \in\tilde  H^{k}_{\foj,\foc_\mu}$, by
Lemma~\ref{theta_forr is in group}. Recall that
$\theta_{\foc_{\mu},p}$ is given by $m\mapsto
f_{\foc_\mu,p}^{-\langle\overline m,\check d_{\rho_\mu}\rangle}$ for
$p\in\foc_\mu$ far away from the origin. Note in addition that for
$m\in P$, $\theta_1(m),\theta_2(m)$ can be written as a sum of terms
$cz^{m'}$ with $-\doublebar{m'}\in \foc_{\mu}$. It then follows that
if $\theta'_{\foc_{\mu}}$ is defined by $m\mapsto \overline{\theta_1}
(f_{\foc_\mu,p})^{-\langle \overline  m,\check d_{\rho_\mu}\rangle}$,
then
\begin{eqnarray*}
\theta_1\circ\theta_{\foc_{\mu}}(m)&=&\overline{\theta_1}
(f_{\foc_\mu,p}^{-\langle \overline m,\check d_{\rho_\mu}\rangle})
\cdot\theta_1(m)\\
&=&\theta'_{\foc_\mu}(m)
\cdot \overline{\theta'_{\foc_{\mu}}}\big(\theta_1(m) \big)
\ =\ \big(\theta'_{\foc_{\mu}}\circ\theta_1\big)(m).
\end{eqnarray*}
Thus $\theta_1\circ\theta_{\foc_{\mu}}\circ\theta_2=
\theta'_{\foc_{\mu}}\circ\theta_1\circ\theta_2$, and since
$\theta_1\circ\theta_2\in\tilde H_{\foj,\foc_\mu}^k$, we can in fact
write $\theta_1\circ\theta_2$ as $m\mapsto g_\mu^{-\langle \overline
m,\check d_{\rho_\mu}\rangle}$ for some $g_\mu$ only involving
monomials $z^m$ with $\langle\overline m,\check
d_{\rho_\mu}\rangle=0$. Thus the log automorphism $m\mapsto
\big(\overline{\theta_1} (f_{\foc_\mu,p}) \cdot
g_\mu\big)^{-\langle\overline m, \check d_{ \rho_\mu} \rangle}$
coincides with $\theta_1\circ \theta_{\foc_{\mu}} \circ\theta_2$.

Now set
\[
\foD:= \foD'_\inc\cup\foD'_\no\cup{\bigcup}_\forr\foD_\forr\cup
\big\{ \overline{\theta_1}(f_{\foc_\mu,p}) \cdot
g_\mu \big\},
\]
where the union is over all rational half-lines $\forr\subseteq \shQ$ that
are not cuts. By construction, $\theta^k_{\foD}=\theta^k_{\foD_l}=1$.

%=======================

\subsubsection{Construction of $\tilde f_i$, $h_\mu$ and 
$\{F_{\mu}\,|\,\mu\in\NS(\tau)\}$.}
\label{par:tilde f_i}
The diagram $\foD$ constructed in \ref{par:construction of foD}
depends on the choices of $\tilde f_1,\ldots,\tilde f_s$, $h_1,\ldots,
h_r$ and $\{F_{\mu}\,|\,\mu\in\NS(\tau)\}$,  and needs not be an
extension of $\foD'$. We now explain how to make these choices so it
is. We will construct a sequence of choices for the $\tilde f_i$'s,
$F_{\mu}$'s, and $h_\mu$'s, $(\tilde f_i^l)$, $(F^l_{\mu})$,
$(h_\mu^l)$, $l=0,1,\ldots$, each yielding via
Construction~\ref{construction foD_l} a scattering diagram
$\foD(l)=\big\{(\forr,m_\forr, c_\forr), f^l_{\foc_\mu}\big\}$ with
the properties
\begin{enumerate}
\item
$\foD(l)$ is equivalent to $\foD'$ modulo $J_l+I_{k-1}$.
\item
For each $\mu$ any monomial $z^{m}$
appearing in $f^l_{\foc_\mu}- f_{\foc_\mu}\mod J_l$ fulfills
$-\doublebar{m}\in\foc_\mu$.
\end{enumerate}

To begin, take $h_\mu^0=0$, $F_{\mu}^0=1$ for $\mu\in\NS(\tau)$, and
$\tilde f^0_i$ to be given by  the condition that $c_{i,m}=0$ unless
$\height(m)=0$ and $\doublebar m\in \foc_\mu$ for some $\mu$ with
$i=i(\rho_\mu)$. Note that $\tilde f^0_i$ is uniquely determined by
the $f_{\rho}$'s. Furthermore, modulo $J_0=I_0$, $\foD(0)$ is
equivalent to the scattering diagram with functions $f_\rho$ and no
rays: indeed, all non-foundational elements in $\foD(0)$ are
irrelevant modulo $I_0$, and the statement for the foundational
elements follows from Lemma~\ref{thing2prop},(1)(i) and
$F_\mu=\kappa_\mu^{-1} (\tilde F_\mu)$ for $\mu\not\in\NS(\tau)$. 

Now suppose we have constructed $\tilde f_i^l$, $F^l_{\mu}$, and
$h^l_{\mu}$ with the desired properties. For each $\mu$, let 
$g^l_{\rho_\mu}\in \kk[P]$ be the sum of terms in
$f^l_{\foc_\mu}-f_{\foc_\mu}$ of the form $cz^{m} \in J_l\setminus
J_{l+1}$ with $\doublebar{m}\in \foc_\mu\setminus \{0\}$. These are
the terms we need to remove from $f_{\foc_\mu}$ to obtain~(2). Denote
by $\tilde g_{\rho}^l\in \kk[\tilde P]$  what we get by replacing each
term $cz^{m}$ in $g^l_{\rho}$ by $cz^{m+ \sum_j \varphi_j(\overline
m)e_j}$. We then take 
\begin{eqnarray*}
\tilde f_i^{l+1}&=&\tilde f_i^{l}-\sum_{\{\rho\,|\, i(\rho)=i\}}
\tilde g_{\rho}^l,\\ 
F_{\mu}^{l+1}&=&F_{\mu}^l-g_{\rho_{\mu}}^l,\quad \mu\in\NS(\tau).
\end{eqnarray*}
At this point, we don't change the $h_{\mu}$'s.
We now look carefully at how the new scattering diagram
$\foD(l+1)$ differs from $\foD(l)$ modulo $J_{l+1}$.

To do so, we first use Lemma~\ref{thing2prop},(2) to see how for
$\mu\not\in\NS(\tau)$, $\tilde F_{\mu}$ and $\tilde G_\mu^\pm$ change.
Modulo ${}_{\rho_\mu}\hspace{-0.5pt} \tilde J_{l+1}^{>k}$, $\tilde
F_\mu$ is replaced by
\[
\tilde F_\mu-\tilde g_{\rho_\mu}^l+
\sum_{\{m\,|\, -\doublebar{\scriptstyle m}\in
\foc_\mu\setminus\{0\} \}} a_{m}z^{m+\sum_j
\varphi_{j,\rho_\mu}(m)e_l},
\]
where $a_{m}\in \kk[\tilde P_{\rho_\mu}]\setminus
{}_{\rho_\mu}\hspace{-0.5pt} \tilde J_0^{>k}$. As
$F_\mu=\kappa_\mu^{-1}(\tilde F_\mu)$, this has the effect, modulo
$J_{l+1}$, of replacing $F_\mu$ with
\[
F_\mu-g_{\rho_\mu}^l+\hbox{terms of the form $cz^{m}$ with
$-\doublebar{m} \in\foc_\mu\setminus\{0\}$}.
\]
Similarly, $\tilde G_\mu^\pm$ are only changed by terms in
${}_{\rho_\mu}\hspace{-0.5pt} \tilde J_l^{>k} \setminus
{}_{\rho_\mu}\hspace{-0.5pt} \tilde J_{l+1}^{>k}$, so $G_\mu^\pm$ are
only changed by terms in $J_l\setminus J_{l+1}$. Tracing through
Construction~\ref{construction foD_l} one sees that this has the
effect of producing $\foD(l+1)$ with the property that modulo
$J_{l+1}$,
\[
f^{l+1}_{\foc_\mu}-f_{\foc_\mu}=\hbox{sum of terms of the
form $cz^{m}$ with $-\doublebar{m}\in\foc_\mu$}.
\]
The same holds for $\mu\in\NS(\tau)$ directly from the construction of
$F_{\mu}^{l+1}$. This yields the desired condition (2). Unfortunately,
we cannot use the uniqueness statement Proposition~\ref{processunique}
yet to deduce~(1) for $l+1$ because $f_{\foc_\mu}^{l+1}- f_{\foc_\mu}$
may contain monomials $z^m$ with $\overline m\in\Lambda_\foj$. Thus
further modification is necessary.

Modify $\foD(l+1)$ to get an auxiliary scattering diagram
$\hat\foD(l+1)$ as follows. Let $\hat f^{l+1}_{\foc_\mu}$ be obtained
from $f_{\foc_\mu}^{l+1}$ by subtracting those terms $c z^m$ in
$f_{\foc_\mu}^{l+1}-f_{\foc_\mu}$ with $\overline m\in
\Lambda_{\foj}$. Note that by induction hypothesis~(1) for $\foD(l)$,
$\hat f^{l+1}_{\foc_\mu}- f_{\foc_\mu}\in J_l+I_{k-1}$. We then
replace each foundational element in $\foD(l+1)$ by replacing
$f^{l+1}_{\foc_\mu}$ with $\hat f^{l+1}_{\foc_\mu}$ to get
$\hat\foD(l+1)$. Since by construction and induction hypothesis~(1),
$\foD(l), \foD(l+1)$ and $\hat \foD(l+1)$ are all equivalent to
$\foD'$ modulo $J_l+ I_{k-1}$, we have
\[
\theta^{k-1}_{\hat \foD(l+1)}= \theta^{k-1}_{\foD(l)}
=\theta^{k-1}_{ \foD'}=1 \mod J_l+I_{k-1}.
\]
Moreover, by construction $\hat f^{l+1}_{\foc_\mu}- f_{\foc_\mu}$
consists only of terms $z^{m}\in J_l\setminus J_{l+1}$ with
$-\doublebar{m}\in \foc_\mu \setminus\{0\}$. On the other hand, using
Proposition~\ref{influence of monomial changes},(1) to compare
$\theta^{k-1}_{\hat \foD(l+1)}$ and $\theta^{k-1}_{\foD(l+1)} =1$, we
see
\[
\theta^{k-1}_{\hat \foD(l+1)}=
\exp \Big(-\sum_\mu \big[(\hat f_{\foc_\mu}^{l+1}
-f^{l+1}_{\foc_\mu})/f_{\rho_\mu}\big]
\partial_{\check d_{\rho_\mu}}\Big)
\mod J_{l+1}+I_{k-1}. 
\]
Hence
Proposition~\ref{processunique} tells us that $\hat \foD(l+1)$ must be
equivalent to $\foD'$ modulo $J_{l+1}+I_{k-1}$. Thus in particular
\[
\theta^{k-1}_{\hat \foD(l+1)}=1 \mod J_{l+1}+I_{k-1}.
\]
By the normalization condition, if $i(\rho)=i(\rho')$, or equivalently
$Z_{\rho}\cap X_\tau=Z_{\rho'} \cap X_\tau$, then $f_\rho=f_{\rho'}
\mod I_0$. 
Thus again modulo $J_{l+1}+I_{k-1}$,
\[
1=\theta^{k-1}_{\hat \foD(l+1)}
=
\exp\bigg(
-\sum_{\mu\in\NS(\tau)}
\big(\hat f^{l+1}_{\foc_{\mu}}-f^{l+1}_{\foc_{\mu}}\big)
\partial_{\check d_{\rho_\mu}}
-\sum_{i=1}^s f_i^{-1}\bigg(\sum_{\{\mu\,|\, i(\rho_\mu)=i\}}
\big(\hat f^{l+1}_{\foc_{\mu}}-f^{l+1}_{\foc_{\mu}}\big)
\partial_{\check d_{\rho_\mu}} \bigg)
\bigg).
\]
Here $f_i:= f_\rho$ for some $\rho$ with $i(\rho)=i$. Now since
$f_1,\ldots,f_s$ are relatively prime modulo $I_0$ by
Definition~\ref{def:locally rigid},(ii), this is only possible if
$f_i$ divides $\sum_{\{\mu| i(\rho_\mu) =i\}}\big(\hat
f^{l+1}_{\foc_{\mu}} -f^{l+1}_{\foc_{\mu}}\big) \partial_{\check
d_{\rho_\mu}}$ modulo $J_{l+1}+I_{k-1}$, for $i=1,\ldots,s$. We now
use changes of $h^l_\mu$ to turn to zero each term for
$\mu\in\NS(\tau)$ and each sum over $\mu$ in the second summation. If
$\mu\in\NS(\tau)$, take 
\[
h^{l+1}_\mu=h^l_{\mu}+(\hat f^{l+1}_{\foc_{\mu}}-f^{l+1}_{\foc_{\mu}}).
\]
If $\mu\not\in\NS(\tau)$, the polygon $\check\Xi_i$ belonging to $\mu$
according to Remark~\ref{Xi_i, varphi_i} is either a line segment or a
triangle. Let $\mu_\nu$, $\nu=1,2$ or $\nu=1,2,3$ be the corresponding
indices, that is, with $i(\rho_{\mu_\nu})=i$. In the case of a line
segment we have $\check d_{\rho_{\mu_1}}= -\check d_{\rho_{\mu_2}}$,
while in the case of a triangle the $\check d_{\rho_{\mu_\nu}}$
generate $\Lambda_\tau^\perp$ as $\QQ$-vector space. In any case we
can write modulo $J_{l+1}+ I_{k-1}$
\begin{eqnarray}\label{def a_{i,nu}}
\sum_{\nu}\big(\hat
f^{l+1}_{\foc_{\mu_\nu}}-f^{l+1}_{\foc_{\mu_\nu}}\big)
\partial_{\check d_{\rho_{\mu_\nu}}}=
f_i \sum_\nu a_{i,\nu} \partial_{\check d_{\rho_{\mu_\nu}}},
\end{eqnarray}
for some $a_{i,\nu}\in \kk[P]$ containing only monomials
$z^m$ with $\overline{m}\in\Lambda_{\foj}$. Now take
\[
h_{\mu_\nu}^{l+1}=h_{\mu_\nu}^l+a_{i,\nu}.
\]

We can now run Construction~\ref{construction foD_l} again, with the
same $\tilde f_i^{l+1}$, $F_\mu^{l+1}$ as previously, but now with the
newly defined $h_\mu^{l+1}$ rather than $h_\mu^l$. Since we only
changed $h_\mu^l$ by terms $cz^m$ in $J_l$ and with $\overline m\in
\Lambda_\foj$, the infinitesimal scattering diagram $\foD_{l+1}$
remains unchanged modulo $J_l$, while modulo $J_{l+1}$ it differs only
on the automorphisms associated to the undirectional s-rays emanating
from $q$ as given by the modification of $h_\mu^l$. In fact, the
definition of the $f_{\foc_\mu,p}$ from (\ref{constr:f_{foc_mu}})
remains unchanged because $\overline{\theta_{\fol_{\mu'}}}$ acts
trivially on $z^m$ if $\overline m\in\Lambda_\foj$, while if
$\overline m\not\in \Lambda_\foj$ then
$\overline{\theta_{\fol_{\mu'}}}(z^m)$ differs from $z^m$ by terms in
$J_{l+1}$.

The effect to the scattering diagram $\foD(l+1)$ is that modulo
$J_{l+1}+ I_{k-1}$, for $\mu\in\NS(\tau)$, $f^{l+1}_{\foc_{\mu}}$ gets
replaced by  $\hat f^{l+1}_{\foc_{\mu}}$. For
$\mu_\nu\not\in\NS(\tau)$, we add $f_i a_{i,\nu}$ to 
$f^{l+1}_{\foc_{\mu_\nu}}$. Hence by the definition of $a_{i,\nu}$
in~(\ref{def a_{i,nu}}) we now obtain
\[
\sum_{\{\mu\,|\, i(\rho_\mu)=i\}}
\big(\hat f^{l+1}_{\foc_{\mu}}-f^{l+1}_{\foc_{\mu}}\big)
\partial_{\check d_{\rho_\mu}}=0\mod J_{l+1}+I_{k-1},
\]
for any $i$. From the condition that each $\check\Xi_i$ is a line
segment or triangle with every edge of unit affine length, this can
only hold if $\hat f^{l+1}_{\foc_{\mu}}-f^{l+1}_{\foc_{\mu}}$ is
independent modulo $J_{l+1}+I_{k-1}$ of the choice of $\mu$ with
$i(\rho_\mu)=i$.

We are now in position to modify $\tilde f_i^{l+1}$ a second time by,
for each term $cz^{m}$ in $\hat f^{l+1}_{\foc_{\mu}}
-f^{l+1}_{\foc_{\mu}}$ for any $\mu$ with $i(\rho_\mu)=i$,
adding 
\[
cz^{m+\sum_j \varphi_j(m)e_j}
\]
to $\tilde f_i^{l+1}$. By applying Lemma~\ref{thing2prop},(2)
and the argument already made, $\foD(l+1)$ is modified so that
now $f^{l+1}_{\foc_\mu}-f_{\foc_\mu}$ contains no terms in
$J_l\setminus J_{l+1}$ of the form $z^m$ with
$-\doublebar{m} \not\in\foc_\mu\setminus\{0\}$. Thus by the same
uniqueness arguments, $\foD(l+1)$ coincides with $\foD'$
modulo $J_{l+1}+ I_{k-1}$. This completes the inductive construction of
$\tilde f_i^{l+1}$ and $h_\mu^l$.
\medskip

Now take $l$ sufficiently large so that $J_l=I_k$. Then
$\theta^k_{\foD(l)}=1\mod I_k$, and $\foD(l)$ is equivalent to $\foD'$
modulo $I_{k}$. The diagram $\foD(l)$ is almost what we want: the
functions $f_{\foc_\mu}$ however may still contain terms of the form
$cz^{m}\in I_{k-1} \setminus I_k$ with $\doublebar{m}
\in\Lambda_{\foj}$, which we do not wish to allow in
Proposition~\ref{scattering proposition}. We simply discard these
terms to get $\foD$; by Proposition~\ref{influence of monomial
changes},(1), $\theta^k_{\foD}$ takes the desired form.
\qed

%===========================================================
%===========================================================
\section{Concluding remarks}

We will end with a number of short remarks and observations 
about our construction.

\begin{remark}
The first point to emphasize is the importance of the normalization
procedure carried out in Step III of the algorithm. Observe that given
$\scrS_{k-1}$ consistent to order $k-1$, we actually can produce many
liftings to obtain a structure $\scrS_k$ consistent to order $k$. We
can do so by modifying the structure $\scrS_k$ produced in our
algorithm as follows. Change all the slabs in a given
$\rho\in\P^{[n-1]}$ in the same way: For each vertex $v'$ of $\rho$
choose $c_{v'}\in \kk$, and  add to $f_{\fob,x}$ for $x\in
\fob\setminus \Delta$ the expression 
\[
D(s_e,\rho,v[x])\cdot
s_e\bigg(
\sum_{v'} c_{v'}t^kz^{m^{\rho}_{v[x]v'}}\bigg),
\]
with $e:v[x]\rightarrow\rho$. Doing so does not destroy consistency
for codimension one joints. If such modifications are made for each
$\rho\in\P^{[n-1]}$, then consistency at codimension two joints is a
cocycle condition which, with proper choices of coefficients $c_{v'}$,
can be satisfied. In the case when $B$ is simple, it turns out this
gives all ``well-behaved'' logarithmic $k$-th order liftings of the
$(k-1)$-st order logarithmic deformation of $X$ specified by
$\scrS_{k-1}$. In fact, one can develop a form of logarithmic
deformation theory for log Calabi-Yau spaces, which is done in
\cite{part II}, which explains what ``well-behaved'' means. In this
context, the set of these well-behaved $k$-th order liftings modulo a
suitable equivalence relation is in fact a vector space of the
expected dimension, defined as $H^1$ of a sheaf of logarithmic 
derivations of $X$. What might seem surprising at first about our
construction  is that we construct a canonical choice of lifting:
normally one expects the set of liftings to be a torsor over this
$H^1$, without a canonical choice of origin.

Of course, in mirror symmetry, there is a natural set of coordinates
on the moduli space of Calabi-Yau varieties near a large complex
structure limit point, namely canonical coordinates. The expectation
is that our canonical choice of lifting makes $t$ into  a canonical
coordinate. We do not wish to make this statement precise here, but
just  illustrate in a simple example why this might arise. Consider a
local three-dimensional example. Suppose we have a vertex $v\in\P$
contained in a two-dimensional monodromy invariant affine subspace, a
plane,  as depicted in Figure \ref{localP2example}. The figure only
shows a neighbourhood of $v$ and only those cells of  $\P$ contained
in the plane. The dotted lines represent the discriminant locus, and
the numbers indicate a choice of representative for $\varphi$ near
$v$; the value given is that on a primitive vector in the direction of
the labelled cell. Assuming all monodromy vectors $m^{\rho}_{v'v}$
appearing in this example are primitive and the gluing data is
trivial, the slab function at the point $x$, as depicted, takes the
form
\[
f_{\fob,x}=1+z^{(1,0,0,0)}+z^{(0,1,0,0)}+z^{(-1,-1,0,1)}+\sum_{k\ge 1} a_kt^k,
\]
where $t=z^{(0,0,0,1)}$. The normalization condition dictates the values
of the coefficients $a_k$, which are easily seen to give the sum
\[
-2t+5t^2-32t^3+286t^4-3038t^5+\cdots.
\]
\begin{figure}
\input{localP2example.pstex_t}
\caption{}
\label{localP2example}
\end{figure}
This can be compared with the mirror of the anti-canonical bundle of
$\PP^2$, as described e.g. in \cite{GbZa}, \S4.2; the extra power
series in $t$ means $t$ is a canonical coordinate for this family.
\qed
\end{remark}

\medskip
The next observation is that our construction is integral in a
precise sense. This may be related to some of the observed
arithmetic properties of mirror phenomena.

\begin{theorem} Let $A\subseteq \kk$ be a subring. Given a
proper, locally rigid, positive, toric log CY-pair given by
open gluing data $s=(s_e)_e$ with $s_e$ taking values in $A^{\times}$ for
all $e$, then the structures $\scrS_k$ produced by our algorithm
are also defined over $A$; i.e. for each slab $\fob\in\scrS_k$ and
wall $\fop\in\scrS_k$,
\begin{align*}
f_{\fob,x}&\in A[P_x]\\
c_{\fop}&\in A.
\end{align*}
\end{theorem}

\proof
We will give a sketch of the argument. One needs to check that at each
step of the algorithm, all coefficients are in $A$. In Step I, we need
to know that in Proposition~\ref{scattering proposition},  if $\foD'$
is defined over $A$, so is $\foD$. To check this, we need to check it
is never necessary to divide by an element of $A$. It is easy to see
this is the case if $\codim\sigma_{\foj}=0$ directly from the proof of
that case. Indeed, in the exponentials which appear in the proof, only
a first order expansion is necessary as terms of the form $z^{m_i}\in
I_{k-1}\setminus I_k$ have square zero modulo $I_k$. Thus no
denominators appear in the expansion of the expressions used in the
proof. Note here it is important that the log automorphisms associated
to walls take the form $\exp(-\log(1+cz^m)\partial_n)$ rather than
$\exp(-cz^m\partial_n)$, to guarantee that no denominators appear in
the automorphisms attached to walls.

When $\codim\sigma_{\foj}=1$, a similar analysis of the argument in
\S4.3 shows the same integrality. However, this is not true of the
argument given for $\codim\sigma_{\foj}=2$, but it is faster in both
cases to argue directly. Once one knows that
Proposition~\ref{scattering proposition} is true, one knows that the
naive algorithm which works for codimension zero joints also works for
the other types of joints, as described at the beginning of \S4.4.1.
We omit the details.

Step II presents no additional problems; the relevant relative
homology groups are zero whether the coefficient ring is $\kk$ or $A$,
and thus we only need to modify slabs by adding terms with
coefficients in $A$. Finally, integrality in Step III requires
understanding the  normalization condition better. Consider the
following situation. Suppose we have $\rho\in \P^{[n-1]}$,
$v\in\rho\subseteq\sigma\in\P_{\max}$, so we have the  set 
\[
E=\{m\in P_{\rho,\sigma}|\overline{m}\in K_v\rho\setminus\{0\}\}.
\]
Suppose $f=f_0+g\in A[P_{\rho,\sigma}]$ such that $f_0$ consists of
all terms with zero $\ord_{\rho}$ appearing in $f$, and $g$ contains
only exponents $m$ with $\overline m\in\Lambda_{\rho}$ and with
$\ord_{\rho} m>0$. Furthermore, assume $f_0$ has constant term $a_0\in
A^{\times}$, all other exponents appearing in $f_0$ are contained in
$E$, and $\tlog_v(f)\in (t^k)$. Then we need to show that
$\tlog_v(f)=at^k\mod t^{k+1}$ for $a\in A$. To do this, we can first
of all replace $f$ by $f/a_0$ (recalling that $a_0\in A^{\times}$)
without changing $\tlog_v f$, so we can assume the constant term is
$1$. We will show how to expand $f$ in an infinite product expansion
\[
f=\prod_{m\in P_{\rho,\sigma}} \big(1+a_mz^m\big)
\]
with $a_m\in A$ and which converges in the completed ring
${\vphantom{R^k}}^v\!\widehat{R}^k_{\id_{\rho},\sigma}$. A sufficient
condition to guarantee this convergence is that for any $\nu$, there
are only a finite number of $m$ with $\overline{m}\not \in\nu E$ with
$a_m\not=0$. Then for any given $\nu$, all but a finite number of the
cross-terms in the expansion of the product are in $\nu E$.

We will construct the infinite product expansion in two steps. First,
write $f=f_0(1+g/f_0)$. This factorization can be performed in
${\vphantom{R^k}}^v\!\widehat{R}^k_{\id_{\rho},\sigma}$, as $f_0$ is
invertible in this ring. We then express both $f_0$ and $1+g/f_0$
as infinite products of the desired form.

To express $f_0$ as an infinite product, we proceed inductively, for
each $\nu\ge 1$ writing $f_0$ as a product $\prod_m (1+a_mz^m)$ with
$m\not\in \nu E$, up to terms in $\nu E$. For $\nu=1$, the product is
taken to be empty. For $\nu>1$, suppose that the product
$\prod_m(1+a_mz^m)$ agrees with $f_0$ up to terms in $(\nu-1) E$.
Then  $f_0-\prod_m (1+a_mz^m)$ contains only a finite number of terms
$\sum b_{m'}z^{m'}$ with $m\in (\nu-1) E\setminus\nu E$. We then can
replace $\prod(1+a_mz^m)$ with 
\[
\prod \big(1+a_mz^m\big)\prod \big(1+b_{m'}z^{m'}\big)
\]
to obtain a product which works for $\nu$.

For $1+g/f_0$, we proceed similarly, but now not all the exponents
occuring are in $\nu E$ for some $\nu$. This time we proceed order by
order. Suppose we have writtten $1+g/f_0=\prod_m(1+a_mz^m)$, an
infinite product defined in the ring ${\vphantom{R^k}}^v\!
\widehat{R}^l_{\id_{\rho},\sigma}$, for some  $0\le l<k$, and we wish
to extend the infinite product to work in the ring
${\vphantom{R^k}}^v\!\widehat{R}^{l+1}_{\id_{\rho},\sigma}$. It is
not difficult to see that after expanding $1+g/f_0$ in any of these
rings, writing $f_0=1+\hbox{terms in $E$}$, that for any given $\nu$,
$1+g/f_0$ contains only a finite number of terms not in $\nu E$. The
same is true of $\prod_m(1+a_mz^m)$. Thus the same is true of
\[
1+g/f_0-\prod_m (1+a_mz^m)=\sum_{m'} b_{m'}z^{m'}
\]
in ${\vphantom{R^k}}^v\!\widehat{R}^{l+1}_{\id_{\rho},\sigma}$,
with $\ord_{\rho} m'=l+1$ for each $m'$. Then replace
$\prod(1+a_mz^m)$ with \[
\prod\big( 1+a_mz^m\big)
\prod \big(1+b_{m'}z^{m'}\big).
\] 
All new cross-terms now have $\ord_{\rho}$
larger than $l+1$. Proceeding for $l$ up to $k$, we obtain the full
expansion.

We now observe that if $f=\prod_m (1+a_mz^m)$ as above, then
$\tlog_v f\in \kk[t]/(t^{k+1})$ coincides with
\[
\sum_m \tlog_v (1+a_mz^m)=
\sum_{\{m\,|\,\overline{m}=0\}} \tlog_v(1+a_mz^m).
\]
Of course, this is true if this is a finite product; for the case of
an infinite product, we observe for any given $f$ and $k$, there is
some $\nu\ge 1$ such that if $f$ and $f'$ agree up to terms in $\nu E$,
then $\tlog_v f=\tlog_v f'$. Hence the infinite product can be truncated
to a finite product without changing the value of $\tlog_v$.

Now note in this construction that as we only took products and
subtracted, all of the $a_m$'s are in the ring $A$, given all the
coefficients of $f$ were. Furthermore, by assumption, $\tlog_v f\in
(t^k)$. Since a term of the form $1+a_lt^l$ appears at most once in
the infinite product expansion of $f$, the only way this can happen is
if the factors of the form $1+a_lt^l$ appearing in the product
expansion have $a_l=0$ for $l<k$. Thus $\tlog_v f=\tlog_v
(1+at^k)=at^k \mod t^{k+1}$ for some $a\in A$, as desired. This shows
integrality in Step III of the algorithm.
\qed

\begin{remark}
Finally, we would like to comment on the dependence on the choice of
discriminant locus $\Delta$ in our construction. In fact, our
construction is independent of this choice. The easiest way to see
this is to run our algorithm on $X_0\times\PP^1$, with intersection
complex $B\times [0,1]$ with the product affine structure. We take the
discriminant locus in $B\times [0,1]$ to be an isotopy between two
choices of discriminant locus $\Delta_0\subseteq B\times \{0\}$ and
$\Delta_1\subseteq B\times \{1\}$, chosen so $\Delta$ contains no
rational points. We then run our algorithm for $B\times [0,1]$, and it
is not difficult to see that the structures $\scrS_k$  on $B\times
[0,1]$ restrict to the structures on $B\times\{0\}$ and $B\times
\{1\}$. In particular, the structures given by the two choices of
discriminant locus in fact lead to the same formal toric
degeneration of CY-pairs.
\end{remark}

%===========================================================
%===========================================================

\end{document}